%% file: thesis.tex
\documentclass[11pt,onecolumn,fleqn]{book}

\input{includes.tex}

\usepackage[slantedGreek]{mathptmx}

\begin{document}

\frontmatter
\frontbackchapterstyle

\input{thesis/titlepage.tex}

\bookmark[page=3,level=0]{Contents}
\tableofcontents

\input{thesis/frontmatter/_frontmatter.tex}

\mainmatter
\mainchapterstyle

\chapter*{Introduction}
    \markboth{Introduction}{}
    \addcontentsline{toc}{chapter}{Introduction}
    \input{thesis/frontmatter/introduction.tex}

\chapter{Background}
    \label{chBackground}
    \input{thesis/ch1-background/_background.tex}

\chapter{Categories of natural models}
    \label{chCategoriesOfNaturalModels}
    \input{thesis/ch2-categories-of-natural-models/_categories-of-natural-models.tex}

\chapter{Polynomials and representability}
    \label{chPolynomialsRepresentability}
    \input{thesis/ch3-polynomials-representability/_polynomials-representability.tex}

\chapter{Natural model semantics}
    \label{chSemantics}
    \input{thesis/ch4-semantics/_semantics.tex}

\frontbackchapterstyle
\chapter{Future work}
    \label{chReflection}
    \input{thesis/ch5-reflection/_reflection.tex}

\backmatter
\newpage
\clearpage
\addcontentsline{toc}{chapter}{Index}
\printindex

\cleardoublepage
\phantomsection
\addcontentsline{toc}{chapter}{Bibliography}
\bibliographystyle{alpha}
\bibliography{thesis/bibliography}
\nocite{*}

\end{document}

%% file: includes.tex
\usepackage{afterpage}
\usepackage{amsmath}
\usepackage{accents}
\usepackage{amsfonts}
\usepackage{amssymb}
\usepackage{amsthm}
\usepackage[titletoc]{appendix}
\usepackage[UKenglish]{babel}
\usepackage{bbm}
\usepackage[numbered]{bookmark}
\usepackage[nameinlink]{cleveref}
\usepackage[yyyymmdd]{datetime}
    
\usepackage{enumerate}
\usepackage[shortlabels]{enumitem}
    \setlist[itemize]{topsep=0pt}
    \setlist[enumerate]{topsep=0pt}
\usepackage{fancyhdr}
\usepackage{float}
\usepackage[paper=letterpaper,marginratio={1:1}]{geometry}
\usepackage{graphicx}
\usepackage[hang,flushmargin]{footmisc}
\usepackage{hyperref}
\usepackage[latin1]{inputenc}
\usepackage{makeidx}
    \makeindex
\usepackage{marginnote}
\usepackage{mathtools}
\usepackage[framemethod=TikZ]{mdframed}
\usepackage{multicol}
\usepackage{pdflscape}
\usepackage{rotating}
\usepackage{stmaryrd}
\usepackage{tikz}
    \tikzset{node distance=3cm and 5cm, auto}
    \usetikzlibrary{arrows}
\usepackage[normalem]{ulem}
\usepackage{thmtools}
\usepackage{tikz-cd}
    \usetikzlibrary{decorations.pathmorphing}
\usepackage{type1cm}
\usepackage{url}
\usepackage{xcolor}
\usepackage{wrapfig}

\makeatletter
\def\thm@space@setup{%
  \thm@preskip=\parskip
  \thm@postskip=0pt
}
\makeatother

\makeatletter
\newcommand\footnoteref[1]{\protected@xdef\@thefnmark{\ref{#1}}\@footnotemark}
\makeatother

\makeatletter
\renewcommand{\chaptermark}[1]{\markboth{\@chapapp\ \thechapter.\ #1}{\rightmark}}
\makeatother

\definecolor{thmcol}{HTML}{002277}
\definecolor{thmbdcol}{HTML}{4488CC}
\definecolor{thmbgcol}{HTML}{F7FCFF}
\definecolor{defcol}{HTML}{772200}
\definecolor{defbdcol}{HTML}{CC8844}
\definecolor{defbgcol}{HTML}{FFFCF7}
\definecolor{defqedcol}{HTML}{993300}
\definecolor{excol}{HTML}{007722}
\definecolor{exbdcol}{HTML}{44CC88}
\definecolor{exbgcol}{HTML}{F7FFFC}
\definecolor{exqedcol}{HTML}{009933}
\definecolor{tipcol}{HTML}{770077}
\definecolor{tipqedcol}{HTML}{990099}
\definecolor{optmarkcol}{HTML}{770077}
\definecolor{notecol}{HTML}{2266AA}
\definecolor{citecol}{HTML}{883388}
\definecolor{refcol}{HTML}{115588}

\hypersetup{
    colorlinks=true,
    linkcolor={refcol},
    citecolor={citecol}
}

\definecolor{wtcol}{HTML}{660066}     
\definecolor{ptcol}{HTML}{006666}     
\definecolor{prcol}{HTML}{338800}     
\definecolor{glcol}{HTML}{CC4400}     
\definecolor{ascol}{HTML}{444444}     
\definecolor{ercol}{HTML}{990000}     
\definecolor{lccol1}{HTML}{3333CC}    
\definecolor{lccol2}{HTML}{000066}    
\definecolor{gridlines}{HTML}{DDDDDD} 
\definecolor{fncol}{HTML}{664466}     

\definecolor{lccol1}{HTML}{3333CC}    
\definecolor{lccol2}{HTML}{000066}    
\definecolor{gridlines}{HTML}{DDDDDD} 
\definecolor{fncol}{HTML}{664466}     
\definecolor{hdcol}{HTML}{444444}     

\newcommand{\optmarksymbol}{{\color{optmarkcol}$\circledast$}}

\declaretheoremstyle[
    spaceabove=6pt,
    spacebelow=6pt,
    headfont={\color{thmcol}\fontfamily{bch}\bfseries},
    bodyfont=\normalfont,
    headformat={{\small\NUMBER}}
]{thesispar}

\declaretheoremstyle[
    spaceabove=6pt,
    spacebelow=6pt,
    headfont={\color{thmcol}\fontfamily{bch}\bfseries},
    headpunct={{\par\nobreak}\\},
    notefont={\itshape\mdseries\color{notecol}\small},
    notebraces={--- }{},
    bodyfont=\normalfont,
    headformat={{\small\NUMBER} {\NAME}\NOTE}
]{thesisthm}

\declaretheorem[style=thesisthm, numberwithin=section, Refname={Theorem,Theorems}, name=Theorem]{theorem}

\declaretheorem[style=thesisthm, numbered=no, refname={theorem,theorems}, Refname={Theorem,Theorems}, name=Theorem]{theorem*}

\declaretheorem[style=thesisthm, sibling=theorem, refname={lemma,lemmas}, Refname={Lemma,Lemmas}, name=Lemma]{lemma}

\declaretheorem[style=thesisthm, numbered=no, refname={lemma,lemmas}, Refname={Lemma,Lemmas}, name=Lemma]{lemma*}

\declaretheorem[style=thesisthm, numbered=no, refname={proposition,propositions}, Refname={Proposition,Propositions}, name=Proposition]{proposition*}

\declaretheorem[style=thesisthm, sibling=theorem, refname={corollary,corollaries}, Refname={Corollary,Corollaries}, name=Corollary]{corollary}

\declaretheorem[style=thesisthm, numbered=no, refname={corollary,corollaries}, Refname={Corollary,Corollaries}, name=Corollary]{corollary*}

\declaretheorem[style=thesisthm, sibling=theorem, refname={conjecture,conjectures}, Refname={Conjecture,Conjectures}, name=Conjecture]{conjecture}

\declaretheorem[style=thesisthm, numbered=no, refname={conjecture,conjectures}, Refname={Conjecture,Conjectures}, name=Proposition]{conjecture*}

\declaretheorem[style=thesisthm, sibling=theorem, refname={definition,definitions}, Refname={Definition,Definitions}, name=Definition]{definition}

\declaretheorem[style=thesisthm, numbered=no, refname={definition,definitions}, Refname={Definition,Definitions}, name=Definition]{definition*}

\declaretheorem[style=thesisthm, sibling=theorem, refname={construction,constructions}, Refname={Construction,Constructions}, name=Construction]{construction}

\declaretheorem[style=thesisthm, numbered=no, refname={construction,constructions}, Refname={Definition,Definitions}, name=Definition]{construction*}

\declaretheorem[style=thesisthm, sibling=theorem, refname={schema,schemata}, Refname={Schema,Schemata}, name=Schema]{schema}

\declaretheorem[style=thesisthm, numbered=no, refname={schema,schemata}, Refname={Definition,Definitions}, name=Definition]{schema*}

\declaretheorem[style=thesisthm, numbered=no, refname={definition,definitions}, Refname={Definition,Definitions}, name=Definition]{predefinition*}

\declaretheorem[style=thesisthm, numbered=no, refname={notation,notations}, Refname={Notation,Notations}, name=Notation]{notation*}

\declaretheorem[style=thesisthm, numbered=no, refname={axiom,axioms}, Refname={Axiom,Axioms}, name=Axiom]{axiom*}

\declaretheorem[style=thesisthm, numbered=no, refname={axioms,axioms}, Refname={Axioms,Axioms}, name=Axioms]{axioms*}

\declaretheorem[style=thesisthm, numbered=no, refname={assumption,assumptions}, Refname={Assumption,Assumptions}, name=Assumption]{assumption*}

\declaretheorem[style=thesisthm, sibling=theorem, refname={example,examples}, Refname={Example,Examples}, name=Example]{example}

\declaretheorem[style=thesisthm, numbered=no, refname={example,examples}, Refname={Example,Examples}, name=Example]{example*}

\declaretheorem[style=thesisthm, numbered=no, refname={exercise,exercises}, Refname={Exercise,Exercises}, name=Exercise]{exercise*}

\declaretheorem[style=thesisthm, numbered=no, refname={discussion,discussions}, Refname={Discussion,Discussions}, name=Discussion]{discussion*}

\declaretheorem[style=thesisthm, numbered=no, refname={problem,problems}, Refname={Problem,Problems}, name=Problem]{problem*}

\declaretheorem[style=thesisthm, sibling=theorem, refname={remark,remarks}, Refname={Remark,Remarks}, name=Remark]{remark}

\newenvironment{verification}{%
\begin{proof}[Verification]%
}{%
\end{proof}%
}

\declaretheorem[style=thesispar, sibling=theorem, Refname={Paragraph,Paragraphs}, name={}]{numbered}

\numberwithin{equation}{section}

\renewcommand{\le}{\leqslant}

\renewcommand{\ge}{\geqslant}

\newcommand{\slice}[1]{/_{\hspace{-1pt}#1}}

\newcommand{\pullback}{\arrow[dr, draw=none, "\text{\Large$\lrcorner$}" description, pos=0.1]}
\newcommand{\pullbackc}[3][0]{\arrow[#2, draw=none, "\text{\rotatebox{#1}{\Large$\lrcorner$}}" description, pos=#3]}
\newcommand{\pullbackright}{\arrow[dl, draw=none, "\text{\Large$\llcorner$}" description, pos=0.1]}

\newcommand{\sqlabel}[2]{\arrow[#1, draw=none, "#2" description, pos=0.5]}

\newcommand{\seqbn}[2]{\left(\left. #1 \ \middle\rvert\ #2 \right.\right)}
\newcommand{\sqseqbn}[2]{\left[\left. #1 \ \middle\rvert\ #2 \right.\right]}
\newcommand{\append}{\mathbin{^\frown}}

\newcommand{\Yon}{\mathsf{y}} 
\newcommand{\op}{^{\mathrm{op}}} 
\newcommand{\catel}[2][]{\int_{#1} #2} 
\newcommand{\Fam}[2][]{\mathrm{Fam}_{#1}(#2)} 
\newcommand{\fisc}[1]{\mathbb{S}(#1)} 

\newcommand{\nmmk}[1]{\accentset{\mbox{\large\bfseries .}}{#1}} 
    \newcommand{\nmmkalt}[1]{\nmmk{#1}}
\newcommand{\nmty}[1]{\mathcal{#1}}         

\newcommand{\nmtm}[1]{\nmmk{\mathcal{#1}}}  
\newcommand{\nmtmalt}[1]{\nmmkalt{\mathcal{#1}}}  

\newcommand{\tree}{_{\mathsf{tree}}}
\newcommand{\nmtytree}[1]{\nmty{#1}\tree}
\newcommand{\nmtmtree}[1]{\nmtm{#1}_{\mathsf{tree}}}
\newcommand{\adjsigma}{_{\upSigma}}
\newcommand{\nmtysigma}[1]{\nmty{#1}\adjsigma}
\newcommand{\nmtmsigma}[1]{\nmtm{#1}_{\upSigma}}

\newcommand{\adjunit}{_{\mathbbm{1}}}

\newcommand{\nmtyunit}[1]{\nmty{#1}\adjunit}
\newcommand{\nmtmunit}[1]{\nmtm{#1}_{\mathbbm{1}}}

\newcommand{\adjt}[1]{_{#1}}
\newcommand{\nmtyadjt}[2]{\nmty{#2}_{#1}}
\newcommand{\nmtmadjt}[2]{\nmtm{#2}_{#1}}

\newcommand{\adjpi}{_{\upPi}}

\newcommand{\nmp}[2][]{\mathsf{p}^{#1}_{#2}} 
\newcommand{\nmq}[2][]{\mathsf{q}^{#1}_{#2}} 
\newcommand{\nmu}[2][]{\mathsf{u}^{#1}_{#2}} 
\newcommand{\nmv}[2][]{\mathsf{v}^{#1}_{#2}} 
\newcommand{\nms}[1]{\mathsf{s}(#1)} 
\newcommand{\nmt}[1]{\mathsf{t}_{#1}} 
\newcommand{\wmto}{\to_{\mathrm{wk}}} 
\newcommand{\nmsw}{\mathsf{sw}} 
\newcommand{\pmto}{\rightharpoondown} 
\newcommand{\cext}{\mathbin{\resizebox{2.5pt}{2.5pt}{$\bullet$}}}
\newcommand{\cextalt}{\mathbin{\resizebox{3.5pt}{3.5pt}{$\circ$}}}
\newcommand{\circdot}{\mathbin{\dot{\circ}}}
\newcommand{\timesdot}{\mathbin{\dot{\times}}}
\newcommand{\cart}{^{\text{cart}}} 
\newcommand{\twocell}[1]{\text{\rotatebox[origin=c]{#1}{$\Rightarrow$}}}
\newcommand{\threecell}[1]{\text{\rotatebox[origin=c]{#1}{$\Rrightarrow$}}}

\newcommand{\ZFC}{{\footnotesize ZFC}}

\newcommand{\upP}{\mathrm{P}}

\newcommand{\thmcite}[1]{{\color{black} \normalfont \cite{#1}}}
\newcommand{\thmcitenote}[2]{{\color{black} \normalfont \cite[#1]{#2}}}

\newenvironment{diagram}{%
\begin{center}
\begin{tikzcd}[row sep=huge, column sep=huge]
}{%
\end{tikzcd}
\end{center}
}

\newenvironment{ldiagram}{%
\begin{center}
\begin{tikzcd}[row sep=large, column sep=large]
}{%
\end{tikzcd}
\end{center}
}

\newenvironment{ndiagram}{%
\begin{center}
\begin{tikzcd}[row sep=normal, column sep=normal]
}{%
\end{tikzcd}
\end{center}
}

\newenvironment{sdiagram}{%
\begin{center}
\begin{tikzcd}[row sep=small, column sep=small]
}{%
\end{tikzcd}
\end{center}
}

\newenvironment{cdiagram}[2]{%
\begin{center}
\begin{tikzcd}[row sep=#1, column sep=#2]
}{%
\end{tikzcd}
\end{center}
}

\let\nsum\sum
\let\nprod\prod
\renewcommand{\sum}{\nsum\limits}
\renewcommand{\prod}{\nprod\limits}

\newcommand{\declpoly}[7]{#1 \xleftarrow{#5} #2 \xrightarrow{#6} #3 \xrightarrow{#7} #4}

\newcommand{\pto}{%
  \mathrel{\ooalign{\hfil$\vcenter{
    \hbox{\scalebox{0.7}{$\scriptscriptstyle|\,$}}}$\hfil\cr$\to$\cr}
  }%
}
\newcommand{\pRightarrow}{%
  \mathrel{\ooalign{\hfil$\vcenter{
    \hbox{\scalebox{0.7}{$\scriptscriptstyle|\,$}}}$\hfil\cr$\Rightarrow$\cr}
  }%
}
\newcommand{\pRrightarrow}{%
  \mathrel{\ooalign{\hfil$\vcenter{
    \hbox{\scalebox{0.7}{$\scriptscriptstyle|\,$}}}$\hfil\cr$\Rrightarrow$\cr}
  }%
}

\definecolor{whiteish}{HTML}{FEFEFE}
\definecolor{lightpurple}{HTML}{CC44FF}

\usepackage{titlesec} 
\definecolor{midgray}{HTML}{777777}
\makeatletter
\newcommand{\frontbackchapterstyle}{%
    \titleformat{\chapter}[display]{\LARGE}{\color{midgray} \@chapapp\ \thechapter}{0pt}{\fontfamily{bch}\Huge\bfseries}
}
\newcommand{\mainchapterstyle}{%
    \titleformat{\chapter}[display]{\thispagestyle{empty}\vfill\LARGE\centering}{\color{midgray} \@chapapp\ \thechapter}{0pt}{\fontfamily{bch}\Huge\bfseries}[\vfill~\newpage]
}
\titleformat{\section}[display]{\Large}{\color{midgray} Section \thesection}{-10pt}{\fontfamily{bch}\LARGE\bfseries}
\makeatother

%% file: thesis/titlepage.tex
\begin{center}
\thispagestyle{empty}

{\Huge \mdseries \fontfamily{bch}\selectfont
Algebraic models of\\
dependent type theory}

by

{\large \sc Clive Newstead}



\vfill

\textit{Thesis submitted in partial fulfilment of the requirements for\\
the degree of Doctor of Philosophy in Mathematical Sciences}

\vfill

Monday 20th August 2018

\vfill

\includegraphics[width=0.3\columnwidth]{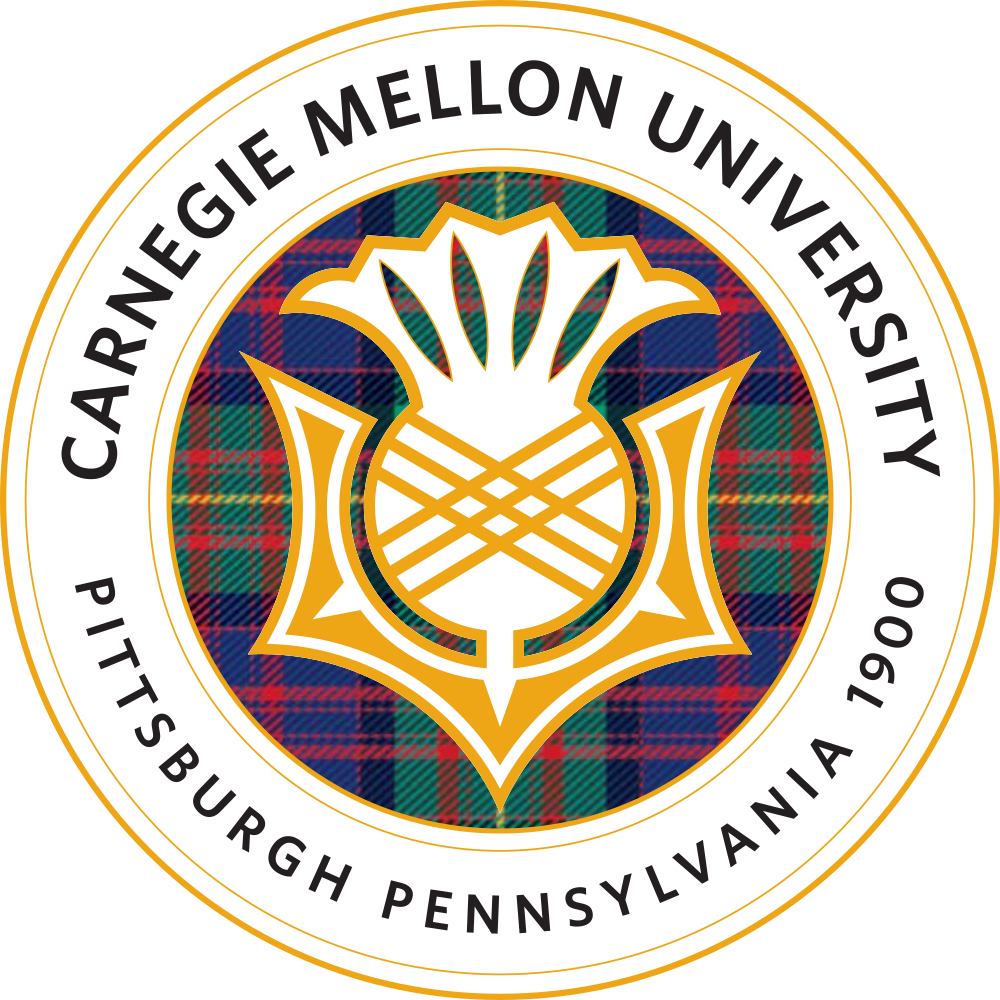}

\vfill

Department of Mathematical Sciences\\
Carnegie Mellon University\\
Pittsburgh, PA

\vfill

\textbf{Doctoral thesis committee}

\begin{tabular}{ccccc}
\textsc{Steve Awodey} && \textsc{Clinton Conley} && \textsc{James Cummings} \\
\textsc{Jonas Frey} && \textsc{Nicola Gambino} && \textsc{Richard Statman}
\end{tabular}
\end{center}

%% file: thesis/frontmatter/_frontmatter.tex
\chapter*{Abstract}
\markboth{Abstract}{}
\addcontentsline{toc}{chapter}{Abstract}
\input{thesis/frontmatter/abstract.tex}

\chapter*{Acknowledgements}
\markboth{Acknowledgements}{}
\addcontentsline{toc}{chapter}{Acknowledgements}
\input{thesis/frontmatter/acknowledgements.tex}

%% file: thesis/frontmatter/abstract.tex
It has been observed \cite{Awodey2016NaturalModels,Fiore2012DiscreteGeneralisedPolynomialFunctors} that the rules governing the essentially algebraic notion of a category with families \cite{Dybjer1995InternalTypeTheory} precisely match those of a representable natural transformation between presheaves. This provides us with a natural, functorial description of essentially algebraic objects which are used to model dependent type theory---following Steve Awodey, we call them \textit{natural models}.

We can view natural models from several different viewpoints, of which we focus on three in this thesis. First, natural models are essentially algebraic, meaning that they can be described by specifying operations between sorts, subject to equational axioms---this allows us to assemble natural models into a category with certain beneficial properties. Second, since natural models are natural transformations between presheaves, they are morphisms in a locally cartesian closed category, meaning that they can be regarded as polynomials \cite{GambinoKock2013PolynomialFunctors}. Third, since natural models admit interpretations of dependent type theory, we can use them to provide a functorial semantics. This thesis develops the theory of natural models in three new directions by viewing them in these three ways.

\textbf{Natural models as essentially algebraic objects.} The first development of the thesis is to bridge the gap between the presentation of natural models as models of an essentially algebraic theory, and the functorial characterisation of natural models as representable natural transformations. We demonstrate that the functorial characterisations of natural models and morphisms thereof align as we hope with the essentially algebraic characterisations.

\textbf{Natural models as polynomials.} The next development is to apply the theory of polynomials in locally cartesian closed categories to natural models. In doing so, we are able to characterise the conditions under which a natural model admits certain type theoretic structure, and under which a natural transformation is representable, entirely in the internal language of a locally cartesian closed category. In particular, we prove that a natural model admits a unit type and dependent sum types if and only if it is a polynomial pseudomonad, that it admits dependent product types if and only if it is a pseudoalgebra, and we prove various facts about the full internal subcategory associated with a natural model.

\textbf{Natural models as models of dependent type theory.} The final development of the thesis is to demonstrate their suitability as a tool for the semantics of dependent type theory. We build the term model of a particularly simple dependent type theory and prove that it satisfies the appropriate universal property, and then we proceed by describing how to turn an arbitrary natural model into one admitting additional type theoretic structure in an algebraically free way.

%% file: thesis/frontmatter/acknowledgements.tex
If I were to give due thanks to everyone who has helped me during my PhD program and my broader development as a mathematician and a scholar, then this section would be the longest of the thesis.

Above all, I am grateful to my doctoral advisor, Steve Awodey, whose influence I cannot understate. I first met him in October 2013 at his house---he was hosting a dinner party following a colloquium in honour of Dana Scott's eighty-first birthday. We discussed category theory and homotopy type theory, and he was kind enough to agree to meet with me to discuss those topics further. The meetings continued and become more focused, and the fruits of our discussions can be found throughout this thesis. I cannot thank him enough for his guidance, patience and generosity.

My thanks extend to my other thesis committee members, who have helped support its progress. In particular, Jonas Frey attended many of my meetings with Steve and provided valuable input on many topics, especially fibrations and locally cartesian closed categories; and Nicola Gambino hosted me at the University of Leeds in March 2017, which catalysed my progress on the work that now constitutes \Cref{secPolynomialPseudomonads}.

I was extremely lucky to be a part of the broader homotopy type theory community at Carnegie Mellon, which includes students, faculty, postdoctoral fellows and visiting scholars spanning three departments. I learnt a great deal from attending the Homotopy Type Theory Seminar for five years, from Bob Harper's course on homotopy type theory in Fall 2013, from Jeremy Avigad's course on interactive theorem proving in Spring 2015, from Jonas Frey's course on categorical logic in Spring 2017, and from countless discussions with my fellow graduate students.

Casting an even wider net, I have learnt much from the international homotopy type theory and category theory communities, and have benefited greatly from discussions and collaborations with many people, particularly Emily Riehl, Chris Kapulkin, Peter Lumsdaine, Ulrik Buchholtz, Pieter Hofstra, Marcelo Fiore, Mike Shulman, Tamara von Glehn, Andr\'{e} Joyal, Pino Rosolini, Guillaume Brunerie, Simon Cho, Cory Knapp and Liang Ze Wong.

Next, I would like to thank the Department of Mathematical Sciences at Carnegie Mellon for supporting me academically and financially throughout my time as a graduate student. My decision to have an advisor in the Department of Philosophy was made seamless by the hard work of the administrative staff and faculty in both departments, and especially Bill Hrusa, Deborah Brandon, Stella Andreoletti, Jeff Moreci, Rosemarie Commisso and Jacqueline DeFazio.

Although the purpose of this thesis is to present my research, of equal importance to me in my mathematical career is my teaching, which has consumed a large part of my time and effort as a graduate student and has been a wonderful experience.

The Department of Mathematical Sciences played an important role in my development as a teacher, offering me opportunities to teach a wide variety of courses both as a teaching assistant and as a course instructor. I was honoured to receive teaching awards from both Carnegie Mellon University and the Mellon College of Science in April 2016 after a departmental nomination, which I am very grateful for.

The opportunity to serve as a teaching assistant for John Mackey in the Fall semesters of 2015 and 2016 transformed me as a teacher, particularly in 2016 when he taught from lecture notes that I had written; I learnt an incredible amount from the experience, and as a result the lecture notes were able to evolve into what I am now calling a `textbook', but has yet to be published.

The Eberly Center for Teaching Excellence and Educational Innovation played a large role in my development as a teacher. Through enrolling in the Future Faculty Program, and then serving as a Graduate Teaching Fellow for three and a half years, I was exposed to the research on teaching and learning. This transformed my own teaching and led to my incorporation of student-centred, evidence-based techniques such as active learning. I would especially like to thank Chad Hershock, Marsha Lovett, Heather Dwyer, Ruth Poproski, Hilary Schuldt, Emily Weiss and Jessica Harrell, as well as the other Graduate Teaching Fellows, for many insightful discussions about teaching and learning over the last five years.

On a more personal note, I would like to thank Bethany, my family, my friends and the Graduate Student Assembly for their roles in making the last five years an enjoyable and fulfilling experience from start to finish.

Finally, I gratefully acknowledge the support of the Air Force Office of Scientific Research through MURI grant FA9550-15-1-0053. Any opinions, findings and conclusions or recommendations expressed in this thesis are my own and do not necessarily reflect the views of the AFOSR.

\begin{flushright}
Clive Newstead\\
August 2018\\
Pittsburgh
\end{flushright}

%% file: thesis/frontmatter/introduction.tex
\subsection*{Outline of the thesis}

In \Cref{chBackground} we provide the fundamental definitions and results underlying the rest of the thesis. We begin with an informal overview of dependent type theory in \Cref{secTypesForTheWorkingMathematician}, followed by a review of polynomials and locally cartesian closed categories in \Cref{secLocallyCartesianClosedCategories} and of presheaves and representability in \Cref{secPresheavesRepresentability}. I do not claim originality for any of the definitions or results in these sections.

\Cref{chCategoriesOfNaturalModels} focuses on natural models in their capacity as models of an essentially algebraic theory. In \Cref{secNaturalModels} we recall the basic definitions and results from \cite{Awodey2016NaturalModels}, before explicitly spelling out the essentially algebraic theories of natural models and of natural models admitting certain type theoretic structure in \Cref{secEssentiallyAlgebraicTheory}, and discussing morphisms of natural models in \Cref{secHomomorphisms}. The goal of this chapter is to provide an explicit demonstration that the convenient, functorial characterisation of natural models as representable natural transformations captures the essentially algebraic notion.

In \Cref{chPolynomialsRepresentability} we enter the realm of locally cartesian closed categories, now viewing natural models in their capacity as polynomials. In \Cref{secInternalCategories} we connect natural models with the theory of internal categories, in particular exploring the properties of their associated full internal subcategory. We build upon this in \Cref{secPolynomialPseudomonads} by using the perspective of internal categories to motivate the definition of a notion of 3-cell in a \textit{tricategory} of polynomials, which allows us to extract a sense in which a natural model admitting a unit type, dependent sum types and dependent product types gives rise to a polynomial pseudomonad and pseudoalgebra. In \Cref{secCharacterisationsOfRepresentability}, we explore the properties possessed by representable natural transformations which can be expressed internally to a locally cartesian closed category.

We sink our teeth into the semantics of dependent type theory in \Cref{chSemantics}, in which we discuss the matter of building the free natural model of a dependent type theory. After introducing the problem and building a free natural model on a very basic type theory in \Cref{secInterpretationsInitiality}, we proceed to discuss how to algebraically freely admit new type theoretic structure to a natural model in \Cref{secExtTerm,secExtType,secExtUnit,secExtSigma}.

A mathematician's work is never done, and this thesis is no exception---in \Cref{chReflection} we discuss some possible directions for future research that are suggested by the work in this thesis.

\subsection*{A remark on foundations}

The official metatheory of this thesis is Zermelo--Fraenkel set theory with the axiom of choice ({\scshape zfc}) together with a fixed \textbf{Grothendieck universe}\index{universe!Grothendieck ---}\label{parGrothendieckUniverses} $\mathfrak{U}$, which is a transitive set containing the von Neumann ordinal $\omega$ and closed under taking power sets and under unions indexed by sets in $\mathfrak{U}$. We will omit reference to $\mathfrak{U}$ by referring to those sets which are elements of $\mathfrak{U}$ as \textbf{small sets}\index{small!{---} set} and to those which are not as \textbf{large sets}. We remark that $\mathfrak{U}$ is itself a model of {\scshape zfc}, and its existence is equivalent (under {\scshape zfc}) to the existence of a strongly inaccessible cardinal. This is one of many solutions to the issues of size arising in category theory---a discussion in far more depth can be found in \cite{Shulman2008SetTheory}---though our results do not depend on which solution is chosen.

\subsection*{Conventions on notation and terminology}

Categories will typically be denoted using calligraphic font $\mathcal{C}, \mathcal{D}, \dots$, but \textit{small}\index{small!{---} category} categories---that is, those whose sets of objects and of morphisms are small---will typically be denoted using blackboard bold font $\mathbb{C}, \mathbb{D}, \dots$.

When working in an $n$-category (strict or otherwise) for $n>1$, the $n$-cells will be denoted by arrows with $n$ horizontal lines. Thus for example a natural transformation from a functor $F$ to a parallel functor $G$ will be denoted by $\varphi : F \Rightarrow G$, and its components by $\varphi_C : F(C) \to G(C)$.

\subsection*{Constructions and verifications}

On many occasions in the thesis, we will define a mathematical object and then prove that it behaves as we say it does. Rather than separating the definition from the theorem proving that the definition makes sense, we combine the two into a single `construction', which reads like a definition of an object, with a subsequent `verification', which reads like a proof that the object we defined behaves as required. See \Cref{cnsPresheafOfTypeTermTrees}, for instance.

\subsection*{Supporting references}

We will assume basic results from category theory and type theory. Useful references for category theory include \cite{MacLane1971CWM} and \cite{Awodey2010CategoryTheory}, references for type theory include \cite{MartinLof1984Intuitionistic} and \cite{HoTTBook2013}, and discussions of categorical models of type theory can be found in \cite{Jacobs1999CategoricalLogic} and \cite{Johnstone2002Elephant}. The results in \Cref{secPolynomialPseudomonads} appear in \cite{AwodeyNewstead2018PolynomialPseudomonads}.

%% file: thesis/ch1-background/_background.tex
\newpage
\input{thesis/ch1-background/types.tex}


\newpage
\input{thesis/ch1-background/lcccs.tex}

\newpage
\input{thesis/ch1-background/representability.tex}

%% file: thesis/ch1-background/types.tex
\section{Dependent type theory}
\label{secTypesForTheWorkingMathematician}

The term \textit{dependent type theory} refers to any one of a number of logical systems derived from those proposed by Per Martin-L\"{o}f in the 1970s (see e.g.~\cite{MartinLof1975DTT} and \cite{MartinLof1984Intuitionistic}), which in turn are descendents of Alonzo Church's $\lambda$-calculus \cite{Church1932Lambda} and, by transitivity, of Bertrand Russell's theory of types \cite{Russell1908Types}.

This section is aimed at a mathematician with a classical mathematical training---that is, first-order logic and \ZFC{} set theory (Zermelo--Fraenkel set theory with the axiom of choice). It intended to be a (very) informal exposition of what dependent type theory is, together with a brief survey of some existing accounts of the semantics of type theory. We will emphasise the similarities and differences between dependent type theory and classical foundations.

Useful references on the syntax and semantics of dependent type theory include \cite{Hofmann1997SyntaxSemantics} and \cite{Jacobs1999CategoricalLogic}.

\subsection*{Types and terms}

The basic objects of dependent type theory are \textit{types} and \textit{terms}. This is in contrast to classical foundations, where all objects are sets. We write $a : A$ to mean that the term $a$ has type $A$. We assume that each term has a unique type---although this assumption is not universally accepted by type theorists, we will need it in order for our notion of `model' (\Cref{defNaturalModel}) to be well-defined.

In some settings, it is helpful for the purpose of intuition to think about a type as being a set, with terms of the type being the elements of the set. In some other settings, it is helpful to think about a type as being a proposition, with the terms of the type being the proofs of the proposition. We will keep this apparent duality between types-as-sets and types-as-propositions, known as the \textit{Curry--Howard correspondence}, in mind.

To illustrate, let $A$ and $B$ be types. We can form their product $A \times B$, whose canonical terms are pairs $\langle a, b \rangle$, where $a : A$ and $b : B$. Under the \textit{types-as-sets} interpretation, we think of $A \times B$ as the cartesian product of $A$ and $B$, whose terms we think of as `ordered pairs of elements'. Under the \textit{types-as-propositions} interpretation, we would think of $A \times B$ as the conjunction of $A$ and $B$, whose terms we think of as `concatenations of proofs': indeed, what is a proof of `$A$ and $B$' if not a proof of $A$ followed by a proof of $B$?

This highlights a key difference between classical foundations and dependent type theory. In classical foundations, we build the theory of sets as a layer on top of first-order logic---the propositions we prove are not themselves the objects of the theory. In dependent type theory, there are just terms and types; we reason about types by constructing terms of new types, which we think about as proofs of propositions. This has the knock-on effect that when we change our theory (say, by adding an axiom), we are also changing the logical system we are working in.

\subsection*{Type dependency, contexts and substitutions}

What sets \textit{dependent} type theory apart from its predecessors is that a type may \textit{depend} on variable terms of other types. For instance, we might consider the type $\mathsf{Vec}_n(\mathbb{R})$ of $n$-dimensional vectors of real numbers, where $n$ is a variable of type $\mathbb{N}$. A list of typed variables that a type may depend on is called a \textit{context}, which is a (possibly empty) list of the form
$$x_1 : A_1, ~ x_2:A_2(x_1), ~ \dots, ~ x_n:A_n(x_1,\dots,x_{n-1})$$
where the parentheses denote the variables the type depends on. We will usually denote contexts by upper-case Greek letters $\Gamma, \Delta, \Theta, \dots$, and we will denote the assertion that $A$ is a type in a context $\Gamma$ by writing $\Gamma \vdash A$, or $\Gamma \vdash A(\vec x)$ if we want to make the variables explicit.

Under the types-as-sets interpretation, a type-in-context $x:A \vdash B(x)$ is interpreted as an $A$-indexed family of sets $(B(x) \mid x \in A)$, or equivalently as a map $B \to A$, where the `set' $B(x)$ corresponds with the preimage of $x \in A$. Under the types-as-propositions interpretation, a type-in-context $x : A \vdash B(x)$ is interpreted as a proposition $B(x)$ depending on a variable term $x : A$, which in turn might be thought of as a hypothesis (together with its proof), or as a variable element of a set.

If $A$ is a type in a context $\Gamma$, we can form the \textit{context extension} of $\Gamma$ by a variable $x$ of type $A$, denoted $\Gamma,~ x:A$; moreover, all contexts can be generated from the empty context by context extension. A type $\vdash A$ in the empty context is called a \textit{basic type}; for example, $\mathbb{N}$ is a basic type.

If $\Gamma \vdash A$ is a type-in-context, we denote by $\Gamma \vdash a : A$ the assertion that $a$ is a term of type $A$ in the presence of the variables $\Gamma$. For example, whenever $\Gamma \vdash A$, it is always the case that $\Gamma,~ x : A \vdash x : A$. We may also assert (definitional) equality of types or of terms, but again relative to a context. Thus $\Gamma \vdash A = B$ asserts that the types-in-context $\Gamma \vdash A$ and $\Gamma \vdash B$ are equal; and $\Gamma \vdash a = a' : A$ asserts that the typed terms-in-context $\Gamma \vdash a : A$ and $\Gamma \vdash a' : A$ are equal. Expressions to the right of the $\vdash$ symbol are called \textit{judgements}.

Given contexts $\Gamma = x_1 : A_1, \dots, x_n : A_n(x_1, \dots, x_{n-1})$ and $\Delta = y_1 : B_1, \dots, y_m : B_m(y_1, \dots, y_{m-1})$, a \textit{substitution} from $\Delta$ into $\Gamma$ is a list of terms
$$\Delta \vdash t_1 : A_1, ~~~ \Delta \vdash t_2 : A_2(t_1), ~~~ \Delta \vdash t_n : A_n(t_1, t_2, \dots, t_{n-1})$$
where $A_i(t_1,\dots,t_{i-1})$ denotes the type obtained by replacing the free variables $x_1,\dots,x_{i-1}$ in $A_i$ by the terms $t_1,\dots,t_{i-1}$, respectively.

The contexts and substitutions (quotiented by provable equality between types and terms) of dependent type theory form a category, called the \textit{category of contexts} of the type theory. Given a substitution $\sigma$ from $\Delta$ to $\Gamma$, write $\Delta \vdash A[\sigma]$ for the result of substituting the variables of $\Gamma$ in a type $\Gamma \vdash A$ according to $\sigma$, and write $\Delta \vdash a[\sigma] : A[\sigma]$ for the result of substituting the variables of $\Gamma$ in a term $\Gamma \vdash a : A$ according to $\sigma$.

\subsection*{Specifying a type}

In \ZFC{} set theory, a set is determined by its elements---this is the content of the axiom of extensionality, which says that two sets with the same elements are equal. In dependent type theory, on the other hand, types are defined according to rules which describe how they interact with other types. In practice, these rules come in four kinds: formation rules, introduction rules, elimination rules and computation rules.
\begin{itemize}
\item \textit{Formation rules} tell us how to build the new type out of old types;
\item \textit{Introduction rules} tell us how to use terms of the old types to obtain terms of the new type;
\item \textit{Elimination rules} tell us how to use the terms of the new type to obtain terms of old types;
\item \textit{Computation rules} tell us how the terms constructed from the introduction and elimination rules interact.
\end{itemize}
These rules are typically specified relative to an arbitrary context $\Gamma$.

The computation rules are further broken down into \textit{$\beta$-reduction} rules, which tell us what happens when we apply an elimination rule after an introduction rule and \textit{$\eta$-expansion} rules, which tell us what happens when we apply an introduction rule after an elimination rule. The $\beta$-reduction and $\eta$-expansion rules can be thought of as the `existence' and `uniqueness' parts, respectively, of universal properties satisfied by the types.

To illustrate, we now proceed by specifying the rules defining the \textit{unit type}, \textit{dependent sum types} and \textit{dependent product types}; these will be of importance to us throughout the thesis.

\begin{definition}[Unit types]
Define the \textbf{unit type} to be the dependent type $\mathbbm{1}$ defined according to the following rules.
\begin{itemize}
\item[(\texttt{$\mathbbm{1}$-F})] $\Gamma \vdash \mathbbm{1}$;
\item[(\texttt{$\mathbbm{1}$-I})] $\Gamma \vdash \star : \mathbbm{1}$;
\item[(\texttt{$\mathbbm{1}$-$\eta$})] If $\Gamma \vdash x : \mathbbm{1}$, then $\Gamma \vdash x = \star : \mathbbm{1}$.
\end{itemize}
\end{definition}

Rule (\texttt{$\mathbbm{1}$-F}) says that $\mathbbm{1}$ is a type in any context; rule (\texttt{$\mathbbm{1}$-I}) says that there is a term $\star$ of type $\mathbbm{1}$ in any context; and rule (\texttt{$\mathbbm{1}$-$\eta$}) says that $\star$ is the unique term of type $\mathbbm{1}$ in any context. There are no elimination or $\beta$-reduction rules for the unit type.

Under the types-as-sets interpretation, we think of the unit type as being a singleton set, whose unique element is $\star$. The $\eta$-expansion rule is what gives us uniqueness of the `element' of $\mathbbm{1}$.

Under the types-as-propositions interpretation, we think of the unit type as being the `true' proposition $\top$. The $\eta$-expansion rule tells us that there is a unique proof of $\top$, which implies that that specifying a proof of a proposition $A$ is equivalent to specifying a proof of $\top \to A$.

\begin{definition}[Dependent sum types]
The \textbf{dependent sum type constructor} $\upSigma$ is defined according to the following rules.
\begin{itemize}
\item[($\upSigma$\texttt{-F})] If $\Gamma \vdash A$ and $\Gamma,\, x:A \vdash B(x)$, then $\Gamma \vdash \sum_{x:A} B(x)$;
\item[($\upSigma$\texttt{-I})] If $\Gamma \vdash a : A$ and $\Gamma \vdash b : B(a)$, then $\Gamma \vdash \langle a, b \rangle : \sum_{x:A} B(x)$;
\item[($\upSigma$\texttt{-E-l})] If $\Gamma \vdash p : \sum_{x:A} B(x)$, then $\Gamma \vdash \mathsf{fst}(p) : A$;
\item[($\upSigma$\texttt{-E-r})] If $\Gamma \vdash p : \sum_{x:A} B(x)$, then $\Gamma \vdash \mathsf{snd}(p) : B(\mathsf{fst}(p))$;
\item[($\upSigma$\texttt{-$\beta$-l})] If $\Gamma \vdash a : A$ and $\Gamma \vdash b : B(a)$, then $\Gamma \vdash \mathsf{fst}(\langle a, b \rangle) = a : A$;
\item[($\upSigma$\texttt{-$\beta$-r})] If $\Gamma \vdash a : A$ and $\Gamma \vdash b : B(a)$, then $\Gamma \vdash \mathsf{snd}(\langle a, b \rangle) = b : B(a)$;
\item[($\upSigma$\texttt{-$\eta$})] If $\Gamma \vdash p : \sum_{x:A} B(x)$, then $\Gamma \vdash p = \langle \mathsf{fst}(p), \mathsf{snd}(p) \rangle : \sum_{x:A} B(x)$.
\end{itemize}
\end{definition}

Under the types-as-sets interpretation, we think of the dependent sum type $\sum_{x:A} B(x)$ as the disjoint union of the $A$-indexed family of sets $(B(x) \mid x \in A)$, with the element $\langle a, b \rangle$ being thought of as the element $b$ in the component of the disjoint union given by the index $a$.

Under the types-as-propositions interpretation, we think of the dependent sum type $\sum_{x:A} B(x)$ as the existentially quantified formula $\exists x : A,~B(x)$, with the proof $\langle a, b \rangle$ being thought of as a specification of a witness $a : A$ together with the proof of $B(a)$.

\begin{definition}[Dependent product types]
The \textbf{dependent product type constructor} $\upPi$ is defined according to the following rules.
\begin{itemize}
\item[($\upPi$\texttt{-F})] If $\Gamma \vdash a : A$ and $\Gamma,\, x : A \vdash B(x)$, then $\Gamma \vdash \prod_{x:A} B(x)$;
\item[($\upPi$\texttt{-I})] If $\Gamma \vdash a : A$ and $\Gamma,\, x : A \vdash b(x) : B(x)$, then $\Gamma \vdash \lambda_{x:A} b(x) : \prod_{x:A} B(x)$;
\item[($\upPi$\texttt{-E})] If $\Gamma \vdash f : \prod_{x:A} B(x)$ and $\Gamma \vdash a : A$, then $\Gamma \vdash \mathsf{app}(f,a) : B(a)$;
\item[($\upPi$\texttt{-$\beta$})] If $\Gamma \vdash a : A$ and $\Gamma, x : A \vdash b(x) : B(x)$, then $\Gamma \vdash \mathsf{app}(\lambda_{x:A} b(x), a) = b(a) : B(a)$;
\item[($\upPi$\texttt{-$\eta$})] If $\Gamma \vdash f : \prod_{x:A} B(x)$, then $\Gamma \vdash f = \lambda_{x:A} \mathsf{app}(f,x) : \prod_{x:A} B(x)$.
\end{itemize}
\end{definition}

Under the types-as-sets interpretation, we think of the dependent product type $\prod_{x:A} B(x)$ as the set of choice functions for the $A$ indexed family of sets $(B(x) \mid x \in A)$; that is, we think of a term $f : \prod_{x : A} B(x)$ as a function $f : A \to \bigcup_{x \in A} B(x)$ such that $\mathsf{app}(f,a)$ ($= f(a)$) $\in B(a)$ for each $a \in A$.

Under the types-as-propositions interpretation, we think of the dependent product type $\prod_{x:A} B(x)$ as the universally quantified formula $\forall x : A,~ B(x)$. A proof $f : \prod_{x:A} B(x)$ is then a family $\mathsf{app}(f,x) : B(x)$ of proofs parametrised by $x : A$.

\subsection*{Proof relevance}

Under the types-as-propositions interpretation, the only way we can assert that a proposition is `true; is by exhibiting a term of the corresponding type. Whereas in first-order logic we can say something like `$\mathbb{N}$ is uncountable', in type theory we would need to first construct a type $A$ asserting (under the types-as-propositions interpretation) that $\mathbb{N}$ is uncountable, and then exhibit a \textit{proof term}, that is a term $\vdash a : A$. This aspect of dependent type theory is known as \textit{proof relevance}, meaning that there is that there is no way to assert the truth of a proposition without also providing a proof---in particular, any proof of an existential statement must provide a witness. This reflects the computational and constructive character of dependent type theory.

\subsection*{Categorical semantics}

We now briefly survey some of the existing notions of categorical model of dependent type theory.

The first notion is that of a \textit{contextual category}, introduced by John Cartmell in his doctoral thesis \cite{Cartmell1978Thesis} and later studied by Vladimir Voevodsky under the name \textit{C-systems} \cite{Voevodsky2016CSystems}.

\begin{definition}[Contextual categories \thmcite{Cartmell1978Thesis,Cartmell1986GAT}]
\index{contextual category}
A \textbf{contextual category} consists of:
\begin{itemize}
\item A small category $\mathbb{C}$ with a terminal object $\diamond$;
\item A tree structure on the objects of $\mathbb{C}$ with root $\diamond$---write $\Gamma \triangleleft A$ to denote the assertion that $\Gamma$ is the parent of $A$ in the tree;
\item For each $\Gamma, A \in \mathrm{ob}(\mathbb{C})$ such that $\Gamma \triangleleft A$, a morphism $\nmp{A} : A \to \Gamma$ in $\mathbb{C}$ and, for each $\sigma : \Delta \to \Gamma$ in $\mathbb{C}$, an object $A[\sigma]$ with $\Delta \triangleleft A[\sigma]$ and a morphism $\sigma \cext A : A[\sigma] \to A$ in $\mathbb{C}$;
\end{itemize}
such that
\begin{enumerate}[(i)]
\item The following square commutes and is a pullback;
\begin{diagram}
A[\sigma]
\arrow[d, "\nmp{A[\sigma]}"']
\arrow[r, "\sigma \cext A"]
&
A
\arrow[d, "\nmp{A}"]
\\
\Delta
\arrow[r, "\sigma"']
&
\Gamma
\end{diagram}
\item $A[\mathrm{id}_{\Gamma}] = A$;
\item $A[\sigma \circ \tau] = A[\sigma][\tau]$ for each $\Theta \xrightarrow{\tau} \Delta \xrightarrow{\sigma} \Gamma$; and
\item $(\sigma \circ \tau) \cext A = (\sigma \cext A) \circ (\tau \cext A[\sigma])$ for each $\Theta \xrightarrow{\tau} \Delta \xrightarrow{\sigma} \Gamma$.
\end{enumerate}
\end{definition}

In Vladimir Voevodsky's \textit{C-systems} approach, the tree structure on $\mathrm{ob}(\mathbb{C})$ is replaced by a grading $(\mathrm{ob}(\mathbb{C})_n)_{n \in \mathbb{N}}$ of the objects of $\mathbb{C}$, together with functions $\mathrm{ob}(\mathbb{C})_{n+1} \to \mathrm{ob}(\mathbb{C})_n$ for each $n \in \mathbb{N}$.

\begin{numbered}
Contextual categories are very close to the syntax of dependent type theory, in the following sense. Viewing $\mathbb{C}$ as the category of contexts and substitutions of a dependent type theory $\mathbb{T}$, the tree structure on $\mathrm{ob}(\mathbb{C})$ gives, for each $\Gamma \in \mathrm{ob}(\mathbb{C})$, a unique factorisation:
$$\Gamma \xrightarrow{!_{\Gamma}} \diamond \qquad = \qquad A_n \xrightarrow{\nmp{A_n}} A_{n-1} \xrightarrow{\nmp{A_{n-1}}} \cdots \xrightarrow{\nmp{A_2}} A_1 \xrightarrow{\nmp{A_1}} \diamond$$
Viewing morphisms $A \to \Gamma$ as dependent types $\Gamma \vdash A$, this tells us that every object $\Gamma$ of $\mathbb{C}$ is built in a finite way from the empty context $\diamond$ by context extension. The terms $\Gamma \vdash a : A$ are then given by sections of $\nmp{A}$. Conditions (i)--(iv) then say that substitution respects typing and respects identity and composition strictly.
\end{numbered}

Another notion of categorical model introduced by John Cartmell in his thesis is that of \textit{categories with attributes}, although presentation we use is due to Andrew Pitts \cite{Pitts2001CategoricalLogic}. 

\begin{definition}[Categories with attributes \thmcite{Cartmell1978Thesis,Pitts2001CategoricalLogic}]
\index{category with attributes}
A \textbf{category with attributes} consists of:
\begin{itemize}
\item A small category $\mathbb{C}$ with a terminal object $\diamond$;
\item For each object $\Gamma$ of $\mathbb{C}$, a set $\mathrm{Ty}(\Gamma)$;
\item For each $\Gamma \in \mathrm{ob}(\mathbb{C})$ and each $A \in \mathrm{Ty}(\Gamma)$, an object $\Gamma \cext A$ of $\mathbb{C}$ and a morphism $\nmp{A} : \Gamma \cext A \to \Gamma$ in $\mathbb{C}$ and, for each $\sigma : \Delta \to \Gamma$ in $\mathbb{C}$, a function $(-)[\sigma] : \mathrm{Ty}(\Gamma) \to \mathrm{Ty}(\Delta)$ and a morphism $\sigma \cext A : \Delta \cext A[\sigma] \to \Gamma \cext A$;
\end{itemize}
such that
\begin{enumerate}[(i)]
\item The following square commutes and is a pullback;
\begin{diagram}
\Delta \cext A[\sigma]
\arrow[d, "\nmp{A[\sigma]}"']
\arrow[r, "\sigma \cext A"]
&
\Gamma \cext A
\arrow[d, "\nmp{A}"]
\\
\Delta
\arrow[r, "\sigma"']
&
\Gamma
\end{diagram}
\item $A[\mathrm{id}_{\Gamma}] = A$;
\item $A[\sigma \circ \tau] = A[\sigma][\tau]$ for each $\Theta \xrightarrow{\tau} \Delta \xrightarrow{\sigma} \Gamma$; and
\item $(\sigma \circ \tau) \cext A = (\sigma \cext A) \circ (\tau \cext A[\sigma])$ for each $\Theta \xrightarrow{\tau} \Delta \xrightarrow{\sigma} \Gamma$.
\end{enumerate}
\end{definition}

\begin{numbered}
We view the elements of $\mathrm{Ty}(\Gamma)$ as depedent types in context $\Gamma$; then the object $\Gamma \cext A$ represents the result of extending the context $\Gamma$ by the type $A$. Note that every contextual category has the structure of a category with attributes: given $\Gamma \in \mathrm{ob}(\mathbb{C})$, take $\mathrm{Ty}(\Gamma) = \{ A \in \mathrm{ob}(\mathbb{C}) \mid \Gamma \triangleleft A \}$, and then define $\Gamma \cext A = A$. The removal of the tree structure on the objects of $\mathbb{C}$ implies that there may be objects that are \textit{not} obtained from the terminal object $\diamond$ by context extension. As such, categories with attributes are further removed from the syntax of dependent type theory. As with contextual categories, terms are interpreted as sections of maps of the form $\nmp{A} : \Gamma \cext A \to \Gamma$.
\end{numbered}

The notion of a \textit{category with families} was introduced by Peter Dybjer in \cite{Dybjer1995InternalTypeTheory}.

\begin{numbered}
Denote by $\mathbf{Fam}$ the category of families of (small) sets. An object of $\mathbf{Fam}$ is a pair $(I, (A_i)_{i \in I})$ consisting of a set $I$ and an $I$-indexed family of sets $(A_i)_{i \in I}$, and a morphism from $(I, (A_i)_{i \in I})$ to $(J, (B_j)_{j \in J})$ is a pair $(f, (g_i)_{i \in I})$ consisting of a function $f : I \to J$ and an $I$-indexed family of functions $(g_i : A_i \to B_{f(i)})_{i \in I}$.
\end{numbered}

\begin{definition}[Categories with families \thmcite{Dybjer1995InternalTypeTheory}]
\label{defCwF}
\index{category with families}
A \textbf{category with families} is a category $\mathbb{C}$ with a distinguished terminal object $\diamond$, together with the following data:
\begin{itemize}
\item A functor $T : \mathbb{C}\op \to \mathbf{Fam}$---we write $T(\Gamma) = (\mathrm{Ty}(\Gamma), \mathrm{Tm}(\Gamma, A)_{A \in \mathrm{Ty}(\Gamma)})$ and denote by $A[\sigma] \in \mathrm{Ty}(\Delta)$ and $a[\sigma] \in \mathrm{Tm}(\Delta, A[\sigma])$ the result of applying $T(\sigma : \Delta \to \Gamma)$ to an element $A \in \mathrm{Ty}(\Gamma)$ and $a \in \mathrm{Tm}(\Gamma, A)$, respectively;
\item For each $\Gamma \in \mathrm{ob}(\mathbb{C})$ and each $A \in \mathrm{Ty}(\Gamma)$, an object $\Gamma \cext A$ of $\mathbb{C}$, a morphism $\nmp{A} : \Gamma \cext A \to \Gamma$ of $\mathbb{C}$ and an element $\nmq{A} \in \mathrm{Ty}(\Gamma \cext A, A[\nmp{A}])$;
\end{itemize}
such that, given any object $\Delta$ of $\mathbb{C}$, morphism $\sigma : \Delta \to \Gamma$ and element $a \in \mathrm{Tm}(\Delta, A[\sigma])$, there is a unique morphism $\langle \sigma, a \rangle : \Delta \to \Gamma \cext A$ such that $\sigma = \nmp{A} \circ \langle \sigma, a \rangle$ and $a = \nmq{A}[\langle \sigma, a \rangle]$.
\end{definition}

\begin{numbered}
As the notation suggests, in a category with families we view the elements of $\mathrm{Ty}(\Gamma)$ as dependent types $\Gamma \vdash A$, and the elements of $\mathrm{Tm}(\Gamma, A)$ as terms $\Gamma \vdash a : A$.
\end{numbered}

The final notion of categorical model that we introduce is that of a \textit{universe category}, introduced by Vladimir Voevodsky \cite{Voevodsky2015UniverseCategory}.

\begin{definition}[Universe categories \thmcite{Voevodsky2015UniverseCategory}]
\index{universe category}
A \textbf{universe category} consists of:
\begin{itemize}
\item A small category $\mathbb{C}$ with a terminal object $\diamond$;
\item A morphism $p : \widetilde{U} \to U$ in $\mathbb{C}$; and
\item A \textbf{universe structure} on $p$---that is, an assignment to each $\Gamma \in \mathrm{ob}(\mathbb{C})$ and each $A : \Gamma \to U$ an object $\Gamma \cext A$ and morphisms $\nmp[\Gamma]{A} : \Gamma \cext A \to \Gamma$ and $\nmq[\Gamma]{A} : \Gamma \cext A \to \widetilde{U}$;
\end{itemize}
such that for each $\Gamma \in \mathrm{ob}(\mathbb{C})$ and each $A : \Gamma \to U$, the following square is a pullback.
\begin{diagram}
\Gamma \cext A
\arrow[d, "{\nmp[\Gamma]{A}}"']
\arrow[r, "{\nmq[\Gamma]{A}}"]
&
\widetilde{U}
\arrow[d, "p"]
\\
\Gamma
\arrow[r, "A"']
&
U
\end{diagram}
\end{definition}

\begin{numbered}
In \cite{Voevodsky2015UniverseCategory}, Vladimir Voevodsky describes how to obtain a C-system from a universe category. As suggested by the notation, we view morphisms $A : \Gamma \to U$ as dependent types $\Gamma \cext A$. The pullback condition tells us that morphisms $a : \Gamma \to \widetilde{U}$ such that $p \circ a = A$ correspond with sections of $\nmp[\Gamma]{A}$, which we can thus think about as terms $\Gamma \vdash a : A$, as we did for contextual categories and categories with attributes.
\end{numbered}

The approach we will use is that of \textit{natural models} \cite{Awodey2016NaturalModels}, which bear similarities with both categories with families and universe categories---they will be defined in \Cref{secNaturalModels} and explored in depth in this thesis.

%% file: thesis/ch1-background/lcccs.tex
\section{Polynomials in locally cartesian closed categories}
\label{secLocallyCartesianClosedCategories}

\begin{definition}[Locally cartesian closed categories]
\label{defLCCC}
\index{category!locally cartesian closed {---}}
\index{cartesian closed!locally {---} category}
A \textbf{locally cartesian closed category} is a category $\mathcal{E}$ with a terminal object $1$ and with all slices $\mathcal{E} \slice{A}$ cartesian closed.
\end{definition}

\begin{numbered}
Locally cartesian closed categories $\mathcal{E}$ are characterised by the fact that every morphism $f : B \to A$ induces a triple of adjoint functors
\begin{diagram}
\mathcal{E} \slice{A}
\arrow[rr, "\upDelta_f" description, ""{name=codtop}, ""'{name=dombot}]
&&
\mathcal{E} \slice{B}
\arrow[ll, bend right, "\upSigma_f"', ""{name=domtop}]
\arrow[ll, bend left, "\upPi_f", ""'{name=codbot}]
\arrow[draw=none, from=domtop, to=codtop, inner sep={1pt}, "\bot" description]
\arrow[draw=none, from=dombot, to=codbot, inner sep={1pt}, "\bot" description]
\end{diagram}
where $\upSigma_f$ is given by postcomposition with $f$ and $\upDelta_f$ is given by pullback along $f$. Since we have adopted the convention that locally cartesian closed categories have a terminal object, it follows that they are cartesian closed and have all finite limits. We emphasise that locally cartesian closed categories are categories with additional \textit{structure}. In particular, given an object $(X,x)$ of $\mathcal{E} \slice{A}$, the functor $\upDelta_f : \mathcal{E}\slice{A} \to \mathcal{E}\slice{B}$ gives a \textit{choice} of pullback $\upDelta_f(x) : \upDelta_f(X) \to B$ of $x : X \to A$ along $f$.
\end{numbered}

\begin{example}
\label{exLCCCs}
Examples of locally cartesian closed categories include the category $\mathbf{Set}$ of sets, the category $\widehat{\mathbb{C}}=\mathbf{Set}^{\mathbb{C}^{\mathrm{op}}}$ of presheaves on a small category $\mathbb{C}$, and more generally, any topos. The category $\mathbf{Cat}$ of categories is \textit{not} locally cartesian closed, even though it is cartesian closed.
\end{example}

\begin{numbered}
\label{parInternalLanguage}
\index{internal language}
Every locally cartesian closed category $\mathcal{E}$ has an \textbf{internal language} \cite{Seely1984LCCCs}, which provides a convenient \textit{syntactic} way of reasoning about the objects and morphisms of $\mathcal{E}$. When reasoning internally, we will view an object $(X, x : X \to A)$ of $\mathcal{E} \slice{A}$ as an $A$-indexed family of objects $\seqbn{X_a}{a \in A}$, and a morphism $(X,x) \to (Y,y)$ in $\mathcal{E} \slice{A}$ as an $A$-indexed family of morphisms $\seqbn{X_a \to Y_a}{a \in A}$. Given a morphism $f : B \to A$ of $\mathcal{E}$, the action on objects of the functors $\upSigma_f$, $\upDelta_f$ and $\upPi_f$ can be described by
$$\upSigma_f \seqbn{Y_b}{b \in B} = \seqbn{\sum_{b \in B_a} X_b}{a \in A}$$
$$\upDelta_f \seqbn{X_a}{a \in A} = \seqbn{X_{f(b)}}{b \in B}$$
$$\upPi_f \seqbn{Y_b}{b \in B} = \seqbn{\prod_{b \in B_a} Y_b}{a \in A}$$
Note that when $\mathcal{E} = \mathbf{Set}$ we really can identify an object $(X,x)$ of $\mathbf{Set} \slice{A}$ as an $A$-indexed family $\seqbn{X_a}{a \in A}$ by defining $X_a = x^{-1}[\{a\}]$ for each $a \in A$. The sum and product operations are realised in this case as the disjoint union and dependent product, respectively.
\end{numbered}

\begin{theorem}[Beck--Chevalley condition]
\label{thmBeckChevalley}
\index{Beck--Chevalley condition}
Let $f,g,u,v$ be morphisms in a locally cartesian closed category $\mathcal{E}$ fitting into the following pullback square.
\begin{ldiagram}
B
\arrow[r, "v"]
\arrow[d, "f"']
\pullback
&
D
\arrow[d, "g"]
\\
A
\arrow[r, "u"']
&
C
\end{ldiagram}
There are natural isomorphisms $\upDelta_g \upSigma_u \cong \upSigma_v \upDelta_f$ and $\upDelta_g \upPi_u \cong \upPi_v \upDelta_f$. \qed
\end{theorem}

\begin{theorem}[Distributivity law \thmcite{Weber2015Polynomials}]
\label{thmDistributivity}
\index{distributivity law}
\index{axiom of choice!type theoretic {---}}
Let $C \xrightarrow{u} B \xrightarrow{f} A$ be morphisms in a locally cartesian closed category $\mathcal{E}$. Construct the following commutative diagram, in which $v = \upPi_f(u)$ is the dependent product of $u$ along $f$, $w = \Delta_f(v)$ is the pullback of $v$ along $f$, and $e$ is the component at $h$ of the counit of the adjunction $\upDelta_f \dashv \upPi_f$.
\begin{diagram}
&
P
\arrow[r, "q"]
\arrow[d, "w" description]
\arrow[dl, "e"']
\pullback
&
D
\arrow[d, "v"]
\\
C
\arrow[r, "u"']
&
B
\arrow[r, "f"']
&
A
\end{diagram}
There is a natural isomorphism $\upPi_f \upSigma_u \cong \upSigma_d \upPi_q \upDelta_e$. \qed
\end{theorem}

\begin{numbered}
In the internal language of $\mathcal{E}$, the Beck--Chevalley conditions say, parametrically in $d \in D$, that
$$\sum_{a \in A_{g(d)}} X_a \cong \sum_{b \in B_d} X_{f(b)} \quad \text{and} \quad \prod_{a \in A_{g(d)}} X_a \cong \prod_{b \in B_d} X_{f(b)}$$
and the distributivity law says, parametrically in $a \in A$, that
$$\prod_{b \in B_a} \sum_{c \in C_b} X_c \cong \sum_{d \in D_a} \prod_{p \in P_d} X_{e(p)} \cong \sum_{m \in \prod_{b \in B_a} C_b}\ \prod_{b \in B_a} X_{m(b)}$$
For this reason, the distributivity law is sometimes referred to as the (\textit{type theoretic}) \textit{axiom of choice}. This is slightly misleading, since although it resembles the axiom of choice, it is a theorem rather than an axiom.
\end{numbered}

\subsection*{Polynomials and polynomial functors}

\begin{definition}[Polynomials \thmcite{GambinoKock2013PolynomialFunctors}]
\label{defPolynomial}
\index{polynomial}
A \textbf{polynomial} $F=(s,f,t)$ in a locally cartesian closed category $\mathcal{E}$ is a diagram of the form
$$\declpoly IBAJsft$$
We say that $F$ is a `polynomial from $I$ to $J$' and write $F : \declpoly IBAJsft$ or just $F : I \pto J$.
\end{definition}

\begin{numbered}
\label{parPolynomialsFromOneToOne}
Polynomials generalise morphisms (by taking $I=J=1$) and spans (by taking $B=A$ and $f=\mathrm{id}_A$). Since most of our attention will be focused on polynomials from $1$ to $1$, we will brazenly blur the distinction between morphisms $f : B \to A$ and polynomials $\declpoly 1BA1{!_{B}}f{!_{A}}$. Beware, though, that composition of polynomials (see \Cref{defCompositionOfPolynomials}) differs from composition of morphisms of $\mathcal{E}$.
\end{numbered}

\begin{definition}[Polynomial functors \thmcite{GambinoKock2013PolynomialFunctors}]
\label{defPolynomialFunctor}
\label{defExtension}
\index{polynomial!functor@{---} functor}
\index{polynomial!extension of a {---}}
The \textbf{extension} of a polynomial $F : \declpoly IBAJsft$ in a locally cartesian closed category $\mathcal{E}$ is the functor $\upP_F = \upSigma_t \upPi_f \upDelta_s : \mathcal{E} \slice{I} \to \mathcal{E} \slice{J}$. Internally, we can define $\upP_F$ by
$$\upP_F \seqbn{X_i}{i \in I} = \seqbn{\sum_{a \in A_j} \prod_{b \in B_a} X_{s(b)}}{j \in J}$$
A \textbf{polynomial functor} is a functor that is naturally isomorphic to the extension of a polynomial.
\end{definition}

When $f : B \to A$ is a morphism of $\mathcal{E}$, we obtain an endofunctor $\upP_f = \upP_{({!}_B,f,{!}_A)} : \mathcal{E} \cong \mathcal{E} \slice{1} \to \mathcal{E} \slice{1} \cong \mathcal{E}$, and this endofunctor is described in the internal language of $\mathcal{E}$ quite simply as
$$\upP_f(X) = \sum_{a \in A} X^{B_a}$$
This explains the use of the term \textit{polynomial}.

We recall the following technical lemma from \cite{Awodey2016NaturalModels}; it will be useful for us later on.

\begin{lemma}[See \thmcitenote{Lemma 5}{Awodey2016NaturalModels}]
\label{lemLemmaFive}
Let $f : B \to A$ be a morphism in a locally cartesian closed category $\mathcal{E}$. There is a natural (in $X$ and in $Y$) correspondence between morphisms $g : Y \to \upP_f(X) = \sum_{a \in A} X^{B_a}$ and pairs $(g_1,g_2)$ of morphisms with $g_1 : Y \to A$ and $g_2 : \upDelta_{g_1}(B) \to X$.
\begin{diagram}
X
&
\upDelta_{g_1}(B)
\arrow[l, "g_2"']
\arrow[r]
\arrow[d, "\upDelta_{g_1}(f)"']
\pullback
&
B
\arrow[d, "f"]
\\
&
Y
\arrow[r, "g_1"']
&
A
\end{diagram}
~\qed
\end{lemma}

The following lemma of a similar flavour will also be useful.

\begin{lemma}
\label{lemLemmaElevenPointFive}
Let $f : B \to A$ be a morphism in a locally cartesian closed category $\mathcal{E}$. There is a natural correspondence between morphisms
$$g : Y \to \sum_{a \in A} \sum_{m \in A^{B_a}} \sum_{b \in B_a} B_{m(b)}$$
and quadruples $(g_1,g_2,g_3,g_4)$ of morphisms, with
\begin{itemize}
\item $g_1 : Y \to A$ in $\mathcal{E}$;
\item $g_2 : \upDelta_{g_1}(B) \to A$ in $\mathcal{E}$;
\item $g_3 : (Y, g_1) \to (B,f)$ in $\mathcal{E} \slice{A}$; and
\item $g_4 : (\upDelta_{g_1}(B), g_2) \to (B, f)$ in $\mathcal{E} \slice{A}$.
\end{itemize}
\end{lemma}

\begin{proof}[Sketch of proof]
The is a direct translation of argument on \cite[pp.~18-19]{Awodey2016NaturalModels} into the more general setting of an arbitrary locally cartesian closed category.
\end{proof}

\begin{definition}[Composition of polynomials \thmcite{GambinoKock2013PolynomialFunctors}]
\label{defCompositionOfPolynomials}
\index{polynomial!composition of {---}s}
The \textbf{polynomial composite} of polynomials $F : \declpoly IBAJsft$ and $G : \declpoly JDCKugv$ in a locally cartesian closed category $\mathcal{E}$ is the polynomial $G \cdot F : \declpoly INMK{s \circ n}{q \circ p}{v \circ w}$ indicated in the following diagram, which is constructed as follows: first take the pullback {\color{blue} (1)}; then form {\color{blue} (2)} from $H \xrightarrow{h} D \xrightarrow{g} C$ as in \Cref{thmDistributivity}; and finally take the pullback {\color{blue} (3)} of $k \circ e$ along $f$.

\begin{ldiagram}
&&
N
\pullbackc{ddr}{0.04}
\arrow[ddl,"n"']
\arrow[r,"p"]
\arrow[ddr,blue,"\text{(3)}" description,draw=none,pos=0.3]
&
L
\pullback
\arrow[r,"q"]
\arrow[d,"e" description]
&
M
\arrow[ddr,"w"]
\arrow[ddl,blue,"\text{(2)}" description,draw=none,pos=0.3]
&
\\
~&~&~&
H
\pullbackc[-45]{dd}{0.05}
\arrow[dr,"h" description]
\arrow[dl,"k" description]
\arrow[dd,blue,"\text{(1)}" description,draw=none]
&~&~&~
\\
&
B
\arrow[r,"f"']
\arrow[dl,"s"']
&
A
\arrow[dr,"t"']
&~&
D
\arrow[dl,"u"]
\arrow[r,"g"']
&
C
\arrow[dr,"v"]
&
\\
I
&~&~&
J
&~&~&
K
\\
\end{ldiagram}
\end{definition}

\begin{numbered}
We will make use of the following explicit descriptions of the objects $H,L,M,N$ in the internal language of $\mathcal{E}$.
\begin{itemize}
\item $H = \sum_{d \in D} A_{u(d)}$;
\item $M = \sum_{c \in C}~\prod_{d \in D_c} A_{u(d)}$;
\item $L = \sum_{(c,m) \in M} D_c$;
\item $N = \sum_{(c,m) \in M}~\sum_{d \in D_c} B_{m(d)}$.
\end{itemize}
The morphisms $e,h,k,n,p,q,w$ are then the appropriate projection morphisms.
\end{numbered}

\begin{theorem}[Extension preserves composition \thmcite{GambinoKock2013PolynomialFunctors}]
\label{thmExtensionPreservesPolynomialComposition}
Let $F : \declpoly IBAJsft$ and $G : \declpoly JDCKugv$ be polynomials in a locally cartesian closed category $\mathcal{E}$. Then $\upP_{G \cdot F} \cong \upP_G \circ \upP_F$, where $\cdot$ represents polynomial composition and $\circ$ represents the usual composition of functors.
\end{theorem}

\begin{proof}
With notation as in \Cref{defCompositionOfPolynomials}, we proceed by calculation.
\begin{align*}
\upP_G \circ \upP_F
&= \upSigma_v \upPi_g \upDelta_u \upSigma_t \upPi_f \upDelta_s
&& \text{by \Cref{defExtension}} \\
&\cong \upSigma_v \upPi_g \upSigma_h \upDelta_k \upPi_f \upDelta_s && \text{by Beck--Chevalley (\Cref{thmBeckChevalley})} \\
&\cong \upSigma_v \upSigma_w \upPi_q \upDelta_e \upDelta_k \upPi_f \upDelta_s && \text{by distributivity (\Cref{thmDistributivity})} \\
&\cong \upSigma_v \upSigma_w \upPi_q \upDelta_{k \circ e} \upPi_f \upDelta_s && \text{by functoriality} \\
&\cong \upSigma_v \upSigma_w \upPi_q \upPi_p \upDelta_n \upDelta_s && \text{by Beck--Chevalley} \\
&\cong \upSigma_{v \circ w} \upPi_{q \circ p} \upDelta_{s \circ n} && \text{by functoriality} \\
&\cong \upP_{G \cdot F} && \text{by \Cref{defExtension,defCompositionOfPolynomials}}
\end{align*}
Each of these isomorphisms is natural and strong.
\end{proof}

\begin{definition}
\label{defMorphismOfPolynomials}
\index{polynomial!morphism of {---}s}
\index{morphism!{---} of polynomials}
Let $F : \declpoly IBAJsft$ and $G : \declpoly IDCJugv$ be polynomials from $I$ to $J$ in $\mathcal{E}$. A \textbf{morphism of polynomials} $\varphi$ from $F$ to $G$ consists of an object $D_{\varphi}$ of $\mathcal{E}$ and a triple $(\varphi_0, \varphi_1, \varphi_2)$ of morphisms in $\mathcal{E}$ fitting into a commutative diagram of the following form, in which the lower square is a pullback:
\begin{diagram}
&
B
\arrow[r,"f"]
\arrow[dl, bend right=20, "s"']
&
A
\arrow[d, equals]
\arrow[dr, bend left=20, "t"]
&
\\
I
&
D_{\varphi}
\pullback
\arrow[r]
\arrow[d,"\varphi_1"']
\arrow[u,"\varphi_2"]
&
A
\arrow[d,"\varphi_0"]
&
J
\\
&
D
\arrow[r,"g"']
\arrow[ul, bend left=20, "u"]
&
C
\arrow[ur, bend right=20, "v"']
&
\end{diagram}
We write $\varphi : F \pRightarrow G$ to denote the assertion that $\varphi$ is a morphism of polynomials from $F$ to $G$.
\end{definition}

Each morphism $\varphi : F \pRightarrow G$ of polynomials induces a strong\footnote{Every polynomial functor has a natural \textit{strength}, and the natural candidate for morphisms between polynomial functors are those natural transformations which are comptable with the strength. See \cite{GambinoKock2013PolynomialFunctors} for more on this.} natural transformation $P_F \Rightarrow P_G$, which we shall by abuse of notation also call $\varphi$, whose component at $\vec X = (X_i \mid i \in I)$ can be expressed in the internal language of $\mathcal{E}$ by
$$(\varphi_{\vec X})_j : \sum_{a \in A_j} \prod_{b \in B_a} X_{s(b)} \to \sum_{c \in C_j} \prod_{d \in D_c} X_{u(b)}; \quad (\varphi_{\vec X})_j(a,t) = (\varphi_0(a), t \cdot (\varphi_2)_a \cdot (\varphi_1)_a^{-1})$$

\begin{definition}
\label{defCartesianMorphism}
A morphism $\varphi : F \pRightarrow G$ is \textbf{cartesian} if $\varphi_2$ is invertible.
\end{definition}

As the name suggests, if $\varphi : F \pRightarrow G$ is a cartesian morphism, then the induced strong natural transformation $P_F \Rightarrow P_G$ is cartesian.

\begin{numbered}
\label{rmkCartesianMorphismIsPullbackSquare}
Every cartesian morphism of polynomials has a unique representation as a commutative diagram of the following form.
\begin{equation} \label{diaCartesianMorphism}
\hspace{72pt}
\begin{tikzcd}[row sep=normal, column sep=huge]
&
B
\arrow[dl, bend right=15, "s"']
\arrow[r, "f"]
\arrow[dd, "\varphi_1"']
\pullbackc{ddr}{0.05}
&
A
\arrow[dd, "\varphi_0"]
\arrow[dr, bend left=15, "t"]
&
\\
I
&&&
J
\\
&
D
\arrow[ul, bend left=15, "u"]
\arrow[r, "g"']
&
C
\arrow[ur, bend right=15, "v"']
&
\end{tikzcd}
\end{equation}

Indeed, if $(\varphi_0,\varphi_1,\varphi_2)$ is cartesian, replacing $\varphi_1$ in the above diagram by $\varphi_1 \circ \varphi_2^{-1}$ yields the desired diagram. Conversely, if $(\varphi_0,\varphi_1)$ are as in the above diagram, then $(\varphi_0,\varphi'_1,\varphi'_2)$ is a cartesian morphism of polynomials, where $\varphi'_1 : \Delta_{\varphi_0}D \to D$ is the chosen pullback of $\varphi_0$ along $g$ and $\varphi'_2 : \Delta_{\varphi_0}D \to B$ is the canonical isomorphism induced by the universal property of pullbacks, as illustrated in the following:
\begin{equation} \label{diaCartesianMorphismFromPullbackSquare}
\begin{tikzcd}[row sep=normal, column sep=normal]
&
B
\arrow[r,"f"]
\arrow[dl, bend right=20, "s"']
\arrow[dd, "\varphi_1"']
\pullback
&
A
\arrow[dr, bend left=20, "t"]
\arrow[dd, "\varphi_0"]
&
&& 
&
B
\arrow[r,"f"]
\arrow[dl, bend right=20, "s"']
&
A
\arrow[d, equals]
\arrow[dr, bend left=20, "t"]
&
\\
I
&&~&
J
&=& 
I
&
\Delta_{\varphi_0}D
\pullbackc{dr}{-0.05}
\arrow[r]
\arrow[d,"\varphi'_1"']
\arrow[u,"\varphi'_2", "\scriptsize\cong"']
&
A
\arrow[d,"\varphi_0"]
&
J
\\
&
D
\arrow[r,"g"']
\arrow[ul, bend left=20, "u"]
&
C
\arrow[ur, bend right=20, "v"']
&
&& 
&
D
\arrow[r,"g"']
\arrow[ul, bend left=20, "u"]
&
C
\arrow[ur, bend right=20, "v"']
&
\end{tikzcd}
\end{equation}

Note that, in general, for each diagram of the form \eqref{diaCartesianMorphism}, there are possibly many cartesian morphisms inducing it. Conversely, there are many potential ways of turning a diagram of the form \eqref{diaCartesianMorphism} into a cartesian morphism. Another possibility would be to take the induced cartesian morphism to be $(\varphi_0,\varphi_1,\mathrm{id}_B)$. Theorem \ref{thmPolyECartTrivial} below implies that these are essentially equivalent.

In particular, when $I=J=1$, we can regard pullback squares as cartesian morphisms in a canonical way.
\end{numbered}

We are now ready to assemble polynomials into a bicategory (and polynomial functors into a $2$-category). In fact, as proved in \cite{GambinoKock2013PolynomialFunctors}, more is true:

\begin{theorem}
\label{thmPolyEBicategory}
Let $\mathcal{E}$ be a locally cartesian closed category.
\begin{enumerate}[(a)]
\item There is a bicategory $\mathbf{Poly}_{\mathcal{E}}$ whose 0-cells are the objects of $\mathcal{E}$, whose 1-cells are polynomials in $\mathcal{E}$, and whose 2-cells are morphisms of polynomials.
\item There is a $2$-category $\mathbf{PolyFun}_{\mathcal{E}}$ whose 0-cells are the slices $\mathcal{E} \slice{I}$ of $\mathcal{E}$, whose 1-cells are polynomial functors, and whose 2-cells are strong natural transformations.
\item Extension defines a biequivalence $\mathrm{Ext} : \mathbf{Poly}_{\mathcal{E}} \xrightarrow{\simeq} \mathbf{PolyFun}_{\mathcal{E}}$.
\item Parts (a)--(c) hold true if we restrict the 1-cells to \textit{cartesian} morphisms of polynomials in $\mathbf{Poly}_{\mathcal{E}}$ and \textit{cartesian} strong natural transformations in $\mathbf{PolyFun}_{\mathcal{E}}$; thus there is a bicategory $\mathbf{Poly}\cart_{\mathcal{E}}$ and a $2$-category $\mathbf{PolyFun}\cart_{\mathcal{E}}$, which are biequivalent.
\end{enumerate}
\end{theorem}

We finish this section with the following technical lemma, which will simplify matters for us greatly down the road as it allows us in most instances to prove results about polynomials in the case when $I=J=1$.

\begin{theorem}
\label{lemPolynomialsFromOneToOne}
For fixed objects $I$ and $J$ of a locally cartesian closed category $\mathcal{E}$, there are full and faithful functors
$$S : \mathbf{Poly}_{\mathcal{E}}(I,J) \to \mathbf{Poly}_{\mathcal{E}\slice{I \times J}}(1,1) \quad \text{and} \quad S\cart : \mathbf{Poly}\cart_{\mathcal{E}}(I,J) \to \mathbf{Poly}\cart_{\mathcal{E} \slice{I \times J}}(1,1)$$
\end{theorem}

\begin{proof}[Proof sketch]
Given a polynomial $F : \declpoly IBAJsft$, define $S(F) = \langle s , f \rangle : B \to I \times A$ over $I \times J$ (considered as a polynomial $1 \pto 1$ in $\mathcal{E} \slice{I \times J}$) as in
\begin{center}
\begin{tikzcd}[row sep=huge, column sep=huge]
B
\arrow[rr, "{\langle s, f \rangle}"]
\arrow[dr, "{\langle s, t \circ f \rangle}"']
&&
I \times A
\arrow[dl, "{\mathrm{id}_I \times t}"]
\\
&
I \times J
&
\end{tikzcd}
\end{center}

Given a morphism of polynomials $\varphi : F \pRightarrow G$, as in
\begin{center}
\begin{tikzcd}[row sep=huge, column sep=huge]
&
B
\arrow[r,"f"]
\arrow[dl, bend right=20, "s"']
&
A
\arrow[d, equals]
\arrow[dr, bend left=20, "t"]
&
\\
I
&
D_{\varphi}
\pullback
\arrow[r]
\arrow[d,"\varphi_1"']
\arrow[u,"\varphi_2"]
&
A
\arrow[d,"\varphi_0"]
&
J
\\
&
D
\arrow[r,"g"']
\arrow[ul, bend left=20, "u"]
&
C
\arrow[ur, bend right=20, "v"']
&
\end{tikzcd}
\end{center}
define $S(\varphi) = (\mathrm{id}_I \times \varphi_0, \varphi_1, \varphi_2) : S(F) \pRightarrow S(G)$, as in the following diagram, where we consider $E$ as an object over $I \times J$ via $\langle s \circ \varphi_2, t \circ f \circ \varphi_2 \rangle : E \to I \times J$.
\begin{center}
\begin{tikzcd}[row sep=huge, column sep=huge]
B
\arrow[r,"{\langle s, f \rangle}"]
&
I \times A
\arrow[d, equals]
\\
D_{\varphi}
\pullback
\arrow[r]
\arrow[d,"\varphi_1"']
\arrow[u,"\varphi_2"]
&
I \times A
\arrow[d,"{\mathrm{id}_I \times \varphi_0}"]
\\
D
\arrow[r,"{\langle u, g \rangle}"']
&
I \times C
\end{tikzcd}
\end{center}

It is easy to see that $\mathrm{id}_I \times \varphi_0$, $\varphi_1$ and $\varphi_2$ are morphisms over $I \times J$ and that the lower square of the above diagram truly is cartesian, so that $S(\varphi)$ is a morphism in $\mathbf{Poly}_{\mathcal{E} \slice{I \times J}}(1,1)$. Verifying functoriality, fullness and faithfulness of $S$ is elementary but tedious.

That $S$ restricts to a full and faithful functor $S\cart : \mathbf{Poly}\cart_{\mathcal{E}}(I,J) \to \mathbf{Poly}\cart_{\mathcal{E}\slice{I \times J}}(1,1)$ is immediate, since $S(\varphi)$ is cartesian if and only if $\varphi_2$ is invertible, which holds if and only if $\varphi$ is cartesian.
\end{proof}

%% file: thesis/ch1-background/representability.tex
\section{Presheaves and representability}
\label{secPresheavesRepresentability}

This section lays out the basic definitions and results concerning presheaves and representability which we will use, normally without citation, in the rest of the thesis. Most of the results in this section can be found in the standard references for category theory, such as \cite{MacLane1971CWM}, \cite{Johnstone2002Elephant} and \cite{Awodey2010CategoryTheory}. They are recalled here because of their fundamental importance to the work to follow.

\begin{definition}[Presheaves]
A \textbf{presheaf}\index{presheaf} on a small category $\mathbb{C}$ is a functor $P : \mathbb{C}\op \to \mathbf{Set}$. The category of all presheaves on $\mathbb{C}$ and natural transformations between them is denoted by $\widehat{\mathbb{C}}$.
\end{definition}

Given an object $A \in \mathrm{ob}(\mathbb{C})$, an element $x \in P(A)$ and a morphism $f : B \to A$, we will write $x[f]$ rather than $P(f)(x)$ when $P$ is understood from context. Note that the rules for $P$ being a functor say precisely that $x[\mathrm{id}_A]=x$ and $x[f][g] = x[f \circ g]$ for all $C \xrightarrow{g} B \xrightarrow{f} A$ and all $x \in P(A)$, so that we might think of a presheaf $P$ as defining a right action of the morphisms of $\mathbb{C}$ on an $\mathrm{ob}(\mathbb{C})$-indexed family of sets.

\begin{definition}[Yoneda embedding]
The \textbf{Yoneda embedding}\index{Yoneda!embedding@{---} embedding} is the functor $\Yon : \mathbb{C} \to \widehat{\mathbb{C}}$ defined on objects by $\Yon(A) = \mathbb{C}({-}, A)$.
\end{definition}

\begin{definition}[Representable presheaves]
\label{defRepresentablePresheaf}
\index{representable!representable presheaf@{---} presheaf}
Let $\mathbb{C}$ be a small category. A presheaf $X : \mathbb{C}\op \to \mathbf{Set}$ is \textbf{representable} if $X \cong \Yon(A)$ for some $A \in \mathrm{ob}(\mathbb{C})$. The object $A$ is called a \textbf{representing object} for $X$.
\end{definition}

\begin{theorem}[Yoneda lemma]
Let $\mathbb{C}$ be a small category. For each presheaf $P$ over $\mathbb{C}$ and each object $A$ of $\mathbb{C}$, there is a bijection $\widehat{\mathbb{C}}(\Yon(A), P) \cong P(A)$. Moreover this bijection is natural in both $A$ and $P$.
\end{theorem}

In light of the Yoneda lemma, we will brazenly and unapologetically identify elements $x \in P(A)$ with natural transformations $x : \Yon(A) \to P$, and we may even use phrases such as `the element $x : \Yon(A) \to P$'. A consequence of the Yoneda lemma is that the Yoneda embedding is full and faithful.

\begin{definition}[Category of elements]
The \textbf{category of elements} of a presheaf $P$ is the category $\catel[\mathbb{C}]{P} = \catel{P}$, whose objects are pairs $(C,x)$ with $C \in \mathrm{ob}(\mathbb{C})$ and $x \in P(C)$, and whose morphisms $f : (C, x) \to (D, y)$ are morphisms $f : C \to D$ in $\mathbb{C}$ such that $y[f]=x$.
\end{definition}

Elementary computations reveal that $\catel[\mathbb{C}] \Yon(A) \cong \mathbb{C} \slice{A}$ for all $A \in \mathrm{ob}(\mathbb{C})$, and that $\widehat{\mathbb{C}} \slice{P} \cong \widehat{\catel[\mathbb{C}]{P}}$ for all $P : \mathbb{C}\op \to \mathbf{Set}$. Combining these results, we see that $\widehat{\mathbb{C} \slice{A}} \cong \widehat{\mathbb{C}} \slice{\Yon(A)}$ for each $A \in \mathrm{ob}(\mathbb{C})$.

Note that there is an evident forgetful functor $\pi : \catel[\mathbb{C}] P \to \mathbb{C}$.

\begin{theorem}[Every presheaf is the colimit of representables]
Let $\mathbb{C}$ be a small category and let $P$ be a presheaf over $\mathbb{C}$. Then $P$ is a colimit of the functor
$$\catel[\mathbb{C}]{P} \xrightarrow{\pi} \mathbb{C} \xrightarrow{\Yon} \widehat{\mathbb{C}}$$
In particular, it is a colimit in $\widehat{\mathbb{C}}$ of presheaves of the form $\Yon(\Gamma)$ for $\Gamma \in \mathrm{ob}(\mathbb{C})$.
\end{theorem}

\begin{theorem}
\label{thmKanExtensions}
Let $\mathbb{C}$ and $\mathbb{D}$ be small categories. Each functor $F : \mathbb{C} \to \mathbb{D}$ induces an adjoint triple
\begin{diagram}
\widehat{\mathbb{C}}
\arrow[r, shift left=1, bend left=20, "F_!", ""'{name=domadj1}]
\arrow[r, shift right=1, bend right=20, "F_*"', ""{name=codadj2}]
&[100pt]
\widehat{\mathbb{D}}
\arrow[l, "F^*" description, ""{name=domadj2}, ""'{name=codadj1}]
\arrow[from={domadj1}, to={codadj1}, draw=none, "\bot" description]
\arrow[from={domadj2}, to={codadj2}, draw=none, "\bot" description]
\end{diagram}
where $F^*$ is given by precomposition by $F$. The functors $F_!$ and $F_*$ are the \textbf{left Kan extension} and \textbf{right Kan extension} operations along $F$, respectively. Explicitly, the functors $F_*$ and $F_!$ can be computed as follows.
$$F_*(X)(D) \cong \widehat{\mathbb{D}}(\Yon(D), F_*(X)) \cong \widehat{\mathbb{C}}(F^*\Yon(D), X) \cong \widehat{\mathbb{C}}(\mathbb{D}(F({-}), \Yon(D)), X)$$
$$F_!(X) \cong \varinjlim\limits_{(A,x) \in \catel{X}} \Yon(FA)$$
In particular, we may choose the values of $F_!$ such that it commutes with the Yoneda embedding $F_! \circ \Yon = \Yon \circ F : \mathbb{C} \to \widehat{\mathbb{D}}$.
\end{theorem}

\begin{lemma}[Adjoint functors lift]
\label{lemAdjointFunctorKanExtension}
Given a pair of functors $R : \mathbb{C} \to \mathbb{D}$ and $L : \mathbb{D} \to \mathbb{C}$. If $L \dashv R$, then $L^* \dashv R^*$ and, therefore, $R_! \cong L^*$.

\begin{diagram}
\widehat{\mathbb{C}}
\arrow[r, bend left=10, shift left=1, "L^*", ""'{name=domadjtop}]
&
\widehat{\mathbb{D}}
\arrow[l, bend left=10, shift left=1, "R^*", ""'{name=codadjtop}]
\arrow[from={domadjtop}, to={codadjtop}, draw=none, "\bot" description]
\\
\mathbb{C}
\arrow[r, bend right=10, shift right=1, "R"', ""{name=codadjbase}]
\arrow[u, hook, dashed, "\Yon"]
&
\mathbb{D}
\arrow[l, bend right=10, shift right=1, "L"', ""{name=domadjbase}]
\arrow[from={domadjbase}, to={codadjbase}, draw=none, "\bot" description]
\arrow[u, hook, dashed, "\Yon"']
\end{diagram}
\end{lemma}

\begin{proof}[Sketch of proof]
Let $\eta : \mathrm{id}_{\mathbb{D}} \Rightarrow RL$ and $\varepsilon : LR \Rightarrow \mathrm{id}_{\mathbb{C}}$ be the unit and counit, respectively, of the adjunction $L \dashv R$. Define $\sigma : \mathrm{id}_{\widehat{\mathbb{C}}} \Rightarrow R^*L^* = (LR)^*$ by letting $(\sigma_P)_C = P(\varepsilon_C) : P(C) \to P(LR(C))$ for all $P : \mathbb{C}^{\mathrm{op}} \to \mathbf{Set}$ and all $C \in \mathrm{ob}(\mathbb{C})$, and define $\tau : L^*R^* = (RL)^* \Rightarrow \mathrm{id}_{\widehat{\mathbb{D}}}$ by letting $(\tau_Q)_D = Q(\eta_D) : Q(RL(D)) \to Q(D)$ for all $Q : \mathbb{D}^{\mathrm{op}} \to \mathbf{Set}$ and all $D \in \mathrm{ob}(\mathbb{D})$. Verifying that $\sigma, \tau$ are well-defined natural transformations forming the unit and counit, respectively, of the adjunction $L^* \dashv R^*$, is elementary. That $L^* \cong R_!$ follows from uniqueness of left adjoints up to natural isomorphism. 
\end{proof}

We conclude this section with the definition of a \textit{representable natural transformation}, which is the fundamental component of a \textit{natural model} \Cref{defNaturalModel}, the main object of study in this thesis.

\begin{definition}[Representable natural transformation]
\label{defRepresentableNaturalTransformation}
\index{representable!natural transformation@{---} natural transformation}
Let $\mathbb{C}$ be a small category and let $X$ and $Y$ be presheaves over $\mathbb{C}$. A natural transformation $f : Y \to X$ is \textbf{representable} if all of its fibres are representable, in the sense that for each $A \in \mathrm{ob}(\mathbb{C})$ and each $x \in X(A)$, there exists $B \in \mathrm{ob}(\mathbb{C})$, $g : B \to A$ in $\mathbb{C}$ and $y \in Y(B)$ such that the following square is a pullback.
\begin{diagram}
\Yon(B)
\arrow[d, "\Yon(g)"']
\arrow[r, "y"]
&
Y
\arrow[d, "f"]
\\
\Yon(A)
\arrow[r, "x"']
&
X
\end{diagram}
\end{definition}

The definition of a representable natural transformation can be found in \cite[Tag 0023]{StacksProject} and has been attributed to Alexander Grothendieck \cite{Awodey2016NaturalModels}.

%% file: thesis/ch2-categories-of-natural-models/_categories-of-natural-models.tex
\newpage
\input{thesis/ch2-categories-of-natural-models/natural-models.tex}

\newpage
\input{thesis/ch2-categories-of-natural-models/essentially-algebraic-theories.tex}

\newpage
\input{thesis/ch2-categories-of-natural-models/homomorphisms.tex}


%% file: thesis/ch2-categories-of-natural-models/natural-models.tex
\section{Natural models}
\label{secNaturalModels}

It was observed independently by Steve Awodey \cite{Awodey2012NotesOnModelsOfTypeTheory,Awodey2016NaturalModels} and Marcelo Fiore \cite{Fiore2012DiscreteGeneralisedPolynomialFunctors} that the notion of a representable natural transformation (\Cref{defRepresentableNaturalTransformation}) captures the type theoretic rules for context extension. We might, therefore, be tempted to take representable natural transformations as our notion of model of dependent type theory and leave it at that. However, we would like to provide an \textit{essentially algebraic} account of the semantics of type theory (\Cref{secEssentiallyAlgebraicTheory})---\Cref{defRepresentableNaturalTransformation} does not quite do this because it posits mere existence, rather than a choice, of data witnessing representability of the natural transformation.

\begin{definition}[Natural models \thmcite{Awodey2016NaturalModels}]
\label{defNaturalModel}
\index{natural model}
A \textbf{natural model} is a small category $\mathbb{C}$ with a distinguished terminal object $\diamond$, presheaves $\nmty{U}$ and $\nmtm{U}$ over $\mathbb{C}$, a natural transformation $p : \nmtm{U} \to \nmty{U}$ and, for each $\Gamma \in \mathrm{ob}(\mathcal{C})$ and $A \in \nmty{U}(\Gamma)$, the following \textbf{representability data}: an object $\Gamma \cext A$ of $\mathcal{C}$, a morphism $\nmp{A} = \nmp[\Gamma]{A} : \Gamma \cext A \to \Gamma$ in $\mathcal{C}$, and an element $\nmq{A} = \nmq[\Gamma]{A} \in \nmtm{U}(\Gamma \cext A)$, such that the following square is a pullback for all such $\Gamma$ and $A$:
\begin{diagram}
\Yon(\Gamma \cext A)
\arrow[r, "\nmq{A}"]
\arrow[d, "\Yon(\nmp{A})"']
&
\nmtm{U}
\arrow[d, "p"]
\\
\Yon(\Gamma)
\arrow[r, "A"']
&
\nmty{U}
\end{diagram}
\end{definition}

\begin{numbered}
As explained in detail in \cite{Awodey2016NaturalModels}, we can informally view a natural model as a model of dependent type theory in the following way. The category $\mathbb{C}$ represents the category of contexts and substitutions, with the terminal object $\diamond$ representing the empty context. For each object $\Gamma$, the set $\nmty{U}(\Gamma)$ represents the set of types in context $\Gamma$, and the set $\nmtm{U}(\Gamma)$ represents the set of terms in context $\Gamma$, with the function $p_{\Gamma} : \nmtm{U}(\Gamma) \to \nmty{U}(\Gamma)$ sending a term to its unique type. (This is where we used the uniqueness of typing, as discussed in \Cref{secTypesForTheWorkingMathematician}.) The action of $\nmty{U}$ and $\nmtm{U}$ on morphisms is that of substitution. Naturality of $p$ says that substitution represents typing, in the sense that if $\Gamma \vdash a : A$ and $\sigma : \Delta \to \Gamma$ is a substitution, then $\Delta \vdash a[\sigma] : A[\sigma]$. The object $\Gamma \cext A$ represents the extension of a context $\Gamma$ by a new variable $x : A$; then $\nmp{A} : \Gamma \cext A \to A$ represents the \textit{weakening} substitution and $\nmq{A} \in \nmtm{U}(\Gamma \cext A)$ represents the new variable $x$. Finally, the universal property of the pullback says that $\Gamma \cext A$ truly does satisfy the syntactic rules for context extension.
\end{numbered}

\begin{numbered}
Under the axiom of choice, every representable natural transformation gives rise to a natural model. Moreover, since representability is defined by a pullback condition, given any $\Gamma \in \mathrm{ob}(\mathbb{C})$ and $A \in \nmty{U}(\Gamma)$, the representability data $(\Gamma \cext A, \nmp{A}, \nmq{A})$ are unique up to canonical isomorphism, in the sense that for any other choice $(\widetilde{\Gamma \cext A}, \widetilde{\nmp{A}}, \widetilde{\nmq{A}})$ of representability data, there is a unique isomorphism $\theta : \Gamma \cext A \to \widetilde{\Gamma \cext A}$ such that $\widetilde{\nmp{A}} \circ \theta = \nmp{A}$ and $\widetilde{\nmq{A}}[\theta] = \nmq{A}$.
\end{numbered}

\begin{numbered}
We will adopt the following notation conventions:
\begin{itemize}
\item $\mathsf{t}_{\Gamma}$ is the unique morphism $\Gamma \to \diamond$ in $\mathbb{C}$;
\item In the internal language of $\widehat{\mathbb{C}}$, write $[A] = \nmtm{U}_A$ for the fibre of $p$ over $A \in \nmty{U}$;
\item Given an object $\Gamma$ of $\mathbb{C}$ and an element $A \in \nmty{U}(\Gamma)$, write $\nmtm{U}(\Gamma; A)$ for the preimage of $p_{\Gamma} : \nmtm{U}(\Gamma) \to \nmty{U}(\Gamma)$ over $A$---we may, further, say `$\Gamma \vdash A$ in $(\mathbb{C}, p)$' to mean that $A \in \nmty{U}(\Gamma)$, and `$\Gamma \vdash a : A$ in $(\mathbb{C}, p)$' to mean that $a \in \nmtm{U}(\Gamma; A)$;
\item Given a morphism $\sigma : \Delta \to \Gamma$ in $\mathbb{C}$, an element $A \in \nmty{U}(\Gamma)$ and an element $a \in \nmtm{U}(\Delta, A[\sigma])$, write $\langle \sigma, a \rangle_A$ for the unique morphism $\Delta \to \Gamma \cext A$ induced by the universal property of pullbacks;
\item For each $a \in \nmtm{U}(\Gamma; A)$, write $\nms{a} = \langle \mathrm{id}_{\Gamma}, a \rangle : \Gamma \to \Gamma \cext A$---note that $\nmp{A} \circ \nms{a} = \mathrm{id}_{\Gamma}$, so that $\nms{a}$ is a section of $\nmp{A}$.
\end{itemize}
\end{numbered}

\begin{numbered}
In order to avoid writing a long list of symbols each time we refer to a natural model, we will typically write just $(\mathbb{C}, p)$, leaving the naming of the remaining data implicit. We will adopt the convention that when we write $(\mathbb{C}, p)$, the additional data is named as in \Cref{defNaturalModel} and, when we write $(\mathbb{D}, q)$, the chosen terminal object of $\mathbb{D}$ is denoted by $\star$, that $q : \nmtmalt{V} \to \nmty{V}$ in $\widehat{\mathbb{D}}$, and that the representability data for a given $\Gamma \in \mathrm{ob}(\mathbb{D})$ and $A \in \nmty{V}(\Gamma)$ is denoted by $\Gamma \cextalt A$, $\nmu{A} : \Gamma \cextalt A \to \Gamma$ and $\nmv{A} \in \nmtmalt{V}(\Gamma)$. Furthermore, internally to $\widehat{\mathbb{D}}$, write $\langle A \rangle = \nmtm{V}_A$ rather than $[A]$ for the fibre of $q$ over $A \in \nmty{V}$.
\end{numbered}

\begin{construction}[Canonical pullback squares]
\label{cnsCanonicalPullbacks}
Let $(\mathbb{C}, p)$ be a natural model. For all $\sigma : \Delta \to \Gamma$ in $\mathbb{C}$ and all $A \in \nmty{U}(\Gamma)$, there is a pullback square
\begin{diagram}
\Delta \cext A[\sigma]
\arrow[r, "\sigma \cext A"]
\arrow[d, "\nmp{A[\sigma]}"']
\pullback
&
\Gamma \cext A
\arrow[d, "\nmp{A}"]
\\
\Delta
\arrow[r, "\sigma"']
&
\Gamma
\end{diagram}
Pullback squares of this form are called \textbf{canonical pullback squares}.
\end{construction}

\begin{verification}
Construct the following diagram using representability data for $(\mathbb{C}, p)$.
\begin{diagram}
\Yon(\Delta \cext A[\sigma])
\arrow[r, dashed]
\arrow[d, "\Yon(\nmp{A[\sigma]})"']
\arrow[rr, bend left=15, "\nmq{A[\sigma]}"]
&
\Yon(\Gamma \cext A)
\arrow[d, "\Yon(\nmp{A})" description]
\arrow[r, "\nmq{A}" description]
\pullback
&
\nmtm{U}
\arrow[d, "p"]
\\
\Yon(\Delta)
\arrow[r, "\Yon(\sigma)" description]
\arrow[rr, bend right=15, "{A[\sigma]}"']
&
\Yon(\Gamma)
\arrow[r, "A" description]
&
\nmty{U}
\end{diagram}
The right and outer squares are pullbacks by representability of $p$. The universal property of the right-hand pullback yields a morphism $\Yon(\Delta \cext A[\sigma]) \to \Yon(\Gamma \cext A)$ as indicated, which is of the form $\Yon(\sigma \cext A)$ for some $\sigma \cext A : \Delta \cext A[\sigma] \to \Gamma \cext A$ in $\mathbb{C}$ since the Yoneda embedding is full and faithful. The left-hand square is a pullback by the two pullbacks lemma, and hence the square in the statement of this construction is a pullback since the Yoneda embedding reflects limits.
\end{verification}

\begin{lemma}
\label{lemPastingCanonicalPullbacks}
Let $(\mathbb{C}, p)$ be a natural model, let $\Theta \xrightarrow{\tau} \Delta \xrightarrow{\sigma} \Gamma$ in $\mathbb{C}$ and let $A \in \nmty{U}(\Gamma)$. With notation as in \Cref{cnsCanonicalPullbacks}, we have
$$(\sigma \circ \tau) \cext A = (\sigma \cext A) \circ (\tau \cext A[\sigma]) : \Theta \cext A[\sigma \circ \tau] \to \Gamma \cext A$$
\end{lemma}

\begin{proof}
This is an immediate consequence of the two pullbacks lemma.
\end{proof}

\Cref{lemPastingCanonicalPullbacks} demonstrates that, in a way that mirrors that of Vladimir Voevodsky's `universes' \cite{Voevodsky2009TypeSystems}, natural models overcome the \textit{coherence problem} for interpreting type theory in a locally cartesian closed category.

\subsection*{Admitting type theoretic structure}

The proofs of \Cref{thmAdmittingUnitType,thmAdmittingSigmaTypes,thmAdmittingPiTypes} can be found in \cite{Awodey2016NaturalModels}.

\begin{theorem}[Admitting a unit type \thmcite{Awodey2016NaturalModels}]
\label{thmAdmittingUnitType}
A natural model $(\mathbb{C}, p)$ admits a unit type if and only if there are morphisms
$$\widehat{\mathbbm{1}} : 1 \to \nmty{U} \qquad \text{and} \qquad \widehat{\star} : 1 \to \nmtm{U}$$
in $\widehat{\mathbb{C}}$ exhibiting $\mathrm{id}_1 : 1 \to 1$ as a pullback of $p$.
\begin{diagram}
1
\arrow[d, equals]
\arrow[r, "\widehat{\star}"]
\pullback
&
\nmtm{U}
\arrow[d, "p"]
\\
1
\arrow[r, "\widehat{\mathbbm{1}}"']
&
\nmty{U}
\end{diagram}
\end{theorem}

\begin{theorem}[Admitting dependent sum types \thmcite{Awodey2016NaturalModels}]
\label{thmAdmittingSigmaTypes}
A natural model $(\mathbb{C}, p)$ admits dependent sum types if and only if there are morphisms
$$\widehat{\upSigma} : \sum_{A : \nmty{U}} \nmty{U}^{[A]} \to \nmty{U} \qquad \text{and} \qquad \widehat{\mathsf{pair}} : \sum_{A : \nmty{U}} \sum_{B : \nmty{U}^{[A]}} \sum_{a:[A]} [B(a)] \to \nmtm{U}$$
in $\widehat{\mathbb{C}}$ exhibiting the projection $\pi : \sum_{A : \nmty{U}} \sum_{B : \nmty{U}^{[A]}} \sum_{a:[A]} [B(a)] \to \sum_{A : \nmty{U}} \nmty{U}^{[A]}$ as a pullback of $p$.
\begin{diagram}
\sum_{A : \nmty{U}} \sum_{B : \nmty{U}^{[A]}} \sum_{a:[A]} [B(a)]
\arrow[d, "\pi"']
\arrow[r, "\widehat{\mathsf{pair}}"]
\pullbackc{dr}{-0.05}
&
\nmtm{U}
\arrow[d, "p"]
\\
\sum_{A : \nmty{U}} \nmty{U}^{[A]}
\arrow[r, "\widehat{\upSigma}"']
&
\nmty{U}
\end{diagram}
Moreover the map $\pi$ is precisely the polynomial composite $p \cdot p$.
\end{theorem}

\begin{theorem}[Admitting dependent product types \thmcite{Awodey2016NaturalModels}]
\label{thmAdmittingPiTypes}
A natural model $(\mathbb{C}, p)$ admits dependent product types if and only if there are morphisms
$$\widehat{\upPi} : \sum_{A : \nmty{U}} \nmty{U}^{[A]} \to \nmty{U} \qquad \text{and} \qquad \widehat{\lambda} : \sum_{A : \nmty{U}} \nmtm{U}^{[A]} \to \nmtm{U}$$
in $\widehat{\mathbb{C}}$ exhibiting $\sum_{A : \nmty{U}} p^{[A]}$ as a pullback of $p$.
\begin{diagram}
\sum_{A : \nmty{U}} \nmtm{U}^{[A]}
\arrow[r, "\widehat{\lambda}"]
\arrow[d, "\sum_{A : \nmty{U}} p^{[A]}"']
\pullbackc{dr}{0}
&
\nmtm{U}
\arrow[d, "p"]
\\
\sum_{A : \nmty{U}} \nmty{U}^{[A]}
\arrow[r, "\widehat{\upPi}"']
&
\nmty{U}
\end{diagram}
Moreover the map $\sum_{A \in \nmty{U}} \nmty{U}^{[A]}$ is precisely the morphism $\upP_p(p)$ obtained by applying the extension $\upP_p$ of $p$ to $p$ itself.
\end{theorem}

Recall (\Cref{rmkCartesianMorphismIsPullbackSquare}) that in a locally cartesian closed category (such as $\widehat{\mathbb{C}}$), pullback squares correspond with cartesian morphisms of polynomials. Therefore we can succinctly rephrase \Cref{thmAdmittingUnitType,thmAdmittingSigmaTypes,thmAdmittingPiTypes} in terms of cartesian morphisms of polynomials.

\begin{corollary}
\label{thmUnitSigmaPiPoly}
Let $(\mathbb{C},p)$ be a natural model.
\begin{enumerate}[(a)]
\item $(\mathbb{C},p)$ admits a unit type if and only if there is a cartesian morphism $\eta : i_1 \pRightarrow p$ in $\mathbf{Poly}_{\widehat{\mathbb{C}}}$;
\item $(\mathbb{C},p)$ admits dependent sum types if and only if there is a cartesian morphism $\mu : p \cdot p \pRightarrow p$ in $\mathbf{Poly}_{\widehat{\mathbb{C}}}$;
\item $(\mathbb{C},p)$ admits dependent product types if and only if there is a cartesian morphism $\zeta : P_p(p) \pRightarrow p$ in $\mathbf{Poly}_{\widehat{\mathbb{C}}}$;
\end{enumerate}
\end{corollary}

Our notation is deliberately suggestive of a \textit{monad} and an \textit{algebra}; exploring this topic further is the subject of \Cref{secPolynomialPseudomonads}.

%% file: thesis/ch2-categories-of-natural-models/essentially-algebraic-theories.tex
\section{The essentially algebraic theory of natural models}
\label{secEssentiallyAlgebraicTheory}
\index{essentially algebraic theory}

A (single-sorted) \textit{algebraic theory} is one which is specified by \textit{operation symbols} $\sigma$ with \textit{arities}, which are natural numbers, subject to conditions which can be expressed as (universally quantified) \textit{equations}. For example, the theory of groups has an operation $\mathsf{unit}$ of arity $0$, an operation $\mathsf{inv}$ of arity $1$, and a binary operation $\mathsf{mult}$ of arity $2$, subject to equations describing associativity, inverse and unit laws. A model $\mathfrak{M}$ of an algebraic theory is then a set $M$ together with functions $\sigma^{\mathfrak{M}} : M^{\alpha(\sigma)} \to M$ for each symbol $\sigma$ whose arity is $\alpha(\sigma)$, such that the functions $\sigma^{\mathfrak{M}}$ satisfy the specified equations; thus a model of the theory of groups is a group.

More generally, given a set $S$ of \textit{sorts}, an $S$-sorted algebraic theory is again specified by operation symbols subject to equations, but now the arities of the operation symbols are sequences of elements of $S$. A symbol $\sigma$ whose arity is $(s_1,\dots,s_n, s)$ (which we suggestively write as $s_1 \times \dots \times s_n \to s$) can be thought of as an $n$-ary operation, whose $i^{\text{th}}$ input has sort $s_i$ and whose value has sort $s$. A model $\mathfrak{M}$ of an $S$-sorted algebraic theory is then given by an $S$-indexed family of sets $\seqbn{M_s}{s \in S}$ with functions $\sigma^{\mathfrak{M}} : M_{s_1} \times \dots \times M_{s_n} \to M_s$ for each operation symbol $\sigma$ of arity $s_1 \times \dots \times s_n \to s$, which satisfy the equations of the theory.

The notion of an ($S$-sorted) \textit{essentially algebraic theory} generalises that of an ($S$-sorted) algebraic theory even further by allowing operations to be \textit{partial}, meaning that operation symbols may be defined only on inputs which satisfy certain equational conditions stated in terms of already-specified operation symbols. An example of such a theory is that of categories, whose sorts are $\mathsf{Obj}$, the sort of objects, and $\mathsf{Mor}$, the sort of morphisms. The composition operation $\mathsf{comp}$ with arity $\mathsf{Mor} \times \mathsf{Mor} \to \mathsf{Mor}$ is partial, since the composite $g \circ f$ of two morphisms on a category is defined only when $\mathsf{cod}(f) = \mathsf{dom}(g)$. The interpretation $\mathsf{comp}^{\mathfrak{M}}$ is then a partial function from $M_{\mathsf{Mor}} \times M_{\mathsf{Mor}}$ to $M_{\mathsf{Mor}}$, whose domain of definition is given by $\{ (g,f) \in M_{\mathsf{Mor}} \mid \mathsf{cod}^{\mathfrak{M}}(f) = \mathsf{dom}^{\mathfrak{M}}(g) \}$.

A precise definition of an essentially algebraic theory, and a model thereof, can be found in \cite{AdamekRosicky1994}.

Our goal in this section is to exhibit the theory of natural models as an essentially algebraic theory $\mathbb{T}_{\mathrm{NM}}$. Once we have done so, natural models will automatically assemble into a category $\mathbf{NM}$, whose objects are the models of $\mathbb{T}_{\mathrm{NM}}$ and whose morphisms are \textit{homomorphisms} of models of $\mathbb{T}_{\mathrm{NM}}$---that is, families of functions between the sorts which commute in the appropriate sense with the operation symbols. The general theory of essentially algebraic categories will then apply to the category $\mathbf{NM}$.

The practical power of natural models comes from their functorial, rather than algebraic, description; as such, the main takeaway of this section is that the functorial description captures the algebraic one, and we will provide a functorial account of the algebraic notion of a homomorphism of models of $\mathbb{T}_{\mathrm{NM}}$ in \Cref{secHomomorphisms}.

\begin{definition}
\label{defEATheoryOfNaturalModels}
\index{natural model!essentially algebraic theory of {---}s}
\index{essentially algebraic theory!of natural models@{---} of natural models}
The \textbf{theory of natural models} is the essentially algebraic theory $\mathbb{T} = \mathbb{T}_{\mathrm{NM}}$ described as follows.

The set $S$ of sorts is $\{ \mathsf{Ctx}, \mathsf{Sub}, \mathsf{Type}, \mathsf{Term} \}$;

\textbf{Note:} In all of what follows, the sorts of the variables are
$$\Delta,\Gamma : \mathsf{Ctx} \qquad \sigma,\tau,\upsilon : \mathsf{Sub} \qquad A,B : \mathsf{Type} \qquad a,b,f,p : \mathsf{Term}$$

The set $\Sigma$ of symbols and their arities is defined in the following table.

\begin{center}
\vspace{-20pt}
\begin{tabular}{l|lclcl|lcl}
\textbf{Name} & \multicolumn{5}{l|}{\textbf{Symbol} (with arity \& sorts)} & \multicolumn{3}{l}{\textbf{Shorthand}} \\ \hline
domain$^{\dagger}$ & $\mathsf{dom}$ & $:$ & $\mathsf{Sub}$ & $\to$ & $\mathsf{Ctx}$ & && \\
codomain$^{\dagger}$ & $\mathsf{cod}$ & $:$ & $\mathsf{Sub}$ & $\to$ & $\mathsf{Ctx}$ & && \\
identity$^{\dagger}$ & $\mathsf{id}$ & $:$ & $\mathsf{Ctx}$ & $\to$ & $\mathsf{Sub}$ & $\mathsf{id}(\Gamma)$ & $=$ & $\mathsf{id}_{\Gamma}$ \\
composition & $\mathsf{comp}$ & $:$ & $\mathsf{Sub} \times \mathsf{Sub}$ & $\to$ & $\mathsf{Sub}$ & $\mathsf{comp}(\sigma,\tau)$ & $=$ & $\tau \circ \sigma$ \\
empty context$^{\dagger}$ & $\mathsf{empty}$ & $:$ & & & $\mathsf{Ctx}$ & $\mathsf{empty}$ & $=$ & $\diamond$ \\
sub$^{\text{n}}$ to empty$^{\dagger}$ & $\mathsf{esub}$ & $:$ & $\mathsf{Ctx}$ & $\to$ & $\mathsf{Sub}$ & $\mathsf{esub}(\Gamma)$ & $=$ & $\mathsf{t}_{\Gamma}$ \\
typing$^{\dagger}$ & $\mathsf{typeof}$ & $:$ & $\mathsf{Term}$ & $\to$ & $\mathsf{Type}$ & && \\
context of types$^{\dagger}$ & $\mathsf{ctxof}_{\mathsf{ty}}$ & $:$ & $\mathsf{Type}$ & $\to$ & $\mathsf{Ctx}$ & $\mathsf{ctxof}_{\mathsf{ty}}(A)$ & $=$ & $\mathsf{ctxof}(A)$ \\
context of terms$^{\dagger}$ & $\mathsf{ctxof}_{\mathsf{tm}}$ & $:$ & $\mathsf{Term}$ & $\to$ & $\mathsf{Ctx}$ & $\mathsf{ctxof}_{\mathsf{tm}}(a)$ & $=$ & $\mathsf{ctxof}(a)$ \\
sub$^{\text{n}}$ on types & $\mathsf{subst}_{\mathsf{ty}}$ & $:$ & $\mathsf{Sub} \times \mathsf{Type}$ & $\to$ & $\mathsf{Type}$ & $\mathsf{subst}_{\mathsf{ty}}(\sigma,A)$ & $=$ & $A[\sigma]$ \\
sub$^{\text{n}}$ on terms & $\mathsf{subst}_{\mathsf{tm}}$ & $:$ & $\mathsf{Sub} \times \mathsf{Term}$ & $\to$ & $\mathsf{Term}$ & $\mathsf{subst}_{\mathsf{tm}}(\sigma,a)$ & $=$ & $a[\sigma]$ \\
context extension$^{\dagger}$ & $\mathsf{cext}$ & $:$ & $\mathsf{Type}$ & $\to$ & $\mathsf{Ctx}$ & $\mathsf{cext}(A)$ & $=$ & $\Gamma \cext A$ $^{\ddagger}$ \\
projection$^{\dagger}$ & $\mathsf{proj}$ & $:$ & $\mathsf{Type}$ & $\to$ & $\mathsf{Sub}$ & $\mathsf{proj}(A)$ & $=$ & $\nmp{A}$ \\
variable$^{\dagger}$ & $\mathsf{var}$ & $:$ & $\mathsf{Type}$ & $\to$ & $\mathsf{Term}$ & $\mathsf{var}(A)$ & $=$ & $\nmq{A}$ \\
induced sub$^{\text{n}}$ & $\mathsf{indsub}$ & $:$ & $\mathsf{Sub} \times \mathsf{Term} \times \mathsf{Type}$ & $\to$ & $\mathsf{Sub}$ & $\mathsf{indsub}(\sigma, a, A)$ & $=$ & $\langle \sigma, a \rangle_A$
\end{tabular}

{\footnotesize $\dagger$ denotes total symbols \\ $\ddagger$ When we write $\Gamma \cext A = \mathsf{cext}(A)$, we are implying that $\Gamma = \mathsf{ctxof}(A)$.}
\end{center}

\newpage
The set $E$ of equations is defined as follows:

\newcounter{axiomcount}
\begin{minipage}[t]{0.48\textwidth}
\begin{itemize}
\item \textbf{Category of contexts}
\begin{enumerate}[(i)]
\item $\mathsf{dom}(\mathsf{id}_{\Gamma}) = \Gamma$
\item $\mathsf{cod}(\mathsf{id}_{\Gamma}) = \Gamma$
\item $\mathsf{dom}(\tau \circ \sigma) = \mathsf{dom}(\sigma)$
\item $\mathsf{cod}(\tau \circ \sigma) = \mathsf{cod}(\tau)$
\item $\sigma \circ \mathsf{id}_{\Gamma} = \sigma$
\item $\mathrm{id}_{\Delta} \circ \tau = \tau$
\item $(\upsilon \circ \tau) \circ \sigma = \upsilon \circ (\tau \circ \sigma)$
\setcounter{axiomcount}{\value{enumi}}
\end{enumerate}

\item \textbf{Empty context is terminal}
\begin{enumerate}[(i)]
\setcounter{enumi}{\value{axiomcount}}
\item $\mathsf{dom}(\mathsf{t}_{\Gamma}) = \Gamma$
\item $\mathsf{cod}(\mathsf{t}_{\Gamma}) = \diamond$
\item $\mathsf{t}_{\Gamma} \circ f = \mathsf{t}_{\Delta}$
\setcounter{axiomcount}{\value{enumi}}
\end{enumerate}

\item \textbf{Presheaf of types}
\begin{enumerate}[(i)]
\setcounter{enumi}{\value{axiomcount}}
\item $A[\mathsf{id}_{\mathsf{ctxof}(A)}] = A$
\item $A[\tau \circ \sigma] = A[\tau][\sigma]$
\item $\mathsf{ctxof}(A[\sigma]) = \mathsf{dom}(\sigma)$
\setcounter{axiomcount}{\value{enumi}}
\end{enumerate}
\end{itemize}
\end{minipage}
\begin{minipage}[t]{0.48\textwidth}
\begin{itemize}
\item \textbf{Presheaf of terms}
\begin{enumerate}[(i)]
\setcounter{enumi}{\value{axiomcount}}
\item $a[\mathsf{id}_{\mathsf{ctxof}(a)}] = a$
\item $a[\tau \circ \sigma] = a[\tau][\sigma]$
\item $\mathsf{ctxof}(a[\sigma]) = \mathsf{dom}(\sigma)$
\setcounter{axiomcount}{\value{enumi}}
\end{enumerate}

\item \textbf{Typing is natural}
\begin{enumerate}[(i)]
\setcounter{enumi}{\value{axiomcount}}
\item $\mathsf{ctxof}(\mathsf{typeof}(a)) = \mathsf{ctxof}(a)$
\item $\mathsf{typeof}(a[\sigma]) = \mathsf{typeof}(a)[\sigma]$
\setcounter{axiomcount}{\value{enumi}}
\end{enumerate}

\item \textbf{Representability}
\begin{enumerate}[(i)]
\setcounter{enumi}{\value{axiomcount}}
\item $\mathsf{dom}(\nmp{A}) = \Gamma \cext A$
\item $\mathsf{cod}(\nmp{A}) = \Gamma$
\item $\mathsf{ctxof}(\nmq{A}) = \Gamma \cext A$
\item $\mathsf{typeof}(\nmq{A}) = A[\nmp{A}]$
\item $\mathsf{dom}(\langle \sigma, a \rangle_A) = \mathsf{dom}(\sigma)$
\item $\mathsf{cod}(\langle \sigma, a \rangle_A) = \mathsf{cod}(\sigma) \cext A$
\item $\nmp{A} \circ \langle \sigma, a \rangle_A = \sigma$
\item $\nmq{A}[\langle \sigma, a \rangle_A] = a$
\item $\langle \nmp{A} \circ \sigma, \nmq{A}[\sigma] \rangle_A = \sigma$
\end{enumerate}
\end{itemize}
\end{minipage}

The domains of definition of the partial symbols are given as follows:
\begin{center}
\begin{tabular}{lcl}
$\mathsf{Def}(\mathsf{comp}(f,g))$ & $=$ & $\{ \mathsf{cod}(f) = \mathsf{dom}(g) \}$ \\
$\mathsf{Def}(\mathsf{subst}_{\mathsf{ty}}(\sigma,A))$ & $=$ & $\{ \mathsf{ctxof}(A) = \mathsf{cod}(\sigma) \}$ \\
$\mathsf{Def}(\mathsf{subst}_{\mathsf{tm}}(\sigma,a))$ & $=$ & $\{ \mathsf{ctxof}(a) = \mathsf{cod}(\sigma) \}$ \\
$\mathsf{Def}(\mathsf{indsub}(\sigma,a,A))$ & $=$ & $\{ \mathsf{cod}(\sigma) = \mathsf{ctxof}(A),\ \mathsf{typeof}(a) = A[\sigma] \}$
\end{tabular}
\end{center}
\end{definition}

\begin{theorem}[Natural models are captured by $\mathbb{T}_{\mathrm{NM}}$]
\label{thmTheoryCapturesDefinition}
Specifying a model of the theory $\mathbb{T}_{\mathrm{NM}}$ (\Cref{defEATheoryOfNaturalModels}) is equivalent to specifying a natural model (\Cref{defNaturalModel}).
\end{theorem}

\begin{proof}
A model $\mathfrak{M}$ of $\mathbb{T}$ consists of four sets $M_{\mathsf{Ctx}}$, $M_{\mathsf{Sub}}$, $M_{\mathsf{Type}}$ and $M_{\mathsf{Term}}$, together with fifteen functions $\mathsf{dom}^{\mathfrak{M}}, \mathsf{cod}^{\mathfrak{M}}, \dots, \mathsf{ind}^{\mathfrak{M}}$ whose domains and codomains are determined by the sorts, arities and domains of definition described in \Cref{defEATheoryOfNaturalModels}.

Equations (i)--(vii) say precisely that the data $(M_{\mathsf{Ctx}}, M_{\mathsf{Sub}}, \mathsf{dom}^{\mathfrak{M}}, \mathsf{cod}^{\mathfrak{M}}, \mathsf{id}^{\mathfrak{M}}, \mathsf{comp}^{\mathfrak{M}})$ defines a (necessarily small) category $\mathbb{C}$. Equations (viii)--(x) say precisely that $\diamond^{\mathfrak{M}}$ is a terminal object of $\mathbb{C}$, with $\mathsf{esub}^{\mathfrak{M}}(\Gamma) : \Gamma \to \diamond^{\mathfrak{M}}$ being the unique morphism from an object $\Gamma$ of $\mathbb{C}$ to $\diamond^{\mathfrak{M}}$. We need not (and do not) specify $\mathsf{esub}^{\mathfrak{M}}$ when defining a natural model, since its unique existence follows from the assertion that $\diamond^{\mathfrak{M}}$ is terminal.

Equations (xi)--(xiii) say that $M_{\mathsf{Type}}$ and the functions $\mathsf{ctxof}_{\mathsf{ty}}^{\mathfrak{M}}$ and $\mathsf{subst}_{\mathsf{ty}}^{\mathfrak{M}}$ together define a presheaf $\nmty{U} : \mathbb{C}\op \to \mathbf{Set}$ in the following way: the set $\nmty{U}(\Gamma)$ is given by $(\mathsf{ctxof}_{\mathsf{ty}}^{\mathfrak{M}})^{-1}(\Gamma) \subseteq M_{\mathsf{Type}}$, and the action of $\nmty{U}$ on morphisms is defined by $\nmty{U}(\sigma) : A \mapsto \mathsf{subst}_{\mathsf{ty}}^{\mathfrak{M}}(\sigma, A)$, with the equations telling us that this action is contravariantly functorial with the correct domain and codomain. Conversely, given a presheaf $\mathcal{U} : \mathbb{C}\op \to \mathbf{Set}$, we can take $M_{\mathsf{Type}}$ to be the disjoint union of the sets $\mathcal{U}(\Gamma)$, with $\mathsf{ctxof}^{\mathfrak{M}}_{\mathsf{ty}}$ given by the projection map to $M_{\mathsf{Ctx}} = \mathrm{ob}(\mathbb{C})$ and the function $\mathsf{subst}^{\mathfrak{M}}_{\mathsf{ty}}$ given by $(\sigma, A) \mapsto \nmty{U}(\sigma)(A)$. Likewise, equations (xiv)--(xvi) say that specifying $M_{\mathsf{Term}}$ and the functions $\mathsf{ctxof}_{\mathsf{tm}}^{\mathfrak{M}}$ and $\mathsf{subst}_{\mathsf{tm}}^{\mathfrak{M}}$ is equivalent to defining a presheaf $\nmtm{U} : \mathbb{C}\op \to \mathbf{Set}$.

Equations (xvii)--(xviii) say that $\mathsf{typeof}^{\mathfrak{M}}$ defines a natural transformation $p : \nmtm{U} \to \nmty{U}$. Indeed, equation (xvii) says that the restriction of $\mathsf{typeof}^{\mathfrak{M}}$ to $\nmtm{U}(\Gamma) = (\mathsf{ctxof}_{\mathsf{tm}}^{\mathfrak{M}})^{-1}(\Gamma)$ has image contained in $\nmty{U}(\Gamma) = (\mathsf{ctxof}_{\mathsf{ty}}^{\mathfrak{M}})^{-1}(\Gamma)$, so that we obtain a function $p_{\Gamma} : \nmtm{U}(\Gamma) \to \nmty{U}(\Gamma)$; and equation (xvii) says that the naturality squares for $p$ commute for each morphism $\sigma : \Delta \to \Gamma$ in $\mathbb{C}$.

Equations (xix)--(xxvii) say precisely that for each $\Gamma \in \mathrm{ob}(\mathbb{C})$ and $A \in \nmty{U}(\Gamma)$, the data $\mathsf{cext}^{\mathfrak{M}}(\Gamma, A)$, $\mathsf{proj}^{\mathfrak{M}}$ and $\mathsf{var}^{\mathfrak{M}}$ exhibit $p$ as a representable natural transformation. Indeed, equations (xix)--(xxi) say that these data have the required types; equation (xxii) says that the following square commutes;
\begin{diagram}
\Yon(\mathsf{cext}^{\mathfrak{M}}(\Gamma, A))
\arrow[r, "\mathsf{var}^{\mathfrak{M}}(A)"]
\arrow[d, "\Yon(\mathsf{proj}^{\mathfrak{M}}(A))"']
&
\nmtm{U}
\arrow[d, "p"]
\\
\Yon(\Gamma)
\arrow[r, "A"']
&
\nmty{U}
\end{diagram}
and equations (xxiii)--(xvii) exhibit the square as a pullback, with $\mathsf{indsub}^{\mathfrak{M}}$ giving the morphisms induced from the universal property of pullbacks. Note that we need not (and do not) specify $\mathsf{indsub}^{\mathfrak{M}}$ when defining a natural model, since its unique existence follows from the universal property.
\end{proof}

\begin{numbered}
Categories with families (\Cref{defCwF}) can also be described as models of $\mathbb{T}_{\mathrm{NM}}$---that is, natural models and categories with families are different (but equivalent) presentations of the same essentially algebraic theory.
\end{numbered}

\subsection*{Type theoretic structure as essentially algebraic structure}

\begin{definition}[Theory of natural models with a distinguished set of basic types]
\index{essentially algebraic theory!of natural models basic types@{---} of natural models admitting a set of basic types}
Given a set $I$, the \textbf{theory of natural models admitting an $I$-indexed set of basic types} is the theory $\mathbb{T}_{(\mathsf{ty}_i)_{i \in I}}$ extending $\mathbb{T}$ by adding a new total symbol $\mathsf{btype}_i : \mathsf{Type}$ for each $i \in I$, together with the equation $\mathsf{ctxof}(\mathsf{btype}_i) = \diamond$ for each $i \in I$.
\end{definition}

Evidently we have the following characterisation of natural models with a distinguished set of basic types.

\begin{theorem}
Specifying a model $\mathfrak{M}$ of the theory $\mathbb{T}_{(\mathsf{ty}_i)_{i \in I}}$ is equivalent to specifying a natural model $(\mathbb{C}, p)$ together with an $I$-indexed set $\{ O_i \mid i \in I \} \subseteq \nmty{U}(\diamond)$. \qed
\end{theorem}

\begin{definition}[Theory of natural models with a distinguished set of terms of basic types]\index{essentially algebraic theory!of natural models basic terms@{---} of natural models admitting a set of terms}
Given a set $J$, the \textbf{theory of natural models admitting a $J$-indexed set of terms of basic types} is the theory $\mathbb{T}_{(\mathsf{tm}_j)_{j \in I}}$ extending $\mathbb{T}$ by adding a new total symbol $\mathsf{bterm}_j : \mathsf{Term}$ for each $j \in J$, together with the equation $\mathsf{ctxof}(\mathsf{bterm}_j) = \diamond$ for each $j \in J$.
\end{definition}

\begin{theorem}
Specifying a model $\mathfrak{M}$ of the theory $\mathbb{T}_{(\mathsf{tm}_j)_{j \in J}}$ is equivalent to specifying a natural model $(\mathbb{C}, p)$ together with an $J$-indexed set $\{ o_j \mid j \in J \} \subseteq \nmtm{U}(\diamond)$. \qed
\end{theorem}

Note that a natural model admitting a $J$-indexed set of basic terms automatically has the structure of a natural model admitting a $J$-indexed set of basic types by taking $\{ p_{\diamond}(o_j) \mid j \in J \}$ to be the distinguished set of basic types. As a result, there is a forgetful functor from the category of natural models admitting a $J$-indexed set of basic types to the category of natural models admitting a $J$-indexed set of basic types.

\begin{definition}[Theory of natural models admitting a unit type]
\label{defEATAdjUnit}
\index{essentially algebraic theory!of natural models unit@{---} of natural models admitting a unit type}
The \textbf{theory of natural models admitting a unit type} is the essentially algebraic theory $\mathbb{T}\adjunit$ extending $\mathbb{T}$ as follows.
\begin{itemize}
\item The new symbols and their arities are indicated in the following table.
\begin{center}
\begin{tabular}{l|lclcl|lcl}
\textbf{Name} & \multicolumn{5}{l|}{\textbf{Symbol} (with arity \& sorts)} & \multicolumn{3}{l}{\textbf{Shorthand}} \\ \hline
unit type$^{\dagger}$ & $\mathsf{unit}$ & $:$ &  &  & $\mathsf{Type}$ & && \\
term of unit type$^{\dagger}$ & $\mathsf{star}$ & $:$ &  &  & $\mathsf{Term}$ & && 
\end{tabular}
\end{center}
\item The new equations are as follows.
\begin{enumerate}[(i)]
\item $\mathsf{ctxof}(\mathsf{unit}) = \diamond$
\item $\mathsf{typeof}(\mathsf{star}) = \mathsf{unit}$
\item $\nmp{\mathsf{unit}} = \nmt{\diamond \cext \mathsf{unit}}$
\item $\nmq{\mathsf{unit}} = \mathsf{star}[\nmt{\diamond \cext \mathsf{unit}}]$
\end{enumerate}
\item There are no additional domains of definition to specify, since both new symbols are total.
\end{itemize}
\end{definition}

\begin{theorem}
\label{thmEATAdjUnit}
Specifying a model $\mathfrak{M}$ of $\mathbb{T}\adjunit$ is equivalent to specifying a natural model $(\mathbb{C}, p)$ together with elements $\widehat{\mathbbm{1}} \in \nmty{U}(\diamond)$ and $\widehat{\star} \in \nmtm{U}(\diamond)$ such that the following square is a pullback.
\begin{diagram}
\Yon(\diamond)
\arrow[d, equals]
\arrow[r, "\widehat{\star}"]
\pullback
&
\nmtm{U}
\arrow[d, "p"]
\\
\Yon(\diamond)
\arrow[r, "\widehat{\mathbbm{1}}"']
&
\nmty{U}
\end{diagram}
\end{theorem}

\begin{proof}
Take $\widehat{\mathbbm{1}} = \mathsf{unit}^{\mathfrak{M}}$ and $\widehat{\star} = \mathsf{star}^{\mathfrak{M}}$. Equations (i) and (ii) say that these have the correct types and that the square in the statement of the theorem commutes. Equation (iii) is redundant since the codomain of $\nmp{\mathbbm{1}}$ is terminal.

Equation (iv) is equivalent to the assertion that the square is a pullback. Indeed, suppose (iv) holds, and let $\Gamma \in \mathrm{ob}(\mathbb{C})$ and $a \in \nmtm{U}(\Gamma; \widehat{\mathbbm{1}}[\nmt{\Gamma}])$.
\begin{diagram}
\Yon(\Gamma)
\arrow[ddr, bend right=15, "\Yon(\nmt{\Gamma})"']
\arrow[drr, bend left=15, "a"]
\arrow[dr, dashed, "\Yon(\nmt{\Gamma})" description]
&[-20pt]&
\\[-20pt]
&
\Yon(\diamond)
\arrow[d, equals]
\arrow[r, "\widehat{\star}"]
&
\nmtm{U}
\arrow[d, "p"]
\\
&
\Yon(\diamond)
\arrow[r, "\widehat{\mathbbm{1}}"']
&
\nmty{U}
\end{diagram}
Then we have
\begin{align*}
a &= \nmq{\mathsf{unit}}[\langle \nmt{\Gamma}, a \rangle_{\mathsf{unit}}] && \text{by \Cref{defEATheoryOfNaturalModels}(xxvi)} \\
&= \mathsf{star}[\nmt{\diamond \cext \mathsf{unit}}][\langle \nmt{\Gamma}, a \rangle_{\mathsf{unit}}] && \text{by \Cref{defEATAdjUnit}(iv)} \\
&= \mathsf{star}[\nmp{\mathsf{unit}}][\langle \nmt{\Gamma}, a \rangle_{\mathsf{unit}}] && \text{by \Cref{defEATAdjUnit}(iii)} \\
&= \mathsf{star}[\nmp{\mathsf{unit}} \circ \langle \nmt{\Gamma}, a \rangle_{\mathsf{unit}}] && \text{by \Cref{defEATheoryOfNaturalModels}(xv)} \\
&= \mathsf{star}[\nmt{\Gamma}] && \text{by \Cref{defEATheoryOfNaturalModels}(xxv)}
\end{align*}
Uniqueness of $\nmt{\Gamma}$ is immediate from the fact that its codomain is terminal.

Conversely, if the square is a pullback, then take $\Gamma = \diamond \cext \widehat{\mathbbm{1}}$ and $a = \nmq{\widehat{\mathbbm{1}}}$ in the above. Since $\nmt{\diamond \cext \widehat{\mathbbm{1}}}$ is the morphism induced by the universal property of the pullback, we have $\widehat{\star}[\nmt{\diamond \cext \widehat{\mathbbm{1}}}] = \nmq{\widehat{\mathbbm{1}}}$, so that (iv) holds.
\end{proof}

\begin{definition}[Theory of natural models admitting dependent sum types]
\index{essentially algebraic theory!of natural models dependent sum@{---} of natural models admitting dependent sum types}
The \textbf{theory of natural models admitting dependent sum types} is the essentially algebraic theory $\mathbb{T}\adjsigma$ extending $\mathbb{T}$ as follows.
\begin{itemize}
\item The new symbols and their arities are indicated in the following table.
\begin{center}
\begin{tabular}{l|lclcl|lcl}
\textbf{Name} & \multicolumn{5}{l|}{\textbf{Symbol} (with arity \& sorts)} & \multicolumn{3}{l}{\textbf{Shorthand}} \\ \hline
dependent sum type & $\mathsf{sigma}$ & $:$ & $\mathsf{Type} \times \mathsf{Type}$ & $\to$ & $\mathsf{Type}$ & && \\
pairing & $\mathsf{pair}$ & $:$ & $\mathsf{Type} \times \mathsf{Type} \times \mathsf{Term} \times \mathsf{Term}$ & $\to$ & $\mathsf{Term}$ & && \\
first projection & $\mathsf{fst}$ & $:$ & $\mathsf{Type} \times \mathsf{Type} \times \mathsf{Term}$ & $\to$ & $\mathsf{Term}$ & && \\
second projection & $\mathsf{snd}$ & $:$ & $\mathsf{Type} \times \mathsf{Type} \times \mathsf{Term}$ & $\to$ & $\mathsf{Term}$ & &&
\end{tabular}
\end{center}
\item The new equations are as follows:
\setcounter{axiomcount}{0}
\begin{itemize}
\item \textbf{Dependent sum type-former}
\begin{enumerate}[(i)]
\item $\mathsf{ctxof}(\mathsf{sigma}(A,B)) = \mathsf{ctxof}(A)$
\item $\mathsf{sigma}(A,B)[\sigma] = \mathsf{sigma}(A[\sigma], B[\sigma \cext A])$
\setcounter{axiomcount}{\value{enumi}}
\end{enumerate}
\item \textbf{Pairing term-former}
\begin{enumerate}[(i)]
\setcounter{enumi}{\value{axiomcount}}
\item $\mathsf{typeof}(\mathsf{pair}_{A,B}(a,b)) = \mathsf{sigma}(A,B)$
\item $\mathsf{pair}_{A,B}(a,b)[\sigma] = \mathsf{pair}_{A[\sigma],B[\sigma \cext A]}(a[\sigma], b[\sigma])$
\setcounter{axiomcount}{\value{enumi}}
\end{enumerate}
\item \textbf{First and second projections}
\begin{enumerate}[(i)]
\setcounter{enumi}{\value{axiomcount}}
\item $\mathsf{typeof}(\mathsf{fst}_{A,B}(p)) = A$
\item $\mathsf{fst}_{A,B}(p)[\sigma] = \mathsf{fst}_{A[\sigma], B[\sigma \cext A]}(p[\sigma])$
\item $\mathsf{typeof}(\mathsf{snd}_{A,B}(p)) = B[\langle \mathsf{id}_{\mathsf{ctxof}(A)}, \mathsf{fst}_{A,B}(p) \rangle_A]$
\item $\mathsf{snd}_{A,B}(p)[\sigma] = \mathsf{snd}_{A[\sigma], B[\sigma \cext A]}(p[\sigma])$
\setcounter{axiomcount}{\value{enumi}}
\end{enumerate}
\item \textbf{Computation rules}
\begin{enumerate}[(i)]
\setcounter{enumi}{\value{axiomcount}}
\item $\mathsf{fst}_{A,B}(\mathsf{pair}(a,b)) = a$
\item $\mathsf{snd}_{A,B}(\mathsf{pair}(a,b)) = b$
\item $\mathsf{pair}_{A,B}(\mathsf{fst}_{A,B}(p),\mathsf{snd}_{A,B}(p)) = p$
\end{enumerate}
\end{itemize}
\item The domains of definition of the partial symbols are given as follows.
\begin{center}
\begin{tabular}{lcl}
$\mathsf{Def}(\mathsf{sigma}(A,B))$ & $=$ & $\{ \mathsf{ctxof}(A) \cext A = \mathsf{ctxof}(B) \}$ \\
$\mathsf{Def}(\mathsf{pair}_{A,B}(a,b))$ & $=$ & $\{ \mathsf{ctxof}(A) \cext A = \mathsf{ctxof}(B),\ \mathsf{typeof}(a)=A,\ \mathsf{typeof}(b) = B[\langle \mathrm{id}_{\mathsf{ctxof}(A)}, a \rangle_A] \}$ \\
$\mathsf{Def}(\mathsf{fst}_{A,B}(p))$ & $=$ & $\{ \mathsf{ctxof}(A) \cext A = \mathsf{ctxof}(B),\ \mathsf{typeof}(p) = \mathsf{sigma}(A,B) \}$ \\
$\mathsf{Def}(\mathsf{snd}_{A,B}(p))$ & $=$ & $\{ \mathsf{ctxof}(A) \cext A = \mathsf{ctxof}(B),\ \mathsf{typeof}(p) = \mathsf{sigma}(A,B) \}$
\end{tabular}
\end{center}
\end{itemize}
\end{definition}

\begin{theorem}
\label{thmEATAdjSigma}
Specifying a model $\mathfrak{M}$ of $\mathbb{T}\adjsigma$ is equivalent to specifying a natural model $(\mathbb{C}, p)$ together with natural transformations $\widehat{\upSigma}$ and $\widehat{\mathsf{pair}}$ such that the following diagram is a pullback.
\begin{diagram}
\sum_{A : \nmty{U}} \sum_{B : \nmty{U}^{[A]}} \sum_{a:[A]} [B(a)]
\arrow[d, "\pi"']
\arrow[r, "\widehat{\mathsf{pair}}"]
\pullbackc{dr}{-0.05}
&
\nmtm{U}
\arrow[d, "p"]
\\
\sum_{A : \nmty{U}} \nmty{U}^{[A]}
\arrow[r, "\widehat{\upSigma}"']
&
\nmty{U}
\end{diagram}
\end{theorem}

\begin{proof}
Suppose $\mathfrak{M}$ is a model of $\mathbb{T} \adjsigma$ with underlying natural model $(\mathbb{C}, p)$. By \Cref{lemLemmaFive} there is a natural bijection
$$\xi_{\Gamma} : \left( \sum_{A \in \nmty{U}} \nmty{U}^{[A]} \right)(\Gamma) \cong \sum_{A \in \nmty{U}(\Gamma)} \nmty{U}(\Gamma \cext A)$$
so define $\widehat{\upSigma}_{\Gamma}$ to be the composite with $\xi_{\Gamma}$ of the restriction of $\mathsf{sigma}^{\mathfrak{M}}$ to $\sum_{A \in \nmty{U}(\Gamma)} \nmty{U}(\Gamma \cext A)$. Equation (i) ensures that the image of $\widehat{\upSigma}_{\Gamma}$ is contained in $\nmty{U}(\Gamma)$, so that the functions $\widehat{\upSigma}_{\Gamma}$ have the correct codomains, and equation (ii) tells us that $\widehat{\upSigma}$ is natural.

Likewise, by \Cref{lemLemmaElevenPointFive} there is a natural bijection
$$\zeta_{\Gamma} : \left( \sum_{A \in \nmty{U}} \sum_{B \in \nmty{U}^{[A]}} \sum_{a \in [A]} [B(a)] \right)(\Gamma) \cong \sum_{A \in \nmty{U}(\Gamma)} \sum_{B \in \nmty{U}(\Gamma \cext A)} \sum_{a \in \nmtm{U}(\Gamma; A)} \nmtm{U}(\Gamma; B[\nms{a}])$$
which allows us to translate between $\mathsf{pair}^{\mathfrak{M}}$ and $\widehat{\mathsf{pair}}$; equation (iii) tells us that the components of $\widehat{\mathsf{pair}}$ have the correct types and that the square commutes, and equation (iv) tells us that it is a natural transformation.

Equations (v)--(viii) then describe the morphisms induced by the universal property of the pullbacks, as indicated in the following diagram, in which we have $A \in \nmty{U}(\Gamma)$, $B \in \nmty{U}(\Gamma \cext A)$ and $p \in \nmtm{U}(\Gamma; \widehat{\upSigma}_{\Gamma}(A,B))$.
\begin{diagram}
\Yon(\Gamma)
\arrow[ddr, bend right=30, "{(A,B)}"']
\arrow[drr, bend left=20, "p"]
\arrow[dr, dashed, "{(A,B,\mathsf{fst}^{\mathfrak{M}}(A,B,p),\mathsf{snd}^{\mathfrak{M}}(A,B,p))}" description]
&[50pt]&
\\
&
\sum_{A : \nmty{U}} \sum_{B : \nmty{U}^{[A]}} \sum_{a:[A]} [B(a)]
\arrow[d, "\pi"']
\arrow[r, "\widehat{\mathsf{pair}}"]
&
\nmtm{U}
\arrow[d, "p"]
\\
&
\sum_{A : \nmty{U}} \nmty{U}^{[A]}
\arrow[r, "\widehat{\upSigma}"']
&
\nmty{U}
\end{diagram}
Equations (ix) and (x) then say that the dashed morphism makes the required triangles commute, and equation (xi) says that it is the unique such morphism.
\end{proof}

\begin{definition}[Theory of natural models admitting dependent product types]
\index{essentially algebraic theory!of natural models dependent product@{---} of natural models admitting dependent product types}
The \textbf{theory of natural models admitting dependent product types} is the essentially algebraic theory $\mathbb{T}\adjpi$ extending $\mathbb{T}$ as follows.
\begin{itemize}
\item The new symbols and their arities are indicated in the following table.
\begin{center}
\begin{tabular}{l|lclcl}
\textbf{Name} & \multicolumn{5}{l}{\textbf{Symbol} (with arity \& sorts)} \\ \hline
dependent product type & $\mathsf{pi}$ & $:$ & $\mathsf{Type} \times \mathsf{Type}$ & $\to$ & $\mathsf{Type}$ \\
$\lambda$-abstraction & $\mathsf{lambda}$ & $:$ & $\mathsf{Type} \times \mathsf{Type} \times \mathsf{Term}$ & $\to$ & $\mathsf{Term}$ \\
application & $\mathsf{app}$ & $:$ & $\mathsf{Type} \times \mathsf{Type} \times \mathsf{Term} \times \mathsf{Term}$ & $\to$ & $\mathsf{Term}$ \\
\end{tabular}
\end{center}
\item The new equations are as follows.
\setcounter{axiomcount}{0}
\begin{itemize}
\item \textbf{Dependent product type-former}
\begin{enumerate}[(i)]
\item $\mathsf{ctxof}(\mathsf{pi}(A,B)) = \mathsf{ctxof}(A)$
\item $\mathsf{pi}(A,B)[\sigma] = \mathsf{pi}(A[\sigma], B[\sigma \cext A])$
\setcounter{axiomcount}{\value{enumi}}
\end{enumerate}
\item \textbf{$\lambda$-abstraction term-former}
\begin{enumerate}[(i)]
\setcounter{enumi}{\value{axiomcount}}
\item $\mathsf{typeof}(\mathsf{lambda}_{A,B}(b)) = \mathsf{pi}(A,B)$
\item $\mathsf{lambda}_{A,B}(b)[\sigma] = \mathsf{lambda}_{A[\sigma],B[\sigma \cext A]}(b[\sigma \cext A])$
\setcounter{axiomcount}{\value{enumi}}
\end{enumerate}
\item \textbf{Application}
\begin{enumerate}[(i)]
\setcounter{enumi}{\value{axiomcount}}
\item $\mathsf{typeof}(\mathsf{app}_{A,B}(f,a)) = B[\langle \mathsf{id}_{\mathsf{ctxof}(A)}, a \rangle_A]$
\item $\mathsf{app}_{A,B}(f,a)[\sigma] = \mathsf{app}_{A[\sigma], B[\sigma \cext A]}(f[\sigma], a[\sigma])$
\setcounter{axiomcount}{\value{enumi}}
\end{enumerate}
\item \textbf{Computation rules}
\begin{enumerate}[(i)]
\setcounter{enumi}{\value{axiomcount}}
\item $\mathsf{app}_{A,B}(\mathsf{lambda}_{A,B}(b),a) = b[\langle \mathsf{id}_{\mathsf{ctxof}(A)}, a \rangle_A]$
\item $\mathsf{lambda}_{A,B}(\mathsf{app}_{A,B}(f[\nmp{A}], \nmq{A})) = f$
\end{enumerate}
\end{itemize}
\item The domains of definition of the partial symbols are given as follows.
\begin{center}
\begin{tabular}{lcl}
$\mathsf{Def}(\mathsf{pi}(A,B))$ & $=$ & $\{ \mathsf{ctxof}(A) \cext A = \mathsf{ctxof}(B) \}$ \\
$\mathsf{Def}(\mathsf{lambda}_{A,B}(b))$ & $=$ & $\{ \mathsf{ctxof}(A) \cext A = \mathsf{ctxof}(B),\ \mathsf{typeof}(b) = B \}$ \\
$\mathsf{Def}(\mathsf{app}_{A,B}(f,a))$ & $=$ & $\{ \mathsf{ctxof}(A) \cext A = \mathsf{ctxof}(B),\ \mathsf{typeof}(f) = \mathsf{pi}(A,B),\ \mathsf{typeof}(a)=A \}$ \\
\end{tabular}
\end{center}
\end{itemize}
\end{definition}

\begin{theorem}
\label{thmEATAdjPi}
Specifying a model $\mathfrak{M}$ of $\mathbb{T}\adjpi$ is equivalent to specifying a natural model $(\mathbb{C}, p)$ together with natural transformations $\widehat{\upPi}$ and $\widehat{\lambda}$ such that the following square is a pullback.
\begin{diagram}
\sum_{A : \nmty{U}} \nmtm{U}^{[A]}
\arrow[r, "\widehat{\lambda}"]
\arrow[d, "\sum_{A : \nmty{U}} p^{[A]}"']
\pullbackc{dr}{0}
&
\nmtm{U}
\arrow[d, "p"]
\\
\sum_{A : \nmty{U}} \nmty{U}^{[A]}
\arrow[r, "\widehat{\upPi}"']
&
\nmty{U}
\end{diagram}
\end{theorem}

\begin{proof}
Suppose $\mathfrak{M}$ is a model of $\mathbb{T}\adjpi$ with underlying natural model $(\mathbb{C}, p)$. Again using \Cref{lemLemmaFive} there is a natural bijection
$$\xi_{\Gamma} : \left( \sum_{A \in \nmty{U}} \nmty{U}^{[A]} \right)(\Gamma) \cong \sum_{A \in \nmty{U}(\Gamma)} \nmty{U}(\Gamma \cext A)$$
so define $\widehat{\upPi}_{\Gamma}$ to be the composite with $\xi_{\Gamma}$ of the restriction of $\mathsf{pi}^{\mathfrak{M}}$ to $\sum_{A \in \nmty{U}(\Gamma)} \nmty{U}(\Gamma \cext A)$. Equation (i) ensures that the image of $\widehat{\upPi}_{\Gamma}$ is contained in $\nmty{U}(\Gamma)$, so that the functions $\widehat{\upPi}_{\Gamma}$ have the correct codomains, and equation (ii) tells us that $\widehat{\upPi}$ is natural.

Likewise, by \Cref{lemLemmaFive} again, there is a natural bijection
$$\zeta_{\Gamma} : \left(\sum_{A \in \nmty{U}} \nmtm{U}^{[A]} \right)(\Gamma) \cong \sum_{A \in \nmty{U}(\Gamma)} \nmtm{U}(\Gamma \cext A)$$
which allows us to translate between $\mathsf{lambda}^{\mathfrak{M}}$ and $\widehat{\lambda}$; equation (iii) tells us that the components of $\widehat{\lambda}$ have the correct types and that the square commutes, and equation (iv) tells us that it is a natural transformation.

Equations (v)--(vi) then describe the morphisms induced by the universal property of the pullback indicated in the following diagram, in which we have $A \in \nmty{U}(\Gamma)$, $B \in \nmty{U}(\Gamma \cext A)$ and $f \in \nmtm{U}(\Gamma; \widehat{\upPi}_{\Gamma}(A,B))$.

\begin{diagram}
\Yon(\Gamma)
\arrow[drr, bend left=20, "f"]
\arrow[ddr, bend right=15, "{(A,B)}"']
\arrow[dr, dashed, "{(A, \lambda x:A.\mathsf{app}^{\mathfrak{M}}(f,x))}" description]
&[30pt]&
\\[-20pt]
&
\sum_{A \in \nmty{U}} \nmtm{U}^{[A]}
\arrow[d, "\sum_{A \in \nmty{U}} \nmty{U}^{[A]}" description]
\arrow[r, "\widehat{\lambda}"]
&
\nmtm{U}
\arrow[d, "p"]
\\
&
\sum_{A \in \nmty{U}} \nmty{U}^{[A]}
\arrow[r, "\widehat{\upPi}"']
&
\nmty{U}
\end{diagram}

Equation (vii) then say that the dashed morphism makes the required triangles commute, and equation (viii) says that it is the unique such morphism.
\end{proof}

\begin{numbered}
Although we have only discussed adding type theoretic structure to a `pure' natural model, we can combine the theories above in a modular way in order to add structure to an already structured natural model. For example, a natural model admitting a unit type \textit{and} dependent sum types is a natural model equipped with data making it both a model of $\mathbb{T}\adjunit$ and of $\mathbb{T}\adjsigma$, or equivalently a natural model equipped with data satisfying the hypotheses of both \Cref{thmEATAdjUnit} and \Cref{thmEATAdjSigma}.
\end{numbered}

All of the theories discussed above are finitary essentially algebraic theories, except possibly for the theories of natural models admitting an $I$-indexed family of basic types or a $J$-indexed family of basic terms, which are $\mathrm{max}\{|I|,|J|\}$-ary (and finitary when $I$ and $J$ are finite).

\begin{definition}[Locally presentable categories \thmcite{AdamekRosicky1994}]
Let $\lambda$ be a regular cardinal and let $\mathcal{C}$ be a category. An object $A$ of $\mathcal{C}$ is \textbf{$\lambda$-presentable} if the functor $\mathcal{C}(A,{-}) : \mathcal{C} \to \mathbf{Set}$ preserves $\lambda$-directed colimits. The category $\mathcal{C}$ is \textbf{locally $\lambda$-presentable} if it is cocomplete and has a small set $\mathcal{A} \subseteq \mathrm{ob}(\mathcal{C})$ of $\lambda$-presentable objects such that every object of $\mathcal{C}$ is a $\lambda$-directed colimit of objects in $\mathcal{A}$.
\end{definition}

It is known \cite[Theorem 3.36]{AdamekRosicky1994} that for an regular cardinal $\lambda$, a category is a model of a $\lambda$-ary essentially algebraic theory if and only if it is locally $\lambda$-presentable. As such all of our categories of suitably structured natural models will satisfy the properties enjoyed by locally presentable categories more generally, such as cocompleteness.

In particular, given a dependent type theory $\mathbb{T}$, the category $\mathbf{NM}_{\mathbb{T}}$ of natural models admitting the type theoretic structure in $\mathbb{T}$ has an initial object $(\mathbb{C}_{\mathbb{T}}, p_{\mathbb{T}})$. A long term goal for future work is to prove that $(\mathbb{C}_{\mathbb{T}}, p_{\mathbb{T}})$ can be described as the \textit{term model} \Cref{cnsFreeNaturalModelKappaBasicTypes} of $\mathbb{T}$, for an arbitrary type theory $\mathbb{T}$.

Moreover, given a subtheory $\mathbb{T}' \subseteq \mathbb{T}$, we obtain from the essentially algebraic character of natural models a forgetful functor $U : \mathbf{NM}_{\mathbb{T}} \to \mathbf{NM}_{\mathbb{T}'}$, which has a left adjoint $F : \mathbf{NM}_{\mathbb{T}'} \to \mathbf{NM}_{\mathbb{T}}$. Given a natural model $(\mathbb{C}, p)$ admitting the type theoretic structure of $\mathbb{T}'$, the object $F(\mathbb{C}, p)$ of $\mathbf{NM}_{\mathbb{T}}$ is then the \textit{free} natural model on $\mathbb{C}$ which supports the type theoretic structure of $\mathbb{T}$. In \Cref{chSemantics}, we will describe some examples explicitly.

%% file: thesis/ch2-categories-of-natural-models/homomorphisms.tex
\section{Morphisms of natural models}
\label{secHomomorphisms}

Since a natural model is a model of an essentially algebraic theory, there is a canonical notion of \textit{homomorphism of natural models}, namely maps between the sorts commuting with the structure.

\begin{definition}[Morphisms of natural models]
\label{defMorphismOfNMAsEAT}
\index{morphism!of natural models@{---} of natural models}
Let $(\mathbb{C}, p)$ and $(\mathbb{D}, q)$ be natural models. A \textbf{morphism of natural models} is a homomorphism from $(\mathbb{C}, p)$ to $(\mathbb{D}, q)$, with $(\mathbb{C}, p)$ and $(\mathbb{D}, q)$ considered as models of the essentially algebraic theory $\mathbb{T}_{\mathrm{NM}}$ (\Cref{defEATheoryOfNaturalModels}). The category of all natural models and morphisms between them is denoted by $\mathbf{NM}$.
\end{definition}

Explicitly, \Cref{defMorphismOfNMAsEAT} says that a morphism of natural models is a quadruple of functions $(F_{\mathsf{Ctx}}, F_{\mathsf{Sub}}, F_{\mathsf{Term}}, F_{\mathsf{Type}})$ such that:
\begin{itemize}
\item $(F_{\mathsf{Ctx}}, F_{\mathsf{Sub}})$ defines a functor $\mathbb{C} \to \mathbb{D}$ strictly preserving distinguished terminal objects;
\item $F_{\mathsf{Term}} : \sum_{\Gamma \in \mathrm{ob}(\mathbb{C})} \nmty{U}(\Gamma) \to \sum_{\Gamma \in \mathrm{ob}(\mathbb{D})} \nmty{V}(\Gamma)$ respecting contexts and substitutions;
\item $F_{\mathsf{Type}} : \sum_{\Gamma \in \mathrm{ob}(\mathbb{C})} \nmtm{U}(\Gamma) \to \sum_{\Gamma \in \mathrm{ob}(\mathbb{D})} \nmtmalt{V}(\Gamma)$ respecting contexts and substitutions;
\end{itemize}
and such that the representability data and typing are preserved, in the sense that for each $\Gamma \in \mathrm{ob}(\mathbb{C})$ and each $A \in \nmty{U}(\Gamma)$, we have
$$F_{\mathsf{Ctx}}(\Gamma \cext A) = F_{\mathsf{Ctx}} \cextalt FA, \quad F_{\mathsf{Sub}}(\nmp{A}) = \nmu{F_{\mathsf{Type}}(A)} \quad \text{and} \quad F_{\mathsf{Term}}(\nmq{A}) = \nmv{F_{\mathsf{Type}}(A)}$$
and $F_{\mathsf{Term}}(a) \in \nmtmalt{V}(F_{\mathsf{Ctx}}(\Gamma), F_{\mathsf{Type}}(A))$ for all $a \in \nmtm{U}(\Gamma, A)$. Identity and composition of natural models are then given by the corresponding identity functions and composites of functions, respectively.

Note that the structure specified for models of $\mathbb{T}_{\mathrm{NM}}$ but not for natural models, that is the substitutions $\nmt{\Gamma} : \Gamma \to \diamond$ and $\langle \sigma, a \rangle_A : \Delta \to \Gamma \cext A$, are preserved automatically provided the above conditions hold.

The goal of this section is to provide characterisations of morphisms of natural models which will be more convenient for our purposes.

\begin{theorem}[Morphisms of natural models via categories of elements]
\label{thmMorphismsOfNMFromCatEl}
Let $(\mathbb{C}, p)$ and $(\mathbb{D}, q)$ be natural models. Specifying a morphism $F : (\mathbb{C}, p) \to (\mathbb{D}, q)$ is equivalent to specifying a functor $F_0 : \mathbb{C} \to \mathbb{D}$ preserving distinguished terminal objects, together with functors
$$F_1 : \catel[\mathbb{C}]{\nmty{U}} \to \catel[\mathbb{D}]{\nmty{V}} \quad \text{and} \quad F_2 : \catel[\mathbb{C}]{\nmtm{U}} \to \catel[\mathbb{D}]{\nmtmalt{V}}$$
such that
\begin{itemize}
\item The following diagram of categories and functors commutes;
\begin{diagram}
\catel[\mathbb{C}]{\nmtm{U}}
\arrow[rrr, "F_2"]
\arrow[dr, "\catel p" description]
\arrow[ddr, bend right, "\pi"']
&[-20pt]&&[-20pt]
\catel[\mathbb{D}]{\nmtmalt{V}}
\arrow[dl, "\catel q" description]
\arrow[ddl, bend left, "\pi"]
\\[-20pt]
&
\catel[\mathbb{C}]{\nmty{U}}
\arrow[r, "F_1" description]
\arrow[d, "\pi" description]
&
\catel[\mathbb{D}]{\nmty{V}}
\arrow[d, "\pi" description]
&
\\
&
\mathbb{C}
\arrow[r, "F_0"']
&
\mathbb{D}
&
\end{diagram}
\item $F_1$ and $F_2$ respect the adjunctions $\catel p \dashv p^*$ and $\catel q \dashv q^*$, in the sense that $F_2 \circ p^* = q^* \circ F_1$ and, letting $(\eta, \varepsilon)$ and $(\sigma, \tau)$ be the respective (unit, counit) pairs, the following pasting diagrams commute.

\begin{diagram}
\int_{\mathbb{C}} \nmty{U}
\arrow[r, bend left=20, "p \circ p^*", ""'{name=domeps, inner sep=1pt}]
\arrow[r, bend right=20, "\mathrm{id}"', ""{name=codeps, inner sep=1pt}]
\arrow[Rightarrow, "\varepsilon", from=domeps, to=codeps]
\arrow[d, "F_1"']
&
\int_{\mathbb{C}} \nmty{U}
\arrow[d, "F_1"]
&
\int_{\mathbb{C}} \nmtm{U}
\arrow[r, bend left=20, "\mathrm{id}", ""'{name=dometa, inner sep=1pt}]
\arrow[r, bend right=20, "q^* \circ q"', ""{name=codeta, inner sep=1pt}]
\arrow[Rightarrow, "\eta", from=dometa, to=codeta]
\arrow[d, "F_2"']
&
\int_{\mathbb{C}} \nmtm{U}
\arrow[d, "F_2"]
\\
\int_{\mathbb{D}} \nmty{V}
\arrow[r, bend left=20, "q \circ q^*", ""'{name=domtau, inner sep=1pt}]
\arrow[r, bend right=20, "\mathrm{id}"', ""{name=codtau, inner sep=1pt}]
\arrow[Rightarrow, "\tau", from=domtau, to=codtau]
&
\int_{\mathbb{D}} \nmty{V}
&
\int_{\mathbb{D}} \nmtmalt{V}
\arrow[r, bend left=20, "\mathrm{id}", ""'{name=domsig, inner sep=1pt}]
\arrow[r, bend right=20, "q^* \circ q"', ""{name=codsig, inner sep=1pt}]
\arrow[Rightarrow, "\sigma", from=domsig, to=codsig]
&
\int_{\mathbb{D}} \nmtmalt{V}
\end{diagram}
\end{itemize}
\end{theorem}

\begin{proof}
Evidently the functor $F_0$ is precisely the functor determined by $F_{\mathsf{Ctx}}$ and $F_{\mathsf{Sub}}$. The functor $F_1$ determines and is determined by $F_{\mathsf{Type}}$; explicitly, on objects $(\Gamma, A)$ we have $F_1(\Gamma, A) = (F_{\mathsf{Ctx}}(\Gamma), F_{\mathsf{Type}}(A))$ and on morphisms $\sigma : \Delta \to \Gamma$ by $F_1(\sigma) = F_{\mathsf{Sub}}(\sigma)$. The fact that $F_1$ is a functor is immediate from the facts that $F_0$ is a functor and that $F_{\mathsf{Type}}$ respects contexts and substitutions. The correspondence between $F_2$ and $F_{\mathsf{Term}}$ is obtained likewise.

Commutativity of the bottom and outer squares is equivalent to the assertion that $F_1$ and $F_2$ agree with $F_0$ on their first components, which is immediate from their construction. Commutativity of the top square is equivalent to the assertion that $F_{\mathsf{Type}}$ and $F_{\mathsf{Term}}$ respect typing.

That $F_1$ and $F_2$ respect the adjunctions $\catel{p} \dashv p^*$ and $\catel{q} \dashv q^*$ is equivalent to preservation of representability data. Indeed, the equation $F_2 \circ p^* = q^* \circ F_1$ says that $F_{\mathsf{Ctx}}(\Gamma \cext A) = F_{\mathsf{Ctx}}(\Gamma) \cext F_{\mathsf{Type}}(A)$ and $F_{\mathsf{Term}}(\nmq{A}) = \nmv{F_{\mathsf{Type}}(A)}$; preservation of counits says that $F_{\mathsf{Sub}}(\nmp{A}) = \nmu{F_{\mathsf{Type}}(A)}$. Preservation of units corresponds with preservation of induced substitutions $\langle \sigma, a \rangle_A$.
\end{proof}

The characterisation of morphisms of natural models as a triple of functors given by \Cref{thmMorphismsOfNMFromCatEl} is useful because it is easy to mediate between this description and the explicit description of homomorphisms of models of $\mathbb{T}_{\mathrm{NM}}$, and composition is simply given by composition of functors. A drawback of this description, however, is that when working with natural models directly, it becomes cumbersome to construct categories of elements and keep track of units and counits.

Our next characterisation of morphisms of natural models allows us to work more directly with the representable natural transformations themselves. Recall that every functor $F : \mathbb{C} \to \mathbb{D}$ between small categories induces an adjoint triple $F_! \dashv F^* \dashv F_*$ between the corresponding categories of presheaves, where $F^*$ is given by precomposition by $F$.
\begin{diagram}
\widehat{\mathbb{C}}
\arrow[r, bend left=15, shift left=1, "F_!", ""'{name=domadj1}]
\arrow[r, bend right=15, shift right=1, "F_*"', ""{name=codadj2}]
&[100pt]
\widehat{\mathbb{D}}
\arrow[l, "F^*" description, ""'{name=codadj1}, ""{name=domadj2}]
\arrow[from=domadj1, to=codadj1, draw=none, "\bot" description]
\arrow[from=domadj2, to=codadj2, draw=none, "\bot" description]
\end{diagram}
A morphism of natural models $F : (\mathbb{C}, p) \to (\mathbb{D}, q)$ then induces functions $\nmty{U}(\Gamma) \to \nmty{V}(F\Gamma) = F^*\nmty{V}(\Gamma)$ and $\nmtm{U}(\Gamma) \to \nmtmalt{V}(F\Gamma) = F^*\nmtmalt{V}(\Gamma)$, given by the restrictions of $F_{\mathsf{Type}}$ and $F_{\mathsf{Term}}$, respectively. We will see (\Cref{thmMorphismsOfNMFunctorially}) that these functions are the components of natural transformations $\varphi : \nmty{U} \to F^*\nmty{V}$ and $\nmmk{\varphi} : \nmtm{U} \to F^*\nmtmalt{V}$, respectively, which correspond under the adjunction $F_! \dashv F^*$ with natural transformations $F_!\nmty{U} \to \nmty{V}$ and $F_!\nmtm{U} \to \nmtmalt{V}$.

The choice of whether to work with the left adjoint $F_!$ or the right adjoint $F^*$ is largely arbitrary, and we will translate between the two freely. The \textit{left adjoint convention} will be convenient in some settings because $F_!$ commutes with Yoneda embeddings; meanwhile, the \textit{right adjoint convention} will be convenient in other settings because $F^*$ can be described explicitly with ease.

\begin{definition}[Premorphisms of natural models]
\label{defPremorphism}
\index{premorphism of natural models}
Let $(\mathbb{C}, p)$ and $(\mathbb{D}, q)$ be natural models. A \textbf{premorphism of natural models} from $(\mathbb{C}, p)$ to $(\mathbb{D}, q)$ is a triple $(F, \varphi, \nmmk{\varphi})$ consisting of a functor $F : \mathbb{C} \to \mathbb{D}$ preserving distinguished terminal objects, and natural transformations $\varphi : F_!\nmty{U} \to \nmty{V}$ and $\nmmk{\varphi} : F_!\nmtm{U} \to \nmtm{V}$ satisfying $q \circ \nmmk{\varphi} = \varphi \circ F_!p$.
\begin{center}
\begin{tikzcd}[row sep=huge, column sep=huge]
F_!\nmtm{U}
\arrow[r, "\nmmk{\varphi}"]
\arrow[d, "F_!p"']
\sqlabel{dr}{(\star)}
&
\nmtm{V}
\arrow[d, "q"]
\\
F_!\nmty{U}
\arrow[r, "\varphi"']
&
\nmty{V}
\end{tikzcd}
\end{center}
\end{definition}

\begin{numbered}
\label{parMorphismShorthand}
We will write $F : (\mathbb{C}, p) \pmto (\mathbb{D}, q)$ to denote the assertion that $F = (F, \varphi, \nmmk{\varphi})$ is a premorphism from $(\mathbb{C}, p)$ to $(\mathbb{D}, q)$. Given $\Gamma \in \mathbb{C}$, $A \in \nmty{U}(\Gamma)$ and $a \in \nmtm{U}(\Gamma)$, write $FA$ for the element of $\nmty{V}(F\Gamma)$ given by composite
$$\mathsf{y}(F\Gamma) = F_!\mathsf{y}(\Gamma) \xrightarrow{F_!A} F_!\nmty{U} \xrightarrow{\varphi} \nmty{V}$$
and write $Fa$ for the element of $\nmtm{V}(F\Gamma)$ given by the composite
$$\mathsf{y}(F\Gamma) = F_!\mathsf{y}(\Gamma) \xrightarrow{F_!a} F_!\nmtm{U} \xrightarrow{\nmmk{\varphi}} \nmtm{V}$$
Note that commutativity of $(\star)$ implies that if $\Gamma \vdash a : A$ in $(\mathbb{C}, p)$, then $F\Gamma \vdash Fa : FA$ in $(\mathbb{D}, q)$, since
$$q_{F\Gamma}(Fa) = q_{F\Gamma} \circ \nmmk{\varphi} \circ F_!a = \varphi \circ F_!p_{\Gamma} \circ F_!a = \varphi \circ F_!(p_{\Gamma}(a)) = \varphi \circ F_!A = FA$$
\end{numbered}

\begin{numbered}
\label{numRightHandedConvention}
In \Cref{defPremorphism} we adopted the left adjoint convention. Under the right adjoint convention, a premorphism from $(\mathbb{C}, p)$ to $(\mathbb{D}, q)$ consists of a functor $F : \mathbb{C} \to \mathbb{D}$ preserving distinguished terminal objects and natural transformations $\varphi : \nmty{U} \to F^*\nmty{V}$ and $\nmmk{\varphi} : \nmtm{U} \to F^*\nmtmalt{V}$ satisfying $F^*q \circ \nmmk{\varphi} = \varphi \circ p$. Then given $\Gamma \in \mathrm{ob}(\mathbb{C})$, $A \in \nmty{U}(\Gamma)$ and $a \in \nmtm{U}(\Gamma; A)$, the elements $FA$ and $Fa$ described in \Cref{parMorphismShorthand} are given by $FA = \varphi_{\Gamma}(A)$ and $Fa = \nmmk{\varphi}_{\Gamma}(a)$.
\end{numbered}

\begin{lemma}[Lax preservation of context extension]
\label{lemMediationMorphism}
Let $F : (\mathbb{C}, p) \rightharpoondown (\mathbb{D}, q)$ be a premorphism of natural models. For each $\Gamma \in \mathrm{ob}(\mathbb{C})$ and $A \in \nmty{U}(\Gamma)$, there is a unique morphism $\tau_A : F(\Gamma \cext A) \to F\Gamma \cext FA$ in $\mathbb{D}$ such that $\nmu{FA} \circ \tau_A = F\nmp{A}$ and $\nmv{FA}[\tau_A] = F\nmq{A}$.
\end{lemma}

\begin{proof}
We obtain $\tau_A$ as the morphism $\langle F\nmp{A}, F\nmq{A} \rangle_{FA} : F(\Gamma \cext A) \to F\Gamma \cext FA$ in $\mathbb{D}$.
\begin{diagram}
\Yon(F(\Gamma \cext A))
\arrow[ddr, bend right=15, "\Yon(F\nmp{A})"']
\arrow[drr, bend left=15, "\Yon(F\nmq{A})"]
\arrow[dr, dashed, "\tau_A" description]
&[-20pt]&
\\[-20pt]
&
\Yon(F\Gamma \cext FA)
\arrow[r, "\nmv{FA}"]
\arrow[d, "\Yon(\nmu{FA})"']
\pullback
&
\nmtmalt{V}
\arrow[d, "q"]
\\
&
\Yon(F\Gamma)
\arrow[r, "FA"']
&
\nmty{V}
\end{diagram}
To see that the outer square truly does commute, note that it is obtained as the outer square of the following diagram.
\begin{diagram}
\Yon (F(\Gamma \cext A))
\arrow[r, equals]
\arrow[d, "\Yon(F\nmp{A})"']
&[-30pt]
F_! \Yon(\Gamma \cext A)
\arrow[r, "F_! \nmq{A}"]
\arrow[d, "F_! \Yon(\nmp{A})" description]
&
F_! \nmtm{U}
\arrow[r, "\nmmk{\varphi}"]
\arrow[d, "F_!p" description]
&
\nmtm{V}
\arrow[d, "q"]
\\
\Yon(F\Gamma)
\arrow[r, equals]
&
F_!\Yon(\Gamma)
\arrow[r, "F_!A"']
&
F_! \nmty{U}
\arrow[r, "\varphi"']
&
\nmty{V}
\end{diagram}
The left-hand square commutes since $F_!$ commutes with Yoneda embeddings; the middle square commutes since it is the result of applying $F_!$ to the pullback square exhibiting $\Yon(\nmp{A})$ as a pullback of $p$; and the right-hand square commutes since $(F, \varphi, \nmmk{\varphi})$ is a premorphism of natural models.
\end{proof}

\begin{theorem}[Functorial characterisation of morphisms of natural models]
\label{thmMorphismsOfNMFunctorially}
\index{morphism!of natural models@{---} of natural models}
Let $(\mathbb{C}, p)$ and $(\mathbb{D}, q)$ be natural models. Specifying a morphism of natural models from $(\mathbb{C}, p)$ to $(\mathbb{D}, q)$ is equivalent to specifying a premorphism $(F, \varphi, \nmmk{\varphi}) : (\mathbb{C}, p) \rightharpoondown (\mathbb{D}, q)$ such that $F$ respects context extension (in the sense that $F(\Gamma \cext A) = F\Gamma \cext FA$ for each $\Gamma \in \mathrm{ob}(\mathbb{C})$ and $A \in \nmty{U}(\Gamma)$), and such that the morphisms $\tau_A : F(\Gamma \cext A) \to F\Gamma \cext FA$ of $\mathbb{D}$ are all identity morphisms.
\end{theorem}

\begin{proof}
Specifying a premorphism $(F, \varphi, \nmmk{\varphi})$ is equivalent to specifying a homomorphism of the theory of a category with a terminal object and a natural transformation between presheaves (equations (i)--(xviii) of \Cref{defEATheoryOfNaturalModels}). To see this, note that specifying a functor $F : \mathbb{C} \to \mathbb{D}$ is equivalent to specifying the pair $(F_{\mathsf{Ctx}}, F_{\mathsf{Sub}})$.

We will use the right adjoint convention (see \Cref{numRightHandedConvention}) for the natural transformations $\varphi, \nmmk{\varphi}$. For each $\Gamma \in \mathrm{ob}(\mathbb{C})$, the component $\varphi_{\Gamma} : \nmty{U}(\Gamma) \to \nmty{V}(F\Gamma)$ corresponds with the appropriate restriction of $F_{\mathsf{Type}}$, and likewise the component $\nmmk{\varphi}_{\Gamma} : \nmtm{U}(\Gamma) \to \nmtmalt{V}(F\Gamma)$ corresponds with the appropriate restriction of $F_{\mathsf{Term}}$. That $\varphi$ and $\nmmk{\varphi}$ are natural corresponds with the fact that they respect substitutions, and that the square ($\star$) commutes corresponds with the fact that they respect typing.

Finally, given $\Gamma \in \mathrm{ob}(\mathbb{C})$ and $A \in \nmty{U}(\Gamma)$, then note that $F\nmp{A} = \nmu{FA}$ and $F\nmq{A} = \nmq{FA}$ if and only if $\tau_A = \mathrm{id}_{F(\Gamma \cext A)}$, since by the universal property of pullbacks, $\tau_A$ is the unique morphism satisfying $\nmu{FA} \circ \tau_A = F\nmp{A}$ and $\nmv{FA}[\tau_A] = F\nmq{A}$.
\end{proof}

We will use the characterisation of morphisms of natural models given in \Cref{thmMorphismsOfNMFunctorially} almost exclusively in our proofs, so from now on we will typically use the term `morphism of natural models' to mean a premorphism preserving representability data.

\begin{numbered}
Let $(F, \varphi, \nmmk{\varphi}) : (\mathbb{C}, p) \to (\mathbb{D}, q)$ and $(G, \psi, \nmmk{\psi}) : (\mathbb{D}, q) \to (\mathbb{E}, r)$ be morphisms of natural models.
\begin{enumerate}[(a)]
\item Under the left adjoint convention, the composite $(G, \psi, \nmmk{\psi}) \circ (F, \varphi, \nmmk{\varphi}) : (\mathbb{C}, p) \to (\mathbb{E}, r)$ is given by $(G \circ F, \psi \circ G_!\varphi \circ \mu, \nmmk{\psi} \circ G_!\nmmk{\varphi} \circ \nmmk{\mu})$, where $\mu : (GF)_!\nmty{U} \to G_!F_!\nmty{U}$ and $\nmmk{\mu} : (GF)_!\nmtm{U} \to G_!F_!\nmtm{U}$ are the canonical isomorphisms (\Cref{thmKanExtensions}).
\item Under the right adjoint convention, the composite $(G, \psi, \nmmk{\psi}) \circ (F, \varphi, \nmmk{\varphi}) : (\mathbb{C}, p) \to (\mathbb{E}, r)$ is given by $(G \circ F, F^*\psi \circ \varphi, F^*\nmmk{\psi} \circ \nmmk{\varphi})$. 
\end{enumerate}
\end{numbered}

\subsection*{Weak morphisms of natural models}

\Cref{thmMorphismsOfNMFunctorially} suggests that we may obtain a weaker notion of morphism of natural models by weakening the requirement that the morphism $\tau_A : F(\Gamma \cext A) \to F\Gamma \cext FA$ be an identity morphism to the requirement that it be an isomorphism.

\begin{definition}[Weak morphisms of natural models]
\label{defWeakMorphismOfNM}
\index{morphism!of natural models weak@{---} of natural models (weak)}
Let $(\mathbb{C}, p)$ and $(\mathbb{D}, q)$ be natural models. A \textbf{weak morphism of natural models} from $(\mathbb{C}, p)$ to $(\mathbb{D}, q)$ is a premorphism $(F, \varphi, \nmmk{\varphi})$ such that the morphisms $\tau_A : F(\Gamma \cext A) \to F\Gamma \cext FA$ are isomorphisms for each $\Gamma \in \mathrm{ob}(\mathbb{C})$ and each $A \in \nmty{U}(\Gamma)$.
\end{definition}

\begin{numbered}
Write $F : (\mathbb{C}, p) \to (\mathbb{D}, q)$ to denote the assertion that $F = (F, \varphi, \nmmk{\varphi})$ is a (strict) morphism of natural models, and write $F : (\mathbb{C}, p) \wmto (\mathbb{D}, q)$ to denote the assertion that $F = (F, \varphi, \nmmk{\varphi}, \sigma, (\tau_A))$ is a weak morphism of natural models. We denote the category of natural models and weak morphisms by $\mathbf{NM}^{\mathrm{wk}}$. Note that there is an embedding $\mathbf{NM} \hookrightarrow \mathbf{NM}^{\mathrm{wk}}$ obtained by taking the coherence isomorphisms to be identities.
\end{numbered}

We will now explore some ways of characterising weak morphisms of natural models.

\begin{definition}[Preservation of canonical pullback squares]
A premorphism of natural models $F : (\mathbb{C}, p) \pmto (\mathbb{D}, q)$ (\textbf{weakly}) \textbf{preserves canonical pullback squares} (\Cref{cnsCanonicalPullbacks}) if, for each $\sigma : \Delta \to \Gamma$ in $\mathbb{C}$ and each $A \in \nmty{U}(\Gamma)$, the following square is a pullback
\begin{center}
\begin{tikzcd}[row sep=huge, column sep=huge]
F(\Delta \cext A[\sigma])
\arrow[r, "F(\sigma \cext A)"]
\arrow[d, "F\nmp{A\sigma}"']
&
F(\Gamma \cext A)
\arrow[d, "F\nmp{A}"{name=domeq}]
\\
F\Delta
\arrow[r, "F\sigma"']
&
F\Gamma
\end{tikzcd}
\end{center}
We say $F$ \textbf{strictly preserves canonical pullback squares} if, additionally, we have
\begin{center}
\begin{tikzcd}[row sep=huge, column sep=large]
F(\Delta \cext A[\sigma])
\arrow[r, "F(\sigma \cext A)"]
\arrow[d, "F\nmp{A\sigma}"']
&
F(\Gamma \cext A)
\arrow[d, "F\nmp{A}"{name=domeq}]
&
F\Delta \cextalt FA[F\sigma]
\arrow[r, "F\sigma \cextalt FA"]
\arrow[d, "\nmu{FA[F\sigma]}"'{name=codeq}]
\pullback
\arrow[from=domeq, to=codeq, draw=none, "=" description]
&
F\Gamma \cextalt FA
\arrow[d, "\nmu{FA}"]
\\
F\Delta
\arrow[r, "F\sigma"']
&
F\Gamma
&
F\Delta
\arrow[r, "F\sigma"']
&
F\Gamma
\end{tikzcd}
\end{center}
\end{definition}

\begin{lemma}
\label{lemCoverOfNaturalModel}
Let $(\mathbb{C}, p)$ be a natural model. The map $p : \nmtm{U} \to \nmty{U}$ is covered by
$$\sum_{(\Gamma,A) \in \int \nmty{U}} \nmp[\Gamma]{A} : \sum_{(\Gamma,A) \in \int \nmty{U}} \Yon (\Gamma \cext A) \longrightarrow \sum_{(\Gamma,A) \in \int \nmty{U}} \Yon(\Gamma)$$
\end{lemma}

\begin{proof}
Since $\nmty{U}$ is a colimit of representable presheaves, indexed by $\int_{\mathbb{C}} \nmty{U}$, we have a cover
$$\rho = \sqseqbn{A}{(\Gamma, A) \in \textstyle \int_{\mathbb{C}} \nmty{U}} : \sum_{(\Gamma,A) \in \int \nmty{U}} \Yon(\Gamma) \twoheadrightarrow \nmty{U}$$
Representability of $p$ yields the following pullback square
\begin{center}
\begin{tikzcd}[row sep=huge, column sep=huge]
\sum_{(\Gamma,A) \in \int \nmty{U}} \Yon (\Gamma \cext A)
\arrow[r, "\nmmk{\rho}"]
\arrow[d, "{\sum \Yon \nmp[\Gamma]{A}}"']
\pullbackc{dr}{-0.1}
&
\nmtm{U}
\arrow[d, "p"]
\\
\sum_{(\Gamma,A) \in \int \nmty{U}} \Yon(\Gamma)
\arrow[r, two heads, "\rho"']
&
\nmty{U}
\end{tikzcd}
\end{center}
where $\nmmk{\rho} = \sqseqbn{\nmq[\Gamma]{A}}{(\Gamma,A) \in \textstyle \int_{\mathbb{C}} \nmty{U}}$. Since $\rho$ is a regular epimorphism, so is $\nmmk{\rho}$.
\end{proof}

\begin{lemma}
\label{lemWeakMorphismIffPreservesCanonicalPullbacks}
Let $F = (F, \varphi, \nmmk{\varphi}) : (\mathbb{C}, p) \pmto (\mathbb{D}, q)$ be a premorphism of natural models. Then $F$ is a weak morphism of natural models if and only if $F$ weakly preserves canonical pullback squares.
\end{lemma}

\begin{proof}
Given $\sigma : \Delta \to \Gamma$ in $\mathbb{C}$ and $A \in \nmty{U}(\Gamma)$, consider the following diagram:

\begin{center}
\begin{tikzcd}[row sep=large, column sep=large]
F\Delta \cextalt FA[F\sigma]
\arrow[rrr, "F\sigma \cextalt FA"]
\arrow[dddr, bend right=15, "\nmp{FA[F\sigma]}"']
&&&
F\Gamma \cextalt FA
\arrow[dddl, bend left=15, "\nmp{FA}"]
\\
&
F(\Delta \cext A[\sigma])
\arrow[r, "F(\sigma \cext A)"]
\arrow[dd, "F\nmp{A[\sigma]}" description]
\arrow[ul, "\tau_{A[\sigma]}" description]
&
F(\Gamma \cext A)
\arrow[dd, "F\nmp{A}" description]
\arrow[ur, "\tau_A" description]
&
\\
~&~&~&~
\\
&
F\Delta
\arrow[r, "F\sigma"']
&
F\Gamma
&
\end{tikzcd}
\end{center}
The outer square is a pullback by representability of $q$, and the fact that $(F,\varphi,\nmmk{\varphi})$ is a premorphism of natural models yields the morphisms $\tau_A$ and $\tau_{A[\sigma]}$ as indicated and making the diagram commute.

Now $F$ is a weak morphism of natural models if and only if $\tau_A$ and $\tau_{A[\sigma]}$ are coherent isomorphisms. If they are coherent isomorphisms, then the inner square is a pullback, so that $F$ preserves display pullbacks.

Conversely, if $F$ preserves display pullbacks, then the inner and outer squares are both pullbacks. This implies that $\tau_A$ and $\tau_{A[\sigma]}$ are the canonical isomorphisms induced by the universal property of pullbacks, and hence that they satisfy the coherence laws necessary for $F$ to be a weak morphism of natural models.
\end{proof}

\begin{lemma}
\label{lemWeakMorphismImpliesPullbackSquare}
Let $F : (\mathbb{C}, p) \to (\mathbb{D}, q)$ be a weak morphism of natural models. Then the following square is a pullback.
\begin{diagram}
F_!\nmtm{U}
\arrow[r, "\nmmk{\varphi}"]
\arrow[d, "F_!p"']
\sqlabel{dr}{(\star)}
&
\nmtm{V}
\arrow[d, "q"]
\\
F_!\nmty{U}
\arrow[r, "\varphi"']
&
\nmty{V}
\end{diagram}
\end{lemma}

\begin{proof}
By \Cref{lemWeakMorphismIffPreservesCanonicalPullbacks} it suffices to assume that $F$ preserves canonical pullback squares. Since $F_!$ is a left adjoint, it preserves coproducts and regular epimorphisms, and furthermore $F_!$ commutes with the Yoneda embedding. Thus by applying $F_!$ to the pullback square in the proof of \Cref{lemCoverOfNaturalModel} and pasting $(\star)$ on the right, we obtain the following diagram:

\begin{center}
\begin{tikzcd}[row sep=huge, column sep=huge]
\sum_{(\Gamma,A) \in \int \nmty{U}} \Yon F(\Gamma \cext A)
\arrow[r, two heads, "F_!\nmmk{\rho}"]
\arrow[d, "{\sum \Yon F\nmp[\Gamma]{A}}"']
&
F_!\nmtm{U}
\arrow[r, "\nmmk{\varphi}"]
\arrow[d, "F_!p" description]
\sqlabel{dr}{(\star)}
&
\nmtm{V}
\arrow[d, "q"]
\\
\sum_{(\Gamma,A) \in \int \nmty{U}} \Yon(F\Gamma)
\arrow[r, two heads, "F_!\rho"']
&
F_!\nmty{U}
\arrow[r, "\varphi"']
&
\nmty{V}
\end{tikzcd}
\end{center}

The outer square is a pullback since $F$ preserves canonical pullback squares. Since $F_!\rho$ and $F_!\nmmk{\rho}$ are regular epimorphisms, it suffices to prove that the left-hand square is a pullback, for which, in turn, it suffices to check this on the individual components $(\Gamma,A)$ of the left-hand vertical coproduct.

To this end, fix $\Gamma \in \mathbb{C}$ and $A \in \nmty{U}(\Gamma)$. It suffices to check the universal property of pullbacks on representables, so let $D \in \mathbb{D}$ and let $d \in (F_! \nmtm{U})(D)$ and $\delta : D \to F\Gamma$ such that $F_!A \circ \Yon \delta = F_!p \circ d$.

\begin{center}
\begin{tikzcd}[row sep=normal, column sep=huge]
\Yon D
\arrow[drr, bend left=15, "d"]
\arrow[dddr, bend right=15, "\Yon \delta"']
&&
\\
&
\Yon F(\Gamma \cext A)
\arrow[r]
\arrow[dd, "\Yon F\nmp{A}" description]
&
F_! \nmtm{U}
\arrow[dd, "F_!p"]
\\
~&~&~
\\
&
\Yon F\Gamma
\arrow[r]
&
F_! \nmty{U}
\end{tikzcd}
\end{center}

Since $F_! \nmtm{U}$ is covered by representables indexed by $\int_{\mathbb{C}} \nmty{U}$ (\Cref{lemCoverOfNaturalModel}), there is some $\Delta \in \mathbb{C}$, $B \in \nmty{U}(\Delta)$ and $d' : D \to F(\Delta \cext B)$ such that $d$ factors through $\Yon d'$. Writing $\delta' = F \nmp{B} \circ d' : D \to F\Delta$, we obtain the following commutative diagram.

\begin{center}
\begin{tikzcd}[row sep=small, column sep=large]
\Yon D
\arrow[ddddrr, lightpurple, "\Yon \delta'"']
\arrow[ddddddr, lightpurple, bend right=10, "\Yon \delta"']
\arrow[drr, crossing over, "\Yon d'" description]
\arrow[ddrrrr, bend left, "d"]
&&&&
\\
&&
\Yon F(\Delta \cext B)
\arrow[drr]
\arrow[ddd, pos=0.3, "\Yon F\nmp{B}" description]
&&
\\
&&&&
F_! \nmtm{U}
\arrow[ddd, "F_!p"]
\\
&
\Yon F(\Gamma \cext A)
\arrow[urrr, crossing over]
\arrow[ddd, "\Yon F\nmp{A}" description]
&&&
\\
&&
\Yon F\Delta
\arrow[drr, "F_!B"]
&&
\\
&&&&
F_! \nmty{U}
\\
&
\Yon F\Gamma
\arrow[urrr, "F_!A"']
&&&
\end{tikzcd}
\end{center}

Now note that, in the category $\mathbf{Set}$, we have
$$(F_! \nmty{U})(D) \cong \left( \sum_{(\Gamma,A) \in \int \nmty{U}} \mathbb{D}(D, F\Gamma) \middle) \right/ {\sim}$$
where $\sim$ is the equivalence relation determined by the characterisation of colimits as coequalisers, and the component at $(\Gamma,A)$ of the quotient map $\nsum_{(\Gamma,A) \in \int \nmty{U}} \mathbb{D} (D, F\Gamma) \to (F_! \nmty{U})(D)$ is given by $(F_!A)_D : \mathbb{D}(D, F\Gamma) \to (F_! \nmty{U})(D)$.

Now identifying maps $\Yon D \to F_! \nmty{U}$ with elements of $(F_! \nmty{U})(D)$ (as usual, by the Yoneda lemma), we have
$$(F_!A)_D(\delta) = F_!A \circ \Yon \delta = F_!p \circ d = F_!B \circ \Yon F\nmp{B} \circ \Yon d' = F_!B \circ \Yon \delta' = (F_!B)_D(\delta')$$
and so $\delta$ and $\delta'$ are in the same $\sim$-equivalence class. It follows that there exists a zigzag of morphisms in $\int_{\mathbb{C}} \nmty{U}$ connecting $\delta$ with $\delta'$, in the sense that there exist:
\begin{itemize}
\item A natural number $n \ge 1$;
\item Objects $(\Theta_i, C_i) \in \int_{\mathbb{C}} \nmty{U}$ for all $0 \le i \le 2n$, with $(\Theta_0,C_0)=(\Gamma,A)$ and $(\Theta_{2n},C_{2n}) = (\Delta,B)$;
\item Morphisms $\Theta_{2i} \xleftarrow{\theta_{2i}} \Theta_{2i+1} \xrightarrow{\theta_{2i+1}} \Theta_{2i+2}$ for all $0 \le i < n$, such that $C_{2i+1}[\theta_{2i}] = C_{2i}$ and $C_{2i+1} [\theta_{2i+1}] = C_{2i+2}$ for all $i$; and
\item Morphisms $\delta_i : D \to F\Theta_i$ for all $0 \le i \le 2n$, with $\delta_0=\delta$ and $\delta_{2n} = \delta'$, such that $F\theta_{2i} \circ \delta_{2i+1} = \delta_{2i}$ and $F\theta_{2i+1} \circ \delta_{2i+1} = \delta_{2i+2}$ for all $0 \le i < n$;
\end{itemize}

By the assumption that $F$ preserves canonical pullback squares, and since the Yoneda embedding preserves limits, there is a zigzag of pullbacks connecting $\Yon F \nmp{A}$ with $\Yon F \nmp{B}$ as indicated by squiggly arrows in the following diagram.

\begin{center}
\begin{tikzcd}[row sep=small, column sep=large]
\Yon D
\arrow[ddddrr, lightpurple, "\Yon \delta'"']
\arrow[ddddddr, lightpurple, bend right=10, "\Yon \delta"']
\arrow[drr, crossing over, "\Yon d'" description]
\arrow[ddrrrr, bend left, "d"]
&&&&
\\
&&
\Yon F(\Delta \cext B)
\arrow[drr]
\arrow[ddd, pos=0.3, "\Yon F\nmp{B}" description]
\arrow[ddl, crossing over, leftrightsquigarrow]
&&
\\
&&&&
F_! \nmtm{U}
\arrow[ddd, "F_!p"]
\\
&
\Yon F(\Gamma \cext A)
\arrow[urrr, crossing over]
\arrow[ddd, "\Yon F\nmp{A}" description]
&&&
\\
&&
\Yon F\Delta
\arrow[drr, "F_!B"]
\arrow[ddl, leftrightsquigarrow]
&&
\\
&&&&
F_! \nmty{U}
\\
&
\Yon F\Gamma
\arrow[urrr, "F_!A"']
&&&
\end{tikzcd}
\end{center}

The detail of this zigzag of pullbacks is illustrated in the following commutative diagram.

\resizebox{\textwidth}{!}{%
\begin{tikzcd}[row sep={60pt}, column sep={10pt}, ampersand replacement=\&]
\&\&\&
\Yon D
\arrow[ddrrr, lightpurple, bend right=15, pos=0.65, "\Yon \delta_{2n}" description]
\arrow[ddrr, lightpurple, bend right=15, pos=0.8, "\Yon \delta_{2n-1}" description]
\arrow[ddr, lightpurple, bend right=15, pos=0.85, "\Yon \delta_{2n-2}" description]
\arrow[ddl, lightpurple,  bend left=15, pos=0.85, "\Yon \delta_{2}" description]
\arrow[ddll, lightpurple,  bend left=15, pos=0.8, "\Yon \delta_{1}" description]
\arrow[ddlll, lightpurple,  bend left=15, pos=0.65, "\Yon \delta_{0}" description]
\arrow[drrr, crossing over, "\Yon d'", bend left=15]
\&\&\&
\\
\Yon F(\Gamma \cext A)
\arrow[d, crossing over, "\Yon F\nmp{A}"']
\&
\Yon F(\Theta_1 \cext C_1)
\arrow[d, crossing over, "\Yon F\nmp{C_1}" description]
\arrow[l, crossing over]
\arrow[r, crossing over]
\pullback
\pullbackright
\&
\Yon F(\Theta_2 \cext C_2)
\arrow[d, crossing over, "\Yon F\nmp{C_2}" description]
\arrow[rr, crossing over, leftrightarrow, "\ \cdots\ " description]
\&
\&
\Yon F(\Theta_{2n-2} \cext C_{2n-2})
\arrow[d, crossing over, "\Yon F\nmp{C_{2n-2}}" description]
\&
\Yon F(\Theta_{2n-1} \cext C_{2n-1})
\arrow[d, crossing over, "\Yon F\nmp{C_{2n-1}}" description]
\arrow[l, crossing over]
\arrow[r, crossing over]
\pullback
\pullbackright
\&
\Yon F(\Delta \cext B)
\arrow[d, crossing over, "\Yon F\nmp{B}"]
\\[50pt]
\Yon F\Gamma
\&
\Yon F\Theta_1
\arrow[l]
\arrow[r]
\&
\Yon F\Theta_2
\arrow[rr, leftrightarrow, "\ \cdots\ " description]
\&\&
\Yon F\Theta_{2n-2}
\&
\Yon F\Theta_{2n-1}
\arrow[l]
\arrow[r]
\&
\Yon F\Delta
\\
\&\&\&
F_! \nmty{U}
\arrow[ulll, leftarrow, bend left=15, "F_!A"]
\arrow[ull, leftarrow, bend left=15, "F_!C_1" description]
\arrow[ul, leftarrow, bend left=15, "F_!C_2" description]
\arrow[ur, leftarrow, bend right=15, "F_!C_{2n-2}" description]
\arrow[urr, leftarrow, bend right=15, "F_!C_{2n-1}" description]
\arrow[urrr, leftarrow, bend right=15, "F_!B"']
\&\&\&
\end{tikzcd}
}


We proceed by induction. Let $d_{2n} = d' : D \to F(\Delta \cext B)$. The universal property of the rightmost pullback square yields a unique morphism $d_{2n-1} : D \to \Theta_{2n-1} \cext C_{2n-1}$ making the appropriate triangles commute, and composing with $F(\theta_{2n-1} \cext C_{2n-1})$, we obtain a morphism $d_{2n-2} : D \to F(\Theta_{2n-2} \to C_{2n-2})$ which satisfies $F \nmp{C_{2n-2}} \circ d_{2n-2} = \delta_{2n-2} : D \to F\Theta_{2n-2}$.

Repeating this process $n-1$ more times yields sequence of morphisms $d_i : D \to F(\Theta_i \cext C_i)$ making the required triangles commute; in particular, we obtain a morphism $d_0 : D \to F(\Gamma \cext A)$ satisfying $F \nmp{A} = \delta_0$ and $F_! \nmq{A} \circ \Yon d_0 = d$, as illustrated in the following diagram.

\resizebox{\textwidth}{!}{%
\begin{tikzcd}[row sep={60pt}, column sep={10pt}, ampersand replacement=\&]
\&\&\&
\Yon D
\arrow[ddrrr, lightpurple, bend right=15, pos=0.65, "\Yon \delta_{2n}" description]
\arrow[ddrr, lightpurple, bend right=15, pos=0.8, "\Yon \delta_{2n-1}" description]
\arrow[ddr, lightpurple, bend right=15, pos=0.85, "\Yon \delta_{2n-2}" description]
\arrow[ddl, lightpurple,  bend left=15, pos=0.85, "\Yon \delta_{2}" description]
\arrow[ddll, lightpurple,  bend left=15, pos=0.8, "\Yon \delta_{1}" description]
\arrow[ddlll, lightpurple,  bend left=15, pos=0.65, "\Yon \delta_{0}" description]
\arrow[drrr, crossing over, "\Yon d_{2n}", bend left=15]
\arrow[drr, crossing over, "\Yon d_{2n-1}" description, bend left=15]
\arrow[dr, crossing over, "\Yon d_{2n-2}" description, bend left=15]
\arrow[dl, crossing over, "\Yon d_2" description, bend right=15]
\arrow[dll, crossing over, "\Yon d_1" description, bend right=15]
\arrow[dlll, crossing over, dashed, blue, "\Yon d_0"', bend right=15]
\&\&\&
\\
\Yon F(\Gamma \cext A)
\arrow[d, crossing over, "\Yon F\nmp{A}"']
\&
\Yon F(\Theta_1 \cext C_1)
\arrow[d, crossing over, "\Yon F\nmp{C_1}" description]
\arrow[l, crossing over]
\arrow[r, crossing over]
\pullback
\pullbackright
\&
\Yon F(\Theta_2 \cext C_2)
\arrow[d, crossing over, "\Yon F\nmp{C_2}" description]
\arrow[rr, crossing over, leftrightarrow, "\ \cdots\ " description]
\&
\&
\Yon F(\Theta_{2n-2} \cext C_{2n-2})
\arrow[d, crossing over, "\Yon F\nmp{C_{2n-2}}" description]
\&
\Yon F(\Theta_{2n-1} \cext C_{2n-1})
\arrow[d, crossing over, "\Yon F\nmp{C_{2n-1}}" description]
\arrow[l, crossing over]
\arrow[r, crossing over]
\pullback
\pullbackright
\&
\Yon F(\Delta \cext B)
\arrow[d, crossing over, "\Yon F\nmp{B}"]
\\[50pt]
\Yon F\Gamma
\&
\Yon F\Theta_1
\arrow[l]
\arrow[r]
\&
\Yon F\Theta_2
\arrow[rr, leftrightarrow, "\ \cdots\ " description]
\&\&
\Yon F\Theta_{2n-2}
\&
\Yon F\Theta_{2n-1}
\arrow[l]
\arrow[r]
\&
\Yon F\Delta
\\
\&\&\&
F_! \nmty{U}
\arrow[ulll, leftarrow, bend left=15, "F_!A"]
\arrow[ull, leftarrow, bend left=15, "F_!C_1" description]
\arrow[ul, leftarrow, bend left=15, "F_!C_2" description]
\arrow[ur, leftarrow, bend right=15, "F_!C_{2n-2}" description]
\arrow[urr, leftarrow, bend right=15, "F_!C_{2n-1}" description]
\arrow[urrr, leftarrow, bend right=15, "F_!B"']
\&\&\&
\end{tikzcd}
}

To see that $d_0$ is the unique such morphism, suppose $\widehat{d}_0$ were another. Repeating the above process from left to right (rather than from right to left) yields a sequence of morphisms $\widehat{d}_i : D \to F(\Theta_i \cext C_i)$ induced from the fundamental property of pullbacks. It follows that $\widehat{d}_{2n} = d_{2n}$, and hence working again from right to left, we obtain $\widehat{d}_i = d_i$ for all $0 \le i \le 2n$. In particular, $\widehat{d}_0=d_0$, so that $(\star)$ is a pullback square.
\end{proof}

Proving or refuting the converse to \Cref{lemWeakMorphismImpliesPullbackSquare} is a topic for future work. It would be convenient if it were true---for example, it would imply that weak morphisms of natural models $(\mathbb{C}, p) \to (\mathbb{C}, q)$ over a common base category $\mathbb{C}$ correspond with cartesian morphisms of polynomials $p \pRightarrow q$ (\Cref{defCartesianMorphism}).

\subsection*{Preservation of type theoretic structure}

Having found a convenient way to describe morphisms of natural models, we now extend that description to natural models admitting extra structure. Again, we start with the essentially algebraic notion.

\begin{definition}[Morphisms of natural models preserving type theoretic structure]
\label{defMorphismOfNMPreservingStructure}
Let $S \subseteq \{ \seqbn{\mathsf{ty}_i}{i \in I}, \seqbn{\mathsf{tm}_j}{j \in J}, \mathbf{0}, \mathbbm{1}, \upSigma, \upPi, \dots \}$ be a set describing additional type theoretic structure that a natural model may possess (\Cref{secEssentiallyAlgebraicTheory}) and let $(\mathbb{C}, p)$ and $(\mathbb{D}, q)$ be natural models admitting the type theoretic structure in $S$. A \textbf{morphism of natural models preserving $S$} from $(\mathbb{C}, p)$ and $(\mathbb{D}, q)$ is a homomorphism from $(\mathbb{C}, p)$ to $(\mathbb{D}, q)$ considered as models of the essentially algebraic theory $\mathbb{T}_S$. Write $\mathbf{NM}_S$ for the category of natural models admitting structure from $S$.
\end{definition}

If $S' \subseteq S$ then there is a forgetful functor $\mathbf{NM}_S \to \mathbf{NM}_{S'}$. In particular, there is a forgetful functor $\mathbf{NM}_S \to \mathbf{NM}$ for any set $S$ of additional type theoretic structure that a natural model may possess. With this in mind, a morphism on $\mathbf{NM}_S$ has an underlying morphism of natural models $(F, \varphi, \nmmk{\varphi})$ in the sense of \Cref{thmMorphismsOfNMFunctorially}, so may be described as a triple $(F, \varphi, \nmmk{\varphi})$ satisfying some properties.

If the additional structure consists only of basic types, terms, an empty type or a unit type, then it is immediately clear that a morphism of natural models admitting this structure is simply one which maps the distinguished types and terms of its domain to the corresponding distinguished types and terms of its codomain.

It remains to describe preservation of dependent sum types and dependent product types. For this, we note the following two technical lemmas concerning polynomials and presheaves.

\begin{lemma}[Precomposition by a functor preserves polynomial composition and application]
Let $F : \mathbb{C} \to \mathbb{D}$ be a functor between small categories, and let $f : B \to A$ and $g : D \to C$ be morphisms in $\widehat{\mathbb{D}}$. Then there are isomorphisms in $\widehat{\mathbb{C}}^{\to}$
$$F^*(g \cdot f) \cong F^*(g) \cdot F^*(f) \quad \text{and} \quad F^*(\upP_g(f)) \cong \upP_{F^*(g)}(F^*(f))$$
\end{lemma}

\begin{proof}[Sketch of proof]
To see that $F^*(g \cdot f) \cong F^*(g) \cdot F^*(f)$, note that the functor $F^*$ preserves limits as it is a right adjoint, and the polynomial composite $g \cdot f$ is obtained by composing a particular pullback of $f$ with a particular pullback of $g$ (\Cref{secPolynomialPseudomonads}). To see that $F^*(\upP_g(f)) \cong \upP_{F^*(g)}(F^*(f))$, note that the following diagram commutes up to isomorphism.
\begin{diagram}
1
\arrow[r, "f"]
\arrow[dr, "F^*(f)"']
&
\widehat{\mathbb{D}}
\arrow[r, "\upDelta_D"]
\arrow[d, "F^*" description]
&
\widehat{\mathbb{D}} \slice{D}
\arrow[r, "\upPi_g"]
\arrow[d, "F^*" description]
&
\widehat{\mathbb{D}} \slice{C}
\arrow[r, "\upSigma_C"]
\arrow[d, "F^*" description]
&
\widehat{\mathbb{D}}
\arrow[d, "F^*" description]
\\
&
\widehat{\mathbb{C}}
\arrow[r, "\upDelta_{F^*(D)}"']
&
\widehat{\mathbb{C}} \slice{F^*(D)}
\arrow[r, "\upPi_{F^*(g)}"']
&
\widehat{\mathbb{C}} \slice{F^*(C)}
\arrow[r, "\upSigma_{F^*(C)}"']
&
\widehat{\mathbb{C}}
\end{diagram}
The composite from the top left to bottom right along the top gives $F^*(\upP_g(f))$, and the composite along the bottom gives $\upP_{F^*(g)}(F^*(f))$.
\end{proof}

\begin{lemma}[Polynomial composition preserves commutative squares]
Let $\mathcal{E}$ be a locally cartesian closed category. Polynomial composition extends to a functor
$$(-) \cdot (-) : \mathcal{E}^{\to} \times \mathcal{E}^{\to} \to \mathcal{E}^{\to}$$
\end{lemma}

\begin{proof}[Sketch of proof]
Given morphisms $f,f',g,g'$ in $\mathcal{E}$, chase commutative squares $f \to f'$ and $g \to g'$ through the construction of the polynomial composites $g \cdot f$ and $g' \cdot f'$ \Cref{defCompositionOfPolynomials} using the universal properties of pullbacks and of dependent products. Functoriality follows from uniqueness in the universal property.
\end{proof}

\begin{theorem}[Functorial description of preservation of dependent sum types]
Let $(\mathbb{C}, p)$ and $(\mathbb{D}, q)$ be natural models admitting dependent sum types. A morphism of natural models $F : (\mathbb{C}, p) \to (\mathbb{D}, q)$ preserves dependent sum types if and only if the following diagram in $\widehat{\mathbb{C}}^{\to}$ commutes.
\begin{diagram}
p \cdot p
\arrow[rr, "{(\widehat{\upSigma}, \widehat{\mathsf{pair}})}"]
\arrow[d, "{(\varphi, \nmmk{\varphi})} \cdot {(\varphi, \nmmk{\varphi})}"']
&&
p
\arrow[d, "{(\varphi, \nmmk{\varphi})}"]
\\
F^*q \cdot F^*q
\arrow[r, "\cong"']
&
F^*(q \cdot q)
\arrow[r, "{F^*(\widehat{\upSigma}, \widehat{\mathsf{pair}})}"']
&
F^*q
\end{diagram}
\end{theorem}

\begin{proof}
The morphism $(\varphi, \nmmk{\varphi}) \cdot (\varphi, \nmmk{\varphi})$ can be expressed in the internal language of $\widehat{\mathbb{C}}$ via the following commutative square
\begin{diagram}
\sum_{A \in \nmty{U}} \sum_{B \in \nmty{U}^{[A]}} \sum_{a \in [A]} [B(a)]
\arrow[d, "p \cdot p"']
\arrow[r, "{((\varphi, \nmmk{\varphi}) \cdot (\varphi, \nmmk{\varphi}))_1}"]
&
\sum_{A' \in F^*\nmty{V}} \sum_{B' \in (F^*\nmty{V})^{\langle A'}} \sum_{a' \in \langle A' \rangle} \langle B(a) \rangle
\arrow[d, "F^*q \cdot F^*q"]
\\
\sum_{A \in \nmty{U}} \nmty{U}^{[A]}
\arrow[r, "{((\varphi, \nmmk{\varphi}) \cdot (\varphi, \nmmk{\varphi}))_0}"']
&
\sum_{A' \in \nmty{V}(F\Gamma)} \nmty{V}(F\Gamma \cext A')
\end{diagram}
The component at $\Gamma$ of this diagram, after application of \Cref{lemLemmaFive,lemLemmaElevenPointFive}, is given by the following commutative square of sets and functions.
\begin{diagram}
\sum_{A \in \nmty{U}(\Gamma)} \sum_{B \in \nmty{U}(\Gamma \cext A)} \sum_{a \in \nmtm{U}(\Gamma; A)} \nmtm{U}(\Gamma; B(a))
\arrow[d, "\pi"']
\arrow[r]
&
\sum_{A' \in \nmty{V}(F\Gamma)} \sum_{B' \in \nmty{V}(F\Gamma \cextalt A')} \sum_{a' \in \nmtmalt{V}(F\Gamma; A')} \nmtmalt{V}(F\Gamma; B'(a'))
\arrow[d, "\pi"]
\\
\sum_{A \in \nmty{U}(\Gamma)} \nmty{U}(\Gamma \cext A)
\arrow[r]
&
\sum_{A' \in F^*\nmty{V}} (F^*\nmty{V})^{\langle A \rangle}
\end{diagram}
The morphisms at the top and bottom apply $F$ to each component.

As such, the original square commutes if and only if for all $\Gamma \in \mathrm{ob}(\mathbb{C})$, $A \in \nmty{U}(\Gamma)$, $B \in \nmty{U}(\Gamma \cext A)$, $a \in \nmtm{U}(\Gamma; A)$ and $b \in \nmtm{U}(\Gamma; B(a))$, we have
$$F\widehat{\upSigma}(A,B) = \widehat{\upSigma}(FA,FB) \quad \text{and} \quad F\widehat{\mathsf{pair}}(A,B,a,b) = \widehat{\mathsf{pair}}(FA,FB,Fa,Fb)$$

This says that the square in the statement of the theorem commutes if and only if $F$ preserves dependent sum types, as required.
\end{proof}

%% file: thesis/ch3-polynomials-representability/_polynomials-representability.tex
\newpage
\input{thesis/ch3-polynomials-representability/internal-categories.tex}

\newpage
\input{thesis/ch3-polynomials-representability/polynomial-pseudomonads.tex}

\newpage
\input{thesis/ch3-polynomials-representability/representability.tex}

%% file: thesis/ch3-polynomials-representability/internal-categories.tex
\section{Internal categories}
\label{secInternalCategories}
\index{category!internal {---}}
\index{internal!category@{---} category}

\begin{construction}[Full internal subcategories \thmcite{Jacobs1999CategoricalLogic}]
\label{cnsFullInternalSubcategory}
\index{full internal subcategory}
Let $f : B \to A$ be a morphism in a locally cartesian closed category $\mathcal{E}$. The \textbf{full internal subcategory} of $\mathcal{E}$ induced by $f$ is the internal category $\fisc{f}$ of $\mathcal{E}$ defined by
\begin{itemize}
\item The object of objects $\fisc{f}_0$ is simply $A$;
\item The object of morphisms $\fisc{f}_1$ together with the pair $\partial = \langle \partial_0, \partial_1 \rangle : \fisc{f}_1 \to A \times A$ is given as an object of $\mathcal{E} \slice{A \times A}$ by taking the exponential $f_2^{f_1}$, where $f_1=\upDelta_{\pi_1}(f) : B_1 \to A \times A$ and $f_2 = \upDelta_{\pi_2}(f) : B_2 \to A \times A$ are the pullbacks of $f$ along the projections $A \xleftarrow{\pi_1} A \times A \xrightarrow{\pi_2} A$.
\begin{diagram}
B
\arrow[d, "f"']
&
B_1
\arrow[l]
\arrow[dr, "f_1"']
\pullbackright
&[-50pt]
&[-50pt]
B_2
\arrow[r]
\arrow[dl, "f_2"]
\pullback
&
B
\arrow[d, "f"]
\\
A
&&
A \times A
\arrow[ll, "\pi_1"]
\arrow[rr, "\pi_2"']
&&
A
\end{diagram}
In the internal language of $\mathcal{E}$, we have $\fisc{f}_1 = \sum_{a,a' \in A} B_{a'}^{B_a}$;
\item The identities morphism $A \to \sum_{a,a' \in A} B_{a'}^{B_a}$ is given by $a \mapsto \langle a, a, \mathrm{id}_{B_a} \rangle$;
\item The composition morphism
$$\sum_{a,a',a'' \in A} B_{a''}^{B_{a'}} \times B_{a'}^{B_a} \to \sum_{a,a'' \in A} B_{a''}^{B_a}$$
is given by internal composition in $\mathcal{E}$ in the expected way.
\end{itemize}
\end{construction}

\begin{numbered}
More generally, a \textit{full internal subcategory} of a locally cartesian closed category $\mathcal{E}$ is an internal category $\mathbb{S}$ of $\mathcal{E}$ for which there is a full and faithful fibred functor over $\mathcal{E}$ from the externalisation of $\mathbb{C}$ to $\mathcal{E}^{\to}$. Every full internal subcategory is isomorphic to one of the form $\fisc{f}$. This matter is discussed in great detail in \cite{Jacobs1999CategoricalLogic}, and we will explore the externalisation of full internal subcategories of $\widehat{\mathbb{C}}$ in \Cref{secCharacterisationsOfRepresentability}.
\end{numbered}

\begin{definition}[Associated full internal subcategory of a natural model]
\index{full internal subcategory!{---} associated with a natural model}
Let $(\mathbb{C}, p)$ be a natural model. Its \textbf{associated full internal subcategory} is the full internal subcategory $\mathbb{U} = \fisc{p}$ of $\widehat{\mathbb{C}}$ induced by $p$. Explicitly, we have $\mathbb{U}_0 = \nmty{U}$ and $\mathbb{U}_1 = \sum_{A, B \in \nmty{U}} [B]^{[A]}$, so that we can think of $\mathbb{U}$ as a category whose objects are types and whose morphisms are functions between types.
\end{definition}

In the same vein as \Cref{lemLemmaFive} and \Cref{lemLemmaElevenPointFive}, the following lemma will allow us to work more concretely with full internal subcategories.

\begin{lemma}
\label{lemMorphismsOfAssociatedFISC}
Let $f : B \to A$ be a morphism in a locally cartesian closed category $\mathcal{E}$. Morphisms $x : X \to \sum_{a,a' \in A} B_{a'}^{B_a}$ in $\mathcal{E}$ correspond naturally with triples $(x_1,x_2,\widetilde{x})$, where $x_1, x_2 : X \to A$ in $\mathcal{E}$ and $\widetilde{x} : \upDelta_{x_1}(f) \to \upDelta_{x_2}(f)$ in $\mathcal{E} \slice{X}$.
\end{lemma}

\begin{proof}
First note that a morphism $x : X \to \sum_{a,a' \in A} B_{a'}^{B_a}$ gives rise to a morphism $(x_1,x_2) \to f_2^{f_1}$ in $\mathcal{E} \slice{A \times A}$, where $x_1$ and $x_2$ are the composites of $x$ with the respective projections $\sum_{a,a' \in A} B_{a'}^{B_a} \to A$.
\begin{diagram}
X
\arrow[rr, "x"]
\arrow[dr, "{(x_1,x_2)}"']
&&
\sum_{a,a' \in A} B_{a'}^{B_a}
\arrow[dl, "f_2^{f_1}"]
\\
&
A \times A
\end{diagram}
Now $f_2^{f_1} = \upPi_{f_1} \upDelta_{f_1}(f_2)$, so under the adjunction $\upDelta_{f_1} \dashv \upPi_{f_1}$, a morphism $x : (x_1,x_2) \to f_2^{f_1}$ corresponds with a morphism $x' : \upDelta_{f_1}(x_1,x_2) \to \upDelta_{f_1}(f_2)$ in $\mathcal{E} \slice{B_1}$. Hence it suffices to show that morphisms $x' : \upDelta_{f_1}(x_1,x_2) \to \upDelta_{f_1}(f_2)$ in $\mathcal{E} \slice{B_2}$ correspond with morphisms $\widetilde{x} : \upDelta_{x_1}(f) \to \upDelta_{x_2}(f)$ in $\mathcal{E} \slice{X}$.

Consider now following diagram, in which all three squares are pullbacks---the fact that the bottom and front squares are pullbacks follows from the fact that $\upDelta_f(x_i) \cong \upDelta_{f_i}(x_1,x_2)$ for $i=1,2$. Our goal is to prove that there is a correspondence between morphisms $x'$ and morphisms $\widetilde{x}$ fitting into the diagram as indicated with dashed arrows.
\begin{ldiagram}
&&&
B_3
\arrow[dr, "\upDelta_{f_2}(f_1)"]
\arrow[dd, pos=0.7, "\upDelta_{f_1}(f_2)" description]
&&
\\
&&
\upDelta_{x_2}(B)
\arrow[rr, pos=0.3, crossing over, "\upDelta_f(x_2)"]
&&
B_2
\arrow[dd, "f_2"]
&
\\
&
\upDelta_{x_1}(B)
\arrow[dr, "\upDelta_{x_1}(f)"']
\arrow[rr, pos=0.3, "\upDelta_f(x_1)" description]
\arrow[ur, dashed, blue, bend left=15, "\widetilde{x}" description]
\arrow[uurr, dashed, red, bend left=30, "x'"]
&&
B_1
\arrow[dr, "f_1" description]
&&
\\
&&
X
\arrow[uu, leftarrow, crossing over, pos=0.7, "\upDelta_{x_2}(f)" description]
\arrow[rr, "{(x_1,x_2)}"']
&~&
A \times A
&
\end{ldiagram}

Given $\widetilde{x}$, we obtain $x'$ from the universal property of the right pullback square as the pair $x' = (\upDelta_f(x_2) \circ \widetilde{x}, \upDelta_f(x_1))$. Conversely, given $x'$, we obtain $\widetilde{x}$ from the universal property of the front pullback as the pair $\widetilde{x} = (\upDelta_{f_2}(f_1) \circ x', \upDelta_{x_1}(f))$. That the assignments $\widetilde{x} \mapsto x'$ and $x' \mapsto \widetilde{x}$ are mutually inverse follows immediately from the `uniqueness' part of the universal property.
\end{proof}

\begin{numbered}
As a result of \Cref{lemMorphismsOfAssociatedFISC}, when reasoning internally about a full internal subcategory $\fisc{f : B \to A}$ of a locally cartesian closed category $\mathcal{E}$, we can view the object of objects $\fisc{f}_0$ as an $A$-indexed family $\seqbn{a}{a \in A}$ and the object of morphisms $\fisc{f}_1$ as an $A \times A$-indexed family $\seqbn{f : B_a \to B_{a'}}{a,a' \in A}$. This allows us to reason internally to $\mathcal{E}$ about full internal subcategories, as well as internal functors and internal natural transformations between them, much like we reason about categories, functors and natural transformations externally.
\end{numbered}

\begin{lemma}
\label{lemMorphismsOfAssociatedFISCNM}
Let $(\mathbb{C}, p)$ be a natural model with associated full internal subcategory $\mathbb{U}$. For each object $\Gamma$ of $\mathbb{C}$, there is a bijection
$$\mathbb{U}_1(\Gamma) \overset{\cong}{\longrightarrow} \sum_{A, B \in \nmty{U}(\Gamma)} \mathbb{C} \slice{\Gamma} \left( (\Gamma \cext A, \nmp{A}), (\Gamma \cext B, \nmp{B}) \right)$$
which is natural in $\Gamma$.

Explicitly, given $g \in \mathbb{U}_1(\Gamma)$ and $\sigma : \upDelta \to \Gamma$ in $\mathbb{C}$, if $g$ corresponds with a triple $(A,B,h)$, where $A,B \in \nmty{U}(\Gamma)$ and $h : \Gamma \cext A \to \Gamma \cext B$ over $\Gamma$, then $g[\sigma]$ corresponds with $(A[\sigma], B[\sigma], h[\sigma])$, where $h[\sigma] : \upDelta \cext A[\sigma] \to \upDelta \cext B[\sigma]$ over $\upDelta$, as indicated with a dashed arrow in the following diagram in which the front and back squares are canonical pullback squares (\Cref{cnsCanonicalPullbacks}).

\begin{cdiagram}{normal}{large}
&&[15pt]
\upDelta \cext B[\sigma]
\arrow[dddl, pos=0.6, "\nmp{B[\sigma]}"]
\arrow[rrr, "\sigma \cext B"]
&&&[15pt]
\Gamma \cext B
\arrow[dddl, "\nmp{B}"]
\\
\upDelta \cext A[\sigma]
\arrow[ddr, "\nmp{A[\sigma]}"']
\arrow[urr, dashed, "{h[\sigma]}"]
\arrow[rrr, crossing over, pos=0.7, "\sigma \cext A"]
&&&
\Gamma \cext A
\arrow[ddr, "\nmp{A}"']
\arrow[urr, "h"]
&&
\\
~&~&~&~&~&~
\\
&
\upDelta
\arrow[rrr, "\sigma"']
&&&
\Gamma
&
\end{cdiagram}
\end{lemma}

\begin{proof}
Apply \Cref{lemMorphismsOfAssociatedFISC} with $f = p$ and $X = \Yon(\Gamma)$.
\end{proof}

\begin{numbered}
\Cref{lemMorphismsOfAssociatedFISCNM} proves that the associated full internal subcategory of a natural model is equivalent to \textit{context-indexed family of types} \cite[Proposition 1]{ClairambaultDybjer2014Biequivalence}, the latter regarded as a $\mathbb{C}$-indexed category rather than a category internal to $\widehat{\mathbb{C}}$, although these notions are equivalent.
\end{numbered}

\begin{construction}[Cartesian morphisms of polynomials induce full and faithful internal functors]
\label{cnsCartesianMorphismsInduceInternalFunctors}
Let $f : B \to A$ and $g : D \to C$ be morphisms in a locally cartesian closed category $\mathcal{E}$. Given a cartesian morphism of polynomials $\varphi : f \pRightarrow g$, let $\fisc{\varphi} : \fisc{f} \to \fisc{g}$ be the full and faithful internal functor defined in the internal language of $\mathcal{E}$ as follows.
\begin{itemize}
\item $\fisc{\varphi}_0(a) = \varphi_0(a)$ for $a \in A$; and
\item $\fisc{\varphi}_1(m : B_a \to B_{a'}) = \varphi_{a'} \circdot k \circdot \varphi_a^{-1})$ for $a,a' \in A$;
\end{itemize}
The assignment $\varphi \mapsto \fisc{\varphi}$ extends to a functor $\mathbb{S} : \mathbf{Poly}_{\mathcal{E}}\cart(1,1) \to \mathbf{Cat}(\mathcal{E})$.
\end{construction}

\begin{verification}
We work internally to $\mathcal{E}$. To see that $\mathbb{S}(\varphi)$ defines a functor, note that for $a,a',a'' \in A$ and $m : B_a \to B_{a'}$ and $n : B_{a'} \to B_{a''}$, we have
$$\fisc{\varphi}(n) \circ \fisc{\varphi}(m) = \varphi_{a''} \circ n \circ \varphi_{a'}^{-1} \circ \varphi_{a'} \circ m \circ \varphi_a = \varphi_{a''} \circ n \circ m \circ \varphi_a^{-1} = \fisc{\varphi}(n \circ m)$$
and evidently $\fisc{\varphi}(\mathrm{id}_{B_a}) = \mathrm{id}_{D_{\varphi_0(a)}}$.

To see that the assignment $\varphi \mapsto \fisc{\varphi}$ is functorial, note that evidently $\fisc{\mathrm{id}_f} = \mathrm{id}_{\fisc{f}}$ for each $f : B \to A$ in $\mathcal{E}$, and given $\varphi : f \pRightarrow g$ and $\psi : g \pRightarrow h$, we have
$$\fisc{\psi \circ \varphi}_0 = (\psi \circ \varphi)_0 = \psi_0 \circ \varphi_0 = \fisc{\psi}_0 \circ \fisc{\varphi}_0$$
and for $m : B_a \to B_{a'}$ we have
$$\fisc{\psi \circ \varphi}_1(m) = (\psi \circ \varphi)_{a'} \circ m \circ (\psi \circ \varphi)_a^{-1} = \psi_{\varphi_0(a')} \circ \varphi_{a'} \circ m \circ \varphi_a^{-1} \circ \psi_{\varphi_0(a)}^{-1} = \fisc{\psi}_1(\fisc{\varphi_1}(m))$$
as required.
\end{verification}

In the following, given an object $A$ of a locally cartesian closed category $\mathcal{E}$, we will write $|A|$ for the discrete internal category on $\mathcal{E}$, whose object of objects and of morphisms are both $A$ and with domain, codomain, identities and composition morphisms given by $\mathrm{id}_A$.

\begin{definition}[Internally cartesian closed categories \thmcite{Jacobs1999CategoricalLogic}]
\index{internal!cartesian closed@{---}ly cartesian closed category}
Let $\mathbb{A}$ be an internal category of a locally cartesian closed category $\mathcal{E}$.
\begin{itemize}
\item $\mathbb{A}$ has an \textbf{internal terminal object} if the terminal internal functor ${!} : \mathbb{A} \to |1|$ has an internal right adjoint $\mathsf{tobj} : |1| \to \mathbb{A}$.
\item $\mathbb{A}$ has \textbf{internal binary products} if the diagonal internal functor $\upDelta : \mathbb{A} \to \mathbb{A} \times \mathbb{A}$ has an internal right adjoint.
\item $\mathbb{A}$ with internal binary products has \textbf{internal exponentials} if the internal functor
$$\mathsf{prod} : |\mathbb{A}_0| \times \mathbb{A} \to |\mathbb{A}_0| \times \mathbb{A}$$
has a right adjoint $\mathsf{exp}$, where $\mathsf{prod}$ is defined as follows and where $\timesdot : \mathbb{A} \times \mathbb{A} \to \mathbb{A}$ is the internal binary product functor.
\begin{itemize}
\item $\mathsf{prod}_0 : \mathbb{A}_0 \times \mathbb{A}_0 \to \mathbb{A}_0 \times \mathbb{A}_0$ is defined by $\mathsf{prod}_0(A,B) = (A, A \timesdot B)$; and
\item $\mathsf{prod}_1 : \mathbb{A}_0 \times \mathbb{A}_1 \to \mathbb{A}_0 \times \mathbb{A}_1$ is defined by $\mathsf{prod}_1(A,f) = (A, \mathrm{id}_A \timesdot f)$;
\end{itemize}
\end{itemize}
We say $\mathbb{A}$ is \textbf{internally cartesian closed} if it has an internal terminal object, internal binary products and internal exponentials.
\end{definition}

\begin{theorem}[Cartesian closure of the associated full internal subcategory]
\label{thmCartesianClosedInternalCategory}
Let $(\mathbb{C}, p)$ be a natural model and let $\mathbb{U} = \fisc{p}$ be its associated full internal subcategory. If $(\mathbb{C}, p)$ admits a unit type, dependent sum types and dependent product types, then $\mathbb{U}$ is internally cartesian closed.
\end{theorem}

\begin{proof}
Let $\eta = (\widehat{\mathbbm{1}}, \widehat{\star}) : i_1 \pRightarrow p$, $\mu = (\widehat{\upSigma}, \widehat{\mathsf{pair}}) : p \cdot p \pRightarrow p$ and $\zeta = (\widehat{\upPi}, \widehat{\lambda}) : \upP_p(p) \pRightarrow p$ be the cartesian morphisms of polynomials arising from the unit, dependent sum and dependent product type structure for $(\mathbb{C}, p)$, as in \Cref{thmUnitSigmaPiPoly}.

To see that $\mathbb{U}$ has an internal terminal object, let $\mathsf{tobj} = \fisc{\eta} : \fisc{\Yon(\diamond)} = |1| \to \fisc{p} = \mathbb{U}$. Now for $A \in \nmty{U}$ and $x \in 1$, there is a natural correspondence between morphisms
$$1_{!A} \to 1_{x} \qquad \text{and} \qquad [A] \to [\mathsf{tobj}(x)]$$
Indeed, the only morphism $1_{!A} \to 1_{x}$ in $|1|$ is $\mathrm{id}_1$, and the only morphism $[A] \to [\widehat{\mathbbm{1}}]$ is $\lambda a.~\widehat{\star}$. So ${!} \dashv \mathsf{tobj}$ as required.

To see that $\mathbb{U}$ has internal products, we will define an internal functor
$$\timesdot : \mathbb{U} \times \mathbb{U} \to \mathbb{U}$$
on objects by $A \timesdot B = \widehat{\upSigma}(A, \underline{B})$, where $\underline{B} = \lambda x.B \in \nmty{U}^{[A]}$, and on morphisms $(f,g) \in [C]^{[A]} \times [D]^{[B]}$ by letting $f \timesdot g : [A \timesdot B] \to [C \timesdot D]$ be defined by
$$(f \timesdot g)(\langle a, b \rangle) = \widehat{\mathsf{pair}}(f(a), g(b))$$
To see that $\upDelta \dashv {\timesdot}$, note that for $A,B,C \in \nmty{U}$ there is a map
$$[A]^{[C]} \times [B]^{[C]} \to [A \timesdot B]^{[C]}$$
given by $(f, g) \mapsto \lambda c.~ \widehat{\mathsf{pair}}(f(c), g(c))$ and with an inverse given by
$$h \mapsto (\lambda c.\mathsf{fst}(h(c)), ~ \lambda c.\mathsf{snd}(h(c)))$$
Hence $\upDelta \dashv {\timesdot}$ as required.

To see that $\mathbb{U}$ has internal exponentials, define an internal functor
$$\mathsf{exp} : |\nmty{U}| \times \mathbb{U} \to |\nmty{U}| \times \mathbb{U}$$
on objects by $\mathsf{exp}_0(A,B) = (A, \widehat{\upPi}(A, \underline{B}))$ and on morphisms $(A, B) \xrightarrow{(\mathrm{id}_A, f)} (A, C)$ by letting
$$\mathsf{exp}(f) : \widehat{\upPi}(A,\underline B) \to \widehat{\upPi}(A, \underline C)$$
be defined by $\mathsf{exp}(f)(t) = \widehat{\lambda}(A, t \circ f)$. To see that $\mathsf{prod} \dashv \mathsf{exp}$, note that for $A,B,C \in \nmty{U}$ there is a map
$$|\nmty{U}|(A,C) \times [D]^{[A \timesdot B]} \to |\nmty{U}|(A,C) \times [\widehat{\upPi}_C D]$$
which is trivial when $A \ne C$ and is defined by the usual currying and uncurrying correspondence when $A=C$. This proves that $\mathsf{prod} \dashv \mathsf{exp}$.
\end{proof}

We will use the following construction in \Cref{thmInternalRightAdjoint} in order to characterise when a natural model admits dependent sum types.

\begin{construction}
Given full internal subcategories $\mathbb{A} = \fisc{B \xrightarrow{f} A}$ and $\mathbb{I} = \fisc{J \xrightarrow{g} I}$ of a locally cartesian closed category $\mathcal{E}$, there is an internal category $\Fam[\mathbb{I}]{\mathbb{A}}$ of $\mathcal{E}$ defined as follows.
\begin{itemize}
\item $\Fam[\mathbb{I}]{\mathbb{A}}_0 = \sum_{i \in I} A^{J_i}$ ($= \upP_{\alpha}(A)$);
\item $\Fam[\mathbb{I}]{\mathbb{A}}_1 = \sum_{(i,a),(i',a')} \sum_{\alpha \in J_{i'}^{J_i}} \prod_{j \in J_i} B_{a'(\alpha(j))}^{B_{a(j)}}$;
\item The domain and codomain morphisms $\mathsf{dom},\mathsf{cod} : \Fam[\mathbb{I}]{\mathbb{A}}_1 \to \Fam[\mathbb{I}]{\mathbb{A}}_0$ are given by the evident projections;
\item The identities morphism $\mathsf{ids} : \Fam[\mathbb{I}]{\mathbb{A}}_0 \to \Fam[\mathbb{I}]{\mathbb{A}}_1$ is defined in the internal language of $\mathcal{E}$ by
$$\mathsf{ids}(i,a) = \Big( (i,a), ~ (i,a), ~ \lambda j.j, ~ \lambda j.\lambda b.b \Big)$$
\item Composition is given by internal composition in $\mathcal{E}$; explicitly, the object of composable pairs of morphisms is given by
$$\Fam[\mathbb{I}]{\mathbb{A}}_2 = \sum_{(i,a),(i',a'),(i'',a'')} ~ \sum_{\beta : J_{i''}^{J_{i'}}} ~ \sum_{\alpha : J_{i'}^{J_i}} ~ \prod_{j \in J_i} ~ B_{a''(\beta(\alpha(j)))}^{B_{a'(\alpha(j))}} \times B_{a'(\alpha(j))}^{B_{a(j)}}$$
and the composition morphism $\mathsf{comp} : \Fam[\mathbb{I}]{\mathbb{A}}_2 \to \Fam[\mathbb{I}]{\mathbb{A}}_1$ is given by
$$\mathsf{comp} \bigg((i,a),(i',a'),(i'',a''), \beta, \alpha, \lambda j.(b'_j,b_j) \Big) = \Big((i,a), (i'',a''), \beta \circdot \alpha, \lambda j.(b'_j \circdot b_j) \Big)$$
\end{itemize}
\end{construction}

\begin{verification}
The fact that $\mathbb{A}$ and $\mathbb{I}$ allow us to check the required equations using the internal language of $\mathcal{E}$; but these equations are exactly the ones that demonstrate that the regular $\mathrm{Fam}$ construction defines a category (see e.g.\
\cite{Jacobs1999CategoricalLogic}).
\end{verification}

\begin{construction}
Let $\mathbb{A}$ be an full internal subcategory with an internal terminal object $\mathbf{1} \in A$. By anology with the diagonal functor $\upDelta : \mathbb{A} \to \mathbb{A} \times \mathbb{A}$, define an internal functor $\widetilde{\upDelta} : \mathbb{A} \to \Fam[\mathbb{A}]{\mathbb{A}}$ defined internally on objects by $\widetilde{\upDelta}(a) = (a, \underline{\mathbf{1}})$ and on morphisms $f : B_{a} \to B_{a'}$ by $\widetilde{\upDelta}(f) = (f : B_a \to B_{a'}, \mathrm{id}_{B_{\mathbf{1}}} : B_{\mathbf{1}} \to B_{\mathbf{1}})$. \qed
\end{construction}

\begin{theorem}
\label{thmInternalRightAdjoint}
Let $(\mathbb{C}, p)$ be a natural model admitting a unit type and let $\mathbb{U} = \fisc{p}$ be its associated full internal subcategory. Then $(\mathbb{C}, p)$ admits dependent sum types if and only if the internal functor $\widetilde{\upDelta} : \mathbb{U} \to \Fam[\mathbb{U}]{\mathbb{U}}$ has an internal right adjoint.
\end{theorem}

\begin{proof}
First recall (\Cref{thmUnitSigmaPiPoly}) that a natural model $(\mathbb{C}, p)$ admits dependent sum types if and only if there exists a cartesian morphism $(\widehat{\mathsf{pair}}, \widehat{\upSigma}) : p \cdot p \pRightarrow p$ of polynomials in $\widehat{\mathbb{C}}$, as indicated in the following pullback square:
\begin{center}
\begin{tikzcd}[row sep=huge, column sep=huge]
\sum_{A,B} \sum_{a:[A]} [B(a)]
\arrow[r, "\widehat{\mathsf{pair}}"]
\arrow[d, "p \cdot p"']
\pullbackc{dr}{-0.05}
&
\nmtm{U}
\arrow[d, "p"]
\\
\sum_{A:\nmty{U}} \nmty{U}^{[A]}
\arrow[r, "\widehat{\upSigma}"']
&
\nmty{U}
\end{tikzcd}
\end{center}

Write $\widehat{\mathsf{pair}}(A,B,a,b) = \langle a,b \rangle$ and $\widehat{\upSigma}(A,B) = \upSigma_A B$.

First suppose that $(\mathbb{C}, p)$ admits dependent sum types, and define $\mathsf{sigma} : \Fam[\mathbb{U}]{\mathbb{U}} \to \mathbb{U}$ as follows:
\begin{itemize}
\item $\mathsf{sigma}_0 = \widehat{\upSigma} : \sum_{A : \nmty{U}} \nmty{U}^{[A]} \to \nmty{U}$, so that $\mathsf{sigma}_0(A,B) = \upSigma_A B$ for each $(A,B) \in (\Fam[\mathbb{U}]{\mathbb{U}})_0$.
\item For $(A,B),(C,D) \in \sum_{A:\nmty{U}} \nmty{U}^{[A]}$ and $(f,g) \in \Fam[\mathbb{U}]{\mathbb{U}}((A,B),(C,D))$, define
$$\mathsf{sigma}_1(f,g) = \lambda p. \langle f(p.0), g_{p.0}(p.1) \rangle : [\upSigma_A B] \to [\upSigma_C D]$$
\end{itemize}

First note that $\mathsf{sigma}$ is an internal functor. That it respects identities is evident; to see that it respects composition, note that
\begin{align*}
\mathsf{sigma}_1((f',g') \circ (f,g))
&= \mathsf{sigma_1}(f' \circ f, \lambda a.~ g'_{f(a)} \circ g_a) && \\
&= \lambda p.~ \langle f'(f(p.0)), g'_{f(p.0)}(g_{p.0}(p.1)) \rangle && \\
&= \lambda p.~ \mathsf{sigma}_1(f',g')(\langle f(p.0), g_{p.0}(p.1) \rangle) && \\
&= \lambda p.~ \mathsf{sigma}_1(f',g')(\mathsf{sigma}_1(f,g)(p)) && \\
&= \mathsf{sigma}_1(f',g') \circ \mathsf{sigma}_1(f,g)
\end{align*}

To see that $\widetilde{\upDelta} \dashv \mathsf{sigma}$, let $C \in \nmty{U}$ and let $(A,B) \in \sum_{A \in \nmty{U}} \nmty{U}^{[A]}$. We obtain a map
$$\Fam[\mathbb{U}]{\mathbb{U}}((C,\underline{\mathbf{1}}), (A,B)) \to \mathbb{U}(\upSigma_C \underline{\mathbf{1}}, \upSigma_A B)$$
via $(f,g) \mapsto \lambda \langle c, \star \rangle.~ (f(c), g_c(\star))$; and we obtain a map
$$\mathbb{U}(\upSigma_C \underline{\mathbf{1}}, \upSigma_A B) \to \Fam[\mathbb{U}]{\mathbb{U}}((C,\underline{\mathbf{1}}), (A,B))$$
via $h \mapsto (\lambda c.~ h(c).0, ~ \lambda c.~ \lambda x : \mathbf{1}.~ h(c).1 \rangle)$.

These maps are mutually inverse, and so $\widetilde{\upDelta} \dashv \mathsf{sigma}$ as required.

Conversely, suppose now that $\upDelta$ has an internal right adjoint $\mathsf{sigma} : \Fam[\mathbb{U}]{\mathbb{U}} \to \mathbb{U}$. Define $\widehat{\upSigma} : \sum_{A : \nmty{U}} \nmty{U}^{[A]} \to \nmty{U}$. Then $\mathsf{sigma}_1$ gives for each $(A,B),(C,D) \in \sum_{A \in \nmty{U}} \nmty{U}^{[A]}$ a map
$$\mathsf{sigma}_1 : \sum_{f : [C]^{[A]}} \prod_{a \in [A]} [D(f(a))]^{[B(a)]} \to [\upSigma_C D]^{[\upSigma_A B]}$$
Given $a \in [A]$ and $b \in [B(a)]$, let $f_a : [\mathbf{1}] \to [A]$ be given by $f_a(\star) = a$ and let $g_b : [\underline{\mathbf{1}}(a)] = [\mathbf{1}] \to [B(a)]$ be given by $g_b(\star) = b$. Then
$$\mathsf{sigma}_1(f_a,g_b) : [\upSigma_{\mathbf{1}} \underline{\mathbf{1}}] \to [\upSigma_A B]$$
Define $\widehat{\mathsf{pair}}(A,B,a,b) = \mathsf{sigma}_1(f_a,g_b)(\langle \star, \star \rangle)$.

By construction, these are maps of the appropriate sorts, and $p \circ \widehat{\mathsf{pair}} = \widehat{\upSigma} \circ (p \cdot p)$. To see that the desired square is a pullback, note that the fibre of $p \cdot p$ over $(A,B)$ is mapped bijectively via $\widehat{\mathsf{pair}}$ to the fibre of $p$ over $\upSigma_A B$.
\end{proof}

In future work, we hope to find a result analogous to \Cref{thmInternalRightAdjoint} which characterises when a natural model admits dependent product types in terms of an internal adjunction.

%% file: thesis/ch3-polynomials-representability/polynomial-pseudomonads.tex
\section{Polynomial pseudomonads}
\label{secPolynomialPseudomonads}
\index{pseudomonad}
\index{polynomial!pseudomonad@{---} pseudomonad}

\begin{definition}
\label{defPolynomialMonad}
\index{polynomial!monad@{---} monad}
A \textbf{polynomial monad} is a monad in the bicategory $\mathbf{Poly}\cart_{\mathcal{E}}$. Specifically, a polynomial monad is a quadruple $\mathbb{P} = (I,p,\eta,\mu)$ consisting of an object $I$ of $\mathcal{E}$, a polynomial $p : I \pto I$ in $\mathcal{E}$ and cartesian morphisms of polynomials $\eta : i_1 \pRightarrow p$ and $\mu : p \cdot p \pRightarrow p$, satisfying the usual monad axioms, namely
$$\mu \circ (\mu \cdot p) = \mu \circ (p \cdot \mu) \quad \text{and} \quad \mu \circ (\eta \cdot p) = \mathrm{id}_p = \mu \circ (p \cdot \eta)$$
\end{definition}

\begin{remark}
\label{rmkPolynomialMonads}
What is usually (e.g.\ \cite{GambinoKock2013PolynomialFunctors}) meant by a \textit{polynomial monad} is a monad $(P,\eta,\mu)$ on a slice $\mathcal{E}\slice{I}$ of $\mathcal{E}$, with $P : \mathcal{E}\slice{I} \to \mathcal{E}\slice{I}$ a polynomial functor and $\eta,\mu$ cartesian natural transformations; equivalently, this is a monad in the $2$-category $\mathbf{PolyFun}\cart_{\mathcal{E}}$. We recover this notion from \Cref{defPolynomialMonad} by applying the extension bifunctor $\mathbf{Poly}\cart_{\mathcal{E}} \to \mathbf{PolyFun}\cart_{\mathcal{E}}$. Furthermore, every polynomial monad in the usual sense is the extension of a polynomial monad in the sense of \Cref{defPolynomialMonad}.
\end{remark}

Recall \Cref{thmUnitSigmaPiPoly}, which says that a natural model $(\mathbb{C}, p)$ admits a unit type if and only if there is a morphism $\eta : i_1 \pRightarrow p$ in $\mathbf{Poly}\cart_{\widehat{\mathbb{C}}}$, admits dependent sum types if and only if there is a cartesian morphism $\mu : p \cdot p \pRightarrow p$ in $\mathbf{Poly}\cart_{\widehat{\mathbb{C}}}$, and admits $(\mathbb{C},p)$ admits dependent product types if and only if there is a cartesian morphism $\zeta : P_p(p) \pRightarrow p$ in $\mathbf{Poly}\cart_{\widehat{\mathbb{C}}}$. It is natural to ask whether $(\Yon(\diamond),p,\eta,\mu)$ is a polynomial monad in the sense of \Cref{defPolynomialMonad}, and that $(p,\zeta)$ is an algebra for this monad in a suitable sense, but unfortunately, this turns out to be false. For example, consider the monad unit laws $\mu \circ (\eta \cdot p) = \mathrm{id}_p = \mu \circ (p \cdot \eta)$---they state precisely that the following equations of pasting diagrams hold:

\begin{ndiagram}
\dot{\mathcal{U}}
\arrow[d, "p"']
\arrow[r, "(\eta \cdot p)_1"]
&
\sum_{A,B} \sum_{a:[A]} [B(a)]
\arrow[d, "p \cdot p" description]
\arrow[r, "\mu_1"]
&
\dot{\mathcal{U}}
\arrow[d, "p" name={domequalsone}]
&[-15pt]
\dot{\mathcal{U}}
\arrow[d,"p"' name={codequalsone}]
\arrow[r,equals]
\arrow[from={domequalsone},to={codequalsone},"=" description,draw=none]
&
\dot{\mathcal{U}}
\arrow[d,"p" name={domequalstwo}]
&[-15pt]
\dot{\mathcal{U}}
\arrow[d, "p"' name={codequalstwo}]
\arrow[r, "(p \cdot \eta)_1"]
\arrow[from={domequalstwo},to={codequalstwo},"=" description,draw=none]
&
\sum_{A,B} \sum_{a:[A]} [B(a)]
\arrow[d, "p \cdot p" description]
\arrow[r, "\mu_1"]
&
\dot{\mathcal{U}}
\arrow[d, "p"]
\\
\mathcal{U}
\arrow[r, "(\eta \cdot p)_0"']
&
\sum_{A:\mathcal{U}} \mathcal{U}^{[A]}
\arrow[r, "\mu_0"']
&
\mathcal{U}
&
\mathcal{U}
\arrow[r,equals]
&
\mathcal{U}
&
\mathcal{U}
\arrow[r, "(p \cdot \eta)_0"']
&
\sum_{A:\mathcal{U}} \mathcal{U}^{[A]}
\arrow[r, "\mu_0"']
&
\mathcal{U}
\end{ndiagram}

However, the monad laws do not hold strictly in general. Indeed, in the internal language of $\widehat{\mathbb{C}}$, we have
$$(\mu \circ (\eta \cdot p))_0(A) = \sum_{x : A} \mathbf{1} = A \times \mathbf{1} \quad \text{and} \quad (\mu \circ (p \cdot \eta))_0(A) = \sum_{x : \mathbf{1}} A = \mathbf{1} \times A$$
But in type theory, the types $A \times \mathbf{1}$, $A$ and $\mathbf{1} \times A$ are not generally \textit{equal}, although there are canonical isomorphisms between them. We therefore cannot, in general, expect the monad laws to hold strictly, for instance if the natural model is one arising from the syntax of dependent type theory (\Cref{schTermModel}). However, it is still reasonable to expect this structure to satisfy the laws of a \textit{pseudomonad}.

Much as monads naturally live in bicategories, pseudomonads naturally live in \textit{tricategories} \cite{Marmolejo1999Pseudomonads,Lack2000Pseudomonads}. To define the notion of a polynomial pseudomonad, we therefore need to endow the bicategory $\mathbf{Poly}\cart_{\mathcal{E}}$ with $3$-cells turning it into a tricategory.

\subsection*{A tricategory of polynomials}

In general, tricategories are fiddly, with lots of coherence data to worry about \cite{GordonPowerStreet1995Tricategories,Gurski2013CoherenceThreeCategories}---fortunately for us, our situation is simplified by the fact that composition of 2-cells of polynomials is strict, so that the $3$-cells turn the hom categories $\mathbf{Poly}\cart`_{\mathcal{E}}(I,J)$ into $2$-categories, rather than bicategories. The emerging structure is that of a \textit{$\mathbf{2Cat}$-enriched bicategory}.

\begin{definition}[$\mathbf{2Cat}$-enriched bicategories]
\label{def2CatEnrichedBicategory}
\index{2Cat-enriched bicategory@$\mathbf{2Cat}$-enriched bicategory}
A \textbf{$\mathbf{2Cat}$-enriched bicategory} $\mathfrak{B}$ consists of:
\begin{itemize}
\item A set $\mathfrak{B}_0$, whose elements we call the \textbf{0-cells} of $\mathfrak{B}$;
\item For all 0-cells $I,J$, a 2-category $\mathfrak{B}(I,J)$, whose 0-cells, 1-cells and 2-cells we call the \textbf{1-cells}, \textbf{2-cells} and \textbf{3-cells} of $\mathfrak{B}$, respectively;
\item For all 0-cells $I,J,K$, a 2-functor $\circ_{I,J,K} : \mathfrak{B}(J,K) \times \mathfrak{B}(I,J) \to \mathfrak{B}(I,K)$, which we call the \textbf{composition} 2-functor;
\item For all 0-cells $I$, a 2-functor $\iota_I : \mathbf{1} \to \mathfrak{B}(I,I)$, which we call the \textbf{identity} 2-functor, where $\mathbf{1}$ is the terminal 2-category;
\item For all 0-cells $I,J,K,L$, a 2-natural isomorphism
\begin{center}
\begin{tikzcd}[column sep=huge, row sep=huge]
\mathfrak{B}(K,L) \times \mathfrak{B}(J,K) \times \mathfrak{B}(I,J)
\arrow[r, "\circ_{J,K,L} \times \mathrm{id}"]
\arrow[d, "\mathrm{id} \times \circ_{I,J,K}"']
&
\mathfrak{B}(J,L) \times \mathfrak{B}(I,J)
\arrow[d, "\circ_{I,J,L}"]
\arrow[dl, draw=none, pos=0.5, "\phantom{\alpha_{I,J,K,L}}\ \twocell{225}\ \alpha_{I,J,K,L}" description]
\\
\mathfrak{B}(K,L) \times \mathfrak{B}(I,K)
\arrow[r, "\circ_{I,K,L}"']
&
\mathfrak{B}(I,L)
\end{tikzcd}

\end{center}
called the \textbf{associator};

\item For all 0-cells $I,J$, 2-natural isomorphisms

\begin{center}
\begin{tikzcd}[row sep=huge, column sep={20pt}]
\mathfrak{B}(I,J) \times \mathbf{1}
\arrow[r, "\mathrm{id} \times \iota_I"]
\arrow[dr, "\cong"']
&
\mathfrak{B}(I,J) \times \mathfrak{B}(I,I)
\arrow[d, "\circ_{I,I,J}"]
\arrow[dl, draw=none, pos=0.2, "\lambda_{I,J} \twocell{225} \phantom{\lambda_{I,J}}" description]
\\
~
&
\mathfrak{B}(I,J)
&
\end{tikzcd}
\begin{tikzcd}[row sep=huge, column sep={20pt}]
\mathfrak{B}(J,J) \times \mathfrak{B}(I,J)
\arrow[d, "\circ_{I,J,J}"']
\arrow[dr, draw=none, pos=0.2, "\rho_{I,J} \twocell{315} \phantom{\rho_{I,J}}" description]
&
\mathbf{1} \times \mathfrak{B}(I,J)
\arrow[l, "\iota_J \times \mathrm{id}"']
\arrow[dl, "\cong"]
\\
\mathfrak{B}(I,J)
&
~
\end{tikzcd}

\end{center}
called the \textbf{left unitor} and \textbf{right unitor}, respectively.
\end{itemize}
such that for all compatible 1-cells $I \xrightarrow{f} J \xrightarrow{g} K \xrightarrow{h} L \xrightarrow{k} M$, the following diagrams commute:
\begin{center}
\begin{tikzcd}[row sep=huge, column sep={0pt}]
((k \circ h) \circ g) \circ f
\arrow[rr, Rightarrow, "\alpha_{I,J,K,M}"]
\arrow[dr, Rightarrow, "\alpha_{J,K,L,M} \circ f"']
&&
(k \circ h) \circ (g \circ f)
\arrow[rr, Rightarrow, "\alpha_{I,K,L,M}"]
&&
k \circ (h \circ (g \circ f))
\\
&
(k \circ (h \circ g)) \circ f
\arrow[rr, Rightarrow, "\alpha_{I,J,L,M}"']
&&
k \circ ((h \circ g) \circ f)
\arrow[ur, Rightarrow, "k \circ \alpha_{I,J,K,L}"']
&
\end{tikzcd}

\end{center}

\begin{center}
\begin{tikzcd}[row sep=huge, column sep=huge]
(g \circ \iota_J) \circ f
\arrow[rr, Rightarrow, "\alpha_{I,J,J,K}"]
\arrow[dr, Rightarrow, "\lambda_{J,K} \circ f"']
&&
g \circ (\iota_J \circ f)
\arrow[dl, Rightarrow, "g \circ \rho_{I,J}"]
\\
&
g \circ f
&
\end{tikzcd}

\end{center}
\end{definition}

Every 3-category is trivially a $\mathbf{2Cat}$-enriched bicategory, and every $\mathbf{2Cat}$-enriched bicategory is a tricategory. Every $\mathbf{2Cat}$-enriched bicategory has an underlying bicategory, obtained by forgetting the 3-cells, and every bicategory can be equipped with the structure of a $\mathbf{2Cat}$-enriched bicategory by taking only identities as 3-cells. An equivalent viewpoint is that $\mathbf{2Cat}$-enriched bicategories are tricategories, whose hom-bicategories are $2$-categories and whose coherence isomorphisms in the top dimension are identities.

Connections between polynomials and $\mathbf{2Cat}$-enriched bicategories have been studied in different but related settings by Tamara von Glehn \cite{vonGlehn2015Thesis} and by Mark Weber \cite{Weber2015Polynomials} (the latter referring to them as `2-bicategories').

In order to motivate our definition of $3$-cells, recall that \Cref{cnsCartesianMorphismsInduceInternalFunctors} yields a functor
$$\mathbb{S} : \mathbf{Poly}\cart_{\mathcal{E}}(1,1) \to \mathbf{Cat}(\mathcal{E})$$
However, $\mathbf{Cat}(\mathcal{E})$ has the structure of a 2-category, so it is therefore reasonable to expect that when we equip $\mathbf{Poly}_{\mathcal{E}}$ with 3-cells, the functor $\mathbb{S}$ should extend to a 2-functor. In particular, any 3-cell between cartesian morphisms of polynomials should induce an internal natural transformation between the induced internal functors. However, since the association of internal functors to morphisms of polynomials works only for \textit{cartesian} morphisms of polynomials, we cannot simply take internal natural transformations as the 3-cells of $\mathbf{Poly}_{\mathcal{E}}$. \Cref{lemHalfInternalNT,lemCartesianAdjustmentsAreInternalNT} provide a correspondence between internal natural transformations $\fisc{\varphi} \Rightarrow \fisc{\psi}$ and particular morphisms of $\mathcal{E}$ in a way that generalises to the case when $\varphi$ and $\psi$ are not required to be cartesian.

\begin{lemma}
\label{lemHalfInternalNT}
Let $f : B \to A$ and $g : D \to C$ be polynomials in a locally cartesian closed category $\mathcal{E}$ and let $\varphi, \psi : f \pRightarrow g$ be cartesian morphisms of polynomials. There is a bijection between the set of morphisms $\alpha : \Delta_{\varphi_0} D \to \Delta_{\psi_0} D$ in $\mathcal{E} \slice{A}$ and the set of morphisms $\widehat{\alpha} : A \to \fisc{g}_1$ in $\mathcal{E} \slice{C \times C}$, as indicated by dashed arrows in the following diagrams, where $\varphi_2, \psi_2$ are canonical isomorphisms induced by the universal property of pullbacks.

\begin{center}
\begin{tikzcd}[row sep=small, column sep=small]
\Delta_{\varphi} D
\arrow[rr, dashed, blue, "\alpha"]
\arrow[d, "\varphi_2"', "{\scriptsize \cong}"]
&&
\Delta_{\psi_0} D
\arrow[d, "\psi_2", "{\scriptsize \cong}"']
\\
B
\arrow[dr, "f"', bend right]
&&
B
\arrow[dl, "f", bend left]
\\
&
A
&
\end{tikzcd}
\hspace{30pt}
\begin{tikzcd}[row sep=huge, column sep=huge]
&
\fisc{g}_1
\arrow[d, "\partial"]
\\
A
\arrow[ur, dashed, blue, "\widehat{\alpha}"]
\arrow[r, "{\langle \varphi_0, \psi_0 \rangle}"']
&
C \times C
\end{tikzcd}
\end{center}
\end{lemma}

\begin{proof}
Given $\alpha : \Delta_{\varphi_0} D \to \Delta_{\psi_0} D$ in $\mathcal{E} \slice{A}$, the exponential transpose of $\alpha$ in $\mathcal{E} \slice{A}$ is, as a morphism in $\mathcal{E}$, a section $\overline{\alpha} : A \to H$ of the projection $H \to A$, where $H = \sum_{a \in A} D_{\psi_0(a)}^{D_{\varphi_0(a)}}$. This projection is precisely the pullback of $\fisc{g}_1 \to C \times C$ along $\langle \varphi_0, \psi_0 \rangle$, as  illustrated in the following diagram:
\begin{center}
\begin{tikzcd}[row sep=huge, column sep=huge]
H
\arrow[d]
\arrow[r]
\pullback
&
\fisc{g}_1
\arrow[d, "\partial"]
\\
A
\arrow[r, "{\langle \varphi_0, \psi_0 \rangle}"']
\arrow[u, bend left, dashed, blue, "\overline{\alpha}"]
&
C \times C
\end{tikzcd}

\end{center}
But sections of the pullback correspond with diagonal fillers $\widehat{\alpha} : A \to \fisc{g}_1$ of the pullback square. This is as required, since such a filler making the lower triangle commute makes the upper triangle commute automatically. This concludes the proof of (a).
\end{proof}

\begin{lemma}
\label{lemCartesianAdjustmentsAreInternalNT}
Let $f : B \to A$ and $g : D \to C$ be polynomials in a locally cartesian closed category $\mathcal{E}$, let $\varphi,\psi : f \pRightarrow g$ be cartesian morphisms of polynomials, and let $\alpha,\widehat{\alpha}$ be as in \Cref{lemHalfInternalNT}. The following are equivalent:
\begin{enumerate}[(i)]
\item $\widehat{\alpha}$ is an internal natural transformation $\fisc{\varphi} \Rightarrow \fisc{\psi}$;
\item In the internal language of $\mathcal{E}$, we have $\fisc{\psi}(k) \circ \alpha_a = \alpha_{a'} \circ \fisc{\varphi}(k)$ for $a,a' \in A$ and $k \in B_{a'}^{B_a}$;
\item In the internal language of $\mathcal{E}$, we have $\gamma_{a'} \circ k = k \circ \gamma_a$ for $a,a' \in A$ and $k \in B_{a'}^{B_a}$, where $\gamma = \psi_2 \circ \alpha \circ \varphi_2^{-1} : B \to B$;
\item $\alpha$ is a morphism in $\mathcal{E} \slice{B}$, i.e.\ $\psi_2 \circ \alpha = \varphi_2$.
\end{enumerate}
\end{lemma}

\begin{proof}
We prove (i)$\Leftrightarrow$(ii)$\Leftrightarrow$(iii)$\Leftrightarrow$(iv).

\begin{itemize}
\item[(i)$\Leftrightarrow$(ii)] In light of \Cref{lemHalfInternalNT}, this is just a translation into the internal language of $\mathcal{E}$ of the definition of an internal natural transformation.

\item[(ii)$\Leftrightarrow$(iii)] Consider the following `internal' diagram, parametrised by $a,a' \in A$ and $k \in B_{a'}^{B_a}$.
\begin{diagram}
B_a
\arrow[d, "k"']
\arrow[r, "(\varphi_2)_a^{-1}"]
&
D_{\varphi(a)}
\arrow[d, "\fisc{\varphi}(k)" description]
\arrow[r, "\alpha_a"]
&
D_{\psi(a)}
\arrow[d, "\fisc{\psi}(k)" description]
\arrow[r, "(\psi_2)_a"]
&
B_a
\arrow[d, "k"]
\\
B_{a'}
\arrow[r, "(\varphi_2)_{a'}^{-1}"']
&
D_{\varphi(a')}
\arrow[r, "\alpha_{a'}"']
&
D_{\psi(a')}
\arrow[r, "(\psi_2)_{a'}"']
&
B_{a'}
\end{diagram}

The left- and right-hand squares commute by functoriality of $\fisc{\varphi}$ and $\fisc{\psi}$.
The centre square commutes if and only if (ii) holds, and the outer square commutes if and only if (iii) holds. But the centre square commutes if and only if the outer square commutes.

\item[(iii)$\Leftrightarrow$(iv)] Let $a \in A$ and $b \in B_a$, and let $k \in B_a^{B_a}$ be the constant (internal) function with value $b$. If (iii) holds, then
$$\gamma_a(b) = \gamma_a(k(b)) = k(\gamma_{a}(b)) = b$$
so that $\seqbn{\gamma_a = \mathrm{id}_{B_a}}{a \in A}$ holds. But this says precisely that $\gamma=\mathrm{id}_B$, and hence $\psi_2 \circ \alpha = \varphi_2$. The converse (iv)$\Rightarrow$(iii) is immediate.
\end{itemize}
\end{proof}

\begin{definition}
\label{defAdjustment}
\index{adjustment}
Let $F : \declpoly IBAJsft$ and $G : \declpoly IDCJugv$ be polynomials and let $\varphi,\psi : F \pRightarrow G$ be morphisms of polynomials, as in:
\begin{center}
\begin{tikzcd}[row sep=normal, column sep=normal]
&
B
\arrow[r,"f"]
\arrow[dl, bend right=20, "s"']
&
A
\arrow[d, equals]
\arrow[dr, bend left=20, "t"]
&
& 
&
B
\arrow[r,"f"]
\arrow[dl, bend right=20, "s"']
&
A
\arrow[d, equals]
\arrow[dr, bend left=20, "t"]
&
\\
I
&
D_{\varphi}
\pullback
\arrow[r]
\arrow[d,"\varphi_1"']
\arrow[u,"\varphi_2"]
&
A
\arrow[d,"\varphi_0"]
&
J
& 
I
&
D_{\psi}
\pullback
\arrow[r]
\arrow[d,"\psi_1"']
\arrow[u,"\psi_2"]
&
A
\arrow[d,"\psi_0"]
&
J
\\
&
D
\arrow[r,"g"']
\arrow[ul, bend left=20, "u"]
&
C
\arrow[ur, bend right=20, "v"']
&
& 
&
D
\arrow[r,"g"']
\arrow[ul, bend left=20, "u"]
&
C
\arrow[ur, bend right=20, "v"']
&
\end{tikzcd}

\end{center}
An \textbf{adjustment} $\alpha$ from $\varphi$ to $\psi$, denoted $\alpha : \varphi \pRrightarrow \psi$, is a morphism $\alpha : D_{\varphi} \to D_{\psi}$ over $B$:
\begin{center}
\begin{tikzcd}[row sep=huge, column sep=huge]
D_{\varphi}
\arrow[rr,"\alpha"]
\arrow[dr,"\varphi_2"']
&&
D_{\psi}
\arrow[dl,"\psi_2"]
\\
&
B
&
\end{tikzcd}

\end{center}
\end{definition}

\begin{numbered}
\Cref{lemCartesianAdjustmentsAreInternalNT} tells us that, when $\varphi$ and $\psi$ are cartesian, adjustments $\alpha : \varphi \pRrightarrow \psi$ can equivalently be described as internal natural transformations $\widehat{\alpha} : \varphi \Rightarrow \psi$.
\end{numbered}

We can now, at least, state the following conjecture.

\begin{conjecture}
\label{conjPolyETricategory}
There is a $\mathbf{2Cat}$-enriched bicategory $\mathfrak{Poly}_{\mathcal{E}}$, whose underlying bicategory is $\mathbf{Poly}_{\mathcal{E}}$ and whose 3-cells are adjustments.
\end{conjecture}

Unfortunately, the details required to fully prove \Cref{conjPolyETricategory} turned out to be somewhat laborious and, since its full force is not required for our main results, we have left the task of verifying these details for future work. Our progress so far is outlined in \Cref{lemPolyEHom2Categories} and \Cref{rmkCompleteProofOfConjecture}, and we prove the analogous result with attention restricted to \textit{cartesian} morphisms of polynomials in \Cref{thmPolyECartTrivial}.

\begin{lemma}
\label{lemPolyEHom2Categories}
Let $I$ and $J$ be objects in a locally cartesian closed category $\mathcal{E}$. There is a 2-category $\mathfrak{Poly}_{\mathcal{E}}(I,J)$ whose underlying category is $\mathbf{Poly}_{\mathcal{E}}(I,J)$ and whose 2-cells are adjustments.
\end{lemma}

\begin{proof}
Given polynomials $F,G : I \pto J$, the category $\mathfrak{Poly}_{\mathcal{E}}(I,J)(F,G)$ has morphisms of polynomials $F \pRightarrow G$ as its objects and adjustments as its morphisms, with identity and composition inherited from $\mathcal{E}\slice{B}$.

Given a polynomial $F : \declpoly IBAJsft$, we have an evident functor $\mathbf{1} \to \mathfrak{Poly}_{\mathcal{E}}(I,J)(F,F)$ picking out the identity morphism $F \pRightarrow F$ and the identity adjustment on this morphism.

Let $F,G,H : I \pto J$ be polynomials. The composition functor
$$c : \mathbf{Poly}_{\mathcal{E}}(I,J)(G,H) \times \mathbf{Poly}_{\mathcal{E}}(I,J)(F,G) \to \mathbf{Poly}_{\mathcal{E}}(I,J)(F,H)$$
is defined as follows. The composite $c(\psi,\varphi)$ of $\varphi : F \pRightarrow G$ and $\psi : G \pRightarrow H$ is defined using a pullback construction, as defined in \cite[3.9]{GambinoKock2013PolynomialFunctors}---in particular, the morphism $(\psi \circ \varphi)_2 : D_{\psi \circ \varphi} \to B$ is induced by the universal property of pullbacks. This yields, for each pair of adjustments $\alpha : \varphi \pRrightarrow \varphi'$ and $\beta : \psi \pRrightarrow \psi'$, a unique morphism $D_{\psi \circ \varphi} \to D_{\psi' \circ \varphi'}$ in $\mathcal{E}$ induced by the universal property of pullbacks, which is an adjustment since it makes the required triangle in $\mathcal{E}\slice{B}$ commute. We take this morphism to be $c(\beta,\alpha)$. Functoriality of $c$ is then immediate from the universal property of pullbacks.

It can be easily verified that this data satisfies the required identity and associativity axioms. Thus we have a 2-category.
\end{proof}

\begin{numbered}
\label{rmkCompleteProofOfConjecture}
In order to prove \Cref{conjPolyETricategory} in its entirety, it remains to define the coherence 2-natural isomorphisms $\alpha,\lambda,\rho$, as described in \Cref{def2CatEnrichedBicategory}, and verify that the required diagrams commute.

To give the reader an idea of the flavour of this task, we present some progress towards defining the associator 2-natural transformation $\alpha$. For each quadruple of objects $I,J,K,L$ of $\mathcal{E}$, this must assign to each triple of polynomials $I \overset{F}{\pto} J \overset{G}{\pto} K \overset{H}{\pto} L$ a morphism of polynomials $\alpha_{F,G,H} : (H \cdot G) \cdot F \pRightarrow H \cdot (G \cdot F)$ and, to each triple of morphisms of polynomials
$$\varphi : F \pRightarrow F', \quad \chi : G \pRightarrow G', \quad \psi : H \pRightarrow H'
$$
an adjustment
$$\alpha_{\varphi,\chi,\psi} : \psi \cdot (\chi \cdot \varphi) \circ \alpha_{F,G,H} \pRrightarrow \alpha_{F',G',H'} \circ (\psi \cdot \chi) \cdot \varphi : (F \cdot G) \cdot H \pRightarrow F' \cdot (G' \cdot H')
$$
which satisfy naturality laws and behave well with respect to composition and identity.

Restricting to the case $I=J=K=L=1$, let $f : B \to A$, $g : D \to C$ and $h : F \to E$ be morphisms of $\mathcal{E}$, considered as polynomials $1 \pto 1$ as usual. We will construct an invertible (and hence cartesian) morphism of polynomials $\alpha_{f,g,h} : (h \cdot g) \cdot f \pRightarrow h \cdot (g \cdot f)$. Such a morphism must fit into the following pullback square:
\begin{center}
\begin{tikzcd}[row sep=huge, column sep=huge]
\sum_{e,n,q} \sum_{f \in F_e} \sum_{d \in D_{n(f)}} B_{q(f,d)}
\arrow[d, "(h \cdot g) \cdot f"']
\arrow[r, "{(\alpha_{f,g,h})_1}"]
\pullbackc{dr}{-0.1}
&
\sum_{e,p} \sum_{f \in F_e} \sum_{d \in D_{c_f}} B_{m_f(d)}
\arrow[d, "h \cdot (g \cdot f)"]
\\
\sum_{e \in E} \sum_{n \in C^{F_e}} \prod_{f \in F_e} \prod_{d \in D_{n(f)}} A
\arrow[r, "{(\alpha_{f,g,h})_0}"']
&
\sum_{e \in E} \prod_{f \in F_e} \sum_{c \in C} \prod_{d \in D_c} A
\end{tikzcd}

\end{center}
In the above, we have overloaded the letter $f$, which is ambiguous between the morphism $f : B \to A$ of $\mathcal{E}$ and an internal `element' $f \in F_e$; and we have written $p(f)=(c_f,m_f)$ for $p \in \prod_{f \in F_e} \sum_{c \in C} \prod_{d \in D_c} A$ and $f \in F_e$.

The isomorphism $(\alpha_{f,g,h})_0$ is given by applying the type theoretic axiom of choice to exchange the middle $\upSigma\upPi$. Specifically, we have
$$(\alpha_{f,g,h})_0(e,n,q) = (e, \lambda f. \langle n(f), q(f) \rangle)
$$
The isomorphism $(\alpha_{f,g,h})_1$ acts trivially; that is, we have
$$(\alpha_{f,g,h})_1(e,n,q,f,d,b) = ((\alpha_{f,g,h})_0(e,n,q),f,d,b)
$$

We suspect that the definition of $\alpha_{\varphi,\chi,\psi}$ will also be an instance of the type theoretic axiom of choice. From this, it will be an exercise in symbolic manipulations to check that the `Mac Lane pentagon' commutes.
\end{numbered}

The situation in which we restrict our attention to cartesian morphisms of polynomials is greatly simplified by the following lemma, allowing us to prove \Cref{conjPolyETricategory} for this case in \Cref{thmPolyECartTrivial}.

\begin{lemma}
\label{lemTriviality}
Let $\varphi$ and $\psi$ be morphisms of polynomials. If $\psi$ is cartesian then there is a unique adjustment from $\varphi$ to $\psi$.
\end{lemma}

\begin{proof}
When $\psi$ is cartesian, the morphism $\psi_2$ is invertible, so that $\alpha = \psi_2^{-1} \circ \varphi_2$ is the only morphism making the required triangle commute.
\end{proof}

From \Cref{thmPolyEBicategory}(d) and \Cref{lemTriviality}, we immediately obtain the following theorem.

\begin{theorem}
\label{thmPolyECartTrivial}
There is a $\mathbf{2Cat}$-enriched bicategory $\mathfrak{Poly}\cart_{\mathcal{E}}$ whose underlying bicategory is $\mathbf{Poly}\cart_{\mathcal{E}}$ and whose hom 2-categories $\mathfrak{Poly}\cart_{\mathcal{E}}(I,J)$ are locally codiscrete for all objects $I,J$ of $\mathcal{E}$.
\end{theorem}

\begin{proof}
The description of the $\mathbf{2Cat}$-enriched bicategory data is described in \Cref{rmkCompleteProofOfConjecture}. The coherence data is uniquely defined and satisfies the required equations by \Cref{lemTriviality}.
\end{proof}

Before moving on, we extend \Cref{lemPolynomialsFromOneToOne} to our tricategorical setting.

\begin{lemma}
\label{lemPolynomialsFromOneToOneTwoCategorical}
For fixed objects $I$ and $J$ of a locally cartesian closed category $\mathcal{E}$, there are full and faithful 2-functors
$$S : \mathfrak{Poly}_{\mathcal{E}}(I,J) \to \mathfrak{Poly}_{\mathcal{E}\slice{I \times J}}(1,1) \quad \text{and} \quad S\cart : \mathfrak{Poly}\cart_{\mathcal{E}}(I,J) \to \mathfrak{Poly}\cart_{\mathcal{E} \slice{I \times J}}(1,1)$$
\end{lemma}

\begin{proof}
Let $F : \declpoly IBAJsft$ and $G : \declpoly IDCJugv$ be polynomials $I \pto J$, and let $\varphi,\psi$ be morphisms of polynomials $F \pRightarrow G$. An adjustment $\alpha : \varphi \pRrightarrow \psi$ is simply a morphism $\alpha : \varphi_2 \to \psi_2$ in $\mathcal{E}\slice{B}$. Since $S(\varphi)_2 = \varphi_2$ and $S(\psi)_2 = \psi_2$, an adjustment $S(\varphi) \pRrightarrow S(\psi)$ is a morphism $\varphi_2 \to \psi_2$ in $(\mathcal{E} \slice{I \times J}) \slice{\langle s, t \circ f \rangle} \cong \mathcal{E} \slice{B}$. So we can take $S$ to be the identity on adjustments. This trivially extends the functors $S$ and $S\cart$ of \Cref{lemPolynomialsFromOneToOne} to full and faithful 2-functors.
\end{proof}

\begin{theorem}
\label{thmPolynomialsYieldInternalCategories2Functor}
Fix objects $I$ and $J$ in a locally cartesian closed category $\mathcal{E}$. There is a locally full and faithful 2-functor
$$\mathbb{A}_{(-)} : \mathfrak{Poly}\cart_{\mathcal{E}}(I,J) \to \mathbf{Cat}(\mathcal{E} \slice{I \times J})
$$
whose underlying 1-functor is as in \Cref{cnsCartesianMorphismsInduceInternalFunctors}.
\end{theorem}

\begin{proof}
Let $\varphi,\psi : F \pRightarrow G$ be cartesian morphisms of polynomials $I \pto J$. We proved in \Cref{lemCartesianAdjustmentsAreInternalNT} that adjustments $\alpha : \varphi \pRrightarrow \psi$ correspond bijectively with internal natural transformations $\widehat{\alpha} : \fisc{\varphi} \Rightarrow \fisc{\psi}$. Moreover, by \Cref{lemTriviality}, there is a unique internal natural transformation $\fisc{\varphi} \Rightarrow \fisc{\psi}$. As such, defining $\mathbb{A}_{\alpha} = \widehat{\alpha}$ for all adjustments $\alpha$, we automatically obtain a 2-functor, which is locally full and faithful since the hom-sets
$$\mathfrak{Poly}\cart_{\mathcal{E}}(I,J)(F,G)(\varphi,\psi) \quad \text{and} \quad \mathbf{Cat}(\mathcal{E}\slice{I \times J})(\fisc{f},\fisc{g})(\fisc{\varphi},\fisc{\psi})$$
are both singletons.
\end{proof}

\subsection*{Polynomial pseudomonads}

We are now ready to define the notion of a polynomial pseudomonad. First, we recall the definition of a pseudomonad in a $\mathbf{2Cat}$-enriched bicategory (in fact, the definition works just fine in an arbitrary tricategory).

\begin{definition}
\label{defPseudomonadIn2CatEnrichedBicategory}
\index{pseudomonad}
Let $\mathfrak{B}$ be a $\mathbf{2Cat}$-enriched bicategory. A \textbf{pseudomonad} $\mathbb{T}$ in $\mathfrak{B}$ consists of:
\begin{itemize}
\item A 0-cell $I$ of $\mathfrak{B}$;
\item A 1-cell $t : I \to I$;
\item 2-cells $\eta : \mathrm{id}_I \Rightarrow t$ and $\mu : t \cdot t \Rightarrow t$, called the \textbf{unit} and \textbf{multiplication} of the pseudomonad, respectively;
\item Invertible 3-cells $\alpha, \lambda, \rho$, called the \textbf{associator}, \textbf{left unitor} and \textbf{right unitor} of the pseudomonad, respectively, as in
\begin{center}
\begin{tikzcd}[row sep=huge, column sep=huge]
t \cdot t \cdot t
\arrow[r, Rightarrow, "t \cdot \mu"]
\arrow[d, Rightarrow, "\mu \cdot t"']
&
t \cdot t
\arrow[d, Rightarrow, "\mu"]
\arrow[dl, draw=none, pos=0.5, "\alpha \threecell{225}" description]
&
t
\arrow[r, Rightarrow, "t \cdot \eta"]
\arrow[dr, Rightarrow, "\mathrm{id}_t"']
&
t \cdot t
\arrow[d, Rightarrow, "\mu" description]
\arrow[dl, draw=none, pos=0.25, "\lambda \threecell{225}" description]
\arrow[dr, draw=none, pos=0.25, "\rho \threecell{315}" description]
&
t
\arrow[l, Rightarrow, "\eta \cdot t"']
\arrow[dl, Rightarrow, "\mathrm{id}_t"]
\\
t \cdot t
\arrow[r, Rightarrow, "\mu"']
&
t
&
~
&
t
&
~
\end{tikzcd}

\end{center}
\end{itemize}
such that the following equations of pasting diagrams hold:

\begin{center}
\begin{tikzcd}[row sep=huge, column sep={40pt}]
t \cdot t \cdot t \cdot t
\arrow[r, Rightarrow, "t \cdot t \cdot \mu"]
\arrow[d, Rightarrow, "\mu \cdot t \cdot t"']
\arrow[dr, Rightarrow, "t \cdot \mu \cdot t" description]
&
t \cdot t \cdot t
\arrow[dr, Rightarrow, "t \cdot \mu"]
\arrow[d, draw=none, pos=0.5, "t \cdot \alpha\ \threecell{270} \phantom{T\alpha}" description]
&
&
t \cdot t \cdot t \cdot t
\arrow[r, Rightarrow, "t \cdot t \cdot \mu"]
\arrow[d, Rightarrow, "\mu \cdot t \cdot t"']
\arrow[dr, draw=none, pos=0.5, "\cong" description]
&
t \cdot t \cdot t
\arrow[dr, Rightarrow, "t \cdot \mu"]
\arrow[d, Rightarrow, "\mu \cdot t" description]
&
\\
t \cdot t \cdot t
\arrow[dr, Rightarrow, "\mu \cdot t"']
\arrow[r, draw=none, pos=0.5, "\threecell{0}"', "\alpha \cdot t"]
&
t \cdot t \cdot t
\arrow[r, Rightarrow, "t \cdot \mu" description]
\arrow[d, Rightarrow, "\mu \cdot t" description]
\arrow[dr, draw=none, pos=0.5, "\alpha \threecell{315} \phantom{\alpha}" description]
&
t \cdot t
\arrow[d, Rightarrow, "\mu"]
\arrow[r, draw=none, "=" description]
&
t \cdot t \cdot t
\arrow[r, Rightarrow, "t \cdot \mu" description]
\arrow[dr, Rightarrow, "\mu \cdot t"']
&
t \cdot t
\arrow[dr, Rightarrow, "\mu" description]
\arrow[r, draw=none, pos=0.5, "\alpha", "\threecell{0}"']
\arrow[d, draw=none, pos=0.5, "\alpha \threecell{270} \phantom{\alpha}" description]
&
t \cdot t
\arrow[d, Rightarrow, "\mu" description]
\\
&
t \cdot t
\arrow[r, Rightarrow, "\mu"']
&
t
&
&
t \cdot t
\arrow[r, Rightarrow, "\mu"']
&
t
\end{tikzcd}

\end{center}

\begin{center}
\begin{tikzcd}[row sep=huge, column sep={40pt}]
t \cdot t \cdot t
\arrow[r, Rightarrow, "t \cdot \mu"]
\arrow[dr, Rightarrow, "\mu \cdot t" description]
&
t \cdot t
\arrow[dr, Rightarrow, "\mu"{name=domequals}]
\arrow[d, draw=none, pos=0.5, "\alpha \threecell{270} \phantom{\alpha}" description]
&
&
t \cdot t \cdot t
\arrow[r, Rightarrow, "t \cdot \mu"]
\arrow[dr, draw=none, pos=0.2, "t \cdot \rho \threecell{315}" description]
&
t \cdot t
\arrow[dr, Rightarrow, "\mu"]
\arrow[d, draw=none, pos=0.5, "=" description]
&
\\
t \cdot t
\arrow[r, Rightarrow, "\mathrm{id}_{t \cdot t}"']
\arrow[u, Rightarrow, "t \cdot \eta \cdot t"]
\arrow[ur, draw=none, pos=0.2, "\lambda \cdot t \threecell{293}" description]
&
t \cdot t
\arrow[r, Rightarrow, "\mu"']
&
t
&
t \cdot t
\arrow[u, Rightarrow, "t \cdot \eta \cdot t"{name=codequals}]
\arrow[from=domequals, to=codequals, draw=none, pos=0.67, "=" description]
\arrow[rr, Rightarrow, "\mu"']
\arrow[ur, Rightarrow, "\mathrm{id}_{t \cdot t}" description]
&
~
&
t
\end{tikzcd}

\end{center}
\end{definition}

\begin{numbered}
\label{rmk2MonadPseumonadOn2Cat}
We reserve the following terminology for particular cases of pseudomonads in $\mathbf{2Cat}$-enriched bicategories:
\begin{itemize}
\item When the 3-cells $\alpha,\lambda,\rho$ are identities, we call $\mathbb{T}$ a \textbf{2-monad} in $\mathfrak{B}$. Note that a 2-monad in $\mathfrak{B}$ restricts to a monad in the underlying bicategory of $\mathfrak{B}$, and that every monad in the underlying bicategory of $\mathfrak{B}$ is automatically a 2-monad in $\mathfrak{B}$.
\item When $\mathfrak{B}=\mathbf{2Cat}$ is the 3-category of 2-categories, 2-functors, pseudo-natural transformations and modifications, and the underlying 0-cell of $\mathbb{T}$ is a 2-category $\mathcal{K}$, we say that $\mathbb{T}$ is a pseudomonad (or 2-monad) \textbf{on} $\mathcal{K}$.
\end{itemize}
\end{numbered}

\begin{definition}
\label{defPolynomialPseudomonad}
\index{polynomial!pseudomonad@{---} pseudomonad}
A \textbf{polynomial 2-monad} (resp.\ \textbf{polynomial pseudomonad}) is a 2-monad (resp.\ pseudomonad) in the $\mathbf{2Cat}$-enriched bicategory $\mathfrak{Poly}\cart_{\mathcal{E}}$. Specifically, a polynomial pseudomonad consists of the following data:
\begin{itemize}
\item An object $I$ of $\mathcal{E}$;
\item A polynomial $p : I \pto I$;
\item Cartesian morphisms of polynomials $\eta : i_I \pRightarrow p$ and $\mu : p \cdot p \pRightarrow p$;
\item Invertible adjustments $\alpha : \mu \circ (p \cdot \mu) \pRrightarrow \mu \circ (\mu \cdot p)$, $\lambda : \mu \circ (\eta \cdot p) \pRrightarrow \mathrm{id}_p$ and $\rho : \mu \circ (p \cdot \eta) \pRrightarrow \mathrm{id}_p$;
\end{itemize}
such that the adjustments $\alpha, \lambda, \rho$ satisfy the coherence axioms of \Cref{defPseudomonadIn2CatEnrichedBicategory}.
\end{definition}

A consequence of \Cref{thmPolyECartTrivial} is that all parallel pairs of cartesian morphisms of polynomials are uniquely isomorphic. It follows that, in this case, simply specifying the \textit{data} for a polynomial monad suffices for defining a polynomial pseudomonad---this is stated precisely in the following lemma, whose proof is immediate.

\begin{lemma}
\label{lemMonadDataIsPseudomonad}
Let $I$ be an object of $\mathcal{E}$, let $p : I \pto I$ be a polynomial and let $\eta : i_I \pRightarrow p$ and $\mu : p \cdot p \pRightarrow p$ be cartesian morphisms of polynomials. Then there are unique adjustments $\alpha,\lambda,\rho$ such that the septuple $\mathbb{P}=(I,p,\eta,\mu,\alpha,\lambda,\rho)$ is a polynomial pseudomonad in $\mathcal{E}$. \qed
\end{lemma}

The next result allows us to lift polynomial 2-monads and polynomial pseudomonads \textit{in} $\mathcal{E}$ to 2-monads and pseudomonads \textit{on} the hom 2-categories of $\mathfrak{Poly}\cart_{\mathcal{E}}$. This will play a key role in identifying the sense in which a natural model $p : \dot{\mathcal{U}} \to \mathcal{U}$ is a pseudoalgebra over the polynomial pseudomonad it induces.

\begin{theorem}
\label{thmPolynomialPseudomonadLifts}
Let $\mathbb{P} = (p,\eta,\mu,\alpha,\lambda,\rho)$ be a polynomial 2-monad (resp.\ pseudomonad) on an object $I$ of a locally cartesian closed category $\mathcal{E}$. Then $\mathbb{P}$ lifts to a 2-monad (resp.\ pseudomonad) $\mathbb{P}^+ = (P, h, m, \dots)$ on $\mathfrak{Poly}\cart_{\mathcal{E}}(I,I)$.
\end{theorem}

\begin{proof}
By \Cref{lemPolynomialsFromOneToOneTwoCategorical}, we may take $I=1$ without loss of generality, so thtat $p$ is just a morphism $p : Y \to X$ in $\mathcal{E}$ and $\eta,\mu$ are pullback squares in $\mathcal{E}$ (cf.\ \Cref{rmkCartesianMorphismIsPullbackSquare}).

For notational simplicity, write $\mathcal{K}$ to denote the 2-category $\mathfrak{Poly}\cart_{\mathcal{E}}(1,1)$. Note $\mathcal{K}$ has as its underlying category the wide subcategory $\mathcal{E}^{\to}_{\text{cart}}$ of $\mathcal{E}^{\to}$ whose morphisms are the pullback squares. Thus the 0-cells of $\mathcal{K}$ are the morphisms of $\mathcal{E}$, the 1-cells of $\mathcal{K}$ are pullback squares in $\mathcal{E}$, and between any two 1-cells there is a unique 2-cell by \Cref{thmPolyECartTrivial}.

First we must define a 2-functor $P : \mathcal{K} \to \mathcal{K}$. Define $P$ on the 0-cells of $\mathcal{K}$ by letting $P(f)=P_p(f)$ for all $f : B \to A$ in $\mathcal{E}$. Given a 1-cell $\varphi : f \pRightarrow g$ of $\mathcal{K}$---that is, a pullback square in $\mathcal{E}$---let $P(\varphi)$ be the result of applying the extension $P_p$ of $p$ to the pullback square defining $\varphi$, as in:
\begin{center}
\begin{tikzcd}[row sep=huge, column sep=huge]
\sum_{x \in X} B^{Y_x}
\arrow[r, "P_p(\varphi_1)"]
\arrow[d, "P_p(f)"']
\pullbackc{dr}{-0.05}
&
\sum_{x \in X} D^{Y_x}
\arrow[d, "P_p(g)"]
\\
\sum_{x \in X} A^{Y_x}
\arrow[r, "P_p(\varphi_0)"']
&
\sum_{x \in X} C^{Y_x}
\end{tikzcd}

\end{center}
Note that $P(\varphi)$ is indeed a pullback square, since polynomial functors preserve all connected limits \cite{GambinoKock2013PolynomialFunctors}. Thus $P(\varphi)$ is a 1-cell from $P(f)$ to $P(g)$ in $\mathcal{K}$.

Now $P$ respects identity 1-cells in $\mathcal{K}$, since if $f : B \to A$ is a 0-cell then
$$P(\mathrm{id}_f)_0 = P_p(\mathrm{id}_B) = \mathrm{id}_{P_p(B)} = (\mathrm{id}_{P(f)})_0
$$
and likewise $P(\mathrm{id}_f)_1 = (\mathrm{id}_{P(f)})_1$; and $P$ respects composition of 2-cells in $\mathcal{K}$, since for $i \in \{0,1\}$ we have
$$P(\psi \circ \varphi)_i = P_p((\psi \circ \varphi)_i) = P_p(\psi_i \circ \varphi_i) = P_p(\psi_i) \circ P_p(\varphi_i) = P(\psi)_i \circ P(\varphi)_i = (P(\psi) \circ P(\varphi))_i
$$
Hence the action of $P$ defines a functor on the underlying category of $\mathcal{K}$.

The fact that $P$ extends to a $2$-functor is trivial: given an adjustment $\alpha : \varphi \pRrightarrow \psi$, there is a unique adjustment $P(\varphi) \pRrightarrow P(\psi)$. We take this to be $P(\alpha)$, and note that the axioms governing identity and composition of 2-cells hold trivially by uniqueness of adjustments.

The pseudo-natural transformations $h : \mathrm{id}_{\mathcal{K}} \Rightarrow P$ and $m : P \circ P \Rightarrow P$ giving the unit and multiplication of $\mathbb{P}^+$ are induced by the unit $\eta : i_1 \pRightarrow p$ and $\mu : p \cdot p \pRightarrow p$ of $\mathbb{P}$. Specifically, define the components $h_f : f \pRightarrow P(f)$ and $m_f : P(P(f)) \pRightarrow P(f)$ at a 0-cell $f : B \to A$ of $\mathcal{K}$ to be the following squares, respectively:
\begin{center}
\begin{tikzcd}[row sep=normal, column sep=huge]
B
\arrow[r, "(P_{\eta})_B"]
\arrow[d, "f"']
\pullback
&
\sum_{x \in X} B^{Y_x}
\arrow[d, "P(f)"]
&
\sum_{(x,t) \in \sum_{x \in X} X^{Y_x}} B^{\left(\sum_{y \in Y_x} Y_{t(y)}\right)}
\arrow[r, "(P_{\mu})_B"]
\arrow[d, "P(P(f))"']
\pullbackc{dr}{-0.2}
&
\sum_{x \in X} B^{Y_x}
\arrow[d, "P(f)"]
\\
A
\arrow[r, "(P_{\eta})_A"']
&
\sum_{x \in X} A^{Y_x}
&
\sum_{(x,t) \in \sum_{x \in X} X^{Y_x}} A^{\left(\sum_{y \in Y_x} Y_{t(y)}\right)}
\arrow[r, "(P_{\mu})_A"']
&
\sum_{x \in X} A^{Y_x}
\end{tikzcd}

\end{center}
Note that these squares commute and are cartesian by naturality and cartesianness of the extensions $P_{\eta},P_{\mu}$ of $\eta,\mu$. That $h$ and $m$ extend to pseudo-natural transformations is immediate from \Cref{thmPolyECartTrivial}: the pseudo-naturality 2-cells in $\mathcal{K}$ are adjustments, so they exist uniquely and satisfy the coherence axioms for pseudo-natural transformations automatically.

If $\mathbb{P}$ is a polynomial 2-monad, it is now easy to verify that the 2-monad laws hold for $\mathbb{P}^+$. If $\mathbb{P}$ is a polynomial pseudomonad, then the pseudomonad laws for $\mathbb{P}^+$ concern existence of and equations between adjustments, hence are trivially true by \Cref{thmPolyECartTrivial}. 
\end{proof}

\begin{definition}
\label{defLiftOfPseudomonad}
Given a polynomial monad (resp.\ pseudomonad) $\mathbb{P}$, the \textbf{lift} of $\mathbb{P}$ is the 2-monad (resp.\ pseudomonad) $\mathbb{P}^+$ as in \Cref{thmPolynomialPseudomonadLifts}.
\end{definition}

\begin{definition}
\label{defPseudoalgebra}
\index{pseudoalgebra}
Let $\mathbb{T} = (T, h, m, \alpha, \lambda, \rho)$ be a pseudomonad on a 2-category $\mathcal{K}$. A \textbf{pseudoalgebra} over $\mathbb{T}$ consists of
\begin{itemize}
\item A 0-cell $A$ of $\mathcal{K}$;
\item A 1-cell $a : T(A) \to A$ in $\mathcal{K}$;
\item Invertible 2-cells $\sigma, \tau$ of $\mathcal{K}$, as in:
\begin{center}
\begin{tikzcd}[row sep=huge, column sep=huge]
T(T(A))
\arrow[r, "T(a)"]
\arrow[d, "m_T"']
&
T(A)
\arrow[d, "a"]
\arrow[dl, draw=none, pos=0.5, "\sigma \twocell{225}" description]
&
A
\arrow[r, "h_A"]
\arrow[dr, "\mathrm{id}_A"']
&
T(A)
\arrow[d, "a"]
\arrow[dl, draw=none, pos=0.25, "\tau \twocell{225}" description]
\\
T(A)
\arrow[r, "a"']
&
A
&
~
&
A
\end{tikzcd}

\end{center}
\end{itemize}
such that the following equations of pasting diagrams hold:

\begin{center}
\begin{tikzcd}[row sep=huge, column sep={40pt}]
T^3A
\arrow[r, "TTa"]
\arrow[d, "m_{TA}"']
\arrow[dr, "Tm_A" description]
&
T^2A
\arrow[dr, "Ta"]
\arrow[d, draw=none, pos=0.5, "T\sigma\ \twocell{270} \phantom{T\sigma}" description]
&
&
T^3A
\arrow[r, "T^2a"]
\arrow[d, "m_{TA}"']
\arrow[dr, draw=none, pos=0.5, "\cong" description]
&
T^2A
\arrow[dr, "Ta"]
\arrow[d, "m_A" description]
&
\\
T^2A
\arrow[dr, "m_A"']
\arrow[r, draw=none, pos=0.5, "\twocell{0}"', "\alpha_A"]
&
T^2A
\arrow[r, "Ta" description]
\arrow[d, "m_A" description]
\arrow[dr, draw=none, pos=0.5, "\sigma \twocell{315} \phantom{\sigma}" description]
&
TA
\arrow[d, "a"]
\arrow[r, draw=none, "=" description]
&
T^2A
\arrow[r, "Ta" description]
\arrow[dr, "m_A"']
&
TA
\arrow[dr, "a" description]
\arrow[r, draw=none, pos=0.5, "\sigma", "\twocell{0}"']
\arrow[d, draw=none, pos=0.5, "\sigma \twocell{270} \phantom{\sigma}" description]
&
TA
\arrow[d, "a" description]
\\
&
TA
\arrow[r, "a"']
&
A
&
&
T(A)
\arrow[r, "a"']
&
A
\end{tikzcd}

\end{center}

\begin{center}
\begin{tikzcd}[row sep=huge, column sep={40pt}]
T^2A
\arrow[r, "Ta"]
\arrow[dr, "m_A" description]
&
TA
\arrow[dr, "a"{name=domequals}]
\arrow[d, draw=none, pos=0.5, "\sigma \twocell{270} \phantom{\sigma}" description]
&
&
T^2A
\arrow[r, "Ta"]
\arrow[dr, draw=none, pos=0.2, "T\tau \twocell{315}" description]
&
TA
\arrow[dr, "a"]
\arrow[d, draw=none, pos=0.5, "=" description]
&
\\
TA
\arrow[r, "\mathrm{id}_{TA}"']
\arrow[u, "Th_A"]
\arrow[ur, draw=none, pos=0.2, "\lambda_A \twocell{293} \phantom{\lambda_A}" description]
&
TA
\arrow[r, "a"']
&
TA
&
TA
\arrow[u, "Th_A"{name=codequals}]
\arrow[from=domequals, to=codequals, draw=none, pos=0.67, "=" description]
\arrow[rr, "a"']
\arrow[ur, "\mathrm{id}_{TA}" description]
&
~
&
A
\end{tikzcd}

\end{center}
\end{definition}

\begin{definition}
\label{defPolynomialPseudoalgebra}
Let $\mathbb{P} = (1, p : Y \to X, \dots)$ be a polynomial pseudomonad in a locally cartesian closed category $\mathcal{E}$. A \textbf{polynomial pseudoalgebra} over $\mathbb{P}$ is a pseudoalgebra over the lift $\mathbb{P}^+$. Specifically, it consists of:
\begin{itemize}
\item A polynomial $f : B \to A$;
\item A cartesian morphism of polynomials $\zeta : P_p(f) \pRightarrow f$;
\item Invertible adjustments $\sigma,\tau$ whose types are as in \Cref{defPseudoalgebra};
\end{itemize}
such that the adjustments $\sigma,\tau$ satisfy the coherence conditions of \Cref{defPseudoalgebra}.
\end{definition}

Much like with polynomial pseudomonads (\Cref{lemMonadDataIsPseudomonad}), merely specifying the \textit{data} for a polynomial pseudoalgebra suffices for the conditions to hold---again, this follows immediately from \Cref{thmPolyECartTrivial}.

\begin{lemma}
\label{lemAlgebraDataIsPseudoalgebra}
Let $\mathbb{P} = (I, p : Y \to X, \dots)$ be a polynomial pseudomonad in a locally cartesian closed category $\mathcal{E}$, let $f : B \to A$ be a polynomial and let $\zeta : P_p(f) \pRightarrow f$ be a morphism of polynomials. Then there are unique adjustments $\sigma,\tau$ making $(f,\zeta,\sigma,\tau)$ into a polynomial pseudoalgebra over $\mathbb{P}$. \qed
\end{lemma}

We are now ready to precisely state the sense in which a natural model admitting a unit type and dependent sum types `almost' gives rise to a polynomial monad, and one admitting dependent product types `almost' gives rise to an algebra over this monad.

\begin{theorem}
\label{thmUnitSigmaIffPolynomialPseudomonad}
Let $(\mathbb{C}, p)$ be a natural model.
\begin{enumerate}[(a)]
\item $(\mathbb{C}, p)$ supports a unit type and dependent sum types if and only if $p$ can be equipped with the structure of a polynomial pseudomonad $\mathbb{P}$ in $\widehat{\mathbb{C}}$.
\item $(\mathbb{C}, p)$ additionally supports dependent product types if and only if $p$ can be equipped with the structure of a polynomial pseudoalgebra over $\mathbb{P}$.
\end{enumerate}
\end{theorem}
\begin{proof}
By \Cref{thmAdmittingUnitType,thmAdmittingSigmaTypes}, $(\mathbb{C}, p)$ supports a unit type and dependent sum types if and only if there exist cartesian morphisms of polynomials $\eta : i_1 \pRightarrow p$ and $\mu : p \cdot p \pRightarrow p$, and by \Cref{thmAdmittingPiTypes}, $(\mathbb{C}, p)$ additionally supports dependent product types if and only if there exists a cartesian morphism of polynomials $\zeta : P_p(p) \pRightarrow p$. By \Cref{lemMonadDataIsPseudomonad,lemAlgebraDataIsPseudoalgebra}, there are unique adjustments turning $(p,\eta,\mu)$ into a polynomial pseudomonad $\mathbb{P}$, and unique adjustments turning $(p,\zeta)$ into a polynomial pseudoalgebra over $\mathbb{P}$.
\end{proof}

\begin{numbered}
\Cref{thmUnitSigmaIffPolynomialPseudomonad} makes a connection between logic and algebra by exhibiting a correspondence between laws concerning dependent sums and dependent products in type theory with laws concerning monads in algebra. Specifically, for $\eta : \iota_1 \pRightarrow p$, $\mu : p \cdot p \pRightarrow p$ and $\zeta : P_p(p) \pRightarrow p$, the pseudomonad and pseudoalgebra isomorphisms in $\mathfrak{Poly}_{\widehat{\mathbb{C}}}\cart$ correspond to certain type isomorphisms as follows:
\begin{center}
\begin{tabular}{c||c|c}
\textbf{Name} & \textbf{Monads and algebras} & \textbf{Type theory} \\ \hline \hline
&& \\
Monad associativity & $\mu \circ (p \cdot \mu) \cong \mu \circ (\mu \cdot p)$ & $\sum_{x:A} \sum_{y:B(x)} C(x,y) \cong \sum_{\langle x,y \rangle : \sum_{x:A} B(x)} C(x,y)$ \\
&& \\ \hline
&& \\
Monad unit (left) & $\mu \circ (p \cdot \eta) \cong \mathrm{id}_p$ & $\sum_{x:A} \mathbf{1} \cong A$ \\ 
&& \\ \hline
&& \\
Monad unit (right) & $\mu \circ (\eta \cdot p) \cong \mathrm{id}_p$ & $\sum_{x:\mathbf{1}} A \cong A$ \\
&& \\ \hline
&& \\
Algebra multiplication & $\zeta \circ (p \cdot \zeta) \cong \zeta \circ (\mu \cdot p)$ & $\prod_{x:A} \prod_{y : B(x)} C(x,y) \cong \prod_{\langle x,y \rangle : \sum_{x:A} B(x)} C(x,y)$ \\
&& \\ \hline
&& \\
Algebra unit & $\zeta \circ (\eta \cdot p) \cong \mathrm{id}_p$ & $\prod_{x:\mathbf{1}} A \cong A$ \\
&&
\end{tabular}
\end{center}
\end{numbered}


%% file: thesis/ch3-polynomials-representability/representability.tex
\section{Representability revisited}
\label{secCharacterisationsOfRepresentability}
\index{representable!natural transformation@{---} natural transformation}

\subsection*{Representability and cocontinuity}

The first goal of this section is to identify a condition for a natural transformation $p : \nmtm{U} \to \nmty{U}$ between presheaves over a small category $\mathbb{C}$ to be representable in terms of cocontinuity of its polynomial extension $\upP_p : \widehat{\mathbb{C}} \to \widehat{\mathbb{C}}$. We will prove that if $p$ is representable, then $\upP_p$ is cocontinuous. For the converse, we will need to assume some conditions on the base category $\mathbb{C}$, namely that it is Cauchy complete and has finite products.

We begin with a discussion of Cauchy completeness and its relation to so-called \textit{tiny} objects in presheaf categories. This matter is confused somewhat by the fact that there are different notions of tininess and, even more confusingly, different words have been used by different authors to refer to the same notion of tininess, and different notions of tininess have been referred to by different authors by the same word! With this in mind, we will first fix our own terminology.

\begin{definition}[Tiny and atomic objects]
\index{tiny object}
\index{atomic object}
Let $\mathcal{E}$ be a locally small, cocomplete, locally cartesian closed category and let $X$ be an object of $\mathcal{E}$. Then:
\begin{enumerate}[(i)]
\item $X$ is \textbf{internally atomic} if $(-)^X : \mathcal{E} \to \mathcal{E}$ has a right adjoint;
\item $X$ is \textbf{internally tiny} if $(-)^X : \mathcal{E} \to \mathcal{E}$ is cocontinuous;
\item $X$ is \textbf{externally atomic} if $\mathcal{E}(X,{-}) : \mathcal{E} \to \mathbf{Set}$ has a right adjoint;
\item $X$ is \textbf{externally tiny} if $\mathcal{E}(X,{-}) : \mathcal{E} \to \mathbf{Set}$ is cocontinuous.
\end{enumerate}
\end{definition}

Condition (i) goes back to William Lawvere, who referred to the right adjoint to the functor $(-)^X$ as the \textit{amazing right adjoint}. An object satisfying (i) is called \textit{infinitesimal} by Lawvere \cite{Lawvere1980Dynamics}, \textit{atomic} by Anders Kock \cite{Kock1981SDG} and \textit{tiny} by David Yetter \cite{Yetter1987RightAdjoints}. An object satisfying condition (iv) is called an \textit{atom} by Marta Bunge \cite{Bunge1965Thesis} and is called \textit{small-projective} by Max Kelly \cite{Kelly1982BasicConcepts}.

\begin{numbered}
When $\mathcal{E} = \widehat{\mathbb{C}}$ for some small category $\mathbb{C}$, it follows from Freyd's adjoint functor theorem that each of the functors $(-)^X$ and $\widehat{\mathbb{C}}(X,{-})$ is cocontinuous if and only if it has a right adjoint. Thus a presheaf is internally atomic (in $\widehat{\mathbb{C}}$) if and only if it is internally tiny, and is externally atomic if and only if it is externally tiny. For this reason, since we will focus on presheaves from now on, we will simply use the term \textit{internally tiny} to refer to conditions (i) and (ii) together, and \textit{externally tiny} to refer to conditions (iii) and (iv) together.
\end{numbered}

\newpage

\begin{definition}[Cauchy complete category, \thmcite{Borceux1986CauchyCompletion}]
\index{Cauchy complete category}
\index{category!Cauchy complete {---}}
A category $\mathcal{C}$ is \textbf{Cauchy complete} if every idempotent in $\mathcal{C}$ splits---that is, if for each $e : A \to A$ in $\mathcal{C}$ such that $e \circ e = e$, there exists a factorisation of $e$ in $\mathcal{C}$ as $A \xrightarrow{r} B \xrightarrow{i} A$ such that $r \circ i = \mathrm{id}_B$.
\end{definition}

\begin{numbered}
\label{parTinyObjects}
We briefly recall some results from \cite[\S{}5.5]{Kelly1982BasicConcepts}, \cite{Borceux1986CauchyCompletion} and \cite{Yetter1987RightAdjoints} concerning tiny objects and Cauchy completions. A presheaf $X : \mathbb{C}\op \to \mathbf{Set}$ over a small category $\mathbb{C}$ is externally tiny in $\widehat{\mathbb{C}}$ if and only if it is a retract of a representable functor. Writing $\bar{\mathbb{C}}$ for the full subcategory of $\widehat{\mathbb{C}}$ whose objects are the retracts of representable functors, we have that $\bar{\mathbb{C}}$ is small and the Yoneda embedding $\Yon : \mathbb{C} \hookrightarrow \widehat{\mathbb{C}}$ factors through the embedding $\mathsf{k} : \mathbb{C} \hookrightarrow \bar{\mathbb{C}}$. The category $\bar{\mathbb{C}}$ is the Cauchy completion of $\mathbb{C}$; moreover, the functor $\mathsf{k}^* : \widehat{\bar{\mathbb{C}}} \to \widehat{\mathbb{C}}$ is an equivalence of categories, and if $\mathbb{C}$ has finite products then so does $\bar{\mathbb{C}}$. A category $\mathbb{C}$ is Cauchy complete if and only if $\mathsf{k}$ itself is an equivalence. In particular, a category $\mathbb{C}$ is Cauchy complete if and only if the externally tiny objects of $\widehat{\mathbb{C}}$ are exactly the representable functors. If $\mathbb{C}$ is Cauchy complete and has a terminal object, then externally tiny objects are internally tiny; the converse holds if $\mathbb{C}$ has finite products.
\end{numbered}

\begin{theorem}[Characterisation of representability by cocontinuity]
\label{thmRepresentabilityAndContinuityOfExtension}
Let $\mathbb{C}$ be a small category and let $p : \nmtm{U} \to \nmty{U}$ be a natural transformation between presheaves over $\mathbb{C}$.
\begin{enumerate}[(a)]
\item If $p$ is representable, then its extension $\upP_p : \widehat{\mathbb{C}} \to \widehat{\mathbb{C}}$ is cocontinuous.
\item If the extension $\upP_p : \widehat{\mathbb{C}} \to \widehat{\mathbb{C}}$ of $p$ is cocontinuous and $\mathbb{C}$ is Cauchy complete with finite products, then $p$ is representable.
\end{enumerate}
\end{theorem}

\begin{proof}
By \Cref{lemLemmaFive}, for each $\Gamma \in \mathrm{ob}(\mathbb{C})$, there is a natural (in $X$ and in $\Gamma$) bijection
$$\upP_p(X)(\Gamma) = \left( \sum_{A \in \nmty{U}} X^{[A]} \right)(\Gamma) \cong \sum_{A \in \nmty{U}(\Gamma)} \widehat{\mathbb{C}}(\upDelta_A(\nmtm{U}), X)$$

If $p$ is representable, then $\upDelta_A(\nmtm{U}) \cong \Yon(\Gamma \cext A)$ for some object $\Gamma \cext A$ of $\mathbb{C}$, so $\upDelta_A(\nmtm{U})$ is externally tiny and we see that $\upP_p$ preserves colimits. This proves (a).

Conversely, if $\upP_p$ preserves colimits then so does $\upPi_p$, so that $\upDelta_A(\nmtm{U})$ is internally tiny. If $\mathbb{C}$ is Cauchy complete and has finite products, then as discussed above we have that $\upDelta_A(\nmtm{U}) \cong \Yon(\Gamma \cext A)$ for some object $\Gamma \cext A$ of $\mathbb{C}$, so that $p$ is representable. This proves (b).
\end{proof}



\begin{corollary}
Let $\mathbb{C}$ be a small category with finite products and let $p : \nmtm{U} \to \nmty{U}$ be a natural transformation between presheaves over $\mathbb{C}$ such that $\upP_p : \widehat{\mathbb{C}} \to \widehat{\mathbb{C}}$ is cocontinuous. By transporting $p$ along the equivalence $\widehat{\mathbb{C}} \simeq \widehat{\bar{\mathbb{C}}}$ discussed in \Cref{parTinyObjects}, we obtain a \textit{representable} natural transformation $\bar p$ between presheaves over $\bar{\mathbb{C}}$.
\end{corollary}

\begin{proof}
Note that $\upP_{\bar p} : \widehat{\bar{\mathbb{C}}} \to \bar{\mathbb{C}}$ is cocontinuous since $\mathsf{k}^* : \widehat{\bar{\mathbb{C}}} \to \widehat{\mathbb{C}}$ is an equivalence. Since $\mathbb{C}$ has finite products, so does $\bar{\mathbb{C}}$, and since $\bar{\mathbb{C}}$ is Cauchy complete, it follows from \Cref{thmRepresentabilityAndContinuityOfExtension} that $\bar p$ is representable.
\end{proof}

\begin{numbered}
In \cite[Theorem 5.26]{Kelly1982BasicConcepts} it is proved that a category $\mathcal{E}$ is equivalent to the category $\widehat{\mathbb{C}}$ of presheaves on a small category $\mathbb{C}$ if and only if $\mathcal{E}$ is cocomplete and there is a small set of tiny objects constituting a strong generator of $\mathcal{E}$. The category $\mathbb{C}$ is obtained as the full subcategory of $\mathcal{E}$ determined by this set of tiny objects.
\end{numbered}

\subsection*{Representability and full internal categories}

\begin{numbered}
\label{parComprehensionCategoryFromFullInternalSubcategory}
We recall from \cite[\S{}4]{Jacobs1993ComprehensionCategories} and \cite[\S{}7]{Jacobs1999CategoricalLogic} some facts about full internal subcategries. Given any morphism $f : B \to A$ of a locally cartesian closed category $\mathcal{E}$, the full internal subcategory $\fisc{f}$ of $\mathcal{E}$ (\Cref{cnsFullInternalSubcategory}) gives rise to a fibration $\mathbb{E}(f) \to \mathcal{E}$ together with a full and faithful cartesian functor $\mathbb{E}(f) \to \mathcal{E}^{\to}$ over $\mathcal{E}$.
\begin{diagram}
\mathbb{E}(f)
\arrow[rr]
\arrow[dr, "\pi"']
&[-20pt]&[-20pt]
\mathcal{E}^{\to}
\arrow[dl, "\mathsf{cod}"]
\\
&
\mathcal{E}
&
\end{diagram}
The category $\mathbb{E}(f)$ is the \textit{externalisation} of $\mathbb{S}(f)$, which can be described as follows.
\begin{itemize}
\item The objects of $\mathbb{E}(f)$ are morphisms $x : X \to \fisc{f}_0 = A$ in $\mathcal{E}$;
\item Given $x : X \to A$ and $y : Y \to A$, a morphism from $x$ to $y$ in $\mathbb{E}(p)$ is a pair $(\sigma, h)$ consisting of a morphism $\sigma : X \to Y$ in $\mathcal{E}$ and a morphism $\upDelta_x(f) \to \upDelta_{y \circ \sigma}(f)$ in $\mathcal{E} \slice{X}$.
\end{itemize}
The fibration $\pi : \mathbb{E}(f) \to \mathcal{E}$ sends each object $x : X \to A$ to its domain $X$ and each morphism $(\sigma, h)$ to its first component $\sigma$; and the cartesian functor $\mathbb{E}(f) \to \mathcal{E}^{\to}$ sends an object $x : X \to A$ of $\mathbb{E}(f)$ to the morphism $\upDelta_f(x) : \upDelta_f(B) \to X$ and a morphism $(\sigma, h) : x \to y$ to the square described by $\sigma$ and $h$.

When $p : \nmtm{U} \to \nmty{U}$ is a natural transformation between presheaves over a small category $\mathbb{C}$, this construction gives rise to a full and faithful fibred functor $\mathbb{E}(p) \to \widehat{\mathbb{C}}^{\to}$ over $\widehat{\mathbb{C}}$. The fibre $\mathbb{E}(p)_{\Yon(\Gamma)}$ over a representable presheaf $\Yon(\Gamma)$ is then exactly the (external) category $\mathbb{S}(p)(\Gamma)$ obtained by applying the data defining $\mathbb{S}(p) \in \mathbf{Cat}(\widehat{\mathbb{C}})$ to the object $\Gamma$ of $\mathbb{C}$. By pulling back $\pi : \mathbb{E}(p) \to \widehat{\mathbb{C}}$ along the Yoneda embedding $\Yon : \mathbb{C} \to \widehat{\mathbb{C}}$, we obtain a fibration $\pi' : \mathbb{E}'(p) \to \mathbb{C}$, where $\mathbb{E}'(p)_{\Gamma} = \mathbb{E}(p)_{\Yon(\Gamma)}$. By abuse of notation, we will simply write $\pi : \mathbb{E}(p) \to \mathbb{C}$ for this fibration.
\end{numbered}

We prove in \Cref{thmRepresentabilityFromFISC} that representability of $p$ can be characterised by the existence of a full and faithful fibred functor (over $\mathbb{C}$) from $\mathbb{E}(p)$ to the full subcategory of $\mathbb{C}^{\to}$ determined by the morphisms in $\mathbb{C}$ classified by $p$ (\Cref{defMorphismClassifiedByNM}). First, we remark that the codomain fibration restricts to this subcategory.

\begin{definition}
\label{defMorphismClassifiedByNM}
\index{morphism!classified@{---} classified by a natural transformation}
Let $p : Y \to X$ be a natural transformation between presheaves over a small category $\mathbb{C}$. We say a morphism $\sigma : \Delta \to \Gamma$ of $\mathbb{C}$ is \textbf{classified by $p$} if $\Yon(\sigma) : \Yon(\Delta) \to \Yon(\Gamma)$ arises as a pullback of $p$ in $\widehat{\mathbb{C}}$. Write $\mathcal{F}_p$ to denote both the set of morphisms of $\mathbb{C}$ classified by $p$, and the corresponding full subcategory of $\mathbb{C}^{\to}$.
\end{definition}

\begin{lemma}
\label{lemClassificationOfMorphismsClassifiedByNM}
Let $(\mathbb{C}, p)$ be a natural model and let $\sigma : \Gamma' \to \Gamma$ in $\mathbb{C}$. Then $\sigma \in \mathcal{F}_p$ if and only if there is an isomorphism $(\Gamma \cext A, \nmp{A}) \to (\Gamma',\sigma)$ in $\mathbb{C} \slice{\Gamma}$ for some $A \in \nmty{U}(\Gamma)$.
\end{lemma}

\begin{proof}
The morphism $\sigma$ is classified by $p$ if and only if there exist $A \in \nmty{U}(\Gamma)$ and $a \in \nmtm{U}(\Delta, A[\sigma])$ such that the following square is a pullback.
\begin{diagram}
\Yon(\Delta)
\arrow[r, "a"]
\arrow[d, "\Yon(\sigma)"']
\pullback
&
\nmtm{U}
\arrow[d, "p"]
\\
\Yon(\Gamma)
\arrow[r, "A"']
&
\nmty{U}
\end{diagram}
By representability of $p$, the natural transformation $\Yon(\nmp{A})$ is also a pullback of $p$, hence the result follows from the universal property of pullbacks together with the fact that the Yoneda embedding is full and faithful.
\end{proof}

\begin{lemma}[Classified morphisms yield a fibration]
Let $\mathbb{C}$ be a small category and $p : \nmtm{U} \to \nmty{U}$ be a natural transformation between presheaves over $\mathbb{C}$. The codomain functor $\mathcal{F}_p \to \mathbb{C}$ is a Grothendieck fibration.
\end{lemma}

\begin{proof}
It suffices to prove that pullbacks of morphisms in $\mathcal{F}_p$ along arbitrary morphisms of $\mathbb{C}$ exist and are in $\mathcal{F}_p$. To this end, let $\gamma : \Gamma' \to \Gamma \in \mathcal{F}_p$ and let $\sigma : \Delta \to \Gamma$ be a morphism in $\mathbb{C}$. By \Cref{lemClassificationOfMorphismsClassifiedByNM}, there is an isomorphism $\theta : (\Gamma \cext A, \nmp{A}) \cong (\Gamma', \sigma)$ in $\mathbb{C} \slice{\Gamma}$. Now consider the following diagram
\begin{diagram}
\Delta \cext A[\sigma]
\arrow[r, "\sigma \cext A"]
\arrow[d, "\nmp{\sigma}"']
\pullback
&
\Gamma \cext A
\arrow[r, "\theta", "\cong"']
\arrow[d, "\nmp{A}" description]
\pullback
&
\Gamma'
\arrow[d, "\gamma"]
\\
\Delta
\arrow[r, "\sigma"']
&
\Gamma
\arrow[r, equals]
&
\Gamma
\end{diagram}
The square on the left is a canonical pullback square (\Cref{cnsCanonicalPullbacks}) and the square on the right is a pullback since it commutes and $\theta$ is an isomorphism, so that the outer square is a pullback. But then $(\sigma, \theta \circ (\sigma \cext A))$ is a cartesian lift of $\sigma$.
\end{proof}

\begin{theorem}[Characterisation of representability from full internal subcategories]
\label{thmRepresentabilityFromFISC}
Let $\mathbb{C}$ be a small category and let $p : \nmtm{U} \to \nmty{U}$ be a natural transformation between presheaves over $\mathbb{C}$. Then $p$ is representable if and only if there is a fibred equivalence $\chi : \mathbb{E}(p) \to \mathcal{F}_p$ over $\mathbb{C}$.
\begin{diagram}
\mathbb{E}(p)
\arrow[dr, "\pi"']
\arrow[rr, "\chi", "\simeq"']
&&
\mathcal{F}_p
\arrow[dl, "\mathsf{cod}"]
\\
&
\mathbb{C}
&
\end{diagram}
\end{theorem}

\begin{proof}
($\Leftarrow$) Suppose there is a fibred equivalence $\chi : \mathbb{E}(p) \to \mathcal{F}_p$, and let $\Gamma \in \mathrm{ob}(\mathbb{C})$ and $A \in \nmty{U}(\Gamma)$. Then $A \in \mathrm{ob}(\mathbb{E}(p)_{\Gamma})$, so that $\chi_{\Gamma}(A) \in (\mathcal{F}_p)_{\Gamma}$. Define $\Gamma \cext A = \mathrm{dom}(\chi_{\Gamma}(A))$ and $\nmp{A} = \chi_{\Gamma}(A)$. Since $\nmp{A} \in \mathcal{F}_p$, there is a morphism $\nmq{A} : \Yon(\Gamma \cext A) \to \nmtm{U}$ making the following square a pullback
\begin{diagram}
\Yon(\Gamma \cext A)
\arrow[r, "\nmq{A}"]
\arrow[d, "\nmp{A}"']
&
\nmtm{U}
\arrow[d, "p"]
\\
\Yon(\Gamma)
\arrow[r, "A"']
&
\nmty{U}
\end{diagram}
This demonstrates that $p$ is representable, and that specifying $\chi$ gives rise to representability data for a natural model $(\mathbb{C}, p)$.

($\Rightarrow$) Suppose $p$ is representable. Representability of $p$ is equivalent to the existence, for each object $\Gamma$ of $\mathbb{C}$ and element $A \in \nmty{U}(\Gamma)$, of a morphism $\nmp{A} : \Gamma \cext A \to \Gamma$ of $\mathbb{C}$ such that $\Yon(\nmp{A})$ is a pullback of $p$. Let the action of $\chi$ on objects choose such a morphism for each pair $(\Gamma, A)$. By \Cref{lemMorphismsOfAssociatedFISC} and by definition of $\mathbb{E}(p)$, morphisms from $A \in \nmty{U}(\Delta)$ to $B \in \nmty{U}(\Gamma)$ in $\mathbb{E}(p)$ correspond naturally with pairs $(\sigma, \tau)$, where $\sigma : \Delta \to \Gamma$ is a morphism of $\mathbb{C}$ and $\tau$ is a morphism from $(\Gamma \cext A[\sigma], \nmp{A[\sigma]})$ to $(\Gamma \cext B, \nmp{A})$ in $\mathbb{C} \slice{\Gamma}$; but this correspondence precisely gives us the action of $\chi$ on morphisms and tells us that it is full and faithful (since the correspondence is bijective) and fibred over $\Gamma$ (since the codomain is respected). Moreover, $\chi$ is essentially surjective: given $\sigma : \Gamma' \to \Gamma \in \mathrm{ob}(\mathcal{F}_p)$, it follows from \Cref{lemClassificationOfMorphismsClassifiedByNM} that for some $A \in \nmty{U}(\Gamma)$ we have $\sigma \cong \chi(\Gamma, A)$ in $\mathcal{F}_p$.
\end{proof}

\begin{corollary}
Specifying a natural model $(\mathbb{C}, p)$ is equivalent to specifying a category $\mathbb{C}$ with a terminal object $\diamond$, a natural transformation $p : \nmtm{U} \to \nmty{U}$ between presheaves over $\mathbb{C}$, and a fibred equivalence $\mathbb{E}(p) \to \mathcal{F}_p$ over $\mathbb{C}$. \qed
\end{corollary}

\subsection*{Closure properties of representable natural transformations}

\begin{theorem}[Closure properties of representable natural transformations]
\label{thmClosurePropertiesOfRepresentability}
Let $\mathbb{C}$ be a small category. The class $\mathcal{R} \subseteq \mathbb{C}^{\to}$ of all representable natural transformations over $\mathbb{C}$ enjoys the following closure properties:
\begin{enumerate}[(a)]
\item $\mathcal{R}$ is closed under composition in $\widehat{\mathbb{C}}$;
\item $\mathcal{R}$ is closed under pullbacks (in $\widehat{\mathbb{C}}$) along arbitrary morphisms of $\widehat{\mathbb{C}}$;
\item $\mathcal{R}$ is closed under polynomial composition in $\widehat{\mathbb{C}}$;
\item $\mathcal{R}$ is closed under isomorphisms in $\widehat{\mathbb{C}}^{\to}$;
\item $\mathcal{R}$ is closed under (small) coproducts in $\widehat{\mathbb{C}}^{\to}$.
\end{enumerate}
\end{theorem}

\begin{proof} ~
\begin{enumerate}[(a)]
\item 
Let $p : Y \to X$ and $q : Z \to Y$ be representable natural transformations and let $C \in \mathrm{ob}(\mathbb{C})$ and $x \in X(C)$. Using representability of $p$ and of $q$, construct the following diagram in which the top and bottom squares are pullbacks.

\begin{diagram}
\Yon(E)
\arrow[r, "z"]
\arrow[d, "g"']
\pullback
&
Z
\arrow[d, "q"]
\\
\Yon(D)
\arrow[r, "y" description]
\arrow[d, "f"']
\pullback
&
Y
\arrow[d, "p"]
\\
\Yon(C)
\arrow[r, "x"']
&
X
\end{diagram}
By the two pullbacks lemma, the outer square is a pullback, so that $p \circ q$ is representable.

\item 
Let $p : Y \to X$ be a representable natural transformation and let $f,g,q$ be natural transformations fitting into the following pullback square.

\begin{diagram}
Y'
\arrow[r, "g"]
\arrow[d, "q"']
\pullback
&
Y
\arrow[d, "p"]
\\
X'
\arrow[r, "f"']
&
X
\end{diagram}

Let $C \in \mathrm{ob}(\mathbb{C})$ and $x \in X'(C)$. Then $f_C(x) \in X(C)$, so there exist $D \in \mathrm{ob}(\mathbb{D})$, $y \in Y(D)$ and $f : D \to C$ in $\mathbb{C}$ making the outer square of the following diagram a pullback.

\begin{diagram}
\Yon(D)
\arrow[d, "\Yon(f)"']
\arrow[r, dashed, "y'" description]
\arrow[rr, bend left=15, "y"]
&
Y'
\arrow[r, "g" description]
\arrow[d, "q"']
\pullback
&
Y
\arrow[d, "p"]
\\
\Yon(C)
\arrow[r, "x"']
&
X'
\arrow[r, "f"']
&
X
\end{diagram}

By the universal property of pullbacks, there is an element $y' \in Y'(D)$ fitting into the diagram as indicated with the dashed morphism, making the left square a pullback by the two pullbacks lemma. Hence $q$ is representable.

\item 
Let $p : Y \to X$ and $q : V \to U$ be representable natural transformations. As described in \Cref{defCompositionOfPolynomials}, the polynomial composite $q \cdot p$ obtained by composing a pullback of $p$ with a pullback of $q$; by parts (a) and (b), it follows that $q \cdot p$ is representable.

\item 

An isomorphism in $\widehat{\mathbb{C}}^{\to}$ is, in particular, a pullback square, so if $p : Y \to X$ is representable and $q : V \to U$ is isomorphic to $p$ in $\widehat{\mathbb{C}}^{\to}$, then $q$ is representable by part (b).

\item 
Let $I$ be a set and let $\seqbn{p_i : Y_i \to X_i}{i \in I}$ be an $I$-indexed family of representable natural transformations. Let $p : Y \to X$ be their coproduct, i.e.
$$p = \sum_{i \in I} p_i : \sum_{i \in I} Y_i \to \sum_{i \in I} X_i$$
Let $C \in \mathrm{ob}(\mathbb{C})$ and let $x \in X(C)$. Then $x = (i,x')$ for some $i \in I$ and $x' \in X_i(C)$. By representability of $p_i$, there exist $D \in \mathrm{ob}(\mathbb{C})$, $y' \in Y_i(D)$ and $f : D \to C$ in $\mathbb{C}$ making the left square in the following diagram a pullback.

\begin{diagram}
\Yon(D)
\arrow[r, "y'"]
\arrow[d, "\Yon(f)"']
\pullback
&
Y_i
\arrow[d, "p_i" description]
\arrow[r, "\iota_i"]
&
Y
\arrow[d, "p"]
\\
\Yon(C)
\arrow[r, "x'" description]
\arrow[rr, bend right=15, "x"']
&
X_i
\arrow[r, "\iota_i" description]
&
X
\end{diagram}
Checking that the outer square is a pullback is a straightforward verification of the universal property of pullbacks. Hence $p$ is representable.
\end{enumerate}
\end{proof}

We know by \Cref{thmClosurePropertiesOfRepresentability}(c) that the composite of two representable natural transformations is representable. By chasing the representability data through the respective proofs that pullbacks and composites of representable natural transformations are representable, we obtain the following construction of the \textit{polynomial composite} of natural models.

\begin{construction}[Polynomial composite of natural models]
\label{cnsPolynomialCompositeOfNaturalModels}
\index{polynomial!composite of natural models@{---} composite of natural models}
\index{natural model!polynomial composite of {---}s}
Fix a small category $\mathbb{C}$. The \textbf{polynomial composite} of natural models $(\mathbb{C}, p)$ and $(\mathbb{C}, q)$ is the natural model $(\mathbb{C}, q \cdot p)$ with representability data given by
\begin{itemize}
\item \textbf{Context extension.} The extension of $\Gamma \in \mathrm{ob}(\mathbb{C})$ by $(A,B) \in \sum_{A : \nmty{V}} \nmty{U}^{\langle A \rangle}$ is given by $(\Gamma \cextalt A) \cext B$;
\item \textbf{Projection.} The projection $(\Gamma \cextalt A) \cext B \to \Gamma$ is given by the composite
$$\nmu{A} \circ \nmp{B} : (\Gamma \cextalt A) \cext B \xrightarrow{\nmp[\Gamma \cextalt A]{B}} \Gamma \cextalt A \xrightarrow{\nmu[\Gamma]{A}} \Gamma$$
\item \textbf{Variable.} The new variable term in context $(\Gamma \cextalt A) \cext B$ is $(A, B, \nmv[\Gamma]{A}, \nmq[\Gamma \cextalt A]{B})$.
\end{itemize}
\end{construction}

\begin{verification}
Using \Cref{defCompositionOfPolynomials}, we can express the polynomial composite $q \cdot p$ as the regular composite $r \circ r'$ of morphisms in $\mathcal{E}$ indicated in the following commutative diagram, in which the unlabelled arrows are the appropriate projection morphisms.

\begin{ldiagram}
&
\sum_{A,B} \sum_{a \in \langle A \rangle} [B(a)]
\arrow[ddl]
\arrow[r, "r'"]
&
\sum_{A,B} \langle A \rangle
\arrow[r, "r"]
\arrow[d]
&
\sum_{A : \nmty{V}} \nmty{U}^{\langle A \rangle}
\arrow[ddr]
&
\\
&&
\nmty{U} \times \nmtmalt{V}
\arrow[dr]
\arrow[dl]
&&
\\
\nmtm{U}
\arrow[r,"p"']
&
\nmty{U}
&&
\nmtmalt{V}
\arrow[r,"q"']
&
\nmty{V}
\end{ldiagram}
Note that the left- and right-hand `squares' are cartesian, exhibiting $r$ as a pullback of $q$ and $r'$ as a pullback of $p$.

We now proceed chase the representability data of $p$ and of $q$ through the proofs that pullbacks and composites of representable natural transformations are representable. To this end, let $\Gamma \in \mathrm{ob}(\mathbb{C})$ and let $(A,B) \in \left( \sum_{A : \nmty{V}} \nmty{U}^{\langle A \rangle} \right)(\Gamma)$.

Since $q$ is representable and the right-hand square of the above diagram is cartesian, we may form the following pasting of pullback squares:
\begin{diagram}
\Yon(\Gamma \cextalt A)
\arrow[r, "{(A, B, \nmv{A})}" description]
\arrow[rr, blue, bend left=15, "\nmv{A}"]
\arrow[d, "\Yon(\nmu{A})" description]
\pullback
&
\sum_{A,B} \langle A \rangle
\arrow[r]
\arrow[d, "r" description]
\pullbackc{dr}{0pt}
&
\nmtmalt{V}
\arrow[d, "q"]
\\
\Yon(\Gamma)
\arrow[r, "{(A,B)}" description]
\arrow[rr, blue, bend right=20, "A"']
&[20pt]
\sum_{A:\nmty{V}} \nmty{U}^{\langle A \rangle}
\arrow[r]
&
\nmty{V}
\end{diagram}

and since the left hand square is cartesian, we may form the following pasting diagram of pullback squares:

\begin{diagram}
\Yon((\Gamma \cextalt A) \cext B)
\arrow[r, "{(A,B,\nmv{A},\nmq{B})}" description]
\arrow[rr, blue, bend left=15, "\nmq{B}"]
\arrow[d, "\Yon(\nmp{B})"']
\pullback
&[20pt]
\sum_{A,B} \sum_{a \in \langle A \rangle} [B(a)]
\arrow[r]
\arrow[d, "r'" description]
\pullbackc{dr}{-0.1}
&
\nmtm{U}
\arrow[d, "p"]
\\
\Yon(\Gamma \cextalt A)
\arrow[r, "{(A, B, \nmv{A})}" description]
\arrow[rr, blue, bend right=20, "A"']
&
\sum_{A,B} \langle A \rangle
\arrow[r]
&
\nmty{U}
\end{diagram}

Pasting the left-hand squares of the previous two diagrams vertically, we obtain the desired pullback square.

\begin{diagram}
\Yon((\Gamma \cextalt A) \cext B)
\arrow[r, "{(A, B, \nmv{A}, \nmq{B})}"]
\arrow[d, "\Yon(\nmu{A} \circ \nmq{B})"']
\pullback
&
\sum_{A,B} \sum_{a \in \langle A \rangle} [B(a)]
\arrow[d, "q \cdot p"]
\\
\Yon(\Gamma)
\arrow[r, "{(A,B)}"']
&
\sum_{A \in \nmty{V}} \nmtm{U}^{\langle A \rangle}
\end{diagram}

This proves that $q \cdot p$ is representable, with representability data as indicated in the statement of this construction, so that $(\mathbb{C}, q \cdot p)$ is a natural model.
\end{verification}

%% file: thesis/ch4-semantics/_semantics.tex
\newpage
\input{thesis/ch4-semantics/interpretations-initiality.tex}

\newpage
\input{thesis/ch4-semantics/ext-term.tex}

\newpage
\input{thesis/ch4-semantics/ext-type.tex}

\newpage
\input{thesis/ch4-semantics/ext-unit.tex}

\newpage
\input{thesis/ch4-semantics/ext-sigma.tex}

%% file: thesis/ch4-semantics/interpretations-initiality.tex
\section{Free natural models}
\label{secInterpretationsInitiality}

In \Cref{secEssentiallyAlgebraicTheory}, we saw that the various theories of natural models equipped with type theoretic structure are essentially algebraic. As we discussed at the end of that section, we can use this fact to apply the machinery of essentially algebraic categories and locally presentable categories to categories of the form $\mathbf{NM}_{\mathbb{T}}$ for a given dependent type theory $\mathbb{T}$---for example, such categories have initial objects (since they are cocomplete), and there are free--forgetful adjunctions between such categories.

The goal of this chapter is to explicitly describe the initial object of $\mathbf{NM}_{\mathbb{T}}$ for a couple of basic dependent type theories $\mathbb{T}$, and to provide an explicit description of some of these left adjoints to forgetful functors.

In this section, we construct the free natural model on a set of basic types---or, more precisely, on an indexed family of elements of $\nmty{U}(\diamond)$. In the subsequent sections, we describe how to freely add type theoretic structure to a natural model without additional type theoretic structure.

\subsection*{A basic free model}

We construct the free model on the theory $\mathbb{T}_{(\mathsf{ty}_i)_{i \in I}}$ of an $I$-indexed family of basic types, where $I$ is a fixed set.

\begin{construction}[Free model of the theory of a family of basic types]
\label{cnsFreeNaturalModelKappaBasicTypes}
The \textbf{free natural model on the theory $\mathbb{T}_{(\mathsf{ty}_i)_{i \in I}}$} is the natural model $(\mathbb{C}_I, p_I)$ given by the following data:
\begin{itemize}
\item The underlying category $\mathbb{C}_{I}$ is $(\mathbf{Fin}\slice{I})^{\mathrm{op}}$.

\item The (only, and hence) chosen terminal object of $\mathbb{C}_{I}$ is $(\varnothing, {!}_I : \varnothing \to I)$.

\item The presheaf of types $\nmty{U}_{I} : \mathbb{C}_{I}^{\mathrm{op}} \to \mathbf{Set}$ is the codomain functor $\mathbf{Fin} \slice{I} \to \mathbf{Set}$, or equivalently the constant functor $\upDelta(I)$ with value $I$. Thus $\nmty{U}_{I}(A, u) = I$ for all $(A, u)$ and $\nmty{U}_{I}(f) = \mathrm{id}_{I}$ for all $f$.

\item The presheaf of terms $\nmtm{U}_I$ is the domain functor $\mathsf{dom} : \mathbf{Fin} \slice{I} \to \mathbf{Set}$. Thus $\nmtm{U}_I(A,u) = A$ for all $(A,u)$ and $\nmtm{U}_I(f) = f$ for all $f$.

\item The natural transformation $p_{I} : \nmtm{U}_{I} \to \nmty{U}_{I}$ is given simply by $(p_I)_{(A,u)} = u : A \to I$.

\item The representability data is defined as follows. Given $(A, u)$ and $j \in I$, we define
\begin{itemize}
\item $(A, u) \cext j = (A + 1, [u,j])$---thus $(A+1)_i = A_i$ if $i \ne j$ and $(A+1)_i = A_j+1$ if $i=j$;
\item $\nmp{j} : (A, u) \cext j \to (A, u)$ in $\mathbb{C}_I$ is given by the left inclusion function $A \to A+1$ in $\mathbf{Fin}\slice{I}$.
\item $\nmq{j} \in \nmtm{U}_{I}((A, u) \cext j) = A+1$ is the added element $\star \in A+1$.
\end{itemize}
\end{itemize}
\end{construction}

\begin{verification}
The only part of the verification that is not immediate is representability of $p_{I}$ as witnessed by the given representability data.

The set-up is as follows. Take $(A, u) \in \mathbb{C}_I$ and $j \in I$. We need to prove that the following square is a pullback
\begin{center}
\begin{tikzcd}[row sep=huge, column sep=huge]
\Yon((A, u) \cext j)
\arrow[d, "\nmp{j}"']
\arrow[r, "\nmq{j}"]
&
\nmtm{U}_{I}
\arrow[d, "p_{I}"]
\\
\Yon(A, u)
\arrow[r, "j"']
&
\nmty{U}_{I}
\end{tikzcd}
\end{center}
It evidently commutes, so it suffices to check the universal property on representables.

To this end, let $(B, v) \in \mathbb{C}_I$, let $f : (B, v) \to (A, u)$ in $\mathbb{C}_{I}$ (so that $f$ is a function $A \to B$ over $I$) and let $b : \Yon(B, v) \to \nmtm{U}_{I}$, and suppose that $p_{I} \circ b = j \circ \Yon(f)$.

\begin{center}
\begin{tikzcd}[row sep=huge, column sep=huge]
\Yon(B, v)
\arrow[ddr, "\Yon(f)"', bend right=20]
\arrow[drr, "b", bend left=20]
\arrow[dr, dashed, blue]
&[-20pt]&
\\[-20pt]
&
\Yon((A, u) \cext j)
\arrow[d, "\Yon(\nmp{j})"']
\arrow[r, "\nmq{j}"]
&
\nmtm{U}_{I}
\arrow[d, "p_{I}"]
\\
&\Yon(A, u)
\arrow[r, "j"']
&
\nmty{U}_{I}
\end{tikzcd}
\end{center}

By the Yoneda lemma, $b$ is an element of $B$, and commutativity of the outer square says that
$$v(b) = (p_I)_{(B,v)}(b) = \nmty{U}_I(f)(j) = j$$

We need to prove that there is a unique $g : (B, v) \to (A, u) \cext j$ in $\mathbb{C}_I$ such that $\nmq{j} \circ \Yon(g) = b$ and $\nmp{j} \circ g = f$. Now:

\begin{itemize}
\item As a morphism in $\mathbb{C}_I$, the map $g$ must be a function $A+1 \to B$ over $I$, which is equivalent to saying that $g = [g', b']$ for some function $g' : A \to B$ over $I$ and some element $b' \in B_j$;
\item The requirement that $\nmp{j} \circ g = f$ in $\mathbb{C}_i$ is equivalent to the requirement that $g'=f$;
\item The requirement that $\nmq{j} \circ \Yon(g) = b$ is equivalent to the requirement $b'=b$.
\end{itemize}

So $g = [f,b] : (B,v) \to (A,u) \cext j$ is the unique morphism satisfying the required conditions. Hence the square is a pullback, and so the representability data of \Cref{cnsFreeNaturalModelKappaBasicTypes} truly does witness representability of $p_{I}$.
\end{verification}

\begin{example}
Take $I=0$. The category $\mathbb{C}_0$ is the terminal category $\mathbf{1}$; the presheaves $\nmty{U}_0, \nmtm{U}_0$ are empty and the natural transformation $p_0 : \nmtm{U}_0 \to \nmty{U}_0$ is the empty natural transformation.
\end{example}

\begin{example}
Take $I=1$. Then $\mathbb{C}_1 \simeq \mathbf{Fin}^{\mathrm{op}}$, which is the free category with finite limits on one object. The presheaf $\nmty{U}_1$ is the constant presheaf $\upDelta(1)$ and the presheaf $\nmtm{U}_1$ is the inclusion $\mathbf{Fin} \hookrightarrow \mathbf{Set}$. The natural transformation $p_1 : \nmtm{U}_1 \to \nmty{U}_1$ is then uniquely determined since $\nmty{U}_1$ is terminal.
\end{example}

We now prove that the term model $(\mathbb{C}_{I}, p_I)$ satisfies the appropriate universal property.

\begin{lemma}[Context extension by a basic type is a product]
\label{lemExtensionByBasicTypeIsProduct}
Let $(\mathbb{C}, p)$ be a natural model. If $\Gamma \in \mathrm{ob}(\mathbb{C})$ and $A \in \nmty{U}(\diamond)$. The span
$$\Gamma \xleftarrow{\nmp{A}} \Gamma \cext A[\nmt{\Gamma}] \xrightarrow{\nmt{\Gamma} \cext A} \diamond \cext A$$
is a product diagram in $\mathbb{C}$.
\end{lemma}

\begin{proof}
Note that the following square is a canonical pullback square (\Cref{cnsCanonicalPullbacks}).

\begin{center}
\begin{tikzcd}[row sep=huge, column sep=huge]
\Gamma \cext A[\nmt{\Gamma}]
\arrow[r, "\nmt{\Gamma} \cext A"]
\arrow[d, "\nmp{A[\nmt{\Gamma}]}"']
&
\diamond \cext A
\arrow[d, "\nmp{A}"]
\\
\Gamma
\arrow[r, "\nmt{\Gamma}"']
&
\diamond
\end{tikzcd}
\end{center}
But $\diamond$ is terminal in $\mathbb{C}$, so this says precisely that $\Gamma \cext A[\nmt{\Gamma}]$ is a product of $\Gamma$ and $\diamond \cext A$, with the required projection morphisms.
\end{proof}

\begin{theorem}[Universal property of the model $(\mathbb{C}_I, p_I)$]
\label{thmFreeNaturalModelBasicTypesInitial}
Let $(\mathbb{C}, p)$ be a natural model and let $\{ O_i \mid i \in I \} \subseteq \nmty{U}(\diamond)$. There is a unique morphism of natural models $F : (\mathbb{C}_I, p_I) \to (\mathbb{C}, p)$ such that $F(i) = O_i$ for all $i \in I$. Hence $(\mathbb{C}_I, p_I)$ is initial in the category $\mathbf{NM}_{(\mathsf{ty}_i)_{i \in I}}$.
\end{theorem}

\begin{proof}
Define the underlying functor $F : \mathbb{C}_{I} \to \mathbb{C}$ on objects by
$$F(A, u) = O_{u(a_0)} \cext \dots \cext O_{u(a_{m-1})}$$

Let $f : (A, u) \to (B, v)$ in $\mathbb{C}_{I}$, so that $f$ is a function $B \to A$ over $I$. By Lemma \ref{lemExtensionByBasicTypeIsProduct} (and an easy induction), the objects $O_{u(a_0)} \cext \dots \cext O_{u(a_{m-1})}$ and $O_{v(b_0)} \cext \dots \cext O_{v(b_{n-1})}$ are the products of their respective component basic types. With this in mind, let
$$F(f) = \langle \pi_{f(b_0)}, \pi_{f(b_1)}, \dots, \pi_{f(b_{n-1})} \rangle : O_{u(a_0)} \cext \dots \cext O_{u(a_{m-1})} \to O_{v(b_0)} \cext \dots \cext O_{v(b_{n-1})}$$
where $\pi_{a_k} : O_{u(a_0)} \cext \dots \cext O_{u(a_{m-1})} \to O_{u(a_k)}$ is the product projection onto the $k^{\text{th}}$ component.

In order to see that $F(f)$ is well-defined, we need the codomain of $\pi_{f(b_{\ell})}$ to be $O_{v(b_{\ell})}$ for each $\ell < n$. To see this, note that for given $\ell < n$ we have $f(b_{\ell}) = a_k$ for some $k<m$, so that the codomain of $\pi_{f(b_{\ell})}$ is $O_{u(a_k)}$. Now $u(a_k) = u(f(b_{\ell})) = v(b_{\ell})$ since $f$ is a morphism over $I$, and hence $O_{u(a_k)} = O_{v(b_{\ell})}$, as required.

That the assignment $f \mapsto F(f)$ is functorial follows from elementary computations using product projections.

The natural transformation $\varphi : \nmty{U}_I \to F^*\nmty{U}$ is defined by
$$\varphi_{(A,u)} : I \to \nmty{U}(O_{u(a_0)} \cext \dots \cext O_{u(a_{m-1})}), \quad i \mapsto O_i$$
and the natural transformation $\nmmk{\varphi} : \nmtm{U}_I \to F^*\nmtm{U}$ is defined by
$$\nmmk{\varphi}_{(A,u)} : A \to \nmtm{U}(O_{u(a_0)} \cext \dots \cext O_{u(a_{m-1})}), \quad a_k \mapsto \nmq{O_{u(a_k)}}$$

Note that, in particular, we have $F(i) = \varphi_{(0,!_I)}( i ) = O_i$, as required.

To see that $\varphi$ and $\nmmk{\varphi}$ are natural, let $f : (A, u) \to (B, v)$ in $\mathbb{C}_I$, so that $f$ is a function $B \to A$ over $I$.
\begin{itemize}
\item The naturality square for $\varphi$ is as follows:
\begin{center}
\begin{tikzcd}[row sep=huge, column sep=huge]
I
\arrow[r, "\varphi_{(A,u)}"]
&
\nmty{U}(O_{u(a_0)} \cext \dots \cext O_{u(a_{m-1})})
\
\\
I
\arrow[r, "\varphi_{(B,v)}"']
\arrow[u, "\mathrm{id}_I"]
&
\nmty{U}(O_{v(b_0)} \cext \dots \cext B_{v(b_{n-1})})
\arrow[u, "{\nmty{U}(\langle \pi_{f(0)}, \dots, \pi_{f(n-1)} \rangle)}"']
\end{tikzcd}
\end{center}
Both composites $\varphi_{(A,u)} \circ \mathrm{id}_I$ and $\nmty{U}(\langle \pi_{f(0)}, \dots, \pi_{f(n-1)} \rangle) \circ \varphi_{(B,v)}$ send $i \in I$ to $O_i \in \nmty{U}(O_{u(a_0)} \cext \dots \cext O_{u(a_{m-1})})$, and so $\varphi$ is indeed natural.
\item The naturality square for $\nmmk{\varphi}$ is as follows:
\begin{center}
\begin{tikzcd}[row sep=huge, column sep=huge]
A
\arrow[r, "\nmmk{\varphi}_{(A,u)}"]
&
\nmtm{U}(O_{u(a_0)} \cext \dots \cext O_{u(a_{m-1})})
\\
B
\arrow[r, "\nmmk{\varphi}_{(B,v)}"']
\arrow[u, "f"]
&
\nmtm{U}(O_{v(b_0)} \cext \dots \cext O_{v(b_{n-1})})
\arrow[u, "{\nmtm{U}(\langle \pi_{f(0)}, \dots, \pi_{f(n-1)} \rangle}"']
\end{tikzcd}
\end{center}
To see that this commutes, let $b \in B$. Then $b=b_{\ell}$ for some $\ell<n$. Let $k<m$ be such that $f(b_{\ell}) = a_k$. It is then evident that both composites send $b$ to $\nmq{O_{u(a_k)}}$, so that $\nmmk{\varphi}$ is natural.
\end{itemize}

It is immediate from its definition that $F_0$ preserves chosen terminal objects.

To see that $F^*(p) \circ \nmmk{\varphi} = \varphi \circ p_I$, note that for all $(A, u) \in \mathbb{C}_I$ we have
\begin{align*}
F^*(p)_{(A,u)} \circ \nmmk{\varphi}_{(A,u)} ( a_k )
&= p_{F(A,u)} ( \nmq{O_{u(a_k)}} ) && \text{definition of $\nmmk{\varphi}$ and $F^*$} \\
&= O_{u(a_k)} && \text{definition of $p$ and $\nmq{O_{u(a_k)}}$} \\
&= \varphi_{(A,u)}( u(a_k) ) && \text{definition of $\varphi$} \\
&= \varphi_{(A,u)} \circ (p_I)_{(A,u)} ( a_k ) && \text{definition of $p_I$} \\
\end{align*}
To see that $F$ preserves the representability data, fix $(A,u) \in \mathbb{C}_I$ and $j \in I$. Then:
\begin{itemize}
\item $F((A,u) \cext j) = O_{u(a_0)} \cext \dots \cext O_{u(a_{m-1})} \cext O_j = F(A,u) \cext F(j)$;
\item $\nmp{F(j)} = \nmp{O_j} : O_{u(a_0)} \cext \dots \cext O_{u(a_{m-1})} \cext O_j \to O_{u(a_0)} \cext \dots \cext O_{u(a_{m-1})}$ is given by projection onto the first $m$ components, which is precisely $F(\nmp{j})$;
\item $F(\nmq{j}) = \nmmk{\varphi}_{(A,u) \cext j}(\star) = \nmq{O_j} = \nmq{F_{\mathrm{ty}}(j)}$.
\end{itemize}

Hence $F$ is a morphism of natural models.

It remains to prove that $F$ is unique. Suppose that $G = (G, \gamma, \nmmk{\gamma}) : (\mathbb{C}_I, p_I) \to (\mathbb{C}, p)$ is another morphism of natural models such that $G(i)=O_i$ for all $i \in I$. It follows by induction on $|A|$ that $G(A,u)=F(A,u)$ for all $(A,u) \in \mathbb{C}_I$. That $G=F$ then follows from the fact that $G$ preserves context extension and all morphisms in $\mathbb{C}_I$ are given by projections.
\end{proof}

\subsection*{A word of warning}

In the next four sections, we describe left adjoints to forgetful functors of the form $\mathbf{NM}_{\mathbb{T}} \to \mathbf{NM}$ when $\mathbb{T}$ is, respectively, the theory of a term of a basic type (\Cref{secExtTerm}), the theory of an (externally) indexed family of basic types (\Cref{secExtType}), the theory of a unit type (\Cref{secExtUnit}), and the theory of dependent sum types (\Cref{secExtSigma}). Before we do so, it is worth pointing out what we do \textit{not} do.

\begin{itemize}
\item We do not prove the initiality of the \textit{term model} of a dependent type theory $\mathbb{T}$ in the category $\mathbf{NM}_{\mathbb{T}}$, which is the natural model built out of the syntax of the theory $\mathbb{T}$. We could do so either by proving that a given term model is isomorphic to the constructed \textit{free} model, or by proving that the term model satisfies the universal property of the free model. For more on this, see the discussion around \Cref{schTermModel}.
\item We do not compose our free constructions. For instance, suppose we are given an arbitrary natural model $(\mathbb{C}, p)$. We could use the work in \Cref{secExtType} to freely adjoin a basic type, and the work in \Cref{secExtSigma} to freely adjoin dependent sum types to the resulting model; or we could first freely adjoin dependent sum types, and then freely adjoin a basic type. The two resulting natural models would, in general, not be isomorphic, since we have implicitly composed with the forgetful functor $\mathbf{NM}_{\mathsf{ty}} \to \mathbf{NM}$ in the first case, and with the forgetful functor $\mathbf{NM}_{\upSigma} \to \mathbf{NM}$ in the second case. In order to resolve this issue, we would need to describe the left adjoint to at least one (preferably both) of the forgetful functors
$$\mathbf{NM}_{\mathsf{ty},\upSigma} \to \mathbf{NM}_{\mathsf{ty}} \quad \text{and} \quad \mathbf{NM}_{\mathsf{ty},\upSigma} \to \mathbf{NM}_{\upSigma}$$
With such adjoint functors described, we would expect the free constructions to commute in the desired way---that is, the composites of free functors
$$\mathbf{NM} \to \mathbf{NM}_{\mathsf{ty}} \to \mathbf{NM}_{\mathsf{ty},\upSigma} \quad \text{and} \quad \mathbf{NM} \to \mathbf{NM}_{\upSigma} \to \mathbf{NM}_{\mathsf{ty},\upSigma}$$
will be naturally isomorphic.
\end{itemize}

We leave the task of overcoming these limitations to future work.

%% file: thesis/ch4-semantics/ext-term.tex
\section{Extending a natural model by a term of a basic type}
\label{secExtTerm}

In a given dependent type theory $\mathbb{T}$, the contexts $\Gamma = x_1 : A_1, x_2: A_2, \dots, x_n:A_n$ satisfy the property that each $A_i$ is a type in context $x_1 : A_1, \dots, x_{i-1} : A_{i-1}$. In particular, $A_1$ is a basic type. If we introduce a new term $o : O$ of a basic type $O$, the contexts $x_1 : A_1, \dots, x_n : A_n$ in the resulting type theory $\mathbb{T}'$ satisfy the property that each $A_i$ may additionally depend on the new term $o:O$. Thus the contexts of $\mathbb{T}'$ are precisely those obtained by taking a context $x:O, x_1 : A_1, \dots, x_n : A_n$ in $\mathbb{T}$ and then (semantically) substituting $o : O$ for $x : O$.

Intuitively speaking, then, given a model $\mathfrak{M}$ of type theory $\mathbb{T}$ such that $\mathfrak{M}$ satisfies the judgement $\vdash o : O$, we can interpret a context of $\mathbb{T}'$ by first interpreting the corresponding context $x : O, x_1 : A_1, \dots, x_n : A_n$ of $\mathbb{T}$ in $\mathfrak{M}$, and then substituting $o$ for $x$ in $\mathfrak{M}$.

Transferring this intuition to a natural model $(\mathbb{C}, p)$, when we freely adjoin a term of a basic type $O \in \nmty{U}(\diamond)$, the contexts of the new natural model `should' look like $O \cext \Gamma$, where $\Gamma$ is an old context. Unfortunately it is not always possible to make sense of the expression $O \cext \Gamma$, unless $\Gamma$ is itself of the form $\diamond \cext A_1 \cext A_2 \cext \dots \cext A_n$. In order to overcome this obstacle, we instead take our new contexts to be those of the form $\Gamma \cext O \cext A_1 \cext \dots \cext A_n$. A technicality we must take care of is that the inclusion from $(\mathbb{C}, p)$ to the new natural model, which sends a context $\Gamma$ to the weakened context $\Gamma \cext O$, must preserve context extension---in order to do this, we only \textit{formally} extend the contexts, and we reduce the contexts to a normal form by pushing the variable $x : O$ as far to the right as possible, so that the formally extended contexts $\Gamma \cext O \cext A$ and $(\Gamma \cext A) \cext O$ become identified whenever $A \in \nmty{U}(\Gamma)$.

\begin{definition}[Swap isomorphisms]
\label{defSwapIsomorphism}
Let $(\mathbb{C}, p)$ be a natural model. For each $\Gamma \in \mathrm{ob}(\mathbb{C})$ and $A,O \in \nmty{U}(\Gamma)$, the \textbf{swap isomorphism}
$$\nmsw_{O,A} : \Gamma \cext O[\nmt{\Gamma}] \cext A[\nmp{O}] \xrightarrow{\cong} \Gamma \cext A \cext O[\nmt{\Gamma \cext A}]$$
is the isomorphism given by the respective canonical pullback squares for $\nmp{A}$ and $\nmp{O}$; we will write $\Gamma \cext O \cext A$ and $\Gamma \cext A \cext O$ to simplify notation.
\begin{diagram}
\Gamma \cext O \cext A
\arrow[ddr, bend right=15, "\nmp{O} \cext A"']
\arrow[drr, bend left=15, "\nmp{A}"]
\arrow[dr, dashed, "\nmsw_{O,A}" description]
&[-20pt]&
\\[-20pt]
&
\Gamma \cext A \cext O
\arrow[r, "\nmp{A} \cext O" description]
\arrow[d, "\nmp{O}"']
&
\Gamma \cext O
\arrow[d, "\nmp{O}"]
\\
&
\Gamma \cext A
\arrow[r, "\nmp{A}"']
&
\Gamma
\end{diagram}
Furthermore, given a list $(A_1, \dots, A_n)$ with $A_{i+1} \in \nmty{U}(\Gamma \cext O \cext A \cext A_1 \cext \dots \cext A_i)$ for each $i<n$, we obtain isomorphisms
$$\nmsw : \Gamma \cext O \cext A \cext A_1 \cext A_2 \cext \dots \cext A_n \xrightarrow{\cong} \Gamma \cext A \cext O \cext A_1 \cext A_2 \cext \dots \cext A_n$$
where again we have suppressed the substitutions. We will also refer to these as swap isomorphisms.
\end{definition}

Note that swap isomorphisms cohere with representability data since they are the isomorphisms induced from the universal property of canonical pullback squares.

\begin{construction}[Category of contexts extended by a term]
\label{cnsCategoryOfContextsExtendedByTerm}
Let $(\mathbb{C}, p)$ be a natural model and let $O \in \nmty{U}(\diamond)$. The \textbf{category of contexts extended by a term $x$ of type $O$} is the category $\mathbb{C}\adjt{x:O}$ defined as follows.
\begin{itemize}
\item The objects of $\mathbb{C}\adjt{x:O}$ are $(n+1)$-tuples $(\Gamma; A_1, \dots, A_n)$, where $n \ge 0$, $\Gamma \in \mathrm{ob}(\mathbb{C})$ and $A_i \in \nmty{U}(\Gamma \cext O[\nmt{\Gamma}] \cext A_1 \cext \dots \cext A_i)$ for each $i < n$, and we additionally identify the list $(\Gamma \cext A; A_1, \dots, A_n)$ with the list $(\Gamma; A[\nmp{O}], A_1, \dots, A_n)$, where we have suppressed the swap isomorphisms.

The object $(\Gamma; \vec A)$ of $\mathbb{C}\adjt{x:O}$ represents the result of extending a context $\Gamma$ first by a variable $x:O$ and then by variables of types which may depend on $x$. Each object $\vec \Gamma$ of $\mathbb{C}\adjt{x:O}$ has a unique expression as an $(n+1)$-tuple $(\Gamma; A_1, \dots, A_n)$ with $n$ minimal, which we call the \textbf{normal form} of $\vec \Gamma$, representing the result of pushing the new variable $x:O$ as far to the right as possible by swap isomorphisms.

Unless otherwise specified, all expressions of objects of $\mathbb{C} \adjt{x:O}$ as lists $(\Gamma; A_1, \dots, A_n)$ will be assumed to be in normal form. Note that if $(\Gamma; A_1, \dots, A_n)$ is in normal form and $n>0$, then $(\Gamma; A_1, \dots, A_n, B)$ is in normal form.

\item A morphism $\sigma : (\Delta; B_1, \dots, B_m) \to (\Gamma; A_1, \dots, A_n)$ in $\mathbb{C}\adjt{x:O}$ is a morphism
$$\sigma : \Delta \cext O[\nmt{\Delta}] \cext B_1 \cext \dots \cext B_m \to \Gamma \cext O[\nmt{\Gamma}] \cext A_1 \cext \dots \cext A_n$$
commuting with the canonical morphisms to $\diamond \cext O$, as indicated in the following diagram.
\begin{diagram}
\Delta \cext O[\nmt{\Delta}] \cext B_1 \cext \dots \cext B_m
\arrow[rr, "\sigma"]
\arrow[d, "\nmp{\vec B}"']
&&
\Gamma \cext O[\nmt{\Gamma}] \cext A_1 \cext \dots \cext A_n
\arrow[d, "\nmp{\vec A}"]
\\
\Delta \cext O[\nmt{\Delta}]
\arrow[dr, "\nmt{\Delta} \cext O"']
\arrow[dd, midgray, "\nmp{O[\nmt{\Delta}]}"']
\pullbackc{dddr, midgray}{0.05}
&&
\Gamma \cext O[\nmt{\Gamma}]
\arrow[dl, "\nmt{\Gamma} \cext O"]
\arrow[dd, midgray, "\nmp{O[\nmt{\Gamma}]}"]
\pullbackc[-90]{dddl, midgray}{0.05}
\\[-10pt]
&
\diamond \cext O
\arrow[dd, midgray, "\nmp{O}" description]
&
\\[-30pt]
{\color{midgray} \Delta}
\arrow[dr, midgray, , "\nmt{\Delta}"']
&&
{\color{midgray} \Gamma}
\arrow[dl, midgray, "\nmt{\Gamma}"]
\\[-10pt]
&
{\color{midgray} \diamond}
&
\end{diagram}
\end{itemize}
There is a full and faithful functor $E : \mathbb{C}\adjt{x:O} \to \mathbb{C} \slice{\diamond \cext O}$ defined on objects by letting $E(\Gamma; \vec A)$ be the composite
$$\Gamma \cext O[\nmt{\Gamma}] \cext \vec A \xrightarrow{\nmp{\vec A}} \Gamma \cext O[\nmt{\Gamma}] \xrightarrow{\nmt{\Gamma} \cext O} \diamond \cext O$$
in $\mathbb{C} \slice{\diamond \cext O}$, and on morphisms by $E(\sigma)=\sigma$, so that $\mathbb{C}\adjt{x:O}$ is equivalent to a full subcategory of $\mathbb{C} \slice{\diamond \cext O}$.
\end{construction}

\begin{verification}
That $\mathbb{C}\adjt{x:O}$ is a category follows immediately from the fact that composition and identity are inherited from $\mathbb{C} \slice{\diamond \cext O}$. The functor $E$ evidently respects domains and codomains, and is full and faithful since it acts as the identity on morphisms. Furthermore, every full and faithful functor corestricts to an equivalence between its domain and its image.
\end{verification}

When clear from context, we will abuse notation by writing `$O$' to refer simultaneously to the element $O \in \nmty{U}(\diamond)$, the object $\diamond \cext O \in \mathbb{C}$ and the elements $O[\nmt{\Gamma}] \in \nmty{U}(\Gamma)$ for all $\Gamma \in \mathrm{ob}(\mathbb{C})$.

\begin{numbered}
We will denote by $\mathbb{C}(O)$ the full subcategory of $\mathbb{C} \slice{O}$ which is the image of $\mathbb{C} \adjt{x:O}$ under $E$. Thus the objects of $\mathbb{C}(O)$ are morphisms of $\mathbb{C}$ of the form $(\nmt{\Gamma} \cext O) \circ \nmp{\vec A}$ for $(\Gamma; \vec A) \in \mathrm{ob}(\mathbb{C} \adjt{x:O})$. Since $E$ is full and faithful, the corestriction $E : \mathbb{C}\adjt{x:O} \to \mathbb{C}(O)$ is an equivalence of categories. The `product with $O$' functor $\upDelta_O : \mathbb{C} \to \mathbb{C} \slice{O}$ given by $\Gamma \mapsto (\Gamma \cext O, \nmt{\Gamma} \cext O : \Gamma \cext O \to O)$ (see \Cref{lemExtensionByBasicTypeIsProduct}) factors through the inclusion $\mathbb{C}(O) \hookrightarrow \mathbb{C} \slice{O}$, and so we obtain an adjunction $\upSigma_O \dashv \upDelta_O : \mathbb{C} \leftrightarrows \mathbb{C}(O)$.
\end{numbered}

\begin{lemma}
\label{lemFactorisationOfProductFunctor}
Let $(\mathbb{C}, p)$ be a natural model and let $O \in \nmty{U}(\diamond)$. The functor $\upDelta_O : \mathbb{C} \to \mathbb{C}(O)$ factors through the functor $E : \mathbb{C}\adjt{x:O} \to \mathbb{C}(O)$.
\begin{diagram}
\mathbb{C}
\arrow[d, "I"']
\arrow[dr, "\upDelta_O"]
&
\\
\mathbb{C} \adjt{x:O}
\arrow[r, "E"']
&
\mathbb{C}(O)
\end{diagram}
Moreover, the functor $I : \mathbb{C} \to \mathbb{C}\adjt{x:O}$ is right adjoint to the composite $\mathbb{C}\adjt{x:O} \xrightarrow{E} \mathbb{C}(O) \xrightarrow{\upSigma_O} \mathbb{C}$.
\end{lemma}

\begin{proof}
Since $\Gamma \xleftarrow{\nmp{O}} \Gamma \cext O \xrightarrow{\nmt{\Gamma} \cext O} O$ is a product diagram in $\mathbb{C}$, we can take $\upDelta_O\Gamma = (\Gamma \cext O, \nmt{\Gamma} \cext O : \Gamma \cext O \to O)$ for each $\Gamma \in \mathbb{C}$. But then $\upDelta_O\Gamma = E(\Gamma)$; so define $I\Gamma = (\Gamma)$ and $I(\sigma : \Delta \to \Gamma) = \sigma \cext O : (\Delta) \to (\Gamma)$, and observe that this defines a functor $\mathbb{C} \to \mathbb{C}\adjt{x:O}$ with $\upDelta_O = E \circ I$, which is well-defined since $(\Delta)$ and $(\Gamma)$ are in normal form.

To see that $\upSigma_O \circ E \dashv I$, observe that we have the following chain of equalities and natural isomorphisms.
\begin{align*}
& \mathbb{C}(\upSigma_O(E(\Delta; \vec B), \Gamma) && \\
&\cong \mathbb{C}(O)(E(\Delta; \vec B), \upDelta\Gamma) && \text{since $\upSigma_O \dashv \upDelta_O$} \\
&= \mathbb{C}(O)(\Delta \cext O \cext B_1 \cext \dots \cext B_m, \nmt{\Delta \cext O} \circ \nmp{\vec B}), (\Gamma \cext O, \nmt{\Gamma} \cext O) && \text{definitions of $\upDelta_O$ and of $E$} \\
&= \mathbb{C}\adjt{x:O}((\Delta; \vec B), (\Gamma)) && \text{definition of morphisms in $\mathbb{C}\adjt{x:O}$} \\
&= \mathbb{C}\adjt{x:O}((\Delta; \vec B), I\Gamma) && \text{definition of $I$}
\end{align*}
\end{proof}

\begin{construction}[Free natural model extended by a term]
\label{cnsFreeNaturalModelExtTerm}
Let $(\mathbb{C}, p)$ be a natural model and let $O \in \nmty{U}(\diamond)$. The \textbf{free natural model extended by a term $x$ of type $O$} is the natural model $(\mathbb{C}\adjt{x:O}, p\adjt{x:O} : \nmtmadjt{x:O}{U} \to \nmtyadjt{x:O}{U})$ defined by the following data. The underlying category is $\mathbb{C}\adjt{x:O}$ (\Cref{cnsCategoryOfContextsExtendedByTerm}) with distinguished terminal object $(\diamond)$. The presheaves $\nmtyadjt{x:O}{U}, \nmtmadjt{x:O}{U} : \mathbb{C}\adjt{x:O}\op \to \mathbf{Set}$ and the natural transformation $p\adjt{x:O} : \nmtmadjt{x:O}{U} \to \nmtyadjt{x:O}{U}$ are obtained from $p : \nmtm{U} \to \nmty{U}$ by precomposing with the composite $\mathbb{C}\adjt{x:O} \xrightarrow{E} \mathbb{C}(O) \xrightarrow{\upSigma_0} \mathbb{C}$. Explicitly, we have
$$\nmtyadjt{x:O}{U}(\Gamma; A_1, \dots, A_n) = \nmty{U}(\Gamma \cext O \cext A_1 \cext \dots \cext A_n)$$
and likewise for $\nmtmadjt{x:O}{U}$, and then
$$(p\adjt{x:O})_{(\Gamma; A_1, \dots, A_n)}(a) = p_{\Gamma \cext O \cext A_1 \cext \dots \cext A_n}(a)$$
for all $(\Gamma; A_1, \dots, A_n) \in \mathrm{ob}(\mathbb{C}\adjt{x:O})$ and all $a \in \nmtmadjt{x:O}{U}(\Gamma; A_1, \dots, A_n)$.

The representability data is defined for $(\Gamma; \vec A) = (\Gamma; A_1, \dots, A_n)$ as follows.
\begin{itemize}
\item Let $(\Gamma; A_1, \dots, A_n) \cext A = (\Gamma; A_1, \dots, A_n, A)$---note that if $n \ge 1$ then this is automatically in normal form, and if $n=0$ and $A = A'[\nmp{O}]$ for some $A' \in \nmty{U}(\Gamma)$, then the normal form is given by $(\Gamma) \cext A = (\Gamma \cext A')$;
\item Let $\nmp{A} : (\Gamma; \vec A, A) \to (\Gamma; \vec A)$ be the usual morphism $\nmp{A} : \Gamma \cext O \cext A_1 \cext \dots \cext A_n \cext A \to \Gamma \cext O \cext A_1 \cext \dots \cext A_n$ in $\mathbb{C}$ (or $\nmp{A} \circ \nmsw^{-1} : \Gamma \cext A' \cext O \to \Gamma \cext O$ in the case discussed above); and
\item Let $\nmq{A} \in \nmtyadjt{x:O}{U}(\Gamma; A_1, \dots, A_n, A) = \nmty{U}(\Gamma \cext O \cext A_1 \cext \dots \cext A_n \cext A)$ be the usual element $\nmq{A}$ (or the element $\nmq{A}[\nmsw^{-1}] \in \nmty{U}(\Gamma \cext A'\cext O)$ in the case discussed above).
\end{itemize}
The distinguished term $x \in \nmtmadjt{x:O}{U}(\diamond\adjt{x:O}; O)$ is given by the element $\nmq{O} \in \nmtm{U}(\diamond \cext O)$.
\end{construction}

\begin{verification}
That $\nmtyadjt{x:O}{U}$ and $\nmtmadjt{x:O}{U}$ are presheaves and that $p\adjt{x:O}$ is a natural transformation are immediate from the fact that they are obtained from $p$ by applying the functor $(\Sigma_O \circ E)^* : \widehat{\mathbb{C}} \to \widehat{\mathbb{C} \adjt{x:O}}$.

To see that $p\adjt{x:O}$ is representable, let $\vec \Gamma = (\Gamma; A_1, \dots, A_n) \in \mathrm{ob}(\mathbb{C} \adjt{x:O})$ and $A \in \nmtyadjt{x:O}{U}(\vec \Gamma)$ and consider the following square in $\widehat{\mathbb{C}\adjt{x:O}}$.
\begin{diagram}
\Yon(\vec \Gamma \cext A)
\arrow[d, "\nmp{A}"']
\arrow[r, "\nmq{A}"]
&
\nmtmadjt{x:O}{U}
\arrow[d, "p\adjt{x:O}"]
\\
\Yon(\vec \Gamma)
\arrow[r, "A"']
&
\nmtyadjt{x:O}{U}
\end{diagram}
Composing with swap isomorphisms if necessary, we can take $\nmp{A}$ and $\nmq{A}$ to be the respective morphism and element of $\mathbb{C}$. To see that the square is a pullback, let $\vec \Delta = (\Delta; B_1, \dots, B_m) \in \mathrm{ob}(\mathbb{C} \adjt{x:O})$ and let $\sigma : \vec \Delta \to \vec \Gamma$ and $a \in \nmtmadjt{x:O}{U}(\vec \Delta; A[\sigma])$.
\begin{diagram}
\Yon(\vec \Delta)
\arrow[ddr, bend right=15, "\sigma"']
\arrow[drr, bend left=15, "a"]
\arrow[dr, dashed]
&[-20pt]
&
\\[-20pt]
&
\Yon(\vec \Gamma \cext A)
\arrow[d, "\nmp{A}"']
\arrow[r, "\nmq{A}"]
&
\nmtmadjt{x:O}{U}
\arrow[d, "p\adjt{x:O}"]
\\
&
\Yon(\vec \Gamma)
\arrow[r, "A"']
&
\nmtyadjt{x:O}{U}
\end{diagram}
Again composing with swap isomorphisms if necessary, we can take $\sigma$ to be a morphism from $\Delta \cext O \cext B_1 \cext \dots \cext B_m$ to $\Gamma \cext O \cext A_1 \cext \dots \cext A_n$ in $\mathbb{C}(O)$ and $a \in \nmtm{U}(\Delta \cext O \cext B_1 \cext \dots \cext B_m; A[\sigma])$. But then by representability of $p$ there is a unique morphism
$$\langle \sigma, a \rangle_A : \Delta \cext O \cext B_1 \cext \dots \cext B_m \to \Gamma \cext O \cext A_1 \cext \dots \cext A_n \cext A$$
in $\mathbb{C}$ such that $\nmp{A} \circ \langle \sigma, a \rangle_A = \sigma$ and $\nmq{A}[\langle \sigma, a \rangle_A] = a$. Moreover, this is a morphism in $\mathbb{C}(O)$ since
\begin{align*}
(\nmt{\Gamma} \cext O) \circ \nmp{\vec A, A} \circ \langle \sigma, a \rangle_A
&= (\nmt{\Gamma} \cext O) \circ \nmp{\vec A} \circ \nmp{A} \circ \langle \sigma, a \rangle_A
&& \text{by definition of $\nmp{\vec A, A}$} \\
&= (\nmt{\Gamma} \cext O) \circ \nmp{\vec A} \circ \sigma && \text{by the universal property of pullbacks} \\
&= (\nmt{\Delta} \cext O) \circ \nmp{\vec B} && \text{since $\sigma$ is a morphism in $\mathbb{C}(O)$}
\end{align*}
So we see that $\langle \sigma, a \rangle_A$, perhaps composed with the relevant swap isomorphisms, is the desired morphism of $\mathbb{C} \adjt{x:O}$. So $p\adjt{x:O}$ is representable.
\end{verification}

\begin{numbered}
Under the equivalence $E : \mathbb{C}\adjt{x:O} \simeq \mathbb{C}(O) \subseteq \mathbb{C}\slice{O}$, the new terminal object $\diamond\adjt{x:O}$ corresponds with the identity morphism $\mathrm{id}_O : O \to O$. The canonical section $\nms{x} : (\diamond) \to (\diamond; O)$ of the new term $x \in \nmtyadjt{x:O}{U}(\diamond\adjt{x:O}; O)$ is then given by the diagonal morphism $\delta_O : (O, \mathrm{id}_O) \to (O \cdot O, \nmp{O})$.
\end{numbered}

\begin{lemma}[Inclusion morphism]
Let $(\mathbb{C}, p)$ be a natural model. The functor $I : \mathbb{C} \to \mathbb{C}\adjt{x:O}$ of \Cref{lemFactorisationOfProductFunctor} extends to a morphism of natural models $(I, \iota, \nmmk{\iota}) : (\mathbb{C}, p) \to (\mathbb{C}\adjt{x:O}, p\adjt{x:O})$.
\end{lemma}

\begin{proof}
In \Cref{cnsFreeNaturalModelExtTerm} we have $p \adjt{x:O} = (\upSigma_O \circ E)^*(p)$. Since $\upSigma_O \circ E \dashv I$, it follows from \Cref{lemAdjointFunctorKanExtension} that $(\upSigma_O \circ E)^* \cong I_!$, so that we can take $I_!(p) = p \adjt{x:O}$. But then we can take
$$\iota = \mathrm{id}_{I_!\nmty{U}} : I_! \nmty{U} \to I_! \nmty{U} = \nmtyadjt{x:O}{U} \quad \text{and} \quad \nmmk{\iota} = \mathrm{id}_{I_! \nmtm{U}} : I_! \nmtm{U} \to I_! \nmtm{U} = \nmtmadjt{x:O}{U}$$
Now note that $(I, \iota, \nmmk{\iota})$ preserves context extension, since for all $\Gamma \in \mathrm{ob}(\mathbb{C})$ and $A \in \nmty{U}(\Gamma)$ we have
$$I\Gamma \cext IA = (\Gamma; A[\nmp{O}]) = (\Gamma \cext A) = I(\Gamma \cext A)$$
by the identification of lists described in \Cref{cnsCategoryOfContextsExtendedByTerm}. The fact that $(I, \iota, \nmmk{\iota})$ is a morphism of natural models now follows trivially from the fact that $\iota$ and $\nmmk{\iota}$ are identity morphisms.
\end{proof}

\begin{theorem}[Extension of a morphism of natural models]
\label{thmFreelyAdjoiningTermFunctorial}
Let $(\mathbb{C}, p)$ be a natural model and let $O \in \nmty{U}(\diamond)$. For each morphism of natural models $F : (\mathbb{C}, p) \to (\mathbb{D}, q)$, there is a morphism of natural models $F_{\adjt{x:O}} : (\mathbb{C} \adjt{x:O}, p \adjt{x:O}) \to (\mathbb{D} \adjt{y : FO}, q \adjt{y : FO})$ such that $F \adjt{x:O} \circ I = I \circ F$ and $F(x)=y \in \nmtmalt{V}(\star; FO)$.
\begin{diagram}
(\mathbb{C}, p)
\arrow[r, "F"]
\arrow[d, "I"']
&
(\mathbb{D}, q)
\arrow[d, "I"]
\\
(\mathbb{C} \adjt{x:O}, p \adjt{x:O})
\arrow[r, dashed, "F \adjt{x:O}"']
&
(\mathbb{D}\adjt{y:FO}, q \adjt{y:FO})
\\[-45pt]
\overset{\text{\rotatebox{90}{$\in$}}}{\underset{x:O}{~}}
\arrow[r, mapsto, shift right=2.83]
&
\overset{\text{\rotatebox{90}{$\in$}}}{\underset{y:FO}{~}}
\end{diagram}
\end{theorem}

We will see in \Cref{corFreelyAdjoiningTermFunctorial} that the assignment $F \mapsto F_{\mathsf{tm}}$ in fact extends to a functor.

\begin{proof}
Let $F = (F, \varphi, \nmmk{\varphi}) : (\mathbb{C}, p) \to (\mathbb{D}, q)$ be a morphism of natural models.

Define the functor $F_{\mathsf{tm}} : \mathbb{C} \adjt{x:O} \to \mathbb{D} \adjt{y:FO}$ on objects by letting
$$F_{\mathsf{tm}}(\Gamma; A_1, \dots, A_n) = (F\Gamma; FA_1, \dots, FA_n)$$
and given a morphism $\sigma : (\Delta; B_1, \dots, B_m) \to (\Gamma; A_1, \dots, A_n)$ in $\mathbb{C} \adjt{x:O}$, let $F_{\mathsf{tm}}(\sigma)$ be the same morphism from $F\Delta \cextalt FO \cextalt FB_1 \cextalt \dots \cextalt FB_m$ to $F\Gamma \cextalt FO \cextalt FA_1 \cextalt \dots \cextalt FA_n$ in $\mathbb{D}$ as is given by $F(\sigma)$ (with $\sigma$ considered as a morphism $\Delta \cext O \cext B_1 \cext \dots \cext B_m \to \Gamma \cext O \cext A_1 \cext \dots \cext A_n$ in $\mathbb{C}$). Then for each $\Gamma \in \mathrm{ob}(\mathbb{C})$ we have
$$F \adjt{x:O} I\Gamma = F \adjt{x:O}(\Gamma) = (F\Gamma) = I(F\Gamma)$$
so that $F \adjt{x:O} \circ I = I \circ F$.

Define $\varphi \adjt{x:O} = I_!(\varphi) : I_!\nmty{U} \Rightarrow I_!\nmty{V}$ and $\nmmk{\varphi} \adjt{x:O} = I_!(\nmmk{\varphi}) : I_!F_! \nmtm{U} \Rightarrow I_!\nmtmalt{V}$. Note that we have
$$I_!F_!p = (I \circ F)_!p = (F \adjt{x:O} \circ I)_!p = (F \adjt{x:O})_! I_!p = F_{\mathsf{tm}} p \adjt{x:O}$$
and $I_!q = q \adjt{y:FO}$, so that $\varphi \adjt{x:O}$ and $\nmmk{\varphi} \adjt{x:O}$ have the correct type. To see that $F_{\mathsf{tm}}$ preserves context extension, note that
\begin{align*}
F_{\mathsf{tm}}\Gamma \cext F_{\mathsf{tm}} A &= (F\Gamma) \cext FA[\nmu{FO}] && \text{definition of $F_{\mathsf{tm}}$ and $\varphi\adjt{x:O}$} \\
&= (F\Gamma; FA[\nmu{FO}]) && \text{context extension in $(\mathbb{C}\adjt{x:O}, p\adjt{x:O})$} \\
&= (F\Gamma; FA[F\nmp{O}]) && \text{$F$ is a morphism of natural models} && \\
&= (F\Gamma \cextalt FA) && \text{normal form} && \\
&= (F(\Gamma \cextalt A)) && \text{$F$ is a morphism of natural models} && \\
&= F_{\mathsf{tm}}(\Gamma \cext A) && \text{definition of $F_{\mathsf{tm}}$}
\end{align*}
Finally note that $F_{\mathsf{tm}}$ preserves the remaining representability data, so that we have a morphism of natural models as required.
\end{proof}

\begin{numbered}
\label{parCanonicalSections}
We are nearly ready to prove the universal property of $(\mathbb{C} \adjt{x:O}, p \adjt{x:O})$. First we must do some acrobatics involving terms of basic types. Given a natural model $(\mathbb{C}, p)$ with a basic type $O \in \nmty{U}(\diamond)$ and a term $o \in \nmtm{U}(\diamond; O)$, we obtain a section $\nms{o} : \diamond \to O$ of the projection $\nmp{O} : O \to \diamond$ from representability of $p$---specifically, we have $\nms{o} = \langle \mathrm{id}_{\diamond}, o \rangle_O$.

Given any object $\Gamma$ of $\mathbb{C}$, this gives rise to a section $\nms{o[\nmt{\Gamma}]} : \Gamma \to \Gamma \cext O$ of $\nmp{O} \cext O : \Gamma \cext O \to \Gamma$; we will just write $\nms{o}$ for $\nms{o}[\nmt{\Gamma}]$. Hence for any object $(\Gamma; \vec A)$ of $\mathbb{C} \adjt{x:O}$, we obtain a section $\nms{o} \cext \vec A : \Gamma \cext \vec A[\nms{o}] \to \Gamma \cext O \cext \vec A$ of $\nmp{O} \cext \vec A : \Gamma \cext O \cext \vec A \to \Gamma \cext A[\nms{O}]$. This is illustrated in the following diagram, in which all four squares are canonical pullbacks and all horizontal composites are identity morphisms.
\begin{diagram}
\Gamma \cext \vec A[\nms{o}]
\pullback
\arrow[r, "\nms{o} \cext \vec A"]
\arrow[d, "\nmp{\vec A[\nms{o}]}"']
&
\Gamma \cext O \cext \vec A
\pullback
\arrow[r, "\nmp{O} \cext \vec A"]
\arrow[d, "\nmp{\vec A}" description]
&
\Gamma \cext \vec A[\nms{o}]
\arrow[d, "\nmp{\vec A[\nms{o}]}"]
\\
\Gamma
\pullback
\arrow[r, "\nms{o}" description]
\arrow[d, "\nmt{\Gamma}"']
&
\Gamma \cext O
\pullback
\arrow[r, "\nmp{O}" description]
\arrow[d, "\nmt{\Gamma} \cext O" description]
&
\Gamma
\arrow[d, "\nmt{\Gamma}"]
\\
\diamond
\arrow[r, "\nms{o}"']
&
O
\arrow[r, "\nmp{O}"']
&
\diamond
\end{diagram}
In particular, the object of $\mathbb{C}$ obtained by pulling back the morphism $E(\Gamma; \vec A) = (\nmt{\Gamma} \cext O) \circ \nmp{\vec A}$ along $\nms{o}$ exists and can be taken to be equal to $\Gamma \cext \vec A[\nms{o}]$. This yields a functor $\upDelta_{\nms{o}} : \mathbb{C}(O) \to \mathbb{C}$.
\end{numbered}

\begin{construction}[Term substitution morphism]
\label{cnsTermSubstitutionMorphism}
Let $(\mathbb{C}, p)$ be a natural model, let $O \in \nmty{U}(\diamond)$ and let $o \in \nmtm{U}(\diamond; O)$. The \textbf{substitution morphism} of $o$ for $x$ is the morphism of natural models $S_o = (S_o, \sigma_o, \nmmk{\sigma}_o) : (\mathbb{C} \adjt{x:O}, p \adjt{x:O}) \to (\mathbb{C}, p)$ satisfying $S_o(x)=o \in \nmtm{U}(\diamond; O)$ and $S_o \circ I = \mathrm{id}_{(\mathbb{C}, p)}$; it is defined as follows.
\begin{itemize}
\item The functor $S_o : \mathbb{C}\adjt{x:O} \to \mathbb{C}$ is the composite
$$\mathbb{C} \adjt{x:O} \xrightarrow{E} \mathbb{C}(O) \xrightarrow{\upDelta_{\nms{o}}} \mathbb{C} \slice{\diamond} \cong \mathbb{C}$$
where $\nms{o} : \diamond \to O$ is as discussed in \Cref{parCanonicalSections}.
\item The natural transformation $\sigma_o : \nmtyadjt{x:O}{U} \to \nmty{U}$ is given by letting $(\sigma_o)_{(\Gamma; \vec A)} : \nmty{U}(\Gamma \cext O \cext \vec A) \to \nmty{U}(\Gamma \cext \vec A[\nms{o}])$ be the function $\nmty{U}(\nms{o} \cext \vec A)$.
\item Likewise, $\nmmk{\sigma}_o : \nmtmadjt{x:O}{U} \to \nmtm{U}$ is defined by $(\nmmk{\sigma}_o)_{(\Gamma; \vec A)} = \nmtm{U}(\nms{o} \cext \vec A)$.
\end{itemize}
\end{construction}

\begin{verification}
Note first that
$$S_o(\diamond) = \upDelta_{\nms{o}}(\diamond \cext O \xrightarrow{\nmp{O}} \diamond) = \diamond$$
so that $S_o$ preserves distinguished terminal objects. Given a morphism $\tau : (\Delta; \vec B) \to (\Gamma; \vec A)$ in $\mathbb{C} \adjt{x:O}$, the corresponding naturality squares for $\sigma_o$ and $\nmmk{\sigma}_o$ are obtained by applying $\nmty{U}$ and $\nmtm{U}$, respectively, to the following diagram $\mathbb{C}$.
\begin{diagram}
\Delta \cext O \cext \vec B
\arrow[r, "\nms{o} \cext \vec B"]
\arrow[d, "\tau"']
&
\Delta \cext \vec B[\nms{o}]
\arrow[d, "\upDelta_{\nms{o}}(\tau)"]
\\
\Gamma \cext O \cext \vec A
\arrow[r, "\nms{o} \cext \vec A"']
&
\Delta \cext \vec A[\nms{o}]
\end{diagram}
These diagrams commute in $\mathbb{C}$, and so the naturality squares commute too.

To see that $S_o$ preserves context extension, let $(\Gamma; \vec A) \in \mathrm{ob}(\mathbb{C} \adjt{x:O})$ and let $A \in \nmtyadjt{x:O}{U}(\Gamma; \vec A) = \nmty{U}(\Gamma \cext \vec A \cext A)$, and note that
\begin{align*}
S_o(\Gamma; \vec A) \cext S_o(A)
&= \Gamma \cext \vec A[\nms{o}] \cext A[\nms{o} \cext \vec A]
&& \text{by definition of $S_o$} \\
&= \Gamma \cext (\vec A \cext A)[\nms{o}] && \text{by our notation convention} \\
&= S_o(\Gamma \cext \vec A \cext A) && \text{by definition of $S_o$}
\end{align*}
And note that we have
$$S_o(\nmp{A}) = \upDelta_{\nms{o}}(E(\nmp{A})) = \nmp{A[\nms{o} \cext \vec A]} = \nmp{S_o(A)}$$
and
$$S_o(\nmq{A}) = \nmq{A}[\nms{o} \cext \vec A] = \nmq{A}[\langle \mathrm{id}_{\Gamma \cext \vec A}, \nmq{A}[\nms{o} \cext \vec A] \rangle_A] = \nmq{A[\nms{o} \cext \vec A]} = \nmq{S_o(A)}$$
so $S_o$ is a morphism of natural models.

To see that $S_o(x) = o$, note that $x = \nmq{O} \in \nmtm{U}(\diamond \cext O)$, so we have
$$S_o(x) = \nmq{O}[\nms{o}] = \nmq{O}[\langle \mathrm{id}_{\diamond}, o \rangle_O] = o$$
as required.

Finally note that for $\Gamma \in \mathrm{ob}(\mathbb{C})$ we have
$$S_o(I(\Gamma)) = S_o(\Gamma) = \upDelta_{\nms{o}}(\Gamma \cext O \xrightarrow{\nmt{\Gamma} \cext O} O) = \Gamma$$
Likewise it is easy to see that $S_o \circ I$ acts as the identity on substitutions, types and terms. So $S_o \circ I = \mathrm{id}_{(\mathbb{C}, p)}$, as required.
\end{verification}

We now have the components needed to prove the universal property of the natural model $(\mathbb{C}\adjt{x:O}, p\adjt{x:O})$.

\begin{theorem}[Universal property of freely extending by a term]
\label{thmAdjoiningTermOfBasicType}
Let $(\mathbb{C}, p)$ be a natural model and let $O \in \nmty{U}(\diamond)$. Given any natural model $(\mathbb{D}, q)$, morphism $F : (\mathbb{C}, p) \to (\mathbb{D}, q)$ and element $o \in \nmtmalt{V}(\star; FO)$, there is a unique morphism of natural models $F^{\sharp} : (\mathbb{C} \adjt{x:O}, p\adjt{x:O}) \to (\mathbb{D}, q)$ such that $F^{\sharp} \circ I = F$ and $F^{\sharp}(x) = o \in \nmtmalt{V}(\star; FO)$.
\begin{diagram}
(\mathbb{C}, p)
\arrow[r, "F"]
\arrow[d, "I"']
&
(\mathbb{D}, q)
&
o \in \nmtmalt{V}(\star; FO)
\\
(\mathbb{C} \adjt{x:O}, p\adjt{x:O})
\arrow[ur, dashed, "F^{\sharp}"']
&
x \in \nmtmadjt{x:O}{U}(\diamond\adjt{x:O}; O)
\arrow[ur, mapsto, dashed]
&
\end{diagram}
\end{theorem}

\begin{proof}
Define $F^{\sharp} = S_o \circ F_{\mathsf{tm}}$, as indicated in the following diagram.
\begin{diagram}
(\mathbb{C}, p)
\arrow[r, "F"]
\arrow[d, "I"']
&
(\mathbb{D}, q)
\arrow[d, "I"]
\arrow[dr, equals]
\\
(\mathbb{C}\adjt{x:O}, p\adjt{x:O})
\arrow[r, "F_{\mathsf{tm}}" description]
\arrow[rr, bend right=15, dashed, "F^{\sharp}"']
&
(\mathbb{D}\adjt{y:FO}, q\adjt{y:FO})
\arrow[r, "S_o" description]
&
(\mathbb{D}, q)
\end{diagram}
Note that $F^{\sharp}$ is a morphism of natural models since it is a composite of morphisms of natural models; it satisfies $F^{\sharp} \circ I = F$, since by \Cref{thmFreelyAdjoiningTermFunctorial,cnsTermSubstitutionMorphism} we have
$$F^{\sharp} \circ I = S_o \circ F_{\mathsf{tm}} \circ I = S_o \circ I \circ F = \mathrm{id}_{(\mathbb{D}, q)} \circ F = F$$
Moreover we have $F^{\sharp}(x) = S_o(F_{\mathsf{tm}}(x)) = S_o(y) = o$, as required.

It remains to prove that $F^{\sharp}$ is unique. To do so, we prove that actions of $F^{\sharp}$ on contexts, substitutions, types and terms are uniquely determined by $(F, \varphi, \nmmk{\varphi})$ and the element $o \in \nmty{V}(\star; FO)$.

For each $(\Gamma; \vec A) \in \mathrm{ob}(\mathbb{C}\adjt{x:O})$, we have
$F^{\sharp}(\Gamma; A_1, \dots, A_n) = F\Gamma \cext F\vec A[\nms{o}]$, so that the action of $F^{\sharp}$ on objects is determined by $F$, $\varphi$ and $o$; likewise on morphisms.

Given $A \in \nmtyadjt{x:O}{U}((\Gamma; \vec A)) = \nmty{U}(\Gamma \cext O \cext \vec A)$, we have $F^{\sharp}(A) = FA[\nms{o} \cext F \vec A]$, so that the action of $F^{\sharp}$ on types is determined by $F$, $\varphi$ and $o$.

Finally, given $a \in \nmtmadjt{x:O}{U}((\Gamma; \vec A); A) = \nmtm{U}(\Gamma \cext O \cext \vec A; A)$, we have $F^{\sharp}a = Fa[\nms{o} \cext F \vec A]$, so that the action of $F^{\sharp}$ on terms is determined by $F$, $\nmmk{\varphi}$ and $o$.
\end{proof}

\begin{numbered}
Although we defined $F^{\sharp}$ in terms of $F_{\mathsf{tm}}$ and $S_o$, we could instead have defined $F^{\sharp}$ directly and recovered $F_{\mathsf{tm}}$ and $S_o$ as instances of morphisms of the form $G^{\sharp}$ for appropriate choices of $G$. Specifically, we can take $F_{\mathsf{tm}} = (I \circ F)^{\sharp}$ and $S_o = (\mathrm{id}_{(\mathbb{D},q)})^{\sharp}$, with the evident choices of distinguished term in each case.
\end{numbered}

\begin{corollary}[Freely extending by a term is functorial]
\label{corFreelyAdjoiningTermFunctorial}
The assignments $(\mathbb{C}, p) \mapsto (\mathbb{C}\adjt{x:O}, p\adjt{x:O})$ and $F \mapsto F_{\mathsf{tm}}$ extend to a functor $(-)_{\mathsf{tm}} : \mathbf{NM}_{\mathsf{ty}} \to \mathbf{NM}_{\mathsf{tm}}$, which is left adjoint to the forgetful functor $\mathbf{NM}_{\mathsf{tm}} \to \mathbf{NM}_{\mathsf{ty}}$. Furthermore, the component at $(\mathbb{C}, p)$ of the unit of this adjunction is $(I, \iota, \nmmk{\iota}) : (\mathbb{C}, p) \to (\mathbb{C}\adjt{x:O}, p\adjt{x:O})$.
\end{corollary}

\begin{proof}
Given natural models $(\mathbb{C}, p)$ and $(\mathbb{D}, q)$, an element $O \in \nmty{U}(\diamond)$, a morphism $F : (\mathbb{C}, p) \to (\mathbb{D}, q)$, note that $F_{\mathsf{tm}} = (I \circ F)^{\sharp}$, where $I : (\mathbb{D}, q) \to (\mathbb{D}\adjt{y:FO}, q\adjt{y:FO})$ and the distinguished term of $(\mathbb{D}\adjt{y:FO}, q\adjt{y:FO})$ is $y = \nmq{FO}$. That $(-)_{\mathsf{tm}}$ is functorial is then immediate from the `uniqueness' part of \Cref{thmAdjoiningTermOfBasicType}, and that it is left adjoint to the forgetful functor with unit as stated is exactly the content of \Cref{thmAdjoiningTermOfBasicType}.
\end{proof}

\begin{corollary}[Free model on a family of basic types and a family of terms]
Let $I$ be an arbitrary set and let $J = \{ j_0, j_1, \dots, j_{n-1} \}$ be a finite set. The free model of the theory of an $I$-indexed family of basic types and a $J$-indexed family of terms of basic types is the natural model $(\mathbb{C}_{I;J}, p_{I;J})$ defined by
$$(\mathbb{C}_{I;J}, p_{I;J}) = (\cdots((\mathbb{C}_{I+J}, p_{I+J})_{x_0:j_0})_{x_1:j_1})\cdots)_{x_{n-1}:j_{n-1}}$$
where $(\mathbb{C}_{I+J}, p_{I+J})$ is the term model on the theory of an $(I+J)$-indexed family of basic types (\Cref{cnsFreeNaturalModelKappaBasicTypes}). In particular, $(\mathbb{C}_{I;J}, p_{I;J})$ is initial in the category $\mathbf{NM}_{(\mathsf{ty}_i)_{i \in I}, (\mathsf{tm}_j)_{j \in J}}$ with an $I$-indexed family of basic types and a $J$-indexed family of terms of basic types.
\end{corollary}

\begin{proof}
As proved in \Cref{thmFreeNaturalModelBasicTypesInitial}, the natural model $(\mathbb{C}_{I+J}, p_{I+J})$ is initial in $\mathbf{NM}_{(\mathsf{ty}_k)_{k \in I+J}}$. The natural model $(\mathbb{C}_{I;J}, p_{I;J})$ is obtained by applying functors of the form $({-})_{\mathsf{tm}}$ finitely many times. Since these functors are left adjoints, they preserve initial objects.
\end{proof}

%% file: thesis/ch4-semantics/ext-type.tex
\section{Extending a natural model by a basic type}
\label{secExtType}

\begin{construction}[Category of contexts extended by a basic type]
\label{cnsCategoryContextsExtendedByType}
Let $(\mathbb{C}, p)$ be a natural model. The \textbf{category of contexts extended by a basic type $X$} is the category $\mathbb{C}_X$ defined as follows.

\begin{itemize}
\item The objects of $\mathbb{C}_X$ are $2(n+1)$-tuples $(\Gamma, k_0, A_1, k_1, \dots, A_n, k_n)$, where $\Gamma \in \mathrm{ob}(\mathbb{C})$, for each $i<n$ we have $A_i \in \nmty{U}(\Gamma \cext A_1 \cext \dots \cext A_n)$ and $k_i \in \mathbb{N}$, and where we identify the lists
$$(\Gamma, 0, A_1, k_1, \dots, A_n, k_n) \quad \text{and} \quad (\Gamma \cext A_1, k_1, \dots, A_n, k_n)$$
Note that every object of $\mathbb{C}_X$ is either of the form $(\Gamma, 0)$ or has a unique representative of the form $(\Gamma, k_0, A_1, k_1, \dots, A_n, k_n)$ with $k_0>0$.

The idea is that the list $(\Gamma, k_0, A_1, k_1, \dots, A_n, k_n)$ should represent the context
$$\Gamma ~ \cext ~ {\underbrace{X \cext \dots \cext X}_{k_0 \text{ copies}}} \cext ~ A_1 ~ \cext ~ {\underbrace{X \cext \dots \cext X}_{k_1 \text{ copies}}} ~ \cext ~ \dots ~ \cext ~ A_n ~ \cext ~ {\underbrace{X \cext \dots \cext X}_{k_n \text{ copies}}}$$

\item A morphism from $(\Delta, \ell_1, B_1, \ell_1, \dots, B_m, \ell_m)$ to $(\Gamma, k_0, A_1, k_1, \dots, A_n, k_n)$ in $\mathbb{C}_X$ is a pair $(\sigma, h)$, where $\sigma : \Delta \cext B_1 \cext \dots \cext B_m \to \Gamma \cext A_1 \cext \dots \cext A_n$ in $\mathbb{C}$ and $h$ is a function from $k_0 + k_1 + \cdots + k_n$ to $\ell_0 + \ell_1 + \cdots + \ell_m$, with identity and composition inherited from $\mathbb{C} \times \mathbf{Fin}\op$.
\end{itemize}

Define functors $I : \mathbb{C} \to \mathbb{C}_X$, $E : \mathbb{C}_X \to \mathbb{C}$ and $G : \mathbb{C}_X \to \mathbf{Fin}\op$ by
\begin{itemize}
\item $I(\Gamma) = (\Gamma, 0)$ and $I(\sigma) = (\sigma, \mathrm{id}_0)$;
\item $E(\Gamma, k_0, A_1, k_1, \dots, A_n, k_n) = \Gamma \cext A_1 \cext \dots \cext A_n$ and $E(\sigma, h) = \sigma$;
\item $G(\Gamma, k_0, A_1, k_1, \dots, A_n, k_n) = k_0 + k_1 + \cdots + k_n$ and $G(\sigma, h) = h$.
\end{itemize}
Then $E \circ I = \mathrm{id}_{\mathbb{C}}$, $G \circ I = 0 = \upDelta(\varnothing)$ (the constant functor whose value is the empty set), and $\langle E, G \rangle : \mathbb{C}_X \to \mathbb{C} \times \mathbf{Fin}\op$ is an equivalence of categories.
\end{construction}

\begin{verification}
Note that the hom sets of $\mathbb{C}_X$ are well-defined under the identification
$$(\Gamma, 0, A_1, k_1, \dots, A_n, k_n) \sim (\Gamma \cext A_1, k_1, \dots, A_n, k_n)$$
and that the associativity and unit laws hold because identity and composition are inherited from $\mathbb{C} \times \mathbf{Fin}\op$. Well-definedness of $I$ is clear, and well-definedness of $E$ and $G$ are immediate from the fact that their action on morphisms is the same as that of the projection functors from $\mathbb{C} \times \mathbf{Fin}\op$ to its components. Evidently $E \circ I = \mathrm{id}_{\mathbb{C}}$ and $G \circ I = 0$.

To see that $\langle E, G \rangle : \mathbb{C}_X \to \mathbb{C} \times \mathbf{Fin}\op$ is an equivalence, note that it is full and faithful since it acts as the identity on morphisms, and it is essentially surjective, since a pair $(\Gamma, k) \in \mathrm{ob}(\mathbb{C} \times \mathbf{Fin}\op)$ is already an object of $\mathbb{C}_X$, and
$$\langle E, G \rangle (\Gamma, k) = (E(\Gamma, k), G(\Gamma, k)) = (\Gamma, k)$$
Hence $\langle E, G \rangle : \mathbb{C}_X \simeq \mathbb{C} \times \mathbf{Fin}\op$, as required.
\end{verification}

\begin{construction}[Free natural model extended by a basic type]
\label{cnsFreeNMExtendedByType}
Let $(\mathbb{C}, p)$ be a natural model. The \textbf{free natural model on $(\mathbb{C}, p)$ extended by a basic type $X$} is the natural model $(\mathbb{C}_X, p_X : \nmtm{U}_X \to \nmty{U}_X)$ defined by the following data. The underlying category is $\mathbb{C}_X$ (\Cref{cnsCategoryContextsExtendedByType}) with distinguished terminal object $\diamond_X = (\diamond, 0)$. The presheaves $\nmty{U}_X, \nmtm{U}_X$ and $p_X : \nmtm{U}_X \to \nmty{U}_X$ are given by
$$p_X = {!} + E^*p : G^*U + E^*\nmtm{U} \to 1 + E^*\nmty{U}$$
where $U$ is the inclusion $\mathbf{Fin} \to \mathbf{Set}$, regarded as an object of $\widehat{\mathbf{Fin}\op}$. Explicitly, we have
\begin{itemize}
\item $\nmty{U}_X(\Gamma, k_0, A_1, k_1, \dots, A_n, k_n) = \{ X \} + \nmty{U}(\Gamma \cext A_1 \cext \dots \cext A_n)$; and
\item $\nmtm{U}_X(\Gamma, k_0, A_1, k_1, \dots, A_n, k_n) = (k_0 + \cdots + k_n) + \nmtm{U}(\Gamma \cext A_1 \cext \dots \cext A_n)$;
\end{itemize}
for all $(\Gamma, k_0, A_1, k_1, \dots, A_n, k_n) \in \mathrm{ob}(\mathbb{C}_X)$, and where we have suggestively written $X$ for the unique element of $1(\vec \Gamma)$.

The representability data is defined as follows. Given $\vec \Gamma = (\Gamma, k_0, A_1, k_1, \dots, A_n, k_n) \in \mathrm{ob}(\mathbb{C}_X)$, an element of $\nmty{U}_X(\vec \Gamma)$ is either $X$ or is some $A \in \nmty{U}_X(\Gamma \cext A_1 \cext \dots \cext A_n)$.

\begin{itemize}
\item Define $\vec \Gamma \cext X = (\Gamma, k_0, A_1, k_1, \dots, A_n, k_n+1)$;
\item The projection $\vec \Gamma \cext X \to \vec \Gamma$ in $\mathbb{C}_X$ is given by the pair $(\mathrm{id}_{\Gamma \cext A_1 \cext \dots \cext A_n}, i)$, where $i : k_0 + \cdots + k_n \hookrightarrow k_0 + \cdots + k_n + 1$ is the inclusion function; and
\item The new variable $\nmtm{U}_X(\vec \Gamma \cext X) = (k_0 + \cdots + k_n + 1) + \nmtm{U}(\Gamma \cext A_1 \cext \dots \cext A_n)$ is element given by the `$+1$' term---identifying natural numbers with the corresponding von Neumann ordinals, we can take this new element to be the natural number $k_0+\dots+k_n$.
\end{itemize}

Given $A \in \nmty{U}(\Gamma \cext A_1 \cext \dots \cext A_n) \subseteq \nmty{U}_X(\vec \Gamma)$:
\begin{itemize}
\item Define $\vec \Gamma \cext A = (\Gamma, k_0, A_1, k_1, \dots, A_n, k_n, A, 0)$;
\item The projection $\vec \Gamma \cext A \to \vec \Gamma$ in $\mathbb{C}_X$ is given by the pair $(\nmp{A}, \mathrm{id}_{k_0 + \cdots + k_n})$, where $\nmp{A} : \Gamma \cext A_1 \cext \dots \cext A_n \cext A \to \Gamma \cext A_1 \cext \dots \cext A_n$ is as in $(\mathbb{C}, p)$;
\item The new variable $\nmtm{U}_X(\vec \Gamma \cext A) = k_0 + \cdots + k_n + \nmtm{U}(\Gamma \cext A_1 \cext \dots \cext A_n \cext A)$ is given by the usual element $\nmq{A} \in \nmtm{U}(\Gamma \cext A_1 \cext \dots \cext A_n)$.
\end{itemize}
The distinguished basic type of $(\mathbb{C}_X, p_X)$ is $X \in \{X\} + \nmty{U}(\diamond) = \nmty{U}_X(\diamond, 0)$.
\end{construction}

\begin{verification}
To see that $(\diamond, 0)$ is terminal, let $\vec \Gamma = (\Gamma, k_0, A_1, k_1, \dots, A_n, k_n) \in \mathrm{ob}(\mathbb{C}_X)$; there is exactly one morphism $\Gamma \cext A_1 \cext \dots \cext A_n \to \diamond$ in $\mathbb{C}$, namely $\nmt{\Gamma \cext A_1 \cext \dots \cext A_n}$, and exactly one function $0 \to k_0 + k_1 + \cdots + k_n$, namely the empty function $\varnothing$, and hence $(\nmt{\Gamma \cext A_1 \cext \dots \cext A_n}, \varnothing)$ is the unique morphism $\vec \Gamma \to (\diamond, 0)$ in $\mathbb{C}_X$. That $\nmty{U}_X$ and $\nmtm{U}_X$ are presheaves and that $p_X$ is natural are immediate from their definitions, so it remains to prove that $p_X$ is representable.

So let $\vec \Gamma = (\Gamma, k_0, A_1, k_1, \dots, A_n, k_n)$. To simplify notation, let $k = G(\vec \Gamma) = k_0 + \cdots + k_n$ and write $\Gamma \cext \vec A$ for $\Gamma \cext A_1 \cext \dots \cext A_n$. We check representability data separately for $X \in 1(\vec \Gamma)$ and for $A \in \nmty{U}(\Gamma \cext A)$.

First consider the following diagram in $\widehat{\mathbb{C}_X}$.
\begin{ldiagram}
\Yon(\vec \Gamma \cext X)
\arrow[r, "k"]
\arrow[d, "{(\mathrm{id}_{\Gamma \cext \vec A}, i)}"']
&
G^*U
\arrow[r, hook]
\arrow[d, "!" description]
&
\nmtm{U}_X
\arrow[d, "p_X"]
\\
\Yon(\vec \Gamma)
\arrow[r, "X"']
&
1
\arrow[r, hook]
&
\nmty{U}_X
\end{ldiagram}
Recall that $i$ denotes the inclusion $k \hookrightarrow k+1$, and $X$ is name we are giving to the unique element of $1(\vec \Gamma)$, and so the diagram evidently commutes. We need to verify that it is a pullback. So take an object $\vec \Delta = (\Delta, \ell_0, B_1, \ell_1, \dots, B_m, \ell_m)$ of $\mathrm{ob}(\mathbb{C}_X)$, a morphism $(\sigma, h) : \vec \Delta \to \vec \Gamma$ and an element $j \in \nmtm{U}_X(\vec \Delta)$, and assume that $(p_X)_{\vec \Delta}(j) = X[(\sigma, h)]$. Since $X : \Yon(\vec \Gamma) \to \nmty{U}_X$ factors through the inclusion $1 \hookrightarrow \nmty{U}_X$, it follows that $j : \Yon(\vec \Delta) \to \nmtm{U}_X$ factors through the inclusion $G^*U \hookrightarrow \nmtm{U}_X$. So it suffices to verify the universal property for the left-hand square. But then this amounts to verifying that there is a unique morphism $\langle \sigma, j \rangle_X$ as indicated with a dashed arrow in the following diagram in $\widehat{\mathbf{Fin}\op}$.

\begin{diagram}
\Yon(\ell)
\arrow[drr, bend left=15, "j"]
\arrow[ddr, bend right=15, "h"']
\arrow[dr, dashed]
&[-20pt]&
\\[-20pt]
&
\Yon(k+1)
\arrow[d, "i" description]
\arrow[r, "k" description]
&
U
\arrow[d, "!"]
\\
&
\Yon(k)
\arrow[r, "X"']
&
1
\end{diagram}
The existence and uniqueness of this morphism follows from representability of $U \to 1$, which is precisely the natural transformation $p_1$ of \Cref{cnsFreeNaturalModelKappaBasicTypes}; in particular, the morphism from $\ell$ to $k+1$ in $\mathbf{Fin}\op$ is the function $[h,j] : k+1 \to \ell$. The morphism $\vec \Delta \to \vec \Gamma \cext A$ in $\widehat{\mathbb{C}_X}$ is then given by $(\sigma, [h,j])$.

Given $A \in \nmty{U}(\Gamma \cext \vec A) \subseteq \nmty{U}_X(\vec \Gamma)$, consider the following diagram in $\widehat{\mathbb{C}_X}$.
\begin{diagram}
\Yon(\vec \Gamma \cext A)
\arrow[d, "{\Yon((\nmp{A}, \mathrm{id}_k))}"']
\arrow[r, "\nmq{A}"]
&
E^*\nmtm{U}
\arrow[d, "E^*p" description]
\arrow[r, hook]
&
\nmtm{U}_X
\arrow[d, "p_X"]
\\
\Yon(\vec \Gamma)
\arrow[r, "A"']
&
E^*\nmty{U}
\arrow[r, hook]
&
\nmty{U}_X
\end{diagram}
The diagram evidently commutes, so we need to verify that it is a pullback. So take an object $\vec \Delta = (\Delta, \ell_0, B_1, \ell_1, \dots, B_m, \ell_m)$ of $\mathrm{ob}(\mathbb{C}_X)$, a morphism $(\sigma, h) : \vec \Delta \to \vec \Gamma$ and an element $a \in \nmtm{U}_X(\vec \Delta)$, and assume that $(p_X)_{\vec \Delta}(a) = A[(\sigma, h)]$. Since $A : \Yon(\vec \Gamma) \to E^*\nmty{U}$ factors through the inclusion $E^*\nmty{U} \hookrightarrow \nmty{U}_X$, we have that $a : \Yon(\vec \Delta) \to \nmtm{U}_X$ factors through the inclusion $a : E^*\nmtm{U} \to \nmtm{U}_X$. So it suffices to verify the universal property for the left-hand square. But then this amounts to verifying that there is a unique morphism $\langle \sigma, a \rangle_A$ as indicated with a dashed arrow in the following diagram in $\widehat{\mathbb{C}}$.

\begin{diagram}
\Yon(\Delta \cext \vec B)
\arrow[drr, bend left=15, "a"]
\arrow[ddr, bend right=15, "\sigma"']
\arrow[dr, dashed]
&[-20pt]&
\\[-20pt]
&
\Yon(\Gamma \cext \vec A \cext A)
\arrow[r, "\nmq{A}" description]
\arrow[d, "\nmp{A}" description]
&
\nmtm{U}
\arrow[d, "p"]
\\
&
\Yon(\Gamma \cext \vec A)
\arrow[r, "A"']
&
\nmty{U}
\end{diagram}

The existence and uniqueness of this morphism follows from representability of $p$; the morphism $\vec \Delta \to \vec \Gamma \cext A$ in $\widehat{\mathbb{C}_X}$ is then given by $(\langle \sigma, a \rangle_A, h)$.

Hence the representability data exhibits $p_X$ as a representable natural transformation, and so $(\mathbb{C}_X, p_X)$ is a natural model, with $X \in \{ X \} \subseteq \nmty{U}_X(\diamond, 0)$ as its distinguished basic type.
\end{verification}

We now work towards verifying that $(\mathbb{C}_X, p_X)$ satisfies the desired universal property.

\begin{lemma}[Inclusion morphism]
\label{lemInclusionNaturalModelType}
Let $(\mathbb{C}, p)$ be a natural model. The functor $I : \mathbb{C} \to \mathbb{C}_X$ of \Cref{cnsCategoryContextsExtendedByType} extends to a morphism of natural models $(I, \iota, \nmmk{\iota}) : (\mathbb{C}, p) \to (\mathbb{C}_X, p_X)$.
\end{lemma}

\begin{proof}
Recall that $\nmty{U}_X = 1 + E^*\nmty{U}$ and $\nmtm{U}_X = G^*U + E^*\nmtm{U}$; since $E \circ I = \mathrm{id}_{\mathbb{C}}$ and $G \circ I = 0$, we have
$$I^*\nmty{U}_X = I^*1 + I^*E^*\nmty{U} = 1 + \nmty{U} \quad \text{and} \quad \nmtm{U}_X = I^*G^*U + I^*E^*\nmtm{U} = 0 + \nmtm{U} {\color{darkgray} (= \nmtm{U})}$$
so that $I^*p_X = {!} + p : \nmtm{U} \to \nmty{U}$. Let $\iota : \nmty{U} \to I^*\nmty{U}_X = 1 + \nmty{U}$ and $\nmtm{\iota} : \nmtm{U} \to I^*\nmtm{U}_X = 0 + \nmtm{U}$ be the respective inclusions. By the identification of objects in \Cref{cnsCategoryContextsExtendedByType}, we have
$$I\Gamma \cext IA = (\Gamma, 0) \cext A = (\Gamma, 0, A, 0) = (\Gamma \cext A, 0) = I(\Gamma \cext A)$$
so that $I$ respects context extension. Moreover we have
$$I(\nmp{A}) = (\nmp{A}, \mathrm{id}_0) \qquad \text{and} \quad I(\nmq{A}) = \nmq{A}$$
so that $I$ respects representability data. Hence $I$ is a morphism of natural models.
\end{proof}

\begin{numbered}
\label{numExtraSwapIsomorphisms}
Let $(\mathbb{C}, p)$ be a natural model admitting a basic type $O \in \nmty{U}(\diamond)$. By iterating swap isomorphisms (\Cref{defSwapIsomorphism}), we obtain isomorphisms
$$\theta : \Gamma ~ \cext ~ \underbrace{O \cext \dots \cext O}_{k_0 \text{ times}} ~ \cext ~ A_1 \cext ~ \underbrace{O \cext \dots \cext O}_{k_1 \text{ times}} ~ \cext ~ \dots ~ \cext A_n ~ \cext ~ \underbrace{O \cext \dots \cext O}_{k_n \text{ times}} ~~ \cong ~~ \Gamma \cext A_1 \cext \dots \cext A_n ~ \cext ~ \underbrace{O \cext \dots \cext O}_{k \text{ times}}$$
for each object $(\Gamma, k_0, A_1, k_1, \dots, A_n, k_n)$ of $\mathbb{C}_X$, where $k = k_0 + k_1 + \cdots + k_n$ and, as usual, we have suppressed projection substitutions. We can choose these isomorphisms such that, for each $j \in k$, the $j^{\text{th}}$ copy of $O$ on the left corresponds with the $j^{\text{th}}$ copy of $O$ on the right, in the sense we can express $\theta$ as a composite of swap isomorphisms containing no swap isomorphisms of the form $\nmsw_{O,O}$. By \Cref{lemExtensionByBasicTypeIsProduct}, then, the object $\Gamma \cext \vec O \cext A_1 \cext \vec O \cext \dots \cext A_n \cext \vec O$ is a product in $\mathbb{C}$ of $\Gamma \cext A_1 \cext \dots \cext A_n$ and $k$ copies of $O$. Furthermore, given another object $(\Delta, \ell_0, B_1, \ell_1, \dots, B_m, \ell_m)$ of $\mathbb{C}_X$, each pair $(\sigma, h)$ consisting of a morphism $\sigma : \Delta \cext B_1 \cext \dots \cext B_m \to \Gamma \cext A_1 \cext \dots \cext A_n$ and a function $h : k \to \ell$ gives rise to a morphism
$$\Delta \cext \vec O \cext B_1 \cext \vec O \cext \dots \cext B_m \cext \vec O \to \Gamma \cext \vec O \cext A_1 \cext \vec O \cext \dots \cext A_n \cext \vec O$$
Explicitly, this morphism is indicated by the dashed arrow in the following diagram.
\begin{diagram}
\Delta \cext \vec O \cext B_1 \cext \vec O \cext \dots \cext B_m \cext \vec O
\arrow[r, dashed]
\arrow[d, "\theta"', "\cong"]
&[50pt]
\Gamma \cext \vec O \cext A_1 \cext \vec O \cext \dots \cext A_n \cext \vec O
\\
\Delta \cext B_1 \cext \dots \cext B_m \times \underbrace{O \times \cdots \times O}_{\ell \text{ times}}
\arrow[r, "{\sigma \times \langle \pi_{h(1)}, \pi_{h(2)}, \dots, \pi_{h(n)} \rangle}"']
&
\Gamma \cext A_1 \cext \dots \cext A_n \times \underbrace{O \times \cdots \times O}_{k \text{ times}}
\arrow[u, "\theta^{-1}"', "\cong"]
\end{diagram}
where $\pi_{j}$ is the projection $O \times \cdots \times O \to O$ onto the $j^{\text{th}}$ coordinate.
\end{numbered}

\begin{theorem}[Universal property of freely extending by a basic type]
\label{thmAdjoinTypeUniversalProperty}
Let $(\mathbb{C}, p)$ be a natural model and let $(\mathbb{D}, q)$ be natural models with distinguished basic type $O \in \nmty{V}(\star)$. For each morphism of natural models $F : (\mathbb{C}, p) \to (\mathbb{D}, q)$, there is a unique morphism of natural models $F^{\sharp} : (\mathbb{C}_X, p_X) \to (\mathbb{D}, q)$ such that $F^{\sharp} \circ I = F$ and $F^{\sharp}X = O$.
\begin{diagram}
(\mathbb{C}, p)
\arrow[r, "F"]
\arrow[d, "I"']
&
(\mathbb{D}, q)
\\
(\mathbb{C}_X, p_X)
\arrow[ur, dashed, "F^{\sharp}"']
&
\end{diagram}
\end{theorem}

\begin{proof}
Let $F : (\mathbb{C}, p) \to (\mathbb{D}, q)$ be a morphism of natural models, and define $F^{\sharp} : (\mathbb{C}_X, p) \to (\mathbb{D}_X, q)$ as follows. The underlying functor $F^{\sharp} : \mathbb{C}_X \to \mathbb{D}$ is defined on objects by
$$F^{\sharp}(\Gamma, k_0, A_1, k_1, \dots, A_n, k_n) = F\Gamma ~ \cext ~ \underbrace{O \cext \dots \cext O}_{k_0 \text{ times}} ~ \cext ~ FA_1 \cext ~ \underbrace{O \cext \dots \cext O}_{k_1 \text{ times}} ~ \cext ~ \dots ~ \cext FA_n ~ \cext ~ \underbrace{O \cext \dots \cext O}_{k_n \text{ times}}$$
and on morphisms $(\sigma, h) : (\Delta, \ell_0, B_1, \ell_1, \dots, B_m, \ell_m) \to (\Gamma, k_0, A_1, k_1, \dots, A_m, k_m)$ by letting $F^{\sharp}(\sigma, h)$ be the morphism
$$F\Delta \cext \vec O \cext FB_1 \cext \vec O \cext \dots \cext FB_m \cext \vec O \to F\Gamma \cext \vec O \cext FA_1 \cext \vec O \cext \dots \cext FA_n \cext \vec O$$
induced by $F\sigma : F\Delta \cext FB_1 \cext \dots \cext FB_m \to F\Gamma \cext FA_1 \cext \dots \cext FA_n$ and $h : n \to m$ as described in \Cref{numExtraSwapIsomorphisms}. Functoriality of $F^{\sharp}$ then follows from functoriality of $F$ and the fact that $F^{\sharp}$ acts by conjugating by isomorphisms.

Given an object $\vec \Gamma = (\Gamma, k_0, A_1, k_1, \dots, A_n, k_n)$, define
$$(\varphi^{\sharp})_{\vec \Gamma} = [O, \varphi_{\Gamma \cext \vec A}] : \{ X \} + \nmty{U}(\Gamma \cext \vec A) \to \nmty{V}(F^{\sharp}\vec \Gamma)$$
Thus we have $F^{\sharp}X = O$ and $F^{\sharp}A = FA$ for each $A \in \nmty{U}(\Gamma \cext \vec A)$. Likewise, define
$$(\nmmk{\varphi}^{\sharp})_{\vec \Gamma} = q + \nmmk{\varphi}_{\Gamma \cext \vec A} : k + \nmty{V}(F^{\sharp} \vec \Gamma)$$
where $q : k \to \nmty{V}(F^{\sharp} \vec \Gamma)$ is defined for $j \in k$ by letting $q(j)$ be the (suitably weakened) element $\nmv{O}$ of $\nmtm{V}(F\vec\Gamma)$ corresponding with the $j^{\text{th}}$ copy of $O$ in $F\vec\Gamma$. Thus we have $F^{\sharp}j = \nmv{O}$ (corresponding with the appropriate copy of $O$), and $F^{\sharp}a = Fa$ for each $a \in \nmtm{U}(\Gamma \cext \vec A; A)$.

To see that $(F^{\sharp}, \varphi^{\sharp}, \nmmk{\varphi}^{\sharp})$ preserves representability data, let $\vec \Gamma = (\Gamma, k_0, A_1, k_1, \dots, A_n, k_n) \in \mathrm{ob}(\mathbb{C}_X)$ and let $A \in \nmty{U}(\Gamma \cext \vec A) \subseteq \nmty{U}_X(\vec \Gamma)$. Then
\begin{align*}
& F^{\sharp}\vec \Gamma \cext F^{\sharp}A && \\
&= (F\Gamma \cext \vec O \cext FA_1 \cext \vec O \cext \dots \cext FA_n \cext \vec O) \cext FA && \text{by definition of $F^{\sharp}$ and $\varphi^{\sharp}$} \\
&= F^{\sharp}(\Gamma, k_0, A_1, k_1, \dots, A_n, k_n, A, 0) && \text{by definition of $F^{\sharp}$} \\
&= F^{\sharp}(\vec \Gamma \cext A) && \text{by definition of context extension in $(\mathbb{C}_X, p_X)$}
\end{align*}
and likewise we have
\begin{align*}
& F^{\sharp}\vec \Gamma \cext F^{\sharp}X && \\
&= (F\Gamma \cext \vec O \cext FA_1 \cext \vec O \cext \dots \cext FA_n \cext \vec O) \cext O \\
&= F^{\sharp}(\Gamma, k_0, A_1, k_1, \dots, A_n, k_n + 1) && \\
&= F^{\sharp}(\vec \Gamma \cext X)
\end{align*}
so $F^{\sharp}$ preserves context extension.

Finally, by the construction of the action of $F^{\sharp}$ on morphisms, we immediately have that $F^{\sharp}(\nmp{A}, \mathrm{id}_k) = \nmu{FA} = \nmu{F^{\sharp}(A)}$, $F^{\sharp}(\mathrm{id}_{\Gamma \cext \vec A}, h) = \nmu{O}$, $F^{\sharp}(\nmq{A}) = \nmv{FA} = \nmv{F^{\sharp}A}$, and $F^{\sharp}(k) = \nmv{O}$. So $F^{\sharp}$ is a morphism of natural models.

We have already remarked that $F^{\sharp}(X) = O$, as required.

To see that $F^{\sharp}$ is unique, we prove that it is determined entirely by $F : (\mathbb{C}_X, p_X) \to (\mathbb{D}, q)$ and the value $F^{\sharp}X$.

For each $\vec \Gamma = (\Gamma, k_0, A_1, k_1, \dots, A_n, k_n) \in \mathrm{ob}(\mathbb{C}_X)$, we have by preservation of context extension that
$$F^{\sharp}(\Gamma, k_0, A_1, k_1, \dots, A_n, k_n) = F^{\sharp}(\Gamma,0) ~ \cext ~ \underbrace{F^{\sharp}X \cext \dots \cext F^{\sharp}X}_{k_0 \text{ times}} ~ \cext ~ F^{\sharp}A_1 \cext ~ \underbrace{F^{\sharp}X \cext \dots \cext F^{\sharp}X}_{k_1 \text{ times}} ~ \cext ~ \dots ~ \cext F^{\sharp}A_n ~ \cext ~ \underbrace{F^{\sharp}X \cext \dots \cext F^{\sharp}X}_{k_n \text{ times}}$$
But $F^{\sharp}(\Gamma,0) = F\Gamma$ and $F^{\sharp}A_i = FA_i$ for each $i$, and so the action of $F^{\sharp}$ on the objects of $\mathbb{C}_X$ is determined by $F$ and $F^{\sharp}X$. The action of $F^{\sharp}$ on morphisms is similarly determined, and hence so is the entire functor $F^{\sharp} : \mathbb{C}_X \to \mathbb{D}$.

Moreover, we have
$$\varphi^{\sharp}_{\vec \Gamma} = F^{\sharp}X + \varphi_{\Gamma \cext A_1 \cext \dots \cext A_n} : \{ X \} + \nmty{U}(\Gamma \cext A_1 \cext \dots \cext A_n) \to \nmty{V}(F^{\sharp}\vec \Gamma)$$
so that $\varphi^{\sharp}$ is determined by $F^{\sharp}X$ and $F$; and
$$\nmmk{\varphi}^{\sharp}_{\vec \Gamma} = q + \varphi_{\Gamma \cext A_1 \cext \dots \cext A_n} : k + \nmtm{U}(\Gamma \cext A_1 \cext \dots \cext A_n) \to \nmtm{V}(F^{\sharp}\vec \Gamma)$$
where $k = k_0 + k_1 + \cdots + k_n$, and $q : k \to \nmtm{V}(F^{\sharp}\vec \Gamma)$ is defined by letting $q(j)$ be the (appropriately weakened) element $\nmv{F^{\sharp}X}$ of $\nmtm{V}(F^{\sharp}\vec\Gamma)$ corresponding with the $j^{\text{th}}$ copy of $F^{\sharp}X$ in $F^{\sharp}\vec \Gamma$.

Hence the entire morphism $(F^{\sharp}, \varphi^{\sharp}, \nmmk{\varphi}^{\sharp})$ is determined by $F$ and $F^{\sharp}X$, as required.
\end{proof}

As a consequence of \Cref{thmAdjoinTypeUniversalProperty}, if $(\mathbb{D}, q)$ is a natural model with distinguished basic type $O \in \nmty{U}(\diamond)$, then specifying a morphism of natural models $(\mathbb{C}_X, p_X) \to (\mathbb{D}, q)$ preserving distinguished basic types is equivalent to specifying a morphism $(\mathbb{C}, p) \to (\mathbb{D}, q)$.

\begin{corollary}[Functoriality of freely extending by a basic type]
Let $(\mathbb{C}, p)$ and $(\mathbb{D}, q)$ be natural models. For every morphism $F : (\mathbb{C}, p) \to (\mathbb{D}, q)$, there is a unique morphism of natural models $F_{\mathsf{ty}} : (\mathbb{C}_X, p_X) \to (\mathbb{D}_Y, q_Y)$ such that $F_{\mathsf{ty}} \circ I = I \circ F$ and $FX=Y \in \nmty{V}_Y(\star_Y)$.
\begin{diagram}
(\mathbb{C}, p)
\arrow[r, "F"]
\arrow[d, "I"']
&
(\mathbb{D}, q)
\arrow[d, "I"]
\\
(\mathbb{C}_X, p_X)
\arrow[r, "F_X"']
&
(\mathbb{D}_Y, q_Y)
\end{diagram}
Moreover, the assignments $(\mathbb{C}, p) \mapsto (\mathbb{C}_X, p_X)$ and $F \mapsto F_{\mathsf{ty}}$ define a functor $(-)_{\mathsf{ty}} : \mathbf{NM} \to \mathbf{NM}_{\mathsf{ty}}$ which is left adjoint to the forgetful functor $\mathbf{NM}_{\mathsf{ty}} \to \mathbf{NM}$, and the component at $(\mathbb{C}, p)$ of the unit of this adjunction is $I : (\mathbb{C}, p) \to (\mathbb{C}_X, p_X)$.
\end{corollary}

\begin{proof}
Define $F_{\mathsf{ty}} = (I \circ F)^{\sharp}$. Then by \Cref{thmAdjoinTypeUniversalProperty} we have that $F_{\mathsf{ty}}$ is the unique morphism of natural models preserving the distinguished basic type and satsifying $F_{\mathsf{ty}} \circ I = (I \circ F)^{\sharp} \circ I = I \circ F$. Functoriality of $(-)_{\mathsf{ty}}$ follows from uniqueness, and the fact that it is left adjoint to the forgetful functor with the unit as described is exactly the content of \Cref{thmAdjoinTypeUniversalProperty}.
\end{proof}

\begin{construction}[Type insertion morphism]
\index{type insertion morphism}
Let $(\mathbb{C}, p)$ be a natural model and let $O \in \nmty{U}(\diamond)$. The \textbf{type insertion morphism} for $O$ is the unique morphism of natural models $S : (\mathbb{C}_X, p_X) \to (\mathbb{C}, p)$ such that $SX=O$ and $S \circ I = \mathrm{id}_{(\mathbb{C}, p)}$.
\end{construction}

\begin{verification}
Take $S = (\mathrm{id}_{(\mathbb{C}, p)})^{\sharp}$, where the distinguished element of $\nmty{U}(\diamond)$ is $O$.
\end{verification}

\begin{corollary}
Let $(\mathbb{C}, p)$ and $(\mathbb{D}, q)$ be natural models and let $O \in \nmty{V}(\star)$. For each morphism of natural models $F : (\mathbb{C}, p) \to (\mathbb{D}, q)$, the morphism $F^{\sharp} : (\mathbb{C}_X, p_X) \to (\mathbb{D}, q)$ of \Cref{thmAdjoinTypeUniversalProperty} factors as $F^{\sharp} = S \circ F_{\mathsf{ty}}$, where $S : (\mathbb{D}_Y, q_Y) \to (\mathbb{D}, q)$ is the type insertion morphism for $O$.
\begin{diagram}
(\mathbb{C}, p)
\arrow[r, "F"]
\arrow[d, "I"']
&
(\mathbb{D}, q)
\arrow[d, "I"]
\arrow[dr, equals]
\\
(\mathbb{C}_X, p_X)
\arrow[r, "F_{\mathsf{ty}}" description]
\arrow[rr, bend right=15, dashed, "F^{\sharp}"']
&
(\mathbb{D}_Y, q_Y)
\arrow[r, "S" description]
&
(\mathbb{D}, q)
\end{diagram}
\end{corollary}

\begin{proof}
Evidently $S \circ F_{\mathsf{ty}}$ is a morphism of natural models which preserves distinguished basic types and extends $F$, so this follows by uniqueness of $F^{\sharp}$.
\end{proof}

\begin{numbered}
The results in this section can be generalised to freely extend a natural model $(\mathbb{C}, p)$ by an $I$-indexed family of basic types $\vec X = \seqbn{X_i}{i \in I}$ for a given index set $I$. If $I$ is finite, we can simply iterate \Cref{cnsCategoryContextsExtendedByType}, but it is in fact possible for index sets of arbitrary cardinality. The new category of contexts $\mathbb{C}_{\vec X}$ is equivalent to $\mathbb{C} \times (\mathbf{Fin} \slice{I})\op$, and the representability data is similarly transported from that of $(\mathbb{C}, p)$ and of $(\mathbb{F}_I, p_I)$ (see \Cref{cnsFreeNaturalModelKappaBasicTypes}). The details of this construction are omitted in this thesis, as they are even more cumbersome.
\end{numbered}

%% file: thesis/ch4-semantics/ext-unit.tex
\section{Extending a natural model by a unit type}
\label{secExtUnit}

\begin{numbered}
To simplify notation in this section, given a natural model $(\mathbb{C}, p)$, we will write $1 = \mathrm{id}_1 : 1 \to 1$ for the identity morphism $\mathrm{id}_{\Yon(\diamond)}$ on the terminal object $\Yon(\diamond)$ of $\widehat{\mathbb{C}}$ and write $\bullet_{\Gamma}$ for the unique element of $1(\Gamma)$ for each $\Gamma \in \mathrm{ob}(\mathbb{C})$; we may just write $\bullet$ if $\Gamma$ can be inferred from context.
\end{numbered}

\begin{theorem}[Representability of $1+p$]
\label{thmOnePlusPIsRepresentable}
Let $(\mathbb{C}, p)$ be a natural model. For each $\Gamma \in \mathrm{ob}(\mathbb{C})$, the following square is a pullback;
\begin{diagram}
\Yon(\Gamma)
\arrow[r, "\bullet_{\Gamma}"]
\arrow[d, equals]
&
1+\nmtm{U}
\arrow[d, "1+p"]
\\
\Yon(\Gamma)
\arrow[r, "\bullet_{\Gamma}"']
&
1+\nmty{U}
\end{diagram}
and for each $A \in \nmty{U}(\Gamma)$, the following square is a pullback.
\begin{diagram}
\Yon(\Gamma \cext A)
\arrow[r, "\nmq{A}"]
\arrow[d, "\Yon(\nmp{A})"']
&
1+\nmtm{U}
\arrow[d, "1+p"]
\\
\Yon(\Gamma)
\arrow[r, "\nmp{A}"']
&
1+\nmty{U}
\end{diagram}
In particular, $1+p$ is representable.
\end{theorem}

\begin{proof}
The terminal natural transformation $1 : 1 \to 1$ is easily seen to be representable---indeed, its pullback along $\Yon(\Gamma) \to 1$ can be taken to be $\mathrm{id}_{\Yon(\Gamma)} : \Yon(\Gamma) \to \Yon(\Gamma)$---so the result follows from \Cref{thmClosurePropertiesOfRepresentability}(e).
\end{proof}

\begin{construction}[Category of contexts with formal unit types]
\label{cnsCategoryOfContextsWithGaps}
Let $(\mathbb{C}, p)$ be a natural model. The \textbf{category of contexts with formal unit types} of $(\mathbb{C}, p)$ is the category $\mathbb{C}\adjunit$ defined as follows.
\begin{itemize}
\item As with \Cref{cnsCategoryContextsExtendedByType}, the objects of $\mathbb{C}\adjunit$ are $2(n+1)$-tuples $(\Gamma, k_0, A_1, k_1, \dots, A_n, k_n)$, where $\Gamma \in \mathrm{ob}(\mathbb{C})$, for each $i < n$ we have $A_i \in \nmty{U}(\Gamma \cext A_1 \cext \dots \cext A_n)$ and $k_i \in \mathbb{N}$, and where we identify the lists $(\Gamma, 0, A_1, k_1, \dots, A_n, k_n)$ and $(\Gamma \cext A_1, k_1, \dots, A_n, k_n)$.

The idea is that the list $(\Gamma, k_0, A_1, k_1, \dots, A_n, k_n)$ should represent the context
$$\Gamma ~ \cext ~ {\underbrace{\mathbbm{1} \cext \dots \cext \mathbbm{1}}_{k_0 \text{ copies}}} \cext ~ A_1 ~ \cext ~ {\underbrace{\mathbbm{1} \cext \dots \cext \mathbbm{1}}_{k_1 \text{ copies}}} ~ \cext ~ \dots ~ \cext ~ A_n ~ \cext ~ {\underbrace{\mathbbm{1} \cext \dots \cext \mathbbm{1}}_{k_n \text{ copies}}}$$

\item A morphism from $(\Delta, \ell_0, B_1, \ell_1, \dots, B_m, \ell_m)$ to $(\Gamma, k_0, A_1, k_1, \dots, A_n, k_n)$ in $\mathbb{C}\adjunit$ is a morphism $\sigma : \Delta \cext B_1 \cext \dots \cext B_m \to \Gamma \cext A_1 \cext \dots \cext A_n$ in $\mathbb{C}$, with identity and composition inherited from $\mathbb{C}$.
\end{itemize}
Define functors $I : \mathbb{C} \to \mathbb{C}\adjunit$ and $E : \mathbb{C}\adjunit \to \mathbb{C}$ by
\begin{itemize}
\item $I(\Gamma) = (\Gamma, 0)$ and $I(\sigma) = \sigma$;
\item $E(\Gamma, k_0, A_1, k_1, \dots, A_n, k_n) = \Gamma \cext A_1 \cext \dots \cext A_n$ and $E(\sigma) = \sigma$.
\end{itemize}
These functors establish an equivalence of categories $\mathbb{C} \simeq \mathbb{C}\adjunit$.
\end{construction}

\begin{verification}
Note that the hom sets of $\mathbb{C}\adjunit$ are well-defined under the identification
$$(\Gamma, 0, A_1, k_1, \dots, A_n, k_n) \sim (\Gamma \cext A_1, k_1, \dots, A_n, k_n)$$
and that the associativity and unit laws hold because composition and identity are inherited from $\mathbb{C}$. Well-definedness of $I$ and $E$ is immediate from the fact that they act trivially on morphisms. Furthermore we have $E \circ I = \mathrm{id}_{\mathbb{C}}$. To see that $I \circ E \cong \mathrm{id}_{\mathbb{C}\adjunit}$, note that for each object $(\Gamma, k_0, A_1, k_1, \dots, A_n, k_n)$ of $\mathbb{C}\adjunit$ we have
$$(I \circ E)(\Gamma, k_0, A_1, k_1, \dots, A_n, k_n) = (\Gamma \cext A_1 \cext \dots \cext A_n, 0)$$
The component at $(\Gamma, k_0, A_1, k_1, \dots, A_n, k_n)$ of the natural isomorphism $I \circ E \cong \mathrm{id}_{\mathbb{C}\adjunit}$ can thus be taken to be the idenitity morphism $\mathrm{id}_{\Gamma \cext A_1 \cext \dots \cext A_n}$, which evidently defines a natural isomorphism. Hence $I$ and $E$ yield an equivalence of categories $\mathbb{C} \simeq \mathbb{C}\adjunit$.
\end{verification}

\begin{construction}[Free natural model admitting a unit type]
Let $(\mathbb{C}, p)$ be a natural model. The \textbf{free natural model admitting a unit type} on $(\mathbb{C}, p)$ is the natural model $(\mathbb{C}\adjunit, p\adjunit)$, where $\mathbb{C}\adjunit$ is as in \Cref{cnsCategoryOfContextsWithGaps} with distinguished terminal object $(\diamond, 0)$, and where the presheaves $\nmtyunit{U}, \nmtmunit{U} : \mathbb{C}\adjunit\op \to \mathbf{Set}$ and the natural transformation $p\adjunit : \nmtmunit{U} \to \nmtyunit{U}$ are obtained from $1+p : 1+\nmtm{U} \to 1+\nmty{U}$ by precopmosing with the functor $E : \mathbb{C}\adjunit \to \mathbb{C}$ (\Cref{cnsCategoryOfContextsWithGaps}).

The representability data is defined for $\vec{\Gamma} = (\Gamma, k_0, A_1, k_1, \dots, A_n, k_n) \in \mathrm{ob}(\mathbb{C}\adjunit)$ as follows.
\begin{itemize}
\item Let $\vec{\Gamma} \cext \bullet = (\Gamma, k_0, A_1, k_1, \dots, A_n, k_n+1)$ and, for each $A \in \nmty{U}(\Gamma \cext A_1 \cext \dots \cext A_n)$, let $\vec{\Gamma} \cext A = (\Gamma, k_0, A_1, k_1, \dots, A_n, k_n, A, 0)$.
\item Let $\nmp{\bullet} : (\Gamma, k_0, A_1, k_1, \dots, A_n, k_n+1) \to (\Gamma, k_0, A_1, k_1, \dots, A_n, k_n)$ be the identity morphism on $\Gamma \cext A_1 \cext \dots \cext A_n$ in $\mathbb{C}$, and let $\nmp{A} : (\Gamma, k_0, A_1, k_1, \dots, A_n, k_n, A, 0) \to (\Gamma, k_0, A_1, k_1, \dots, A_n, k_n)$ in $\mathbb{C}\adjunit$ be the morphism $\nmp{A} : \Gamma \cext A_1 \cext \dots \cext A_n \cext A \to \Gamma \cext A_1 \cext \dots \cext A_n$ in $\mathbb{C}$.
\item Let $\nmq{\bullet} = \bullet_{\Gamma \cext A_1 \cext \dots \cext A_n}$ and let the element $\nmq{A}$ be as in $(\mathbb{C}, p)$.
\end{itemize}

The unit type structure is defined by $\widehat{\mathbbm{1}} = I(\bullet_{\diamond})$ and $\widehat{\star} = I(\bullet_{\diamond})$.
\end{construction}

\begin{verification}
Note that $(\diamond, 0)$ is indeed terminal in $\mathbb{C}\adjunit$, since morphisms $(\Gamma, k_0, A_1, k_1, \dots, A_n, k_n) \to (\diamond, 0)$ in $\mathbb{C}\adjunit$ are exactly morphisms $\Gamma \cext A_1 \cext \dots \cext A_n \to \diamond$ in $\mathbb{C}$, of which there is exactly one.

Next, note that $p\adjunit = E^*(1+p) : E^*(1+\nmtm{U}) \to E^*(1+\nmty{U})$, so that $\nmty{U}\adjunit$ and $\nmtm{U}\adjunit$ are presheaves over $\mathbb{C}\adjunit$ and $p\adjunit$ is a natural transformation. Since $1+p$ is representable (\Cref{thmOnePlusPIsRepresentable}), and since $E$ is an equivalence of categories (\Cref{cnsCategoryOfContextsWithGaps}) sending the described representability data to the maps in the pullback squares witnessing representability of $1+p$, it follows that $p\adjunit$ is representable and the representability data for $p\adjunit$ exhibits $(\mathbb{C}\adjunit, p\adjunit)$ as a natural model.

Finally, consider the following square in $\widehat{\mathbb{C}\adjunit}$.
\begin{diagram}
\Yon(\diamond, 0)
\arrow[r, "\widehat{\star} = I(\bullet)"]
\arrow[d, equals]
&
\nmtmunit{U}
\arrow[d, "p\adjunit"]
\\
\Yon(\diamond, 0)
\arrow[r, "\widehat{\mathbbm{1}} = I(\bullet)"']
&
\nmtyunit{U}
\end{diagram}
It is a pullback, since the corresponding square in $\widehat{\mathbb{C}}$, namely the top square in the statement of \Cref{thmOnePlusPIsRepresentable} with $\Gamma = \diamond$, is a pullback, and $I$ is an equivalence of categories.
\end{verification}

\begin{lemma}[Inclusion morphism]
Let $(\mathbb{C}, p)$ be a natural model. The embedding $I : \mathbb{C} \hookrightarrow \mathbb{C}\adjunit$ extends to a morphism of natural models $(I, \iota, \nmmk{\iota}) : (\mathbb{C}, p) \to (\mathbb{C}\adjunit, p\adjunit)$.
\end{lemma}

\begin{proof}
First note that $I^*\nmtyunit{U} = \nmty{U}$ and $I^*\nmtmunit{U} = \nmtm{U}$, so we can let $\iota$ and $\nmmk{\iota}$ be the respective identity natural transformations. Next note that
$$I(\Gamma \cext A) = (\Gamma \cext A, 0) \overset{\star}{=} (\Gamma, 0, A, 0) = (\Gamma, 0) \cext A = I\Gamma \cext IA$$
where the equation marked $\star$ follows by our identification of presentations of objects of $\mathbb{C}\adjunit$, as described in \Cref{cnsCategoryOfContextsWithGaps}. Hence $I$ preserves context extension. Furthermore, $I\nmp{A}$ is equal as a morphism of $\mathbb{C}$ to $\nmp{A} : \Gamma \cext A \to A$, and $I\nmq{A}$ and $\nmq{A}$ are equal elements of $\nmtm{U}(\Gamma \cext A) \subseteq \nmtmunit{U}(\Gamma \cext A, 0)$, so that $(I, \iota, \nmmk{\iota})$ is indeed a morphism of natural models.
\end{proof}

\begin{lemma}[Extension of a morphism of natural models]
\label{lemAdjUnitFunctor}
For each morphism of natural models $F : (\mathbb{C}, p) \to (\mathbb{D}, q)$, there is a morphism of natural models $F\adjunit : (\mathbb{C}\adjunit, p\adjunit) \to (\mathbb{D}\adjunit, q\adjunit)$ which preserves unit type structure and for which $F\adjunit \circ I = I \circ F$.
\begin{diagram}
(\mathbb{C}, p)
\arrow[r, "F"]
\arrow[d, hook, "I"']
&
(\mathbb{D}, q)
\arrow[d, hook, "I"]
\\
(\mathbb{C}\adjunit, p\adjunit)
\arrow[r, dashed, "F\adjunit"']
&
(\mathbb{D}\adjunit, q\adjunit)
\end{diagram}
\end{lemma}

\begin{proof}
Given a morphism of natural models $F = (F, \varphi, \nmmk{\varphi}) : (\mathbb{C}, p) \to (\mathbb{D}, q)$, define $F\adjunit = (F\adjunit, \varphi\adjunit, \nmmk{\varphi}\adjunit) : (\mathbb{C}\adjunit, p\adjunit) \to (\mathbb{D}\adjunit, q\adjunit)$ as follows.
\begin{itemize}
\item Define $F\adjunit : \mathbb{C}\adjunit \to \mathbb{D}\adjunit$ on objects by
$$F\adjunit(\Gamma, k_0, A_1, k_1, \dots, A_n, k_n) = (F\Gamma, k_0, FA_1, k_1, \dots, FA_n, k_n)$$
and on morphisms by $F\adjunit(\sigma) = F(\sigma)$.
\item Define $\varphi\adjunit = \varphi_E$ and $\nmmk{\varphi}\adjunit = \nmmk{\varphi}_E$; explicitly, given $\vec{\Gamma} = (\Gamma, k_0, A_1, k_1, \dots, A_n, k_n) \in \mathrm{ob}(\mathbb{C}\adjunit)$, we have
$$(\varphi\adjunit)_{\vec{\Gamma}} = 1 + \varphi_{\Gamma \cext A_1 \cext \dots \cext A_n} : \underbrace{1 + \nmtm{U}(\Gamma \cext A_1 \cext \dots \cext A_n)}_{= \nmty{U}\adjunit(\vec{\Gamma})} \to \underbrace{1 + \nmty{V}(F\Gamma \cext FA_1 \cext \dots \cext FA_n)}_{= (F\adjunit^*\nmty{V}\adjunit)(\vec{\Gamma})}$$
and likewise $(\nmmk{\varphi}\adjunit)_{\vec{\Gamma}} = 1 + \nmmk{\varphi}_{\Gamma \cext A_1 \cext \dots \cext A_n}$.
\end{itemize}

To see that $F$ is well-defined, note that a morphism
$$\sigma : (\Delta, \ell_0, B_1, \ell_1, \dots, B_m, \ell_m) \to (\Gamma, k_0, A_1, k_1, \dots, A_n, k_n)$$
in $\widehat{\mathbb{C}}\adjunit$ is a morphism $\sigma : \Delta \cext B_1 \cext \dots \cext B_m \to \Gamma \cext A_1 \cext \dots \cext A_n$ in $\mathbb{C}$. Since $F$ preserves context extension, we have
$$F(\sigma) : F\Delta \cext FB_1 \cext \dots \cext FB_m \to F\Gamma \cext FA_1 \cext \dots \cext FA_n$$
in $\mathbb{C}$, so that $F\adjunit(\sigma) = F(\sigma)$ is a morphism of the appropriate type in $\mathbb{C}\adjunit$. That $F\adjunit$ preserves identity and composition is then immediate from functoriality of $F$, and so $F\adjunit$ is a functor.

To see that $F\adjunit$ preserves distinguished terminal objects, note that
$$F\adjunit(\diamond, 0) = (F\diamond, 0) = (\star, 0)$$
since $F$ preserves distinguished terminal objects.

That $\varphi\adjunit$ and $\nmmk{\varphi}\adjunit$ are natural transformations of the appropriate types and that $F\adjunit^*q\adjunit \circ \varphi\adjunit = \nmmk{\varphi}\adjunit \circ p\adjunit$ is immediate from their definitions.

To see that $F\adjunit$ preserves context extension, take $\vec{\Gamma} = (\Gamma, k_0, A_1, k_1, \dots, A_n, k_n) \in \mathrm{ob}(\mathbb{C}\adjunit)$ and $A \in \nmtyunit{U}(\vec{\Gamma}) = 1 + \nmty{U}(\Gamma \cext A_1 \cext \dots \cext A_n)$. If $A = \bullet$, then
\begin{align*}
F\adjunit(\vec{\Gamma} \cext \bullet) &= F\adjunit(\Gamma, k_0, A_1, k_1, \dots, A_n, k_n + 1) && \text{definition of context extension in $\mathbb{C}\adjunit$} \\
&= (F\Gamma, k_0, FA_1, k_1, \dots, FA_n, k_n + 1) && \text{definition of $F\adjunit$} \\
&= (F\Gamma, k_0, FA_1, k_1, \dots, F_n, k_n) \cext \bullet && \text{definition of context extension in $\mathbb{D}\adjunit$} \\
&= F\adjunit \vec{\Gamma} \cext F\adjunit \bullet && \text{since $F\bullet = (\varphi\adjunit)_{\vec{\Gamma}}(\bullet) = \bullet$}
\end{align*}
and if $A \in \nmty{U}(\Gamma \cext A_1 \cext \dots \cext A_n)$, then $FA \in \nmty{V}(F\Gamma \cext FA_1 \cext \dots \cext FA_n)$, and so
\begin{align*}
F\adjunit(\vec{\Gamma} \cext A) &= F\adjunit(\Gamma, k_0, A_1, k_1, \dots, A_n, k_n, A, 0) && \text{definition of context extension in $\mathbb{C}\adjunit$} \\
&= (F\Gamma, k_0, FA_1, k_1, \dots, FA_n, k_n, FA, 0) && \text{definition of $F\adjunit$} \\
&= (F\Gamma, k_0, FA_1, k_1, \dots, F_n, k_n) \cext FA && \text{definition of context extension in $\mathbb{D}\adjunit$} \\
&= F\adjunit \vec{\Gamma} \cext F\adjunit A && \text{since $F\adjunit A = (\varphi\adjunit)_{\vec{\Gamma}}(A) = FA$}
\end{align*}
so $F\adjunit$ preserves context extension.

Now $F\adjunit \nmp{\bullet}$ and $\nmp{F\adjunit \bullet}$ are equal since they are both equal to the identity morphism on $F\adjunit \vec{\Gamma}$; and $F\adjunit \nmp{A}$ and $\nmp{F\adjunit A}$ are equal since, as morphisms of $\mathbb{C}$, the former is equal to $F\nmp{A}$ and the latter is equal to $\nmp{FA}$, which are equal to each other since $F$ is a morphism of natural models. Likewise $F\adjunit \nmq{\bullet}$ and $\nmq{F\adjunit \bullet}$ are equal to the unique element of $1 \subseteq 1 + \nmty{V}(F\Gamma \cext FA_1 \cext \dots \cext FA_n) = \nmty{V}\adjunit(F\vec{\Gamma} \cext F\bullet)$, and $F\adjunit \nmq{A} = F\nmq{A}$ and $\nmq{F\adjunit A} = \nmq{FA}$ as elements of $\nmty{V}(F\Gamma \cext FA_1 \cext \dots \cext FA_n \cext FA) \subseteq \nmtyunit{V}(F\vec{\Gamma} \cext FA)$, so that $F\adjunit \nmq{A} = \nmq{F\adjunit A}$ since $F$ is a morphism of natural models.

So $F\adjunit$ is a morphism of natural models; moreover, we have already established that $F\adjunit$ preserves the unit type structure.
\end{proof}

\begin{numbered}
\label{parUnitTypeContextExtensionIsomorphisms}
In a natural model $(\mathbb{C}, p)$ admitting a unit type $\mathbbm{1} \in \nmtm{U}(\diamond)$, the morphism $\nmp{1} : \diamond \cext \mathbbm{1} \to \diamond$ is an isomorphism. By induction (and suppressing substitutions), each composite of the form
$$\nmp{\vec{\mathbbm{1}}} : \Gamma \cext \mathbbm{1} \cext \mathbbm{1} \cext \dots \cext \mathbbm{1} \xrightarrow{\nmp{1}} ~ \cdots ~ \xrightarrow{\nmp{1}} \Gamma \cext \mathbbm{1} \cext \mathbbm{1} \xrightarrow{\nmp{1}} \Gamma \cext \mathbbm{1} \xrightarrow{\nmp{1}} \Gamma$$
is an isomorphism for each $\Gamma \in \mathrm{ob}(\mathbb{C})$. It then follows that the substitution
$$\nmp{\vec{\mathbbm{1}}} \cext A : \Gamma \cext \vec{\mathbbm{1}} \cext A[\nmp{\mathbbm{1}}] \to \Gamma \cext A$$
is an isomorphism for all $A \in \nmty{U}(\Gamma)$, since it is obtained by pulling back the isomorphism $\nmp{\vec{\mathbbm{1}}}$ along $\nmp{A}$. To simplify notation, we will write the domain of $\nmp{\vec{\mathbbm{1}}} \vec A$ as $\Gamma \cext \vec{\mathbbm{1}} \cext A$. But now replacing $\Gamma$ by $\Gamma \cext \vec{\mathbbm{1}} \cext A \cext \vec{\mathbbm{1}}$, we see by iterating this process inductively that for each object $(\Gamma, k_0, A_1, k_1, \dots, A_n, k_n)$ of $\mathbb{C}\adjunit$, there is an isomorphism
$$\theta = \theta_{(\Gamma, k_0, A_1, k_1, \dots, A_n, k_n)} : \Gamma \cext \vec{\mathbbm{1}} \cext A_1 \cext \vec{\mathbbm{1}} \cext \dots \cext A_n \cext \vec{\mathbbm{1}} \longrightarrow \Gamma \cext A_1 \cext \dots \cext A_n$$
in $\mathbb{C}$. Hence for each pair of objects $(\Delta, \ell_0, B_1, \ell_1, \dots, B_m, \ell_m)$ and $(\Gamma, k_0, A_1, k_1, \dots, A_n, k_n)$ of $\mathbb{C}\adjunit$ and each morphism $\sigma : \Delta \cext B_1 \cext \dots \cext B_m \to \Gamma \cext A_1 \cext \dots \cext A_n$, we obtain a morphism
$$\sigma' : \Delta \cext \vec{\mathbbm{1}} \cext B_1 \cext \vec{\mathbbm{1}} \cext \dots \cext B_m \cext \vec{\mathbbm{1}} \to \Gamma \cext \vec{\mathbbm{1}} \cext A_1 \cext \vec{\mathbbm{1}} \cext \dots \cext A_n \cext \vec{\mathbbm{1}}$$
defined by $\sigma' = \theta^{-1} \circ \sigma \circ \theta$, and moreover $\sigma'$ is the unique morphism satisfying $\theta \circ \sigma' = \sigma \circ \theta$.
\end{numbered}

\begin{construction}[Unit insertion morphism]
\label{cnsUnitInsertionMorphism}
Let $(\mathbb{C}, p)$ be a natural model admitting a unit type. The \textbf{unit insertion morphism} is the unit type preserving morphism of natural models $N = (N, \nu, \nmmk{\nu}) : (\mathbb{C}\adjunit, p\adjunit) \to (\mathbb{C}, p)$ satisfying $N \circ I = \mathrm{id}_{(\mathbb{C}, p)}$, which is defined as follows.

The functor $N : \mathbb{C} \adjunit \to \mathbb{C}$ is defined on objects by
$$N(\Gamma, k_0, A_1, k_1, \dots, A_n, k_n) = \Gamma \cext \vec{\mathbbm{1}} \cext A_1 \cext \vec{\mathbbm{1}} \cext \dots \cext A_n \cext \vec{\mathbbm{1}}$$
where the $i^{\text{th}}$ instance of $\vec{\mathbbm{1}}$ has length $k_i$. Given a morphism $\sigma : (\Delta, \ell_0, B_1, \ell_1, \dots, B_m, \ell_m) \to (\Gamma, k_0, A_1, k_1, \dots, A_n, k_n)$ in $\mathbb{C}\adjunit$, which is a morphism $\sigma : \Delta \cext B_1 \cext \dots \cext B_m \to \Gamma \cext A_1 \cext \dots \cext A_n$ in $\mathbb{C}$, define
$$N(\sigma) = \theta^{-1} \circ \sigma \circ \theta : \Delta \cext \vec{\mathbbm{1}} \cext B_1 \cext \vec{\mathbbm{1}} \cext \dots \cext B_m \cext \vec{\mathbbm{1}} \to \Gamma \cext \vec{\mathbbm{1}} \cext A_1 \cext \vec{\mathbbm{1}} \cext \dots \cext A_n \cext \vec{\mathbbm{1}}$$
where the symbol $\theta$ refers in each case to the relevant isomorphism as described in \Cref{parUnitTypeContextExtensionIsomorphisms}.

The natural transformation $\nu : \nmtyunit{U} \to N^*\nmty{U}$ is defined by letting the component of $\nu$ at an object $(\Gamma, k_0, A_1, k_1, \dots, A_n, k_n)$ of $\mathbb{C}\adjunit$ be the function
$$\nu_{(\Gamma, k_0, A_1, k_1, \dots, A_n, k_n)} = [\mathbbm{1}, \nmty{U}(\theta)] : 1 + \nmty{U}(\Gamma \cext A_1 \cext \dots \cext A_n) \to \nmty{U}(\Gamma \cext \vec{\mathbbm{1}} \cext A_1 \cext \vec{\mathbbm{1}} \cext \dots \cext A_n \cext \vec{\mathbbm{1}})$$
where $\theta : \Gamma \cext \vec{\mathbbm{1}} \cext A_1 \cext \vec{\mathbbm{1}} \cext \dots \cext A_n \cext \vec{\mathbbm{1}} \longrightarrow \Gamma \cext A_1 \cext \dots \cext A_n$ is the isomorphism described in \Cref{parUnitTypeContextExtensionIsomorphisms}. The natural transformation $\nmmk{\nu} : \nmtmunit{U} \to N^*\nmtm{U}$ is defined likewise.
\end{construction}

\begin{verification}
First note that $N$ defines a functor: it respects identity and composition since it is defined on morphisms by conjugating by isomorphisms. Moreover $N(\diamond, 0) = \diamond$, so that $N$ preserves distinguished terminal objects.

The naturality squares for $\nu$ and $\nmmk{\nu}$ are obtained by applying $1+\nmty{U}$ and $1+\nmtm{U}$, respectively, to squares in $\mathbb{C}$ of the form
\begin{diagram}
\Delta \cext \vec{\mathbbm{1}} \cext B_1 \cext \dots \cext \vec{\mathbbm{1}} \cext B_m \cext \vec{\mathbbm{1}}
\arrow[r, "\theta"]
\arrow[d, "\theta^{-1} \circ \sigma \circ \theta"']
&
\Delta \cext B_1 \cext \dots \cext B_m
\arrow[d, "\sigma"]
\\
\Gamma \cext \vec{\mathbbm{1}} \cext A_1 \cext \dots \cext \vec{\mathbbm{1}} \cext A_n
\arrow[r, "\theta"']
&
\Gamma \cext A_1 \cext \dots \cext A_n
\end{diagram}
These evidently commute in $\mathbb{C}$, and so $\nu, \nmmk{\nu}$ are natural.

That $p \circ \nmmk{\nu} = \nu \circ p\adjunit$ follows from naturality of $p$ and the fact that $p\adjunit = E^*(1+p)$. That $N$ preserves the representability data and unit type structure is evident from the explicit definition given above.
\end{verification}

\begin{theorem}[Universal property of the free natural model admitting a unit type]
\label{thmAdjUnitType}
Let $(\mathbb{C}, p)$ be a natural model, let $(\mathbb{D}, q)$ be a natural model admitting a unit type, and let $F : (\mathbb{C}, p) \to (\mathbb{D}, q)$ be a morphism of natural models. There is a unique unit type structure preserving morphism of natural models $F^{\sharp} : (\mathbb{C}\adjunit, p\adjunit) \to (\mathbb{D}, q)$ such that $F^{\sharp} \circ I = F$.
\begin{diagram}
(\mathbb{C}, p)
\arrow[r, "F"]
\arrow[d, hook, "I"']
&
(\mathbb{D}, q)
\\
(\mathbb{C}\adjunit, p\adjunit)
\arrow[ur, dashed, "F^{\sharp}"']
&
\end{diagram}
\end{theorem}

\begin{proof}
Define $F^{\sharp} = N \circ F\adjunit$, as indicated in the following diagram.
\begin{diagram}
(\mathbb{C}, p)
\arrow[r, "F"]
\arrow[d, hook, "I"']
&
(\mathbb{D}, q)
\arrow[d, hook, "I"]
\arrow[dr, equals]
&
\\
(\mathbb{C}\adjunit, p\adjunit)
\arrow[r, "F\adjunit" description]
\arrow[rr, bend right=15, dashed, "F^{\sharp}"']
&
(\mathbb{D}\adjunit, q\adjunit)
\arrow[r, "N" description]
&
(\mathbb{D}, q)
\end{diagram}
Since $N$ and $F\adjunit$ preserve unit type structure (\Cref{lemAdjUnitFunctor} and \Cref{cnsUnitInsertionMorphism}), so does $F^{\sharp}$. Moreover we have
$$F^{\sharp} \circ I = N \circ F\adjunit \circ I = N \circ I \circ F = F$$
as required.

To see that $F^{\sharp}$ is the unique such morphism, we prove that its actions on contexts, substitutions, types and terms are determined entirely by $(F, \varphi, \nmmk{\varphi})$.

To this end, note that for each $\vec \Gamma = (\Gamma, k_0, A_1, k_1, \dots, A_n, k_n) \in \mathrm{ob}(\mathbb{C}\adjunit)$, we have
$$F^{\sharp} \vec \Gamma = F\Gamma \cextalt \vec{\mathbbm{1}} \cextalt FA_1 \cextalt \vec{\mathbbm{1}} \cextalt \dots \cextalt FA_n \cextalt \vec{\mathbbm{1}}$$
so that the action of $F^{\sharp}$ on objects is determined by that of $F$. Furthemore $F^{\sharp}$ is determined by $F$ on morphisms, since we have
$$F^{\sharp}(\sigma) = N(F\adjunit(\sigma)) = N(F(\sigma)) = \theta \circ F(\sigma) \circ \theta^{-1}$$
with the symbol $\theta$ representing the isomorphisms in $\mathbb{D}$ as discussed in \Cref{parUnitTypeContextExtensionIsomorphisms}.

Now given $A \in \nmty{U}(\Gamma \cext A_1 \cext \dots \cext A_n)$ and $a \in \nmtm{U}(\Gamma \cext A_1 \cext \dots \cext A_n; A)$, we have $F^{\sharp}A = FA$ and $F^{\sharp}a = Fa$; and $F^{\sharp}$ is uniquely determined on the unit type structure since it must preserve it.

Hence $F^{\sharp}$ is the unique unit type preserving morphism satisfying $F^{\sharp} \circ I = I \circ F$.
\end{proof}

\begin{corollary}[Freely extending by a unit type is functorial]
The assignments $(\mathbb{C}, p) \mapsto (\mathbb{C}\adjunit, p\adjunit)$ and $F \mapsto F\adjunit$ determine a functor $(-)\adjunit : \mathbf{NM} \to \mathbf{NM}\adjunit$, which is left adjoint to the forgetful functor $U : \mathbf{NM}\adjunit \to \mathbf{NM}$. Moreover, the component at $(\mathbb{C}, p)$ of the unit of this adjunction is $(I, \iota, \nmmk{\iota}) : (\mathbb{C}, p) \to (\mathbb{C}\adjunit, p\adjunit)$.
\end{corollary}

\begin{proof}
We can recover $F\adjunit$ as $(I \circ F)^{\sharp}$, where $I : (\mathbb{D}, q) \to (\mathbb{D}\adjunit, q\adjunit)$ is the inclusion morphism. As such, functoriality of $(-)\adjunit$ follows from the `uniqueness' part of \Cref{thmAdjUnitType}. That this functor is left adjoint to the forgetful functor $\mathbf{NM}\adjunit \to \mathbf{NM}$ with unit as described is then precisely the content of \Cref{thmAdjUnitType}.
\end{proof}

%% file: thesis/ch4-semantics/ext-sigma.tex
\section{Extending a natural model by dependent sum types}
\label{secExtSigma}

The idea behind freely adjoining dependent sum type structure to a natural model $(\mathbb{C}, p)$ is similar to that of freely adjoining unit type structure. First we modify the representable natural transformation $p$ to obtain a new representable natural transformation which additionally admits dependent sum types, and then we replace the base category $\mathbb{C}$ by an equivalent one that allows formal extensions of objects by dependent sum types.

Given a type $A$ and a dependent type $B$ over $A$, their dependent sum type $\nsum_{x : A} B(x)$ has as terms pairs $\langle a, b \rangle$, where $a : A$ and $b : B(a)$. Given a further dependent type $C$ over $B$, we obtain a type $\nsum_{x : A} \nsum_{y : B(x)} C(x,y)$, whose terms take the form $\langle a, \langle b, c \rangle \rangle$, and a type $\nsum_{\langle x, y \rangle : \nsum_{x:A} B(x)} C(x,y)$, whose terms take the form $\langle \langle a, b \rangle, c \rangle$. More generally, given $n$ types $A_1, A_2, \dots, A_n$, with $A_{i+1}$ depending on $A_i$ for all $i<n$, there is one iterated dependent sum type for each way of parenthesising a list with $n$ elements. As discussed in \Cref{parTrees}, these correspond with particular kinds of \textit{trees}. In order to freely admit dependent sum types, then, we take these trees of types to be our new types, whose terms are trees of terms (\Cref{defTypeTermTree}, \Cref{cnsPresheafOfTypeTermTrees}).

\begin{numbered}
\label{parTrees}
\index{tree}
Given a set $S$, the polynomial functor $\mathbf{Set} \to \mathbf{Set}$ defined on objects by $X \mapsto S + X^2$ has an initial algebra, which we can denote by $\mathrm{Tree}(S)$. The elements of $\mathrm{Tree}(S)$ are \textit{leaf-labelled finite rooted binary trees} with labels from $S$. We can generate the set $\mathrm{Tree}(S)$ inductively by declaring that $a \in \mathrm{Tree}(S)$ for each $a \in S$ and $[T_1,T_2] \in \mathrm{Tree}(S)$ for each $T_1,T_2 \in \mathrm{Tree}(S)$, so that specifying an element $T \in \mathrm{Tree}(S)$ is equivalent to specifying an inhabited list $a_1,a_2,\dots,a_n$ of elements of $S$ together with a parenthesisation of the list. To illustrate, the following leaf-labelled finite rooted binary tree is represented by the parenthesised list $[[[a,b],c],[d,e]]$.
\begin{sdiagram}
&&&&
\bullet
\arrow[dll, dash]
\arrow[drr, dash]
\\
&&
\bullet
\arrow[dl, dash]
\arrow[dr, dash]
&&&&
\bullet
\arrow[dl, dash]
\arrow[dr, dash]
\\
&
\bullet
\arrow[dl, dash]
\arrow[dr, dash]
&&
c
&&
d
&&
e
\\
a
&&
b
\end{sdiagram}
We say two leaf-labelled finite rooted binary trees \textit{have the same shape} if their underlying (unlabelled) trees are isomorphic---in practice, this means that the parenthesised lists have the same parenthesisation but may have different labels.

Write $L(T) = (a_1, a_2, \dots, a_n)$ for the ordered list of the leaves of a tree $T$---more precisely, $L$ is defined inductively by $L(a) = (a)$ and $L([T_1, T_2]) = L(T_1) \append L(T_2)$, where $\append$ is concatenation of sequences. For instance, $L([[[a,b],c],[d,e]])=(a,b,c,d,e)$.
\end{numbered}

\begin{definition}[Type trees and term trees]
\label{defTypeTermTree}
\index{tree!type {---}}
\index{tree!term {---}}
Let $(\mathbb{C}, p)$ be a natural model and let $\Gamma \in \mathrm{ob}(\mathcal{C})$.
\begin{enumerate}[(i)]
\item The set $\nmtytree{U}(\Gamma)$ of (\textbf{dependent}) \textbf{type trees} over $\Gamma$, and the set $\{ \Gamma \cext T \mid T \in \nmtytree{U}(\Gamma) \}$ of extensions of $\Gamma$ by type trees, are defined simultaneously inductively by the following rules.
\begin{itemize}
\item $A \in \nmtytree{U}(\Gamma)$ for each $A \in \nmty{U}(\Gamma)$, and $\Gamma \cext A$ coincides with the regular notion;
\item If $T_1 \in \nmtytree{U}(\Gamma)$, $\Gamma \cext T_1$ is defined, $T_2 = \nmtytree{U}(\Gamma \cext T_1)$ and $\Gamma \cext T_1 \cext T_2$ is defined, then $[T_1, T_2] \in \nmtytree{U}(\Gamma)$ and $\Gamma \cext [T_1, T_2] = \Gamma \cext T_1 \cext T_2$.
\end{itemize}
\item The set $\nmtmtree{U}(\Gamma)$ of (\textbf{dependent}) \textbf{term trees} over $\Gamma$ and the function $(p\tree)_{\Gamma} : \nmtmtree{U}(\Gamma) \to \nmtytree{U}(\Gamma)$ are defined simultaneously inductively by the following rules.
\begin{itemize}
\item $a \in \nmtmtree{U}(\Gamma)$ for each $a \in \nmtm{U}(\Gamma)$ and $(p\tree)_{\Gamma}(a) = p_{\Gamma}(a)$;
\item If $t_1 \in \nmtmtree{U}(\Gamma)$, $T_1 = (p\tree)_{\Gamma}(t_1)$ is defined, $t_2 \in \nmtmtree{U}(\Gamma \cext T_1)$ and $(p\tree)_{\Gamma \cext T_1}(t_2)$ is defined, then $[t_1,t_2] \in \nmtmtree{U}(\Gamma)$ and $(p\tree)_{\Gamma}([t_1,t_2]) = [T_1,T_2]$.
\end{itemize}
Write $\nmtmtree{U}(\Gamma; T)$ for the set of term trees $t$ over $\Gamma$ with $(p\tree)_{\Gamma}(t)=T$.
\end{enumerate}
\end{definition}

The following lemma is useful for working with the definitions of type trees and term trees given in \Cref{defTypeTermTree}.

\begin{lemma}
\label{lemLeavesOfTrees}
Let $(\mathbb{C}, p)$ be a natural model and let $\Gamma \in \mathrm{ob}(\mathbb{C})$.
\begin{enumerate}[(i)]
\item Given a type tree $T$ over $\Gamma$ with $L(T) = (A_1,A_2,\dots,A_n)$, we have $A_{i+1} \in \nmty{U}(\Gamma \cext A_1 \cext \dots \cext A_i)$ for each $i < n$; and
\item Given a term tree $t$ over $\Gamma$ with $(p\tree)_{\Gamma}(t)=T$ and $L(t) = (a_1, a_2, \dots, a_n)$, the trees $t$ and $T$ have the same shape and $a_{i+1} \in \nmtmtree{U}(\Gamma \cext A_1 \cext \dots \cext A_i; A_{i+1})$ for each $i<n$.
\end{enumerate}
\end{lemma}

\begin{proof}
Both proofs are straightforward inductions using the inductive definitions of type trees and term trees.
\end{proof}

In light of \Cref{lemLeavesOfTrees}, given a morphism $\sigma : \Delta \to \Gamma$ in a natural model $(\mathbb{C}, p)$ and a type tree $T$ over $\Gamma$ with $L(T) = (A_1,A_2,\dots,A_n)$, we will write $\sigma \cext T$ for the iterated extension $\sigma \cext A_1 \cext \dots \cext A_n$ of $\sigma$ by the leaves of $T$.

\begin{construction}[Presheaves of type trees and term trees]
\label{cnsPresheafOfTypeTermTrees}
\index{tree!presheaf of type/term {---}s}
\index{presheaf!presheaf of type/term trees@{---} of type/term trees}
Let $(\mathbb{C}, p)$ be a natural model.
\begin{enumerate}[(i)]
\item The \textbf{presheaf of type trees} in $(\mathbb{C}, p)$ is the presheaf $\nmtytree{U} : \mathbb{C}^{\mathrm{op}} \to \mathbf{Set}$ defined on objects as in \Cref{defTypeTermTree}(i) and defined on morphisms $\sigma : \Delta \to \Gamma$ inductively as follows: if $T = A \in \nmty{U}(\Gamma)$, then define $T[\sigma] = A[\sigma]$; and if $T = [T_1, T_2]$, then define $T[\sigma] = [T_1[\sigma], T_2[\sigma \cext T_1]]$.
\item The \textbf{presheaf of term trees} in $(\mathbb{C}, p)$ is the presheaf $\nmtmtree{U} : \mathbb{C}\op \to \mathbf{Set}$ defined on objects as in \Cref{defTypeTermTree}(ii) and defined on morphisms $\sigma : \Delta \to \Gamma$ inductively as follows: if $t = a \in \nmtmtree{U}(\Gamma)$, then define $t[\sigma] = a[\sigma]$; and if $t = [t_1, t_2]$, then define $t[\sigma] = [t_1[\sigma], t_2[\sigma \cext T_1]]$, where $T_1 = (p\tree)_{\Gamma}(t_1)$.
\item The natural transformation $p\tree : \nmtmtree{U}(\Gamma) \to \nmtytree{U}(\Gamma)$ is defined componentwise as in \Cref{defTypeTermTree}(ii).
\end{enumerate}
\end{construction}

\begin{verification}
Most of what needs to be verified is immediate by induction on the trees. To see that $\nmtytree{U}$ is functorial, note that by iterating \Cref{lemPastingCanonicalPullbacks} we have
$$T[\sigma][\tau] = [T_1[\sigma][\tau], (T_2[\sigma \cext T_1])[\tau \cext T_1[\sigma]]] = [T_1[\sigma \circ \tau], T_2[(\sigma \circ \tau) \cext T_1]]$$
and likewise for funtoriality of $\nmtm{U}$.
\end{verification}

\begin{theorem}[Representability of $p\tree$]
\label{thmTreeMapIsRepresentable}
Let $(\mathbb{C}, p)$ be a natural model. For each $\Gamma \in \mathrm{ob}(\mathbb{C})$ and each $T \in \nmtytree{U}(\Gamma)$, the following square is a pullback,
\begin{diagram}
\Yon(\Gamma \cext T)
\arrow[d, "\nmp{T}"']
\arrow[r, "\nmq{T}"]
&
\nmtmtree{U}
\arrow[d, "p\tree"]
\\
\Yon(\Gamma)
\arrow[r, "T"']
&
\nmtytree{U}
\end{diagram}
where $\nmp{T}$ and $\nmq{T}$ are defined inductively by the following two rules.
\begin{itemize}
\item If $T=A \in \nmty{U}(\Gamma)$, then $\nmp{T} = \nmp{A}$ and $\nmq{T} = \nmq{A}$;
\item If $T = [T_1,T_2]$ and $\nmp{T_1}, \nmq{T_1}, \nmp{T_2}, \nmq{T_2}$ are defined, then let $\nmp{T} = \nmp{T_1} \circ \nmp{T_2}$ and $\nmq{T} = [\nmq{T_1}[\nmp{T_2}], \nmq{T_2}[\nmp{T_1[\nmp{T}]}]]$.
\end{itemize}
In particular, $p\tree$ is representable.
\end{theorem}

\begin{proof}
We prove that the square commutes and is a pullback by induction on $T$. When $T=A \in \nmty{U}$, this is immediate, so suppose $T = [T_1,T_2]$ and that we have the following two pullback squares.

\begin{diagram}
\Yon(\Gamma \cext T_1)
\arrow[r, "\nmq{T_1}"]
\arrow[d, "\nmp{T_1}"']
\pullback
&
\nmtmtree{U}
\arrow[d, "p\tree"]
&[-25pt]
\Yon(\Gamma \cext T_1 \cext T_2)
\arrow[d, "\nmp{T_2}"']
\arrow[r, "\nmq{T_2}"]
\pullback
&
\nmtytree{U}
\arrow[d, "p\tree"]
\\
\Yon(\Gamma)
\arrow[r, "T_1"']
&
\nmtytree{U}
&
\Yon(\Gamma \cext T_1)
\arrow[r, "T_2"']
&
\nmtytree{U}
\end{diagram}

First we must prove that the square in the statement of the theorem commutes, which amounts to showing that $\nmq{T} \in \nmtmtree{U}(\Gamma \cext T, T[\nmp{T}])$. Now $\Gamma \cext T = \Gamma \cext T_1 \cext T_2$ and
\begin{align*}
T[\nmp{T}] &= [T_1,T_2][\nmp{T_1}][\nmp{T_2}] && \\
&= [T_1[\nmp{T_1}],\ T_2[\nmp{T_1} \cext T_1]][\nmp{T_2}] && \\
&= [\underbrace{T_1[\nmp{T_1}][\nmp{T_2}}_{= T_1[\nmp{T}]}],\ \underbrace{T_2[\nmp{T_1} \cext T_1][\nmp{T_2} \cext T_1[\nmp{T_1}]]}_{= T_2[\nmp{T} \cext T_1[\nmp{T}]]}]
\end{align*}
and so what we must prove is that
$$\nmq{T_1}[\nmp{T_2}] \in \nmtmtree{U}(\Gamma \cext T_1 \cext T_2; T_1[\nmp{T_1}][\nmp{T_2}])$$
and that
$$\nmq{T_2}[\nmp{T_1[\nmp{T_1}]}] \in \nmtmtree{U}(\Gamma \cext T_1 \cext T_2 \cext T_1[\nmp{T}]; T_2[\nmp{T_1} \cext T_1][\nmp{T_2} \cext T_1[\nmp{T_1}]])$$

The fact that $\nmq{T_1}[\nmp{T_2}] \in \nmtmtree{U}(\Gamma \cext T_1 \cext T_2; T_1[\nmp{T_1}][\nmp{T_2}])$ is immediate from naturality of $p\tree$; note also that the following diagram commutes by definition of $\nmp{T}$.
\begin{diagram}
\Yon(\Gamma \cext T_1 \cext T_2)
\arrow[r, "\Yon(\nmp{T_2})"]
\arrow[dr, dashed, "\Yon(\nmp{T})"']
&
\Yon(\Gamma \cext T_1)
\arrow[r, "\nmq{T_1}"]
\arrow[d, "\Yon(\nmp{T_1})" description]
\pullback
&
\nmtmtree{U}
\arrow[d, "p\tree"]
\\
&
\Yon(\Gamma)
\arrow[r, "T_1"']
&
\nmtytree{U}
\end{diagram}

To see that $\nmq{T_2}[\nmp{T_1[\nmp{T_1}]}] \in \nmtmtree{U}(\Gamma \cext T_1 \cext T_2 \cext T_1[\nmp{T}]; T_2[\nmp{T_1} \cext T_1][\nmp{T_2} \cext T_1[\nmp{T_1}]])$, consider the following diagram.
\begin{diagram}
\Yon(\Gamma \cext T_1 \cext T_2 \cext T_1[\nmp{T_1}][\nmp{T_2}])
\arrow[r, "\Yon(\nmp{T_1[\nmp{T}]})"]
\arrow[d, "{\Yon(\nmp{T_2} \cext T_1[\nmp{T_1}])}"']
\arrow[dr, dashed, "{\Yon(\nmp{T} \cext T_1[\nmp{T}])}" description]
&[50pt]
\Yon(\Gamma \cext T_1 \cext T_2)
\arrow[r, "\nmq{T_2}"]
\arrow[d, "\Yon(\nmp{T_2})" description]
&[30pt]
\nmtmtree{U}
\arrow[d, "p\tree"]
\\
\Yon(\Gamma \cext T_1 \cext T_1[\nmp{T_1}])
\arrow[r, "\Yon(\nmp{T_1} \cext T_1) = \Yon(\nmp{T_1[\nmp{T_1}]})" description]
\arrow[d, "\Yon(\nmp{T_1[\nmp{T_1}]})"']
&
\Yon(\Gamma \cext T_1)
\arrow[d, "\Yon(\nmp{T_1})"]
\arrow[r, "T_2"']
&
\nmtytree{U}
\\
\Yon(\Gamma \cext T_1)
\arrow[r, "\Yon(\nmp{T_1})"']
&
\Yon(\Gamma)
&
\end{diagram}
The top right square commutes by the induction hypothesis, and the top left and bottom left squares commute since they are the results of applying the Yoneda embedding to canonical pullback squares (\Cref{cnsCanonicalPullbacks}). The composite of the top two morphisms represents $\nmq{T_2}[\nmp{T_1[\nmp{T}]}]$, and the fact that this is an element of $\nmtmtree{U}(\Gamma \cext T_1 \cext T_2 \cext T_1[\nmp{T_1}][\nmp{T_2}]; T_2[\nmp{T_1} \cext T_1][\nmp{T_2} \cext T_1[\nmp{T_1}]])$ is exactly the assertion that the pasting of the top two squares commutes.

Hence $\nmq{T} \in \nmtmtree{U}(\Gamma \cext T; T[\nmp{T}])$ as required.

To see that the square in the statement of the theorem is a pullback, let $\sigma : \Delta \to \Gamma$ in $\mathbb{C}$ and let $t \in \nmtmtree{U}(\Delta; T[\sigma])$, as indicated in the outer square of the following diagram.

\begin{diagram}
\Yon(\Delta)
\arrow[drr, bend left=20, "t"]
\arrow[ddr, bend right=20, "\sigma"']
\arrow[dr, dashed, blue, "{\langle \sigma, t \rangle_T}" description]
&[-10pt]&
\\[-15pt]
&
\Yon(\Gamma \cext T)
\arrow[d, "\nmp{T}" description]
\arrow[r, "\nmq{T}" description]
&
\nmtmtree{U}
\arrow[d, "p\tree"]
\\
&
\Yon(\Gamma)
\arrow[r, "T"']
&
\nmtytree{U}
\end{diagram}

Define $\langle \sigma, t \rangle_T = \langle \langle \sigma, t_1 \rangle_{T_1}, t_2 \rangle_{T_2}$. Then
$$\nmp{T} \circ \langle \sigma, t \rangle_T = \nmp{T_1} \circ \nmp{T_2} \circ \langle \langle \sigma, t_1 \rangle_{T_1}, t_2 \rangle_{T_2} = \nmp{T_1} \circ \langle \sigma, t_1 \rangle_{T_1} = \sigma$$
and
\begin{align*}
&\nmq{T}[\langle \sigma, t \rangle_T] && \\
&= [\nmq{T_1}[\nmp{T_2}], \nmq{T_2}[\nmp{T_1[\nmp{T}]}]][\langle \langle \sigma, t_1 \rangle_{T_1}, t_2 \rangle_{T_2}] && \text{unpacking definitions} \\
&= [\nmq{T_1}[\nmp{T_2} \circ \langle \langle \sigma, t_1 \rangle_{T_1}, t_2 \rangle_{T_2}], \nmq{T_2}[\nmp{T_1[\nmp{T}]} \circ \langle \langle \sigma, t_1 \rangle_{T_1}, t_2 \rangle_{T_2} \cext T_1] && \text{by \Cref{cnsPresheafOfTypeTermTrees}(ii)} \\
&= [\nmq{T_1}[\langle \sigma, t_1 \rangle_{T_1}, \nmq{T_2}[\langle \sigma \cext T_1, t_2 \rangle_{T_2}]] && \text{reducing} \\
&= [t_1,t_2] && \text{induction hypothesis} \\
&= t && \text{definition of $t$}
\end{align*}
Uniqueness of $\langle \sigma, t \rangle_T$ then follows from its having been defined by a universal property.
\end{proof}

\begin{construction}[Category of contexts of trees]
\label{cnsCategoryOfContextsOfTrees}
Let $(\mathbb{C}, p)$ be a natural model. The \textbf{category of contexts of trees} of $(\mathbb{C}, p)$ is the category $\mathbb{C}\tree$ defined by
\begin{itemize}
\item \textbf{Objects} are lists $(\Gamma, T_1, \dots, T_k)$, where $k \ge 0$, $\Gamma \in \mathrm{ob}(\mathbb{C})$ and $T_{i+1} \in \nmtytree{U}(\Gamma \cext T_1 \cext \dots \cext T_i)$ for all $i<n$, where we identify $(\Gamma, A, T_1, \dots, T_k)$ with $(\Gamma \cext A, T_1, \dots, T_k)$ for all $\Gamma \in \mathrm{ob}(\mathbb{C})$, $A \in \nmty{U}(\Gamma)$ and type trees $T_1,\dots,T_k$;
\item \textbf{Morphisms.} A morphism $(\Delta, U_1, \dots, U_{\ell}) \to (\Gamma, T_1, \dots, T_k)$ in $\mathbb{C}\tree$ is a morphism $\sigma : \Delta \cext U_1 \cext \dots \cext U_{\ell} \to \Gamma \cext T_1 \cext \dots \cext T_k$, with identity and composition inherited from $\mathbb{C}$.
\end{itemize}
\end{construction}

\begin{verification}
Note that the hom sets of $\mathbb{C}\tree$ are well-defined under the identification
$$(\Gamma, A, T_1, \dots, T_k) \sim (\Gamma \cext A, T_1, \dots, T_k)$$
and that the associativity and unit laws hold because composition and identity are inherited from $\mathbb{C}$.
\end{verification}

\begin{lemma}
\label{lemEquivalenceWithCategoryOfTrees}
Let $(\mathbb{C}, p)$ be a natural model. The assignment $\Gamma \mapsto (\Gamma)$ extends to a full embedding $I : \mathbb{C} \hookrightarrow \mathbb{C}\tree$, and the assignment $(\Gamma, T_1, \dots, T_k) \mapsto \Gamma \cext T_1 \cext \dots \cext T_k$ extends to a full and faithful functor $E : \mathbb{C}\tree \to \mathbb{C}$. Moreover, the pair $(I,E)$ is an equivalence of categories.
\end{lemma}

\begin{proof}
Functoriality of $I$ and $E$ is immediate from the definitions. Furthermore, we have $E \circ I = \mathrm{id}_{\mathbb{C}}$. The natural isomorphism $\varepsilon : I \circ E \to \mathrm{id}_{\mathbb{C}\tree}$ is defined componentwise by letting
$$\varepsilon_{(\Gamma, T_1, \dots, T_k)} : (\Gamma, T_1, \dots, T_k) \to (\Gamma \cext T_1 \cext \dots \cext T_k)$$
in $\mathbb{C}\tree$ be the identity morphism $\Gamma \cext T_1 \cext \dots \cext T_k \to \Gamma \cext T_1 \cext \dots \cext T_k$ in $\mathbb{C}$. Naturality and invertibility of $\varepsilon$ are then trivial since all its components are identity morphisms.
\end{proof}

\begin{construction}[Free admission of dependent sum types]
\label{cnsFreeNMAdmittingSigmaTypes}
Let $(\mathbb{C}, p)$ be a natural model. The \textbf{free natural model admitting dependent sum types} on $(\mathbb{C}, p)$ is defined by the following data. The underlying category is $\mathbb{C}\adjsigma = \mathbb{C}\tree$ (\Cref{cnsCategoryOfContextsOfTrees}) with distinguished terminal object $(\diamond)$. The presheaves $\nmtysigma{U}, \nmtmsigma{U} : \mathbb{C}\adjsigma\op \to \mathbf{Set}$ and the natural transformation $p\adjsigma : \nmtmsigma{U} \to \nmtysigma{U}$ are obtained from $p\tree : \nmtmtree{U} \to \nmtytree{U}$ by precomposing with the functor $E : \mathbb{C}\adjsigma \to \mathbb{C}$ (\Cref{lemEquivalenceWithCategoryOfTrees}).

The representability data is defined for $\vec \Gamma = (\Gamma, T_1, \dots, T_k) \in \mathbb{C}\adjsigma$ and $T \in \nmtysigma{U}(\vec \Gamma)$ as follows.
\begin{itemize}
\item Let $(\Gamma, T_1, \dots, T_k) \cext T = (\Gamma, T_1, \dots, T_k, T)$;
\item Let $\nmp{T} : (\vec \Gamma; T) \to \vec \Gamma$ be the morphism $\nmp{T} : \Gamma \cext T_1 \cext \dots \cext T_k \cext T \to \Gamma \cext T_1 \cext \dots \cext T_k$ in $\mathbb{C}$ defined in the proof of \Cref{thmTreeMapIsRepresentable}.
\item Let $\nmq{T} \in \nmtmsigma{U}(\vec \Gamma, T) = \nmtmtree{U}(\Gamma \cext T_1 \cext \dots \cext T_k \cext T)$ be the element $\nmq{T}$ defined in the proof of \Cref{thmTreeMapIsRepresentable}.
\end{itemize}

The dependent sum type structure is defined as follows.
\begin{itemize}
\item The natural transformation $\widehat{\upSigma} : \sum_{T \in \nmtysigma{U}} \nmtysigma{U}^{(\nmtmsigma{U})_{T}} \to \nmtysigma{U}$ is defined by letting $\widehat{\upSigma}_{\vec \Gamma}$ be the function
$$\sum_{T \in \nmtytree{U}(\vec \Gamma)} \nmtytree{U}(\vec \Gamma \cext T) \xrightarrow{(T,T') \mapsto [T,T']} \nmtysigma{U}(\vec \Gamma)$$
where we have implicitly composed with the natural isomorphism given by \Cref{lemLemmaFive};

\item The natural transformation $\widehat{\mathsf{pair}} : \sum_{T,T'} \sum_{t \in (\nmtmsigma{U})_T} (\nmtmsigma{U})_{T'(t)} \to \nmtm{U}\adjsigma$ is defined by letting $\widehat{\mathsf{pair}}_{\vec \Gamma}$ be the function
$$\sum_{T \in \nmtytree{U}(\vec \Gamma)} \sum_{T' \in \nmtytree{U}(\vec \Gamma \cext T)} \sum_{t \in \nmtmtree{U}(\vec \Gamma; T)} \nmtmtree{U}(\vec \Gamma \cext T; T') \xrightarrow{(T,T',t,t') \mapsto [t,t']} \nmtmsigma{U}(\vec \Gamma)$$
where we have implicitly composed with the natural isomorphism given by \Cref{lemLemmaElevenPointFive}.
\end{itemize}
\end{construction}

\begin{verification}
That $(\diamond)$ is terminal in $\mathbb{C}\adjsigma$ is immediate from the fact that $\diamond$ is terminal in $\mathbb{C}$. We have
$$\nmtysigma{U} = E^*\nmty{U}\tree, \quad \nmtmsigma{U} = E^*\nmtmtree{U} \quad \text{and} \quad p\adjsigma = E^*(p\tree)$$
so that $\nmtysigma{U}$ and $\nmtmsigma{U}$ are presheaves and $p\adjsigma$ is a natural transformation. Representability of $p\adjsigma$ with representability data as defined follows immediately from \Cref{thmTreeMapIsRepresentable} and the fact that $E$ is an equivalence $\mathbb{C} \simeq \mathbb{C}\adjsigma$ which sends the representability data for $p\tree$ to that of $p\adjsigma$ (\Cref{lemEquivalenceWithCategoryOfTrees}).

The functions $\widehat{\upSigma}_{\vec \Gamma}$ and $\widehat{\mathsf{pair}}_{\vec \Gamma}$ respect substitution since $\nmtytree{U}$ and $\nmtmtree{U}$ are natural, so that $\widehat{\upSigma}$ and $\widehat{\mathsf{pair}}$ are natural transformations. Given an object $\vec \Gamma$ of $\mathbb{C}\adjsigma$, consider the following square
\begin{diagram}
\sum_{T \in \nmtytree{U}(\vec \Gamma)} \sum_{T' \in \nmtytree{U}(\vec \Gamma \cext T)} \sum_{t \in \nmtmtree{U}(\vec \Gamma; T)} \nmtmtree{U}(\vec \Gamma \cext T; T')
\arrow[r, "\widehat{\mathsf{pair}}_{\vec \Gamma}"]
\arrow[d, "\pi"']
&
\nmtmsigma{U}(\vec \Gamma)
\arrow[d, "(p\adjsigma)_{\vec \Gamma}"]
\\
\sum_{T \in \nmtytree{U}(\vec \Gamma)} \nmtytree{U}(\vec \Gamma \cext T)
\arrow[r, "\widehat{\upSigma}_{\vec \Gamma}"']
&
\nmtysigma{U}(\vec \Gamma)
\end{diagram}
The square commutes since given $(T,T',t,t')$ we have
$$(p\adjsigma)_{\vec \Gamma} (\widehat{\mathsf{pair}}_{\vec \Gamma}(T,T',t,t')) = (p\tree)_{E(\vec \Gamma)}([t,t']) = [T,T'] = \widehat{\upSigma}_{\vec \Gamma}(T,T') = \widehat{\upSigma}_{\vec \Gamma}(\pi(T,T',t,t'))$$
It is a pullback since the function $\widehat{\mathsf{pair}}_{\vec \Gamma}$ evidently restricts to bijections between the respective fibres of $\pi$ and of $(p\adjsigma)_{\vec \Gamma}$.

Hence the corresponding square in $\widehat{\mathbb{C}\adjsigma}$ is a pullback, so that $(\mathbb{C}\adjsigma, p\adjsigma)$ admits dependent sum types.
\end{verification}

\begin{lemma}[Inclusion morphism]
\label{lemInclusionMorphismSigma}
Let $(\mathbb{C}, p)$ be a natural model. The embedding $I : \mathbb{C} \hookrightarrow \mathbb{C}\adjsigma$ extends to a morphism of natural models $(I, \iota, \nmmk{\iota}) : (\mathbb{C}, p) \to (\mathbb{C}\adjsigma, p\adjsigma)$.
\end{lemma}

\begin{proof}
Note first that since $E \circ I = \mathrm{id}_{\mathbb{C}}$ we have $I^*p\adjsigma = I^*E^*p\tree = p\tree$. Let $\iota : \nmty{U} \hookrightarrow I^*\nmtysigma{U} = \nmtytree{U}$ and $\nmmk{\iota} : \nmtm{U} \hookrightarrow I^*\nmtmsigma{U} = \nmtmtree{U}$ be the respective inclusions, noting that for each $\Gamma$ we have $\nmty{U}(\Gamma) \subseteq \nmtytree{U}(\Gamma)$ and $\nmtm{U}(\Gamma) \subseteq \nmtmtree{U}(\Gamma)$. That these are natural and satisfy $I^*(p\adjsigma) \circ \nmmk{\iota} = \iota \circ p$ is immediate. We must prove that $I$ preserves distinguished terminal objects---which it does by definition---and that $(I, \iota, \nmmk{\iota})$ preserves the representability data. So let $\Gamma \in \mathrm{ob}(\mathbb{C})$ and $A \in \nmty{U}(\Gamma)$.
\begin{itemize}
\item We have $I(\Gamma) \cext I(A) = (\Gamma, A) = (\Gamma \cext A) = I(\Gamma \cext A)$, using the identification of lists described in \Cref{cnsFreeNMAdmittingSigmaTypes};
\item The projection $\nmp{IA} : (\Gamma, A) \to (\Gamma)$ in $\mathbb{C}\adjsigma$ is precisely the morphism $\nmp{A} : \Gamma \cext A \to \Gamma$ in $\mathbb{C}$, so that $I\nmp{A} = \nmp{IA}$;
\item The element $\nmq{IA} \in \nmtmsigma{U}(\Gamma, A)$ is exactly the element $\nmq{A} \in \nmtm{U}(\Gamma \cext A)$, so that $I\nmq{A} = \nmq{IA}$.
\end{itemize}
Hence $(I, \iota, \nmmk{\iota})$ is a morphism of natural models.
\end{proof}

\begin{lemma}[Extension of a morphism of natural models]
\label{lemAdjSumIsFunctorial}
For each morphism of natural models $F : (\mathbb{C}, p) \to (\mathbb{D}, q)$, there is a morphism of natural models $F\adjsigma : (\mathbb{C}\adjsigma, p\adjsigma) \to (\mathbb{D}\adjsigma, q\adjsigma)$ which preserves dependent sum types and for which $F\adjsigma \circ I = I \circ F$.
\begin{diagram}
(\mathbb{C}, p)
\arrow[r, "F"]
\arrow[d, hook, "I"']
&
(\mathbb{D}, q)
\arrow[d, hook, "I"]
\\
(\mathbb{C}\adjsigma, p\adjsigma)
\arrow[r, dashed, "F\adjsigma"']
&
(\mathbb{D}\adjsigma, q\adjsigma)
\end{diagram}
\end{lemma}

\begin{proof}
Let $(F, \varphi, \nmmk{\varphi}) : (\mathbb{C}, p) \to (\mathbb{D}, q)$ be a morphism of natural models, and define $(F\adjsigma, \varphi\adjsigma, \nmmk{\varphi}\adjsigma)$ as follows.
\begin{itemize}
\item The functor $F\adjsigma : \mathbb{C}\adjsigma \to \mathbb{D}\adjsigma$ is defined on objects by $F\adjsigma(\Gamma, T_1, \dots, T_k) = (F\Gamma, FT_1, \dots, FT_k)$ and on morphisms by $F\adjsigma(\sigma) = F(\sigma)$. Note that $F\adjsigma$ respects the identification of lists since $F$ preserves context extension, so this specification is well-defined.
\item The natural transformation $\varphi\adjsigma : \nmtysigma{U} \to F\adjsigma^*\nmtysigma{V}$ is given by defining
$$(\varphi\adjsigma)_{(\Gamma, T_1, \dots, T_k)} : \nmtytree{U}(\Gamma \cext T_1 \cext \dots \cext T_n) \to \nmtytree{V}(F\Gamma \cext FT_1 \cext \dots \cext FT_k)$$
inductively by
$$(\varphi\adjsigma)_{\vec \Gamma}(A) = \varphi_{\Gamma \cext T_1 \cext \dots \cext T_k}(A) \quad \text{and} \quad (\varphi\adjsigma)_{\vec \Gamma}([T,T']) = [(\varphi\adjsigma)_{\vec \Gamma}(T), (\varphi\adjsigma)_{\vec \Gamma \cext T}(T')]$$
\item The natural transformation $\nmmk{\varphi}\adjsigma : \nmtmsigma{U} \to F\adjsigma^*\nmtmsigma{V}$ is defined likewise.
\end{itemize}

Naturality of $\varphi\adjsigma$ and $\nmmk{\varphi}\adjsigma$ then follows from naturality of $\varphi$ and $\nmmk{\varphi}$, as does the fact that $F\adjsigma^*q \circ \nmmk{\varphi} = \varphi \circ p$. This construction further ensures that $F\adjsigma$ preserves dependent sum types.

Note that $F\adjsigma$ preserves context extension, since
\begin{align*}
&  F\adjsigma((\Gamma, T_1, \dots, T_k) \cext T) && \\
&= F\adjsigma(\Gamma, T_1, \dots, T_k, T) && \\
&= (F\Gamma, FT_1, \dots, FT_k, FT) && \\
&= (F\Gamma, FT_1, \dots, FT_k) \cext FT && \\
&= F(\Gamma, T_1, \dots, T_k) \cext FT
\end{align*}
and similarly we see that $F\adjsigma \nmp{T} = \nmp{F\adjsigma T}$ and $F\adjsigma \nmq{T} = \nmq{F\adjsigma T}$.
\end{proof}

\begin{numbered}
\label{parTreeSummationIsomorphism}
By \Cref{cnsPolynomialCompositeOfNaturalModels}, in any natural model $(\mathbb{C}, p)$ admitting dependent sum types, there is for each $\Gamma \in \mathrm{ob}(\mathbb{C})$ and each $A \in \nmty{U}(\Gamma)$ and $B \in \nmty{U}(\Gamma \cext A)$ an isomorphism $\theta : \Gamma \cext \widehat{\upSigma}(A,B) \cong \Gamma \cext A \cext B$ over $\Gamma$.
\begin{diagram}
\Gamma \cext \widehat{\upSigma}(A, B)
\arrow[dr, "{\nmp{\widehat{\upSigma}(A,B)}}"']
\arrow[rr, "\theta", "\cong"']
&&
\Gamma \cext A \cext B
\arrow[dl, "\nmp{A} \circ \nmp{B}"]
\\
&
\Gamma
&
\end{diagram}
in $\mathbb{C} \slice{\Gamma}$. Recalling that $\widehat{\upSigma}(T)$ is defined inductively for $T \in \nmtytree{U}(\Gamma)$ by $\widehat{\upSigma}(A)=A$ and $\widehat{\upSigma}([T,T']) = \widehat{\upSigma}(\widehat{\upSigma}(T), \widehat{\upSigma}(T'))$ (\Cref{thmTreeMapIsRepresentable}), we see by induction that there are isomorphisms $\theta : \Gamma \cext \widehat{\upSigma}(T) \cong \Gamma \cext T$ over $\Gamma$ for each $T \in \nmtytree{U}(\Gamma)$.
\begin{diagram}
\Gamma \cext \widehat{\upSigma}(T)
\arrow[rr, "\theta", "\cong"']
\arrow[dr, "{\nmp{\widehat{\upSigma}(T)}}"']
&&
\Gamma \cext T
\arrow[dl, "\nmp{T}"]
\\
&
\Gamma
&
\end{diagram}
Suppressing substitutions, for each $(\Gamma, T_1, \dots, T_k) \in \mathrm{ob}(\mathbb{C}\adjsigma)$, we obtain an isomorphism
$$\theta : \Gamma \cext \widehat{\upSigma}(T_1) \cext \dots \cext \widehat{\upSigma}(T_k) \xrightarrow{\cong} \Gamma \cext T_1 \cext \dots \cext T_k$$
over $\Gamma$, defined by taking canonical pullbacks. 
\end{numbered}

\begin{construction}[Tree summation morphism]
\label{cnsSumMorphism}
Let $(\mathbb{C}, p)$ be a natural model admitting dependent sum types. The \textbf{tree summation morphism} is the dependent sum type preserving morphism of natural models $S=(S,\sigma,\nmmk{\sigma}) : (\mathbb{C}\adjsigma, p\adjsigma) \to (\mathbb{C}, p)$ satisfying $S \circ I = \mathrm{id}_{(\mathbb{C}, p)}$, which is defined as follows.

The functor $S : \mathbb{C}\adjsigma \to \mathbb{C}$ is defined on objects by
$$S(\Gamma, T_1, \dots, T_k) = \Gamma \cext \widehat{\upSigma}(T_1) \cext \dots \cext \widehat{\upSigma}(T_k)$$
Given a morphism $\tau : (\Delta, U_1, \dots, U_{\ell}) \to (\Gamma, T_1, \dots, T_{\ell})$ in $\mathbb{C}\adjsigma$, define
$$S(\tau) = \theta^{-1} \circ \tau \circ \theta : \Delta \cext \widehat{\upSigma}(U_1) \cext \dots \cext \widehat{\upSigma}(U_k) \to \Gamma \cext \widehat{\upSigma}(T_1) \cext \dots \cext \widehat{\upSigma}(T_{\ell})$$
where the symbol $\theta$ refers in each case to the relevant isomorphism as described in \Cref{parTreeSummationIsomorphism}.

The natural transformation $\sigma : \nmtysigma{U} \to S^*\nmty{U}$ is defined by letting the component of $\sigma$ at an object $(\Gamma, T_1, \dots, T_k)$ of $\mathbb{C}\adjsigma$ be the function
$$\sigma_{\vec \Gamma} = \widehat{\upSigma} \circ \nmtytree{U}(\theta) : \nmtytree{U}(\Gamma \cext T_1 \cext \dots \cext T_k) \to \nmty{U}(\Gamma \cext \widehat{\upSigma}(T_1) \cext \dots \cext \widehat{\upSigma}(T_k))$$
and likewise $\nmmk{\sigma} : \nmtmsigma{U} \to S^*\nmtm{U}$ is defined by $\nmmk{\sigma}_{\vec \Gamma} = \widehat{\mathsf{pair}} \circ \nmtmtree{U}(\theta)$.
\end{construction}

\begin{verification}
First note that $S$ defines a functor: it respects identity and composition since it is defined on morphisms by conjugating by isomorphisms. Moreover $S(\diamond) = \diamond$, so that $S$ preserves distinguished terminal objects. To see that $\sigma$ is natural, let $\tau : (\Delta, \vec U) \to (\Gamma, \vec T)$ in $\mathbb{C}\adjsigma$ and note that
\begin{align*}
S^*\nmty{U}(\tau) \circ \sigma_{\vec \Gamma}
&= \nmty{U}(S\tau) \circ \sigma_{\vec \Gamma} \\
&= \nmty{U}(\theta \circ \tau \circ \theta^{-1}) \circ \nmty{U}(\theta) \circ \widehat{\upSigma} \\
&= \nmty{U}(\theta) \circ \nmty{U}(\tau) \circ \widehat{\upSigma} \\
&= \nmty{U}(\theta) \circ \widehat{\upSigma} \circ \nmtysigma{U}(\tau) && \text{since $\widehat{\upSigma}$ is natural} \\
&= \sigma_{\vec \Delta} \circ \nmtysigma{U}(\tau)
\end{align*}
as required; likewise for $\nmmk{\sigma}$.

That $S^*(p) \circ \nmmk{\sigma} = \sigma \circ p\adjsigma$ and that $S$ preserves representability data and dependent sum types follow immediately from their definitions.
\end{verification}

\begin{theorem}[Universal property of the free natural model admitting dependent sum types]
\label{thmAdjSigmaTypes}
Let $(\mathbb{C}, p)$ be a natural model, let $(\mathbb{D}, q)$ be a natural model admitting dependent sum types, and let $F : (\mathbb{C}, p) \to \mathbb{D}, q)$ be a morphism of natural models. There is a unique dependent sum preserving morphism of natural models $F^{\sharp} : (\mathbb{C}\adjsigma, p\adjsigma) \to (\mathbb{D}, q)$ such that $F^{\sharp} \circ I = F$.

\begin{diagram}
(\mathbb{C}, p)
\arrow[r, "F"]
\arrow[d, hook, "I"']
&
(\mathbb{D}, q)
\\
(\mathbb{C}\adjsigma, p\adjsigma)
\arrow[ur, dashed, "F^{\sharp}"']
&
\end{diagram}
\end{theorem}

\begin{proof}
Define $F^{\sharp} = S \circ F\adjsigma$, as indicated in the following diagram.
\begin{diagram}
(\mathbb{C}, p)
\arrow[r, "F"]
\arrow[d, hook, "I"']
&
(\mathbb{D}, q)
\arrow[d, hook, "I"]
\arrow[dr, equals]
&
\\
(\mathbb{C}\adjsigma, p\adjsigma)
\arrow[r, "F\adjsigma" description]
\arrow[rr, dashed, bend right=15, "F^{\sharp}"']
&
(\mathbb{D}\adjsigma, q\adjsigma)
\arrow[r, "S" description]
&
(\mathbb{D}, q)
\end{diagram}
Since $S$ and $F\adjsigma$ preserve dependent sum types (\Cref{lemAdjSumIsFunctorial} and \Cref{cnsSumMorphism}), so does $F^{\sharp}$. Moreover we have
$$F^{\sharp} \circ I = S \circ F\adjsigma \circ I = S \circ I \circ F = F$$
as required.

To see that $F^{\sharp}$ is unique, we prove that its action on contexts, substitutions, types and terms is determined entirely by $(F, \varphi, \nmmk{\varphi})$.

We proceed by induction. First note that, since $F^{\sharp} \circ I = F$, we have $F^{\sharp}(\Gamma) = F\Gamma$ for each $\Gamma \in \mathrm{ob}(\mathbb{C})$, and $F^{\sharp}A = FA$ and $F^{\sharp}a = Fa$ for each $A \in \nmty{U}(\Gamma) \subseteq \nmtysigma{U}((\Gamma))$ and $a \in \nmtm{U}(\Gamma; A) \subseteq \nmtm{U}((\Gamma); A)$. Now suppose $T=[T_1,T_2] \in \nmtysigma{U}((\Gamma))$ and $t = [t_1,t_2] \in \nmtysigma{U}((\Gamma); T)$, and that the values $F^{\sharp}T_1$, $F^{\sharp}T_2$, $F^{\sharp}t_1$ and $F^{\sharp}t_2$ are uniquely determined by $F$. Since $F^{\sharp}$ preserves the dependent sum type structure of $\mathbb{C}\adjsigma$, we have
$$F^{\sharp}T = \widehat{\upSigma}(F^{\sharp}T_1, F^{\sharp}T_2) \quad \text{and} \quad F^{\sharp}t = \widehat{\mathsf{pair}}(F^{\sharp}t_1, F^{\sharp}t_2)$$
so that the values $F^{\sharp}T$ and $F^{\sharp}t$ are uniquely determined by $F$.

Now let $\vec \Gamma = (\Gamma, T_1, \dots, T_k) \in \mathrm{ob}(\mathbb{C}\adjsigma)$ and suppose that the action of $F^{\sharp}$ on type trees and term trees over $\Gamma \cext T_1 \cext \dots \cext T_k$ is determined by that of $F$. Let $T \in \nmtysigma{U}(\vec \Gamma)$. Since $F^{\sharp}$ preserves context extension, we have $F^{\sharp}(\vec \Gamma \cext T) = F^{\sharp} \vec \Gamma \cext F^{\sharp}T$, which is uniquely determined by $F$ by our induction hypotheses; and then repeating the argument from the previous paragraph demonstrates that the action of $F^{\sharp}$ on type trees and term trees over $\vec \Gamma \cext T$ is uniquely determined by that of $F$.

Hence the entire morphism $F^{\sharp}$ is uniquely determined by $F$, as required.
\end{proof}

\begin{corollary}[Freely extending by dependent sum types is functorial]
The assignments $(\mathbb{C}, p) \mapsto (\mathbb{C}\adjsigma, p\adjsigma)$ and $F \mapsto F\adjsigma$ determine a functor $(-)\adjsigma : \mathbf{NM} \to \mathbf{NM}\adjsigma$, which is left adjoint to the forgetful functor $U : \mathbf{NM}\adjsigma \to \mathbf{NM}$. Moreover, the component at $(\mathbb{C}, p)$ of the unit of this adjunction is $(I, \iota, \nmmk{\iota}) : (\mathbb{C}, p) \to (\mathbb{C}\adjsigma, p\adjsigma)$.
\end{corollary}

\begin{proof}
We can recover $F\adjsigma$ as $(I \circ F)^{\sharp}$, where $I : (\mathbb{D}, q) \to (\mathbb{D}\adjsigma, q\adjsigma)$ is the inclusion morphism. As such, functoriality of $(-)\adjsigma$ follows from the `uniqueness' part of \Cref{thmAdjSigmaTypes}. That this functor is left adjoint to the forgetful functor $\mathbf{NM}\adjsigma \to \mathbf{NM}$ with unit as described is then precisely the content of \Cref{thmAdjSigmaTypes}.
\end{proof}

%

\begin{numbered}
Since $(\mathbb{C}\adjsigma, p\adjsigma)$ admits dependent sum types, we might hope---in presence of a unit type---that the corresponding polynomial pseudomonad (as in \Cref{thmUnitSigmaIffPolynomialPseudomonad}) resembles the \textit{algebraically-free monad} \cite{Kelly1980Monads} on the polynomial endofunctor $\upP_p : \widehat{\mathbb{C}} \to \widehat{\mathbb{C}}$. However, the free dependent sum type structure as described in \Cref{thmUnitSigmaIffPolynomialPseudomonad} does not yield a strict monad in general: if it did, then the type trees $[[A,B],C]$ and $[A,[B,C]]$ would be identified as a result of the associativity axiom.
\end{numbered}

%% file: thesis/ch5-reflection/_reflection.tex
The goal of this brief chapter is to outline some avenues for future research suggested by the work in this thesis.

\fancyhead[LO]{{\color{hdcol}\fontfamily{bch}\itshape \nouppercase{\leftmark}}}
\subsection*{Adjusting adjustments}

The definition of an \textit{adjustment} $\alpha : \varphi \pRrightarrow \psi$ between morphisms of polynomials (\Cref{defAdjustment}) is motivated by the observation that adjustments between cartesian morphisms of polynomials correspond with natural transformations between the full and faithful internal functors induced by those cartesian morphisms. A consequence is that there is at most one adjustment between any parallel pair of cartesian morphisms of polynomials. Unfortunately, we were not able to prove that adjustments between \textit{arbitrary} morphisms of polynomials correspond with anything meaningful, or indeed that they form the 3-cells of a tricategory $\mathfrak{Poly}_{\mathcal{E}}$. So although the definition we provided captures \textit{some} notion of 3-cell, which works for our purposes, a worthwhile goal in future work is to find a more suitable (and likely more general) notion of 3-cell, or to demonstrate that adjustments do in fact form a meaningful notion of 3-cell.

\subsection*{Generalised natural models}

Our work relating natural models with polynomials in \Cref{chPolynomialsRepresentability} was done for the most part without relying on any aspects of $\widehat{\mathbb{C}}$ other than its locally cartesian closed structure. It may therefore be possible to extend the definition of natural model to a more general class of categories, so that a natural model is a morphism $p : \nmtm{U} \to \nmty{U}$ in (say) a cocomplete, locally cartesian closed category $\mathcal{E}$, subject to certain conditions that are equivalent to representability in the case when $\mathcal{E} = \widehat{\mathbb{C}}$ for some small category $\mathbb{C}$. In order to make this definition meaningful, it would need to be established how to interpret the rules governing dependent type theory in such an object.

Along similar lines, recall that the main challenge of \Cref{secPolynomialPseudomonads} was to find a notion of \textit{equivalence} with respect to which a natural model admitting certain type theoretic structure gave rise to a \textit{pseudo}monad, since we discovered that it does not give rise to a strict monad. If we were to define the notion of a natural model within homotopy type theory, say, then it may be the case that the conditions for a natural model to admit a unit type and dependent sum types can now be expressed in terms of (homotopy) pullbacks yielding a monad (up to propositional equality). If this is the case, then the results of \Cref{chPolynomialsRepresentability} could be recast in terms of locally cartesian closed \textit{quasi}categories, with the rules for polynomial monads and algebras holding up to propositional equality.

\subsection*{Free natural models}

In \Cref{secExtSigma} we remarked that the set $\mathrm{Tree}(S)$ of finite rooted binary trees with leaves labelled by elements of a set $S$ is obtained as an initial algebra for the endofunctor $X \mapsto S + X \times X$. Given a natural model $(\mathbb{C}, p)$, it appears that the natural model $(\mathbb{C}\adjsigma, p\adjsigma)$ of \Cref{cnsFreeNMAdmittingSigmaTypes} is an initial algebra for the `endofunctor' $f \mapsto p + f \cdot f$, where $\cdot$ refers to polynomial composition. Indeed, there is a morphism of natural models $(\mathbb{C}, p) \to (\mathbb{C}\adjsigma, p\adjsigma)$ (this is \Cref{lemInclusionMorphismSigma}) and there is a cartesian morphism of polynomials $p\adjsigma \cdot p\adjsigma \pto p\adjsigma$ since $p\adjsigma$ admits dependent sum types (\Cref{cnsFreeNMAdmittingSigmaTypes}, \Cref{thmAdmittingSigmaTypes}), and these morphisms satisfy nice universal properties. The problem is to find the category in which $f \mapsto p + f \cdot f$ defines an endofunctor with respect to which $(\mathbb{C}\adjsigma, p\adjsigma)$ is an initial algebra.

This train of thought could be explored even further. Just as polynomial endofunctors on $\mathbf{Set}$ generalise those of the form $X \mapsto A_0 + A_1 \times X + {\cdots} + A_n \times X^n$, it would be worthwhile to explore whether there is a similar generalisation of endofunctors of the form
$$f \mapsto a_0 + a_1 \cdot f + a_2 \cdot f \cdot f + {~\cdots~} + a_n \cdot f \cdot {\dots} \cdot f$$
where $a_0,\dots,a_n,f$ denote (suitable) morphisms in a (suitable) locally cartesian closed category, and where $\cdot$ is polynomial composition. With such a notion established, it could be applied to dependent type theory to see if a natural model can be freely extended by other kinds of type theoretic structure by taking an initial algebra for such an endofunctor.

Much work remains to be done on the construction of free natural models. For example, it remains an open problem to construct a left adjoint to the forgetful functor $\mathbf{NM}_{\upPi} \to \mathbf{NM}$, thus obtaining the free natural model admitting dependent product types on a given natural model $(\mathbb{C}, p)$. Furthermore, as discussed at the end of \Cref{secInterpretationsInitiality}, it remains open to find a general way to compose these free functors.

\subsection*{Term models and interpretations}

The free natural models studied in \Cref{chSemantics} were \textit{algebraic}, rather than \textit{logical}, constructions. An important task for the future is to define a logical notion of \textit{interpretation} of a theory $\mathbb{T}$ in a natural model $(\mathbb{C}, p)$, and to construct the \textit{term model} on a given dependent type theory $\mathbb{T}$, according to the following schema.

\begin{schema}[Term model of a type theory $\mathbb{T}$]
\label{schTermModel}
Let $\mathbb{T}$ be a dependent type theory. The \textbf{term model} of $\mathbb{T}$ is the natural model $(\mathbb{C}_{\mathbb{T}}, p_{\mathbb{T}} : \nmtm{U}_{\mathbb{T}} \to \nmty{U}_{\mathbb{T}})$ defined as follows.

\begin{itemize}
\item The underlying \textbf{category} $\mathbb{C}_{\mathbb{T}}$ has as its set of objects the quotient $\mathsf{Ctx}_{\mathbb{T}} / {\equiv_{\mathsf{Ctx}}}$, where $\mathsf{Ctx}_{\mathbb{T}}$ is the set of all well-formed contexts in $\mathbb{T}$ and $\equiv_{\mathsf{Ctx}}$ identifies $\Gamma = {x_1:A_1,x_2:A_2, \dots, x_m:A_m}$ and $\Delta = {y_1:a_1,y_2:A_2,\dots,y_n:A_n}$ if and only if $m=n$ and $x_1:A_1, \dots x_{i-1}:A_{i-1} \vdash A_i=B_i$ is provable in $\mathbb{T}$ for each $i$. Given contexts $\Gamma = {x_1:A_1,x_2:A_2, \dots, x_m:A_m}$ and  $\Delta = {y_1:a_1,y_2:A_2,\dots,y_n:A_n}$, the hom set $\mathbb{C}_{\mathbb{T}}([\Delta], [\Gamma])$ is the quotient $\mathsf{Sub}_{\mathbb{T}}(\Delta,\Gamma)/{\equiv_{\mathsf{Sub}}}$, where $\mathsf{Sub}_{\mathbb{T}}(\Delta, \Gamma)$ is the set of all well-formed substitutions $(t_1,t_2,\dots,t_m)$ from $\Delta$ to $\Gamma$ (see \Cref{secTypesForTheWorkingMathematician}), and $\equiv_{\mathsf{Sub}}$ identifies $(t_1,t_2,\dots,t_m)$ with $(t'_1,t'_2,\dots,t'_m)$ whenever $\Delta \vdash t_i=t'_i : A_i(t_1,\dots,t_{i-1})$ is provable in $\mathbb{T}$ for each $i$.
\item The distinguished \textbf{terminal object} of $\mathbb{C}_{\mathbb{T}}$ is the ($\equiv_{\mathsf{Ctx}}$-equivalence class of the) empty context.
\item The \textbf{presheaf of types} $\nmty{U}_{\mathbb{T}}$ is defined on objects by letting $\nmty{U}_{\mathbb{T}}([\Gamma])$ be the quotient $\mathsf{Type}_{\mathbb{T}}(\Gamma) / {\equiv_{\mathsf{Type}}}$, where $\mathsf{Type}_{\mathbb{T}}(\Gamma)$ is the set of well-formed types in context $\Gamma$ and $\equiv_{\mathsf{Type}}$ identifies $A$ with $A'$ whenever $\Gamma \vdash A=A'$ is provable in $\mathbb{T}$; and $\nmty{U}_{\mathbb{T}}$ is defined on morphisms by letting $\nmty{U}_{\mathbb{T}}([\sigma])([A]) = [A[\sigma]]$.
\item Likewise, the \textbf{presheaf of terms} $\nmtm{U}_{\mathbb{T}}$ is defined on objects by letting $\nmtm{U}_{\mathbb{T}}([\Gamma])$ be the quotient $\mathsf{Term}_{\mathbb{T}}(\Gamma) / {{\equiv}_{\mathsf{Term}}}$, where $\mathsf{Term}_{\mathbb{T}}(\Gamma)$ is the set of well-formed terms in context $\Gamma$ and $\equiv_{\mathsf{Term}}$ identifies $a$ with $a'$ whenever the types of $a$ and $a'$ are identified by $\equiv_{\mathsf{Type}}$ and $\Gamma \vdash a = a' : A$ is provable in $\mathbb{T}$.
\item The \textbf{typing natural transformation} $p_{\mathbb{T}} : \nmtm{U}_{\mathbb{T}} \to \nmty{U}_{\mathbb{T}}$ is defined componentwise by letting $(p_{\mathbb{T}})_{[\Gamma]}([a])$ be $[A]$ for the unique $[A] \in \nmty{U}_{\mathbb{T}}([\Gamma])$ such that $\Gamma \vdash a : A$ is provable in $\mathbb{T}$.
\item The \textbf{representability data} for $(\mathbb{C}_{\mathbb{T}}, p_{\mathbb{T}})$ is defined for $[\Gamma] = [x_1:A_1, \dots, x_n:A_n] \in \mathrm{ob}(\mathbb{C}_{\mathbb{T}})$ and $[A] \in \nmty{U}_{\mathbb{T}}([\Gamma])$ as follows.
\begin{itemize}
\item Let $[\Gamma] \cext [A] = [\Gamma, x:A]$, where $x$ is a fresh variable;
\item Let $\nmp{A} = [(x_1,\dots,x_n)] : [\Gamma, x:A] \to [\Gamma]$;
\item Let $\nmq{A} = [x] \in \nmtm{U}_{\mathbb{T}}([\Gamma, x:A], [A])$.
\end{itemize}
\item If the theory $\mathbb{T}$ admits a unit type, dependent sum types, dependent product types, or some combination thereof, define the corresponding structure on $\mathbb{C}_{\mathbb{T}}$ in the evident way. \qed
\end{itemize}
\end{schema}

With this done, it should be the case that the term model is initial in $\mathbf{NM}_{\mathbb{T}}$, and from this it would follow that interpretations of $\mathbb{T}$ in a (suitably structured) natural model $(\mathbb{C}, p)$ correspond with (structure preserving) morphisms of natural models $(\mathbb{C}_{\mathbb{T}}, p_{\mathbb{T}}) \to (\mathbb{C}, p)$. The parenthetical remarks about structure depend the \textit{doctrine} in which the theory $\mathbb{T}$ lives---for example, if $\mathbb{T}$ is a theory in the doctrine of dependent type theories with a unit type and dependent sum types, then the setting for interpretations of $\mathbb{T}$ is the category of natural models admitting a unit type and dependent sum types.